\documentclass[bj]{imsart}

\RequirePackage{amsthm,amsmath,amsfonts,amssymb,amsfonts,amsbsy,mathrsfs,algorithmicx,algorithm,multirow,enumitem}
\RequirePackage[authoryear]{natbib}

\def\T{{\scriptscriptstyle \top}}
\usepackage{color}

\newcommand{\Cov}[0]{\text{Cov}}
\newcommand{\Var}[0]{\text{Var}}

\newcommand{\diag}[0]{\text{diag}}

\newcommand{\calF}[0]{\mathcal{F}}

\newcommand{\calP}[0]{\mathcal{P}}

\newcommand{\R}[0]{\mathbb{R}}

\newcommand{\Z}[0]{\mathbb{Z}}

\newcommand{\bx}{{\mathbf x}}

\newcommand{~}{\quad}

\newtheorem{theorem}{Theorem}
\newtheorem{lemma}{Lemma}
\newtheorem{cy}{Corollary}
\newtheorem{proposition}{Proposition}

\theoremstyle{definition}
\newtheorem{definition}{Definition}
\newtheorem{remark}{{\rm\bf Remark}}

\newtheorem{ass}{Condition}

\renewcommand{\leq}{\leqslant}
\renewcommand{\geq}{\geqslant}

\def\qed{ \hfill $\blacksquare$}

\newcommand{\cA}{\mathcal{A}}

\newcommand{\cI}{\mathcal{I}}
\newcommand{\cJ}{\mathcal{J}}

\newcommand{\bP}{\mathbb{P}}

\DeclareMathOperator{\var}{Var}
\DeclareMathOperator{\cov}{Cov}

\begin{document}

\begin{frontmatter}
\title{Central limit theorems for high dimensional dependent data}
\runtitle{CLTs for high dimensional dependent data}

\begin{aug}
	\author[A,B]{\fnms{Jinyuan} \snm{Chang}\ead[label=e1,mark]{changjinyuan@swufe.edu.cn}},
	\author[C,D]{\fnms{Xiaohui} \snm{Chen}\ead[label=e2,mark]{xiaohuic@usc.edu}} and
	\author[A]{\fnms{Mingcong} \snm{Wu}\ead[label=e3,mark]{wumingcong@smail.swufe.edu.cn}}
	\address[A]{Joint Laboratory of Data Science and Business Intelligence, Southwestern University of Finance and Economics, Chengdu, China, \printead{e1,e3}}
	
	\address[B]{Academy of Mathematics and Systems Science, Chinese Academy of Sciences, Beijing, China}
	
	\address[C]{Department of Mathematics, University of Southern California, Los Angeles, CA, USA, \printead{e2}}
	
	\address[D]{Department of Statistics, University of Illinois at Urbana-Champaign, Champaign, IL, USA}
\end{aug}

\begin{abstract}
Motivated by statistical inference problems in high-dimensional time series data analysis, we first derive non-asymptotic error bounds for Gaussian approximations of sums of high-dimensional dependent random vectors on hyper-rectangles, simple convex sets and sparsely convex sets. We investigate the quantitative effect of temporal dependence on the rates of convergence to a Gaussian random vector over three different dependency frameworks ($\alpha$-mixing, $m$-dependent, and physical dependence measure). In particular, we establish new error bounds under the $\alpha$-mixing framework and derive faster rate over existing results under the physical dependence measure. To implement the proposed results in practical statistical inference problems, we also derive a data-driven parametric bootstrap procedure based on a kernel-type estimator for the long-run covariance matrices. The unified Gaussian and parametric bootstrap approximation results can be used to test mean vectors with combined $\ell^2$ and $\ell^\infty$ type statistics, do change point detection, and construct confidence regions for covariance and precision matrices, all for time series data.
\end{abstract}



\end{frontmatter}

\section{Introduction}

High-dimensional dependent data are frequently encountered in current practical problems of finance, biomedical sciences, geological studies and many more areas. Due to the complicated dependency among different components and nonlinear dynamical behaviors in the series, there have been tremendous challenges in developing principled statistical inference procedures for such data. 
Most existing methods require certain parametric assumptions on the underlying data generation mechanism or structural assumptions on the dependency among different components in order to derive asymptotically pivotal distributions of the involved statistics. Assumptions of this kind are not only difficult to be verified but also often violated in real data. How to derive statistically valid inference procedures that do not rely on specific structural assumptions imposed on the dependency among different components of high-dimensional dependent data has been an urgent demand.

In this paper, we focus on establishing quantitative high-dimensional Central Limit Theorems (CLTs) and related parametric bootstrap approximations for dependent (and possibly non-stationary) data. Let $\mathcal{X}_n = \{X_1,\ldots,X_n\}$ be a sequence of $p$-dimensional dependent random vectors with mean zero, i.e., $\mathbb{E}(X_t)=0$. Write $S_{n,x} = n^{-1/2} \sum_{t=1}^n X_t$. Denote the instantaneous covariance matrix of $X_t$ at time point $t$ by $\Sigma_t=\Cov(X_t)$ and the long-run covariance matrix of $\{X_t\}_{t=1}^n$ by $\Xi=\Cov(S_{n,x})$. Our main goal is to bound
\begin{align}\label{eq:gb}
\rho_n(\mathcal{A}) := \sup_{A\in\mathcal{A}}|\mathbb{P}(S_{n,x} \in A)-\mathbb{P}(G\in A)|\,,
\end{align}
where $G\sim N(0,\Xi)$ and $\mathcal{A}$ is a class of Borel subsets in $\mathbb{R}^p$. Gaussian and parametric bootstrap approximation results over a rich index class $\mathcal{A}$ for large $p$ are fundamental tools in developing downstream statistical inference procedures for a wide spectrum of problems in the high-dimensional setting, including, for example, inference of mean vector, change point detection, structure checking of instantaneous covariance matrix, and testing white noise hypothesis. We refer the readers to Section \ref{sec:apps} for more details of these applications.


When $X_1,\dots,X_n$ are independent random vectors in $\R^p$, the problem of bounding $\rho_n(\mathcal{A})$ for a variety of choices $\cA$ is a classical research topic in probability theory \citep{Petrov_1995,BhattacharyaRao_2010}.
Bounds on $\rho_n(\mathcal{A})$ under mild assumptions yield some useful statistics for a variety of high-dimensional inference problems. For instance, \cite{ChenQin_2010} and \cite{CaiLiuXia_2014} studied the $\ell^2$-type statistic and $\ell^\infty$-type statistic for testing high-dimensional mean vectors, respectively. Asymptotic validity of those test procedures relies on restrictive assumptions such as weak dependence in covariance matrix or sparsity in precision matrix. In contrast, Gaussian approximation results derived in this paper impose no explicit structural assumptions on the component-wise dependence structure, thus allowing to derive associated parametric bootstrap procedures for arbitrary dependence among different components of high-dimensional data. Recent years has witnessed a renewed interest in the accuracy of Gaussian approximation with explicit dependence on the dimension $p$ since such results are particularly useful in modern large-scale statistical inference problems such as change point detection \citep{YuChen_2020,YuChen_2019} and multiple testing for high-dimensional data \citep{ChangZhengZhouZhou_2017, ChangZhouZhouWang_2017}. For isotropic distributions with bounded third moments, \cite{Bentkus_2003} derived a Berry-Esseen type bound $O(p^{7/4} n^{-1/2})$ and $O(p^{3/2} n^{-1/2})$ over the class of convex subsets and Euclidean balls in $\R^p$, respectively. For independent (not necessarily identically distributed) sums, \cite{CCK_2013} considered the problem of approximating maxima of $S_{n,x}$ by its Gaussian analogue and established an error bound that allows the dimension $p$ to grow sub-exponentially fast in the sample size $n$.

Since the seminal work \cite{CCK_2013}, there have been substantial progresses being made in several directions. For instances, generalization of the index set from the max-rectangles to hyper-rectangles with improved rates of convergence can be found in \cite{CCK_2017}, \cite{LopesLinMueller_2020}, \cite{DengZhang_2020},   \cite{Deng_2020}, \cite{KuchibhotlaRinaldo_2020}, \cite{Koike_2019b}, \cite{FangKoike_2020}, \cite{DasLahiri_2020},  \cite{CCK_2020}, \cite{Lopes_2020} and \cite{CCKK_2019}; extension from linear sums to $U$-statistics with nonlinear kernels can be found in \cite{Chen_2018}, \cite{ChenKato_2019}, \cite{SongChenKato_2019}, and \cite{Koike_2019c}; generalization to dependent random vectors over max-rectangles can be found in \cite{ZhangWu_2017}, \cite{ZhangCheng_2018}, and \cite{CCK_2014}. 

In the literature, some popular assumptions imposed on the temporal dependence of the sequence $\mathcal{X}_n$ include: (i) strong-mixing (or $\alpha$-mixing) \citep{Rosenblatt_1956}, (ii) $m$-dependent sequence \citep{HoeffdingRobbins_1948}, and (iii) physical (or functional) dependence measure for casual time series \citep{Wu_2005}. Various CLTs for univariate (or fixed dimensional) dependent data have been developed under these dependence frameworks, see \cite{DoukhanMasssartRio_1994}, \cite{Bradley_2007}, \cite{Wu_2007}, and \cite{BerkesLiuWu_2014}. We remark that there are many other mixing coefficients measuring the temporal dependence of the past and future, among which the $\alpha$-mixing coefficient (see Definition \ref{defn:alpha-mixing-coef} in Section \ref{sec:alpha}) is the weakest one in the literature \citep{Bradley_2005}. In particular, for $p=1$, if the time series has finite third moment, the Komolgorov distance between the normalized random variable $\Xi^{-1/2}S_n$ and the standard univariate Gaussian distribution obeys a nearly optimal Berry-Esseen bound $O(n^{-1/2} \log^2n)$ with geometrically decaying $\alpha$-mixing coefficients \citep{Sunklodas_1984}, or the sharp Berry-Esseen bound $O(n^{-1/2})$ for either $m$-dependent sequence with fixed $m$ \citep{ChenShao_2004} or weakly dependent sequence under the physical dependence framework \citep{Hormann_2009,Jirak_2016}. Note that neither a dependence framework in (i)-(iii) implies the others. Thus there is a pressing call for a unified collection for Gaussian approximation tools under these temporal dependence frameworks for high dimensional dependent data.

Previous related works on high-dimensional CLTs for dependent data in the literature are complementary results for different dependence frameworks on max-rectangles, a subclass of hyper-rectangles. For examples, \cite{ZhangCheng_2018} studied the Gaussian approximation for $m$-dependent sequences with extension to dependent random vectors satisfying a geometric moment contraction condition \citep{WuShao_2004}; \cite{ZhangWu_2017} derived the Gaussian approximation result for causal stationary time series under a polynomial decay of the physical dependence measure; \cite{CCK_2014} studied the validity of a block multiplier bootstrap under the $\beta$-mixing condition. 
All the aforementioned papers are only applicable to approximating the distributions of the $\ell^\infty$-type statistics and not applicable to approximating the distributions of some more general and complicated statistics involved in high-dimensional statistical inference. See Section \ref{sec:apps} for details.

We conclude the introduction by summarizing our main contributions. Specifically, we develop a comprehensive and off-the-shelf probability toolbox containing the explicit rates of convergence of the high-dimensional CLTs for a combination of different index sets (including hyper-rectangles, simple convex sets, and sparsely convex sets) and different dependence frameworks (including $\alpha$-mixing, $m$-dependent, and physical dependence measure). Our error bounds are non-asymptotic in all key parameters, including the sample size $n$ and the data dimension $p$. In particular, our results established under the $\alpha$-mixing framework are new in the literature, while results established under the physical dependence measure improve over existing results. In addition, we provide a parametric bootstrap procedure to implement the proposed results with a kernel-type estimator for the long-run covariance matrix. For both Gaussian and parametric bootstrap approximations, the data dimension $p$ is allowed to grow sub-exponentially fast in the sample size $n$.
The rest of the paper is organized as follows. Section \ref{sec:main_results} presents the error bounds of $\rho_n(\mathcal{A})$ defined as \eqref{eq:gb} with selecting $\mathcal{A}$ as hyper-rectangles, simple convex sets, and sparsely convex sets, respectively. Section \ref{sec:parametric_boostrap} proposes a data-driven parametric bootstrap to approximate the probability $\mathbb{P}(S_{n,x}\in A)$ uniformly over $A\in\mathcal{A}$. Section \ref{sec:apps} discusses how to implement the proposed results in several statistical inference problems of interest. Section \ref{se:pftm1-3} includes the proofs of high-dimensional CLTs on hyper-rectangles presented in Section \ref{sec:hyper}, which provide the backbone for deriving high-dimensional CLTs on simple convex sets and sparsely convex sets stated, respectively, in Sections \ref{sec:simple} and \ref{sec:sparse}. The technical proofs of high-dimensional CLTs on simple convex sets and sparsely convex sets are given in the supplementary material.

\section{High-dimensional central limit theorems}
\label{sec:main_results}

We define some notation first. For any positive integer $m$, we write $[m]:=\{1,\ldots,m\}$. Denote by $I(\cdot)$ the indicator function. For two sequences of positive numbers $\{a_n\}$ and $\{b_n\}$, we write $a_n\lesssim b_n$ or $b_n\gtrsim a_n$ if there exists a universal constant $c>0$ such that $\limsup_{n\rightarrow\infty}a_n/b_n\leq c$. For any two $p$-dimensional vectors $v=(v_1,\ldots,v_p)^\T$ and $u=(u_1,\ldots,u_p)^\T$, $v \leq u$ means that $v_j \leq u_{j}$ for all $j\in[p]$. Given $\alpha>0$, we define the function $\psi_{\alpha}(x):=\exp(x^{\alpha})-1$ for any $x>0$. For a real-valued random variable $\xi$, we define $\|\xi\|_{\psi_{\alpha}}:=\inf[\lambda>0:\mathbb{E}\{\psi_{\alpha}(|\xi|/\lambda)\}\leq 1]$ and write $\xi\in\mathcal{L}^q$ for some $q>0$ if $\|\xi\|_q:=\{\mathbb{E}(|\xi|^q)\}^{1/q}<\infty$.  For a thricely differentiable function $f:\mathbb{R}^p\rightarrow\mathbb{R}$, we write $\partial_{j}f(x)=\partial f(x)/\partial x_{j}$, $\partial_{jk}f(x)=\partial^2f(x)/\partial x_{j}\partial x_{k}$ and $\partial_{jkl}f(x)=\partial^3 f(x)/\partial x_{j}\partial x_{k}\partial x_{l}$ for any $j,k,l\in [p]$. For a $q_1\times q_2$ matrix $B=(b_{i,j})_{q_1\times q_2}$, let $|B|_\infty=\max_{i\in[q_1],j\in[q_2]}|b_{i,j}|$ be the super-norm, and $\|B\|_2=\lambda_{\max}^{1/2}(BB^\T)$ be the spectral norm. Specifically, if $q_2=1$, we use $|B|_0=\sum_{i=1}^{q_1}I(b_{i,1}\neq0)$, $|B|_1=\sum_{i=1}^{q_1}|b_{i,1}|$ and $|B|_2=(\sum_{i=1}^{q_1}b_{i,1}^2)^{1/2}$ to denote the $\ell^0$-norm, $\ell^1$-norm and $\ell^2$-norm of the $q_1$-dimensional vector $B$, respectively.

Recall $S_{n,x}=n^{-1/2}\sum_{t=1}^nX_t$ and $\Xi={\rm Cov}(S_{n,x})$. Let $G\sim N(0,\Xi)$ which is independent of $\mathcal{X}_n=\{X_1,\ldots,X_n\}$. We will first consider in Section \ref{sec:hyper} the upper bounds for
\begin{align}\label{eq:varrhon}
\varrho_n:=\sup_{u\in\mathbb{R}^p,\nu\in[0,1]}|\mathbb{P}(\sqrt{\nu}S_{n,x}+\sqrt{1-\nu}G\leq u)-\mathbb{P}(G\leq u)|
\end{align}
when $\{X_t\}$ is (i) an $\alpha$-mixing sequence, (ii) an $m$-dependent sequence, and (iii) a physical dependence sequence, respectively. Based on such derived upper bounds, we can easily translate them to the upper bounds for $\rho_n(\mathcal{A})$ when $\mathcal{A}$ is selected as the class of all hyper-rectangles in $\mathbb{R}^p$. In Sections \ref{sec:simple} and \ref{sec:sparse}, we will consider the upper bounds for $\rho_n(\mathcal{A})$ when $\mathcal{A}$ is selected as the class of simple convex sets and $s$-sparsely convex sets, respectively.
  Write $X_t=(X_{t,1},\ldots,X_{t,p})^\T$. Throughout the rest of this paper (unless otherwise explicitly stated), we shall focus on the high-dimensional scenario by assuming that $p\geq n^{\kappa}$ for some universal constant $\kappa>0$. Here $\kappa>0$ can be selected as some  sufficiently small constant. Assuming $p\geq n^\kappa$ is a quite mild condition in the literature of high-dimensional data analysis which is not necessary for our theoretical analysis and just used to simplify our presentation. In our theoretical proofs, we need to compare $\log p$ and $\log n$ in lots of places. Without the restriction $p\geq n^\kappa$, some $\log p$ terms in the theoretical results should be replaced by $\log(pn)$.  

\subsection{High-dimensional CLT for hyper-rectangles}\label{sec:hyper}
Let $\mathcal{A}^{{\rm re}}$ be the class of all hyper-rectangles in $\mathbb{R}^p$; that is, $\mathcal{A}^{{\rm re}}$ consists of all sets $A$ of the form   $A=\{(w_1,\ldots,w_p)^\T\in\mathbb{R}^p:a_j\leq w_j\leq b_j\textrm{ for all }j\in[p]\}$ with some $-\infty\leq a_j\leq b_j\leq \infty$. Define $S_{n,\check{x}}=n^{-1/2}\sum_{t=1}^n\check{X}_t$ with $\check{X}_t=(X_t^\T,-X_t^\T)^\T$ and let $\check{G}\sim N(0,\check{\Xi})$ with $\check{\Xi}={\rm Cov}(n^{-1/2}\sum_{t=1}^n\check{X}_t)$. We then have
$$
\rho_n(\mathcal{A}^{\rm re})\leq \sup_{u\in\mathbb{R}^{2p},\nu\in[0,1]}|\mathbb{P}(\sqrt{\nu}S_{n,\check{x}}+\sqrt{1-\nu}\check{G}\leq u)-\mathbb{P}(\check{G}\leq u)|\,,
$$
where the term on the right-hand side is a $(2p)$-dimensional analogue of $\varrho_n$ defined as \eqref{eq:varrhon} over one-sided hyper-rectangles. To derive the convergence rate of $\rho_n(\mathcal{A}^{{\rm re}})$, it suffices to consider that for $\varrho_n$.

\subsubsection{$\alpha$-mixing sequence}\label{sec:alpha}

\begin{definition}[$\alpha$-mixing coefficient]
\label{defn:alpha-mixing-coef}
Let $\{X_t\}$ be a random sequence. Denote by $\mathcal{F}_{-\infty}^u$ and $\mathcal{F}_{u}^\infty$ the $\sigma$-fields generated respectively by $\{X_t\}_{t\leq u}$ and $\{X_t\}_{t\geq u}$. The {\it $\alpha$-mixing coefficient} at lag $k$ of the sequence $\{X_t\}$ is defined as
$$
\alpha_n(k):=\sup_t\sup_{A\in\mathcal{F}_{-\infty}^t, B\in\mathcal{F}_{t+k}^\infty}|\mathbb{P}(AB)-\mathbb{P}(A)\mathbb{P}(B)|\,.
$$
We say the sequence $\{X_t\}$ is {\it $\alpha$-mixing} if $\alpha_n(k) \to 0$ as $k\rightarrow\infty$.
\end{definition}

The {\it long-run variance} of the $j$-th coordinate marginal sequence $\{X_{t,j}\}_{t=1}^n$ is defined as
\begin{align}\label{eq:Vtj}
V_{n,j}={\rm Var}\bigg(\frac{1}{\sqrt{n}}\sum_{t=1}^{n}X_{t,j}\bigg)\,.
\end{align}
To investigate the convergence rate of $\varrho_n$ defined as \eqref{eq:varrhon} for the $\alpha$-mixing sequence $\{X_t\}$, we need the following regularity conditions.

\begin{ass}[Sub-exponential moment]
\label{as:tail}
There exist a sequence of constants $B_{n}\geq 1$ and a universal constant $\gamma_{1}\geq 1$ such that $\|X_{t,j}\|_{\psi_{\gamma_{1}}}\leq B_{n}$ for all $t\in [n]$ and $j\in [p]$.
\end{ass}

\begin{ass}[Decay of $\alpha$-mixing coefficients]
\label{as:alpha-mixing}
There exist some universal constants $K_1>1$, $K_2>0$ and $\gamma_2>0$ such that
$
\alpha_n(k)\leq K_1\exp(-K_2k^{\gamma_2})$ for any $k\geq 1$.
\end{ass}

\begin{ass}[Non-degeneracy]
\label{as:longrun}
There exists a universal constant $K_3>0$ such that $\min_{j\in[p]}V_{n,j}\geq K_3$.
\end{ass}

Since Condition \ref{as:tail} implies that
$\mathbb{E}\{\exp(|X_{t,j}|^{\gamma_{1}}B_{n}^{-\gamma_{1}})\}\leq 2$, it follows from Markov's inequality that $\mathbb{P}(|X_{t,j}|>u) \leq 2 \exp(-u^{\gamma_{1}}B_{n}^{-\gamma_{1}})$ for all $u > 0$. If each $X_{t,j}$ is sub-gaussian, we have $\gamma_1=2$ and $B_n=O(1)$.  On one hand, the sub-exponential moment condition is a widely used condition in the high-dimensional statistics literature as this would generally entail vanishing rates of convergence of the sample mean for mean-zero independent data when the dimension $p$ scales sub-exponentially fast in the sample size $n$ \citep{wainwright2019high}. On the other hand, it can be relaxed to polynomial moment condition (e.g., 3rd moment condition). The trade-off is that under such condition, we can only expect the dimension $p$ scales polynomially fast in the sample size $n$ (i.e., $p = O(n^c)$ for some constant $c > 0$) to obtain vanishing rates of the sample mean for independent data. We can certainly expect that similar rates can be established for temporally dependent data. However, the proof techniques would be similar to the sub-exponential moment case. Moreover, if we replace the sample mean by its self-normalized version, then in the independent data case the studentized mean still has exponential decay tail under 3rd moment condition. It would be an interesting future work to investigate the self-normalization in the high-dimensional time series setting for the Gaussian approximation.

The $\alpha$-mixing assumption is mild in the literature. Causal ARMA processes with continuous innovation distributions are $\alpha$-mixing with exponential decay rates. So are stationary Markov chains satisfying certain conditions. See Section 2.6.1 of \cite{FanYao_2003} and references within. In fact stationary GARCH models with finite second moments and continuous innovation distributions are also $\alpha$-mixing with exponential decay rates. Under certain conditions, VAR processes, multivariate ARCH processes, and multivariate GARCH processes are all $\alpha$-mixing with exponential decay rates;  see \cite{HP09}, \cite{BFS11} and \cite{Wong2020}. The next two examples also satisfy Condition \ref{as:alpha-mixing}.
\begin{itemize}
\item Let $X_t=A_tf_t+\varepsilon_t$, where $A_t$ is a nonrandom loading matrix, $f_t\in\mathbb{R}^r$ is the latent factor with some fixed integer $r$, and $\{\varepsilon_t\}$ is an independent sequence that is also independent of $\{f_t\}$. If $\{f_t\}$ is selected as VAR processes, multivariate ARCH processes, or multivariate GARCH processes, due to that $r$ is fixed and $\{\varepsilon_t\}$ is an independent sequence, we know such defined $\{X_t\}$ satisfies Condition \ref{as:alpha-mixing} under certain conditions imposed on the model of $\{f_t\}$. 

\item Let $X_t=g_t(U_t)$, where $\{U_t\}$ is a $q$-dimensional latent sequence, and $g_t(\cdot):\mathbb{R}^{q}\rightarrow\mathbb{R}^p$ is a Borel function. Here we do not impose any relationship between $p$ and $q$, and allow $q=\infty$. Write $U_t=(U_{t,1},\ldots,U_{t,q})^{\T}$. Assume the sequence $\{U_{t,j}\}$ is $\rho$-mixing with exponential decay rates for each $j\in[q]$. If $\{U_{t,1}\},\ldots,\{U_{t,q}\}$ are independent of each other, Theorem 5.1 of \cite{Bradley_2005} implies $\{U_t\}$ is $\rho$-mixing with exponential decay rates. Due to the relationship between $\rho$-mixing and $\alpha$-mixing, we know such defined $\{X_t\}$ satisfies Condition \ref{as:alpha-mixing}.  
\end{itemize}

Condition \ref{as:longrun} assumes the partial sum $n^{-1/2}\sum_{t=1}^{n}X_{t,j}$ is non-degenerated which is required when we apply Nazarov's inequality \cite[Lemma A.1]{CCK_2017} 
to bound the probability of a Gaussian vector taking values in a small region. When $\{X_{t,j}\}_{t\geq1}$ is stationary, we know $V_{n,j}=\Gamma_j(0)+2\sum_{k=1}^{n-1}(1-kn^{-1})\Gamma_j(k)$, where $\Gamma_j(k)=\cov(X_{1,j},X_{1+k,j})$ is the autocovariance of $\{X_{t,j}\}_{t\geq1}$ at lag $k$. If each component sequence $\{X_{t,j}\}$ is stationary,  Condition \ref{as:longrun} holds if $\Gamma_j(0)+2\sum_{k=1}^\infty \Gamma_j(k)\geq C$ holds for any $j\in[p]$, where $C>0$ is a universal constant. Based on Conditions \ref{as:tail}--\ref{as:longrun}, Theorem \ref{tm:1} gives an upper bound for $\varrho_n$ when the underlying sequence $\{X_t\}$ is $\alpha$-mixing, whose proof is given in Section \ref{sec:pfalphamix}.

\begin{theorem}
[Gaussian approximation for partial sums of the $\alpha$-mixing sequence]
\label{tm:1}
Assume $\{X_t\}$ is an $\alpha$-mixing sequence with $p\geq n^{\kappa}$ for some universal constant $\kappa>0$. Under Conditions {\rm\ref{as:tail}}--{\rm\ref{as:longrun}}, it holds that 
$$
\varrho_n\lesssim \frac{B_{n}^{2/3}(\log p)^{(1+2\gamma_2)/(3\gamma_2)}}{n^{1/9}}+\frac{B_{n}(\log p)^{7/6}}{n^{1/9}}$$
provided that $(\log p)^{3-\gamma_2}=o(n^{\gamma_2/3})$.
\end{theorem}

\begin{remark}[Comparison with existing results under mixing dependence measure]
Appendix B of \cite{CCK_2014} derived the validity of a block multiplier bootstrap (BMB) under the $\beta$-mixing assumption. There are several differences between our Theorem 1 and Theorem B.1 in \cite{CCK_2014}. First, since the $\beta$-mixing assumption implies the $\alpha$-mixing assumption, our Theorem 1 is applicable for wider class of dependent data. Second, Theorem B.1 in \cite{CCK_2014} is proved and stated with the ``large-and-small-blocks" argument, where conditions of their Theorem B.1 involve the ``tuning parameter" of the block sizes. It is empirically known that the performance of BMB is sensitive to the block sizes. Although the BMB procedure given in \cite{CCK_2014} is theoretically valid with suitable divergence rates imposed on the block sizes, how to propose a valid data-driven procedure to select the two involved tuning parameters is unclear in the framework of Gaussian approximation. Thus, it is an undesirable feature of Theorem B.1 in \cite{CCK_2014} to rule out bootstraps without a hard truncation block size to estimate the long-run covariance matrices. In Section \ref{sec:parametric_boostrap}, we consider the kernel-type estimator of \cite{Andrews_1991} to estimate the long-run covariance matrix of $S_{n,x}$, which is more appealing from a practical standview (e.g., with the optimal quadratic spectral kernel and optimal data-driven bandwidth formula). Although the optimal data-driven bandwidth is obtained in the fixed dimensional scenario, extensive numerical studies in \cite{ChangJiangShao_2021} indicate that such formula still works well in high-dimensional setting and the associated performance is quite robust when the bandwidth is selected in a large range. Third, result from \cite{CCK_2014} holds only for max-norm statistics, while our paper derives the convergence rates of Gaussian and parametric bootstrap approximations under much broader classes of index sets for high-dimensional dependent data (see Sections \ref{sec:simple} and \ref{sec:sparse} below) that can be applied to approximate the distributions of more general and complicated statistics used in high-dimensional statistical inference. 
\end{remark}

\begin{remark}
Theorem \ref{tm:1} extends the Gaussian approximation result for independent data in \cite{CCK_2017} to dependent data. When the eigenvalues of $\Xi$ are bounded below from zero (i.e., strongly non-degenerate case), \cite{CCK_2020} derived a nearly optimal rate of convergence for independent data. Our analysis can be adapted with the sharper results from \cite{CCK_2020} to yield an improved error bound in Theorem \ref{tm:1} under stronger conditions.
\end{remark}

\subsubsection{$m$-dependent sequence}
\label{subsec:m-dependent}

Based on the temporal dependency among $\{X_t\}_{t=1}^n$, we can define an undirected graph $G_{n} = (V_{n}, E_{n})$, where $V_{n} = [n]$ is a set of nodes with node $t$ denoting $X_t$, and $E_{n}$ is a set of undirected edges connecting the nodes such that $X_{t}$ and $X_{s}$ are independent whenever $(t,s) \notin E_{n}$. Here we adopt the convention $(t,t)\in E_n$ for any $t\in V_n$. We call such defined $G_n$ is the {\it dependency graph} of the sequence $\{X_t\}_{t=1}^n$. The dependency graph is a flexible model to study CLTs with increasing dependence strength \citep{BaldiRinott_1989} that covers the $m$-dependent sequence as a special case. For any $t\in[n]$, let $\mathcal{N}_{t} = \{s\in V_n: (t,s) \in E_{n}\}$ be the neighbor nodes of node $t$ in $G_n$. Define $D_{n} = \max_{t \in[n]} \sum_{s=1}^{n} I\{(t,s) \in E_{n}\}$ as the maximum degree of the first-degree connections in $G_n$, and $D_{n} ^*= \max_{t \in[n]} \sum_{s=1}^{n} I\{s\in\cup_{\ell\in\mathcal{N}_t}\mathcal{N}_\ell\}$ as the maximum degree of the second-degree connections in $G_n$. Theorem \ref{lem:GA_dependency_graph} gives an upper bound for $\varrho_n$ defined as \eqref{eq:varrhon} based on the maximum degrees $D_n$ and $D^*_n$ of the dependency graph determined by the underlying sequence $\{X_t\}_{t=1}^n$, whose proof is given in Section \ref{sec:5.2}.

\begin{theorem}
[Gaussian approximation for partial sums of a sequence under dependency graph]
\label{lem:GA_dependency_graph}
Assume $p\geq n^\kappa$ for some universal constant $\kappa>0$. Under Conditions {\rm \ref{as:tail}} and {\rm\ref{as:longrun}}, it holds that
$$
\varrho_n \lesssim \frac{B_{n}(D_nD_n^*)^{1/3}(\log p)^{7/6}}{n^{1/6}}\,,
$$
where $D_n$ and $D^*_n$ are the maximum degrees of the first-degree and second-degree of connections in the dependency graph generated by the sequence $\{X_t\}_{t=1}^n$, respectively.
\end{theorem}

If $\{X_{t}\}_{t=1}^n$ is a centered $m$-dependent sequence, i.e., $X_{t}$ and $X_{s}$ are independent for all $|t-s|>m$, then $\{X_t\}_{t=1}^n$ has a dependency graph with $D_{n} = 2m+1$ and $D_n^*=4m+1$. The next corollary states a result for $m$-dependent sequences.

\begin{cy}
[Gaussian approximation for partial sums of an $m$-dependent sequence]
\label{cor:GA_m-dependent_subexp}
Assume $\{X_t\}_{t=1}^n$ is an $m$-dependent sequence with $p\geq n^\kappa$ for some universal constant $\kappa>0$. Under Conditions {\rm \ref{as:tail}} and {\rm\ref{as:longrun}}, it holds that 
$$
\varrho_n \lesssim \frac{B_{n}(m\vee 1)^{2/3}(\log p)^{7/6}}{n^{1/6}}\,.$$
\end{cy}

Since the $0$-dependent sequence reduces to the independent sequence, Corollary \ref{cor:GA_m-dependent_subexp} for $m = 0$ reads $O(B_n n^{-1/6}\log^{7/6}p)$, which has the same sample complexity in $n$ and dimension dependence in $p$ as the independent data case $O(B_n^{1/3}n^{-1/6} \log^{7/6}p)$ up to a moment factor $B_n^{2/3}$ (cf. Proposition 2.1 in \cite{CCK_2017}). The extra cost $B_n^{2/3}$ is due to the argument that we need to decouple the distribution tail and dependence simultaneously. In particular, for the data with $B_n = O(1)$, the rate obtained from our $m$-dependent CLT achieves the CLT rate for independent data derived in \cite{CCK_2017}.

Corollary \ref{cor:GA_m-dependent_subexp} is a stepping stone to study the Gaussian approximation under the physical dependence framework with better rate of convergence than the best known results in \cite{ZhangWu_2017} based on the large-and-small-blocks technique in the weaker temporal dependence regime. See Theorem \ref{prop:GA_weakly-dependent} and the discussions in Section \ref{sec:comparison_physical_dependence} for more details.

\subsubsection{Sequence with physical dependence}\label{sec:physical}

Let $\{\varepsilon_{i}\}_{i \in \mathbb{Z}}$ be a sequence of independent and identically distributed random elements. Consider the (causal) time series model
\begin{equation}
\label{eqn:weakly_dependent_ts_fdp}
X_{t} =f_{t}(\varepsilon_{t},\varepsilon_{t-1},\dots)\,, \quad t\geq1\,,
\end{equation}
where $f_{t}(\cdot)$ is a jointly measurable function taking values in $\R^{p}$ and $\mathbb{E}(X_{t})= 0$. Here $\{\varepsilon_{i}\}_{i \in \Z}$ are {\it innovations} that can be viewed as the input of the non-linear system (\ref{eqn:weakly_dependent_ts_fdp}). Since the data generation mechanism $f_{t}(\cdot)$ may change over time, $X_{t}$ is allowed to be {\it non-stationary}. Non-linear time series of the form (\ref{eqn:weakly_dependent_ts_fdp}) are first introduced in \cite{Wu_2005} for $p=1$ and $f_{t}(\cdot)\equiv f(\cdot)$ for some measurable function $f(\cdot)$ (i.e., stationary univariate time series), and their temporal dependence can be quantified by the {\it functional dependence measure} based on the idea of coupling. In particular, let $\varepsilon'_{i}$ be an independent copy of $\varepsilon_{i}$ and
\begin{align*}
X'_{t,\{m\}} =f_{t}(\varepsilon_{t},\dots,\varepsilon_{t-m+1},\varepsilon'_{t-m},\varepsilon_{t-m-1},\dots)
\end{align*}
be the coupled version of $X_{t}$ at the time lag $m$ with $\varepsilon_{t-m}$ replaced by $\varepsilon'_{t-m}$. By causality, $X'_{t,\{m\}} = X_{t}$ for $m < 0$. Write $X'_{t,\{m\}}=(X'_{t,1,\{m\}},\ldots,X'_{t,p,\{m\}})^\T$. The (uniform) functional dependence measure for the $j$-th coordinate marginal sequence $\{X_{t,j}\}$ is defined as
\begin{align*}
\theta_{m,q,j} = \sup_{t \geq 1} \|X_{t,j} - X'_{t,j,\{m\}}\|_{q}\,, \quad q > 0\,.
\end{align*}
In essence, $\theta_{m,q,j}$ quantifies the uniform impact of coupling on the $j$-th coordinate marginal time series at lag $m$. For any $m\geq0$, write $\Theta_{m,q,j} = \sum_{i=m}^{\infty} \theta_{i,q,j}$. For $\alpha \in (0,\infty)$, define the {\it dependence adjusted norm} introduced in \cite{WuWu_2016} as
\begin{align*}
\|X_{.,j}\|_{q,\alpha} = \sup_{m  \geq 0} (m+1)^{\alpha} \Theta_{m,q,j} \quad \text{and} \quad \|X_{.,j}\|_{\psi_{\nu},\alpha} = \sup_{q \geq 2} q^{-\nu} \|X_{.,j}\|_{q,\alpha}\,,
\end{align*}
whenever the supremums are finite. Define further the aggregated norms as follows:
\begin{align}\label{eq:aggnorm}
\Psi_{q,\alpha} = \max_{j \in[p]} \|X_{.,j}\|_{q,\alpha} \quad \text{and} \quad \Phi_{\psi_{\nu},\alpha} = \max_{j \in[p]} \|X_{.,j}\|_{\psi_{\nu},\alpha}\,.
\end{align}

\begin{theorem}
[Gaussian approximation for maxima of partial sums of time series under functional dependence]
\label{prop:GA_weakly-dependent}
Assume $\{X_t\}$ satisfies the model \eqref{eqn:weakly_dependent_ts_fdp} with $p\geq n^\kappa$ for some universal constant $\kappa>0$. Let $\Phi_{\psi_{\nu},\alpha} < \infty$ for some $\alpha, \nu \in (0,\infty)$. 

{\rm(i)} Under Condition {\rm\ref{as:longrun}}, it holds that 
$$
\varrho_n
\lesssim \frac{\Phi_{\psi_\nu,0}(\log p)^{7/6}}{n^{\alpha/(3+9\alpha)}}+\frac{\Psi_{2,\alpha}^{1/3} \Psi_{2,0}^{1/3}(\log p)^{2/3}}{n^{\alpha/(3+9\alpha)}}+\frac{\Phi_{\psi_\nu,\alpha}(\log p)^{1+\nu}}{n^{\alpha/(1+3\alpha)}}$$
provided that $(\log p)^{\max\{6\nu-1,(5+6\nu)/4\}}= o\{n^{\alpha/(1+3\alpha)}\}$. 

{\rm(ii)} Under Conditions {\rm \ref{as:tail}} and {\rm\ref{as:longrun}},  it holds that
$$
\varrho_n\lesssim \frac{B_{n}(\log p)^{7/6}}{n^{\alpha/(12+6\alpha)}}+\frac{\Psi_{2,\alpha}^{1/3} \Psi_{2,0}^{1/3}(\log p)^{2/3}}{n^{\alpha/(12+6\alpha)}}+\frac{\Phi_{\psi_{\nu},\alpha}(\log p)^{1+\nu}}{n^{\alpha/(4+2\alpha)}}\,.$$
\end{theorem}
The proof of Theorem \ref{prop:GA_weakly-dependent} is given in Section \ref{sec:pflma5}. We remark that the two rates of convergence given in Theorems \ref{prop:GA_weakly-dependent}(i) and \ref{prop:GA_weakly-dependent}(ii) are based on the {\it large-and-small-blocks} and $m$-dependent approximation techniques, respectively. The large-and-small-blocks technique is widely used in time series analysis to approximate the sum of a time series sequence by the sum over its large blocks. It is interesting to note that the large-and-small-blocks technique gives a faster (or slower) rate than the $m$-dependent argument when $0 < \alpha < 3$ (or $\alpha > 3$). In particular, when the temporal dependence is weak (for large values of $\alpha$), the improvement of Theorem \ref{prop:GA_weakly-dependent}(ii) than Theorem \ref{prop:GA_weakly-dependent}(i) is more significant. The intuition is that the large-and-small-blocks technique used to establish Theorem \ref{prop:GA_weakly-dependent}(i) may lose sample size efficiency when the temporal dependence is weak. In such regime, throwing away the data in small blocks may reduce the {\it effective} sample size, while the $m$-dependent approximation directly approximates the sequence $\mathcal{X}_n$ without throwing away data. On the other hand, when the temporal dependence is strong (for small values of $\alpha$), we need to use much larger values of $m$ for constructing an $m$-dependent sequence in Section \ref{subsec:m-dependent}, so the $m$-dependent approximation becomes less effective than throwing a reasonable amount of small blocks to reduce the dependence.

Now combining the two parts of Theorem \ref{prop:GA_weakly-dependent}, we obtain the overall rate of convergence under the physical dependence measure.

\begin{cy}[Overall rate of convergence under physical dependence]
\label{cor:eqn:GA_weakly-dependent_combined}
Assume $\{X_t\}$ satisfies the model \eqref{eqn:weakly_dependent_ts_fdp} with $p\geq n^\kappa$ for some universal constant $\kappa>0$, and $\Phi_{\psi_{\nu},\alpha} < \infty$ for some $\alpha, \nu \in (0,\infty)$. Let $\alpha' = \alpha / \min\{ 1+3\alpha, 4+2\alpha \}$. Under Conditions {\rm \ref{as:tail}} and {\rm\ref{as:longrun}}, it holds that 
$$
\varrho_n
\lesssim \frac{\max\{\Phi_{\psi_\nu,0}, B_n \} (\log p)^{7/6}}{n^{\alpha'/3}}+\frac{\Psi_{2,\alpha}^{1/3} \Psi_{2,0}^{1/3}(\log p)^{2/3}}{n^{\alpha'/3}}+\frac{\Phi_{\psi_\nu,\alpha}(\log p)^{1+\nu}}{n^{\alpha'}}$$
provided that $(\log p)^{\max\{6\nu-1,(5+6\nu)/4\}}= o\{n^{\alpha/(1+3\alpha)}\}$.
\end{cy}

\subsubsection{Comparison with existing result under physical dependence measure}\label{sec:comparison_physical_dependence}
Under the physical dependence and a sub-exponential moment condition, \cite{ZhangWu_2017} derived a Gaussian approximation result for the $\ell^\infty$-norm of normalized sums of a class of stationary time series:
\begin{equation*}
\omega_n := \sup_{u\geq 0}| \bP(|D^{-1} S_{n,x}|_\infty \geq u) - \bP(|D^{-1} G|_\infty \geq u) |\,,
\end{equation*}
where $D= \{\diag(\Xi)\}^{1/2}$ and $G \sim N(0,\Xi)$. Specifically, Theorem 7.4 in \cite{ZhangWu_2017} gives the following error bound: for any $\lambda \in (0,1)$ and $\eta > 0$,
\begin{equation}
\label{eqn:zhang_wu_rate_phystical_dependence}
\omega_n \lesssim f^{\diamond}(\sqrt{n} \eta) + \eta (\log{p})^{1/2} + h\{\lambda, u_m^{\diamond}(\lambda)\} + \pi\{\chi(m,M)\}
\end{equation}
with
\begin{align*}
f^{\diamond}(y) =&\,\, p \exp(-C_\beta y^\beta m^{\alpha\beta}n^{-\beta/2}\Phi_{\psi_\nu,\alpha}^{-\beta}) + p \exp\{ -C_\beta y^\beta (mw)^{-\beta/2} \Phi_{\psi_\nu,0}^{-\beta}\}\,,\\
h\{\lambda, u_m^{\diamond}(\lambda)\} =&\,\, \lambda+w^{-1/8}\max\{\Psi_{3,0}^{3/4},\Psi_{4,0}^{1/2}\}\log^{7/8}(pw\lambda^{-1})\\
&+ w^{-1/2}\max\{  \Phi_{\psi_\nu,0}\log^{1/\beta}(pw\lambda^{-1}), \log^{1/2}(pw\lambda^{-1})\}\log^{3/2}(pw\lambda^{-1})\,,\\
\pi(x) = &\,\,x^{1/3} \max\{1, \, \log^{2/3}(px^{-1})\}\,,\\
\chi(m,M) =&\,\, \Psi_{2,\alpha} \Psi_{2,0}\{m^{-\alpha}+v(M)\} + wmn^{-1}\,,
\end{align*}  where $\beta = 2/(1+2\nu)$, $v(M) = M^{-1}I(\alpha > 1)+(M^{-1} \log{M})I(\alpha = 1)+M^{-\alpha}I(0 < \alpha < 1)$, and $(m, M, w)$ are tuning parameters involved in the ``large and small blocks" technique for deriving \eqref{eqn:zhang_wu_rate_phystical_dependence} with $m$ and $M$ being, respectively, the sizes of small blocks and large blocks satisfying $m=o(M)$, and $w=\lfloor n/(M+m)\rfloor$. They first approximated $S_{n,x}$ by the sum of an $m$-dependent sequence, and then applied the large-and-small-blocks technique to approximate the sum of the $m$-dependent sequence by the sum over large blocks.

To simplify the convergence rate of $\omega_n$ specified in \eqref{eqn:zhang_wu_rate_phystical_dependence}, we assume $\Phi_{\psi_\nu,\alpha} = O(1)$ for some $\alpha,\nu \in(0,\infty)$ and $p\geq n^\kappa$ for some $\kappa>0$. Choose $\lambda = n^{-c_1}, \eta = n^{-c_2}$, $w \asymp n^{c_3}$ and $m\asymp n^{c_4}$ for some constants $c_1,c_2,c_3, c_4>0$. By optimizing $(c_1,c_2,c_3,c_4)$ according to the right-hand side of \eqref{eqn:zhang_wu_rate_phystical_dependence}, we have the following proposition whose proof is given in Section \ref{sec:xic1c4} of the supplementary material.
\begin{proposition}[Rate of convergence under physical dependence in \cite{ZhangWu_2017}]
\label{pn:zhangwu2017}
Assume $\Phi_{\psi_\nu,\alpha} = O(1)$ for some $\alpha,\nu \in(0,\infty)$ and $p\geq n^\kappa$ for some $\kappa>0$. Then the upper bound of $\omega_n$ given in \eqref{eqn:zhang_wu_rate_phystical_dependence} can be simplified as
$$
\omega_n \lesssim \frac{{\rm polylog}(p)}{ n^{\alpha/(3+11\alpha)}}\,,$$
where ${\rm polylog}(p)$ is a polynomial factor of $\log p$.
\end{proposition}
 Contrasting Proposition \ref{pn:zhangwu2017} with Corollary \ref{cor:eqn:GA_weakly-dependent_combined} specialized to max-rectangles, we see that, up to a $\mbox{polylog}(p)$ factor, our rate of convergence reads ${\rm polylog}(p) \cdot n^{-\alpha/[3 \min\{ 1+3\alpha, 4+2\alpha \}]}$, which is uniformly faster than that given in Proposition \ref{pn:zhangwu2017} for all $\alpha > 0$. In other words, our rate has a better sample size dependence than  \cite{ZhangWu_2017}. The reason can be seen that the optimal choice of \cite{ZhangWu_2017} throws away $w \asymp n^{8\alpha/(3+11\alpha)}$ small blocks of size $m\asymp n^{3/(3+11\alpha)}$, which leads a total reduction of $O\{n^{(3+8\alpha)/(3+11\alpha)}\}$ data points in the sample size. In our result, we only throw away $w \asymp n^{2\alpha/(1+3\alpha)}$ small blocks of size $m \asymp n^{1/(1+3\alpha)}$, leading to a total reduction of $O\{n^{(1+2\alpha)/(1+3\alpha)}\}$ data points in the sample size. Moreover, the improvement of our result over \cite{ZhangWu_2017} is more significant for larger values of $\alpha > 3$.

\subsection{High-dimensional CLT for simple convex sets}\label{sec:simple}

In this section, we consider the class of {\it simple convex sets} introduced by \cite{CCK_2017}. Formally, a simple convex set can be well approximated by a convex polytope with a controlled number of facets. Simple convex sets serve an important intermediate step to derive similar error bounds in Gaussian approximation on the class of $s$-sparsely convex sets considered in Section \ref{sec:sparse}. Geometrically, $s$-sparsely convex sets can be represented as an intersection of possibly many convex sets whose indicator functions depend at most on $s$ elements of their coordinates.

For a closed convex set $A\subset\mathbb{R}^p$, we define its support function:
\begin{align*}
\mathcal{S}_{A}:\mathbb{S}^{p-1}\mapsto\mathbb{R}\cup\{\infty\}\,,~~v\mapsto\mathcal{S}_{A}(v):=\sup\{w^{\T}v:w\in A\}\,,
\end{align*}
where $\mathbb{S}^{p-1}$ is the unit sphere in $\R^p$. Specially, if $A^{K}$ is \emph{$K$-generated} (that is, $A^{K}$ is generated by the intersection of $K$ half-spaces), we could characterize $A^{K}$ by its support function for the set $\mathcal{V}(A^K)$ consisting $K$ unit normal vectors outward to the facets of $A^K$:
\begin{align*}
A^{K}=\bigcap_{v\in\mathcal{V}(A^{K})}\{w\in\mathbb{R}^{p}:w^{\T}v\leq\mathcal{S}_{A^{K}}(v)\}\,.
\end{align*}
Moreover, for $\epsilon>0$ and a $K$-generated convex set $A^K$, we also define
\begin{align*}
A^{K,\epsilon}=\bigcap_{v\in\mathcal{V}(A^{K})}\{w\in\mathbb{R}^{p}:w^{\T}v\leq\mathcal{S}_{A^{K}}(v)+\epsilon\}\,.
\end{align*}

\begin{definition}[Simple convex set]
\label{defn:simple_convex_set}
We say $A$ is a \textit{simple convex set}, if there exist two constants $a\geq0$, $d>0$ and an $K$-$generated$ $A^{K}$ satisfying $K\leq (pn)^d$ such that
\begin{align}\label{eq:inclusion}		
A^{K}\subset A\subset A^{K,\epsilon}
\end{align}
with $\epsilon=a/n$. In this case, $A^K$ provides an approximation to $A$ with precision $\epsilon$.
\end{definition}

Let $\mathcal{A}^{\rm si}(a,d)$ be the class of all sets $A$ satisfying \eqref{eq:inclusion} with $K\leq(pn)^d$ and $\epsilon=a/n$. For any $v\in \mathbb{R}^{p}$, define
\begin{align}\label{eq:Vtv}
V_{n}(v)={\rm Var}\bigg(\frac{1}{\sqrt{n}}\sum_{t=1}^{n}v^{\T}X_{t}\bigg)\,.
\end{align}
In the sequel, we shall slightly abuse the notation by using $A^K(\cdot)$ to also denote the operator defined on $\mathcal{A}^{\rm si}(a,d)$ such that $A^K(A)=A^K$ with $A^K$ specified in \eqref{eq:inclusion} for any $A\in\mathcal{A}^{\rm si}(a,d)$. To construct the upper bounds for $\rho_n(\mathcal{A})$ for some $\mathcal{A}\subset\mathcal{A}^{\rm si}(a,d)$, we need the following condition that imposes the moment assumption on $v^\T X_t$ for $v\in\mathcal{V}\{A^K(A)\}$.

\begin{ass}\label{as:ExtendedTail}
	There exist a sequence of constants $B_{n}\geq 1$ and a universal constant $\gamma_{1}\geq 1$ such that $\|v^{\T}X_{t}\|_{\psi_{\gamma_{1}}}\leq B_{n}$ for all $t\in [n]$ and $v\in\mathcal{V}\{A^{K}(A)\}$.
\end{ass}

\subsubsection{$\alpha$-mixing sequence} To obtain an upper bound for $\rho_n(\mathcal{A})$ for some $\mathcal{A}\subset\mathcal{A}^{\rm si}(a,d)$, we need the next condition that requires the long-run variance of the sequence $\{v^\T X_{t}\}_{t=1}^n$ is not degenerated for any $v\in\mathcal{V}\{A^K(A)\}$.

\begin{ass}\label{as:longrun_simple_convex}
	There exists a universal constant $K_{4}>0$ such that $V_{n}(v)\geq K_{4}$ for any $v\in\mathcal{V}\{A^{K}(A)\}$.
\end{ass}

Condition \ref{as:longrun_simple_convex} holds automatically if the smallest eigenvalue of $\Xi={\rm Cov}(n^{-1/2}\sum_{t=1}^nX_t)$ is uniformly bounded away from zero.

\begin{theorem}
[Gaussian approximation for partial sums of the $\alpha$-mixing sequence for simple convex sets]\label{pn:GA_simple_convex}	Assume $\{X_t\}$ is an $\alpha$-mixing sequence with $p\geq n^{\kappa}$ for some universal constant $\kappa>0$ and Condition {\rm\ref{as:alpha-mixing}} being satisfied. Let $\mathcal{A}$ be a subclass of $\mathcal{A}^{\rm si}(a,d)$ such that Conditions {\rm\ref{as:ExtendedTail}} and {\rm\ref{as:longrun_simple_convex}} are satisfied for any $A\in\mathcal{A}$. Then
	$$
		\rho_n(\mathcal{A})\lesssim \frac{a(d\log p)^{1/2}}{n} +\frac{B_{n}^{2/3}(d\log p)^{(1+2\gamma_2)/(3\gamma_2)}}{n^{1/9}}+\frac{B_{n}(d\log p)^{7/6}}{n^{1/9}}$$
	provided that $(d\log p)^{3-\gamma_2}=o(n^{\gamma_2/3})$.
\end{theorem}
The proof of Theorem \ref{pn:GA_simple_convex} is given in Section \ref{sec:pf_simple_convex} of the supplementary material.

\subsubsection{$m$-dependent sequence} Theorem \ref{pn:GA_dependency_graph_simple_convex} gives the high-dimensional CLT for simple convex sets based on the maximum degrees $D_n$ and $D_n^*$ of the dependency graph determined by the underlying sequence $\{X_t\}_{t=1}^n$, whose proof is given in Section \ref{sec:gadpsc} of the supplementary material. See the definitions of dependency graph and its associated maximum degrees in Section \ref{subsec:m-dependent}.
\begin{theorem}
	[Gaussian approximation for partial sums of a sequence under dependency graph for simple convex sets]
	\label{pn:GA_dependency_graph_simple_convex}
	Assume $p\geq n^\kappa$ for some universal constant $\kappa>0$. Let $\mathcal{A}$ be a subclass of $\mathcal{A}^{\rm si}(a,d)$ such that  Conditions {\rm\ref{as:ExtendedTail}} and {\rm\ref{as:longrun_simple_convex}} are satisfied for any $A\in\mathcal{A}$. Then
	$$
	\rho_n(\mathcal{A})\lesssim \frac{a(d\log p)^{1/2}}{n}+\frac{B_{n}(D_{n}D_n^*)^{1/3}(d\log p)^{7/6}}{n^{1/6}}\,,$$
where $D_n$ and $D^*_n$ are the maximum degrees of the first-degree and second-degree of connections in the dependency graph generated by the sequence $\{X_t\}_{t=1}^n$, respectively.
\end{theorem}


The next corollary states a result for $m$-dependent sequences.

\begin{cy}
	\label{cor:GA_m-dependent_subexp_simple_convex}
	Assume $\{X_t\}_{t=1}^n$ is an $m$-dependent sequence with $p\geq n^\kappa$ for some universal constant $\kappa>0$. Let $\mathcal{A}$ be a subclass of $\mathcal{A}^{\rm si}(a,d)$ such that Conditions {\rm\ref{as:ExtendedTail}} and {\rm\ref{as:longrun_simple_convex}} are satisfied for any $A\in\mathcal{A}$. Then
	$$
	\rho_n(\mathcal{A})\lesssim \frac{a(d\log p)^{1/2}}{n}+\frac{B_{n}(m\vee 1)^{2/3}(d\log p)^{7/6}}{n^{1/6}}\,.
	$$
\end{cy}

\subsubsection{Sequence with physical dependence}\label{sec:phydepsim}
For a $K$-generated convex set $A^K$, write $\mathcal{V}(A^K)=\{v_1,\ldots,v_K\}$. Given $\mathcal{V}(A^K)$ and $\{X_t\}$, we can define a new $K$-dimensional sequence $X_t(A^K)=(v_1^\T X_t,\ldots,v_K^\T X_t)^\T$. If $\{X_{t}\}_{t=1}^n$ satisfies model \eqref{eqn:weakly_dependent_ts_fdp} with the jointly measurable function $f_t(\cdot)$, $\{X_t(A^K)\}$ also satisfies \eqref{eqn:weakly_dependent_ts_fdp} with a jointly measurable function $\tilde{f}_t(\cdot)=\{v_1^\T f_t(\cdot),\ldots,v_K^\T f_t(\cdot)\}^\T$. We further define $\Psi_{q,\alpha}(A^K)$ and  $\Phi_{\psi_{\nu},\alpha}(A^K)$ for any $K$-generated convex set $A^K$ in the same manner as $\Psi_{q,\alpha}$ and $\Phi_{\psi_{\nu},\alpha}$ in \eqref{eq:aggnorm} by replacing $\{X_{t}\}$ with $\{X_t(A^K)\}$. Given $\mathcal{A}\subset\mathcal{A}^{{\rm si}}(a,d)$, let
\begin{equation}\label{eq:cobphy}
\Psi_{q,\alpha,\mathcal{A}}=\sup_{A\in\mathcal{A}}\Psi_{q,\alpha}\{A^K(A)\}\quad\text{and}\quad\Phi_{\psi_{\nu},\alpha,\mathcal{A}}=\sup_{A\in\mathcal{A}}\Phi_{\psi_{\nu},\alpha}\{A^K(A)\}\,.
\end{equation}

\begin{theorem}
	[Gaussian approximation for partial sums of time series under functional dependence for simple convex sets]
	\label{pn:GA_weakly-dependent_simple_convex}
	Assume the sequence $\{X_t\}$ satisfies the model \eqref{eqn:weakly_dependent_ts_fdp} with $p\geq n^\kappa$ for some universal constant $\kappa>0$.  Let $\mathcal{A}$ be a subclass of $\mathcal{A}^{\rm si}(a,d)$ such that Condition {\rm\ref{as:longrun_simple_convex}} is satisfied for any $A\in\mathcal{A}$, and $\Phi_{\psi_{\nu},\alpha,\mathcal{A}}< \infty$ for some $\alpha, \nu \in (0,\infty)$.
	
{\rm (i)} It holds that
	\begin{align*}
	\rho_n(\mathcal{A})\lesssim \frac{a(d\log p)^{1/2}}{n}+&\frac{\Phi_{\psi_{\nu},0,\mathcal{A}}(d\log p)^{7/6}}{n^{\alpha/(3+9\alpha)}}\\
	+&\frac{\Psi_{2,\alpha,\mathcal{A}}^{1/3}\Psi_{2,0,\mathcal{A}}^{1/3}(d\log p)^{2/3}}{n^{\alpha/(3+9\alpha)}}+\frac{\Phi_{\psi_\nu,\alpha,\mathcal{A}}(d\log p)^{1+\nu}}{n^{\alpha/(1+3\alpha)}}
	\end{align*}
	provided that $(d\log p)^{\max\{6\nu-1,(5+6\nu)/4\}}=o\{n^{\alpha/(1+3\alpha)}\}$. 
	
	{\rm (ii)} If Condition {\rm\ref{as:ExtendedTail}} is satisfied for any $A\in\mathcal{A}$, it holds that 
\begin{align*}
	\rho_n(\mathcal{A})\lesssim \frac{a(d\log p)^{1/2}}{n} +&\frac{B_{n}(d\log p)^{7/6}}{n^{\alpha/(12+6\alpha)}}\\
	+&\frac{\Psi_{2,\alpha,\mathcal{A}}^{1/3}\Psi_{2,0,\mathcal{A}}^{1/3}(d\log p)^{2/3}}{n^{\alpha/(12+6\alpha)}}+\frac{\Phi_{\psi_{\nu},\alpha,\mathcal{A}}(d\log p)^{1+\nu}}{n^{\alpha/(4+2\alpha)}}\,.
	\end{align*}
\end{theorem}
The proof of Theorem \ref{pn:GA_weakly-dependent_simple_convex} is given in Section \ref{sec:gawdsc} of the supplementary material.

\subsection{High-dimensional CLT for sparsely convex sets}\label{sec:sparse}

We consider sparsely convex sets here, as a generalization of hyper-rectangles, that can be represented as intersections of convex sets whose indicator functions depend only on a small subset of their coordinates.

\begin{definition}[$s$-sparsely convex set]
\label{defn:s-sparsely_convex_set}
For an integer $s>0$, we say $A\subset\mathbb{R}^{p}$ is an \textit{$s$-sparsely convex set} if {\rm(i)} $A$ admits a sparse representation
$
A=\cap_{q=1}^{K_*}A_{q}$ for some positive integer $K_*$ and convex sets $A_{1},\ldots,A_{K_*}\subset\mathbb{R}^{p}$, and {\rm(ii)} the indicator function $I(w\in A_{q})$ depends on at most $s$ components of the vector $w\in\mathbb{R}^{p}$ (which we call the main components of $A_{q}$).
\end{definition}

Denote by $\mathcal{A}^{\rm sp}(s)$ the class of all $s$-sparsely convex sets in $\mathbb{R}^{p}$. In this section, we target on deriving the upper bounds for $\rho_n\{\mathcal{A}^{\rm sp}(s)\}$ when the observed data $\{X_t\}$ are (i) an $\alpha$-mixing sequence, (ii) an $m$-dependent sequence, and (iii) a physical dependence sequence.

\subsubsection{$\alpha$-mixing sequence}

\begin{ass}\label{as:longrun_s-sparsely}
	For $V_{n}(v)$ defined in \eqref{eq:Vtv}, there exists a universal constant $K_{5}>0$ such that $V_{n}(v)\geq K_{5}$	for any $v\in\mathbb{S}^{p-1}$ with $|v|_{0}\leq s$.
\end{ass}

Condition \ref{as:longrun_s-sparsely} holds automatically if the smallest eigenvalue of $\Xi={\rm Cov}(n^{-1/2}\sum_{t=1}^nX_t)$ is uniformly bounded away from zero.

\begin{theorem}
		[Gaussian approximation for partial sums of the $\alpha$-mixing sequence for $s$-sparsely convex sets]
		\label{pn:GA_s-sparsely}
				Assume $\{X_t\}$ is an $\alpha$-mixing sequence with $p\geq n^\kappa$ for some universal constant $\kappa>0$. Under Conditions {\rm\ref{as:tail}}, {\rm\ref{as:alpha-mixing}} and {\rm\ref{as:longrun_s-sparsely}}, it holds that
		$$\rho_n\{\mathcal{A}^{\rm sp}(s)\}\lesssim \frac{B_{n}^{2/3}s^{(2+6\gamma_{2})/(3\gamma_{2})}(\log p)^{(1+2\gamma_2)/(3\gamma_2)}}{n^{1/9}}+\frac{B_{n}s^{10/3}(\log p)^{7/6}}{n^{1/9}}$$
		provided that $(s^2\log p)^{3-\gamma_2}=o(n^{\gamma_2/3})$.
\end{theorem}
The proof of Theorem \ref{pn:GA_s-sparsely} is given in Section \ref{sec:pf_s-sparsely} of the supplementary material.

\subsubsection{$m$-dependent sequence}

\begin{theorem}
	[Gaussian approximation for partial sums of a sequence under dependency graph for $s$-sparsely convex sets]
	\label{pn:GA_dependency_graph_s_sparsely} Assume $p\geq n^\kappa$ for some universal constant $\kappa>0$. Under Conditions {\rm\ref{as:tail}} and {\rm\ref{as:longrun_s-sparsely}}, 
	it holds that
$$
	\rho_n\{\mathcal{A}^{\rm sp}(s)\} \lesssim \frac{s^{10/3}B_{n}(D_nD_n^*)^{1/3}(\log p)^{7/6}}{n^{1/6}}\,,$$
where $D_n$ and $D^*_n$ are the maximum degrees of the first-degree and second-degree of connections in the dependency graph generated by the sequence $\{X_t\}_{t=1}^n$, respectively.
\end{theorem}

The proof of Theorem \ref{pn:GA_dependency_graph_s_sparsely} is given in Section \ref{sec:gas} of the supplementary material. 
The next corollary states a result for $m$-dependent sequences.

\begin{cy}[Gaussian approximation for partial sums of an $m$-dependent sequence for $s$-sparsely convex sets]
	\label{cor:GA_m-dependent_subexp_s-sparsely}
	Assume $\{X_t\}_{t=1}^n$ is an $m$-dependent sequence with $p\geq n^\kappa$ for some universal constant $\kappa>0$. Under Conditions {\rm\ref{as:tail}} and {\rm\ref{as:longrun_s-sparsely}}, it holds that
	$$
	\rho_n\{\mathcal{A}^{\rm sp}(s)\}\lesssim \frac{s^{10/3}(m\vee 1)^{2/3}B_{n}(\log p)^{7/6}}{n^{1/6}}\,.$$
\end{cy}

\subsubsection{Sequence with physical dependence}


For any $c>0$, define $$\Omega_{s,c}=\bigg\{A\in\mathcal{A}^{\rm si}(1,cs^2): \max_{v\in\mathcal{V}\{A^{K}(A)\}}|v|_{0}\leq s\bigg\}\,.$$ Analogous to $\Psi_{q,\alpha,\mathcal{A}}$ and $\Phi_{\psi_{\nu},\alpha,\mathcal{A}}$ defined in \eqref{eq:cobphy}, we also define
$$
\Psi_{q,\alpha,\Omega_{s,c}}=\sup_{A\in\Omega_{s,c}}\Psi_{q,\alpha}\{A^K(A)\}\quad\text{and}\quad \Phi_{\psi_{\nu},\alpha,\Omega_{s,c}}=\sup_{A\in\Omega_{s,c}}\Phi_{\psi_{\nu},\alpha}\{A^K(A)\}$$
with $\Psi_{q,\alpha}(A^K)$ and  $\Phi_{\psi_{\nu},\alpha}(A^K)$ defined in Section \ref{sec:phydepsim}. It then holds that $\Psi_{q,\alpha,\Omega_{s,c}}\geq \Psi_{q,\alpha}$ and $\Phi_{\psi_{\nu},\alpha,\Omega_{s,c}}\geq \Phi_{\psi_{\nu},\alpha}$.

\begin{theorem}
	[Gaussian approximation for partial sums of time series under functional dependence for $s$-sparsely convex sets]
	\label{pn:GA_weakly-dependent_s_sparsely}
	Assume the sequence $\{X_t\}$ satisfies the model \eqref{eqn:weakly_dependent_ts_fdp} with $p\geq n^\kappa$ for some universal constant $\kappa>0$. For some sufficiently large constant $c>0$, let $\Phi_{\psi_{\nu},\alpha,\Omega_{s,c}}< \infty$ with some $\alpha, \nu \in (0,\infty)$.
	
{\rm (i)} Under Condition {\rm\ref{as:longrun_s-sparsely}}, it holds that
		\begin{align*}
		\rho_n\{\mathcal{A}^{\rm sp}(s)\}\lesssim \frac{\Phi_{\psi_{\nu},0,\Omega_{s,c}}(s^2\log p)^{7/6}}{n^{\alpha/(3+9\alpha)}}+&\frac{\Psi_{2,\alpha,\Omega_{s,c}}^{1/3}\Psi_{2,0,\Omega_{s,c}}^{1/3}(s^2\log p)^{2/3}}{n^{\alpha/(3+9\alpha)}}\\
		+&\frac{\Phi_{\psi_{\nu},\alpha,\Omega_{s,c}}(s^2\log p)^{1+\nu}}{n^{\alpha/(1+3\alpha)}}+\frac{\Phi_{\psi_\nu,0}^3s^{17/4}}{n^{\alpha/(2+4\alpha)}}		\end{align*}
		 provided that $(s^2\log p)^{\max\{6\nu-1,(5+6\nu)/4\}}=o\{n^{\alpha/(1+3\alpha)}\}$.
		 
		  {\rm (ii)} Under Conditions {\rm \ref{as:tail}} and {\rm\ref{as:longrun_s-sparsely}}, it holds that 
\begin{align*}
		\rho_n\{\mathcal{A}^{\rm sp}(s)\}\lesssim \frac{s^{10/3}B_{n}(\log p)^{7/6}}{n^{\alpha/(12+6\alpha)}}+&\frac{\Psi_{2,\alpha,\Omega_{s,c}}^{1/3}\Psi_{2,0,\Omega_{s,c}}^{1/3}(s^2\log p)^{2/3}}{n^{\alpha/(12+6\alpha)}}\\
		+&\frac{\Phi_{\psi_{\nu},\alpha,\Omega_{s,c}}(s^2\log p)^{1+\nu}}{n^{\alpha/(4+2\alpha)}}\,.
		\end{align*}
\end{theorem}
The proof of Theorem \ref{pn:GA_weakly-dependent_s_sparsely} is given in Section \ref{sec:gawds} of the supplementary material.

\section{Parametric bootstrap}\label{sec:parametric_boostrap}
In Section \ref{sec:main_results}, we have established the error bounds for $\rho_n(\mathcal{A})$ defined as \eqref{eq:gb} when $\{X_t\}$ is (i) an $\alpha$-mixing sequence, (ii) an $m$-dependent sequence, and (iii) a physical dependence sequence. Since $\mathbb{P}(G\in A)$ with $G\sim N(0,\Xi)$ depends on the unknown long-run covariance matrix $\Xi$, to approximate $\mathbb{P}(S_{n,x}\in A)$ in practice, we need to construct a data-dependent Gaussian analogue $\hat{G}$ of $G$. In this section, we propose a parametric bootstrap procedure to construct $\hat{G}\sim N(0,\hat{\Xi}_{n})$ for some covariance matrix $\hat{\Xi}_n$ that is close to $\Xi$ and establish its theoretical validity.  Define
\begin{align*}
\hat{\rho}_n(\mathcal{A}) := \sup_{A\in\mathcal{A}}|\mathbb{P}(S_{n,x} \in A)-\mathbb{P}(\hat{G}\in A\,|\,\mathcal{X}_n)|
\end{align*}
with $\mathcal{X}_n=\{X_{1},\ldots,X_{n}\}$.  Let $\Delta_{n,r}=|\hat{\Xi}_{n}-\Xi|_{\infty}$ and 
$$
\Delta_{n}(\mathcal{A})=\sup_{A\in \mathcal{A}}\sup_{v_{1},v_{2}\in \mathcal{V}\{A^{K}(A)\}}|v_{1}^{\T}(\hat{\Xi}_{n}\\-\Xi)v_{2}|$$ 
for any $\mathcal{A}\subset\mathcal{A}^{\rm si}(a,d)$.
Theorem \ref{tm:pb_rectangle} establishes primitive error bounds for $\hat{\rho}_n(\mathcal{A})$ when $\mathcal{A}$ is selected as the class of all hyper-rectangles, the class of simple convex sets and the class of $s$-sparsely convex sets, respectively, whose proof is given in Section \ref{pf:pd_rectangle} of the supplementary material.
\begin{theorem}[Rates of convergence for parametric bootstrap]
	\label{tm:pb_rectangle}
	Assume $p\geq n^\kappa$ for some universal constant $\kappa>0$. 
	
{\rm (i)} Under Condition {\rm\ref{as:longrun}}, it holds that
$
	\hat{\rho}_n(\mathcal{A}^{\rm re})\lesssim \rho_n(\mathcal{A}^{\rm re})+\Delta_{n,r}^{1/3}(\log p)^{2/3}$.
	
{\rm (ii)} Let $\mathcal{A}$ be a subclass of $\mathcal{A}^{\rm si}(a,d)$ such that Condition {\rm\ref{as:longrun_simple_convex}} is satisfied for every $A\in\mathcal{A}$. It holds that
$
	\hat{\rho}_n(\mathcal{A})\lesssim \rho_n(\mathcal{A})+  an^{-1}(d\log p)^{1/2}+\Delta_{n}^{1/3}(\mathcal{A})(d\log p)^{2/3}$.	
	
{\rm (iii)} Under Conditions {\rm\ref{as:tail}} and {\rm\ref{as:longrun_s-sparsely}}, it holds that  
	$\hat{\rho}_n\{\mathcal{A}^{\rm sp}(s)\}\lesssim \rho_n\{\mathcal{A}^{\rm sp}(s)\}+s^{2}\Delta_{n,r}^{1/3}(\log p)^{2/3}+\{B_{n}+s(\log p)^{1/2}\}n^{-1}$.
\end{theorem}

\begin{remark}\label{remark:paraboot}
{\rm(i)} For the case of hyper-rectangles, $\Delta_{n,r}=o_{\rm p}\{(\log p)^{-2}\}$ is necessary to guarantee $\hat{\rho}_n(\mathcal{A}^{\rm re})=o_{\rm p}(1)$. {\rm(ii)} For $\mathcal{A}\subset \mathcal{A}^{\rm si}(a,d)$, to make $\hat{\rho}_n(\mathcal{A})=o_{\rm p}(1)$, $\hat{\Xi}_n$ should satisfy $\Delta_n(\mathcal{A})=o_{\rm p}\{(d\log p)^{-2}\}$. Notice that $\Delta_n(\mathcal{A})\leq \Delta_{n,r}\sup_{A\in\mathcal{A}}\sup_{v\in\mathcal{V}\{A^K(A)\}}|v|_1^2$. If the $\ell^1$-norm of the unit normal vectors outward to the facets of $A^K(A)$ is uniformly bounded away from infinity over $A\in\mathcal{A}$,  it suffices to require $\Delta_{n,r}=o_{\rm p}\{(d\log p)^{-2}\}$. {\rm(iii)} For the case of $s$-sparsely convex sets, to make $\hat{\rho}_n\{\mathcal{A}^{\rm sp}(s)\}=o_{\rm p}(1)$, we need to require $\Delta_{n,r}=o_{\rm p}\{(s^3\log p)^{-2}\}$. 
\end{remark}

As we have discussed in Remark \ref{remark:paraboot}, the validity of our proposed parametric bootstrap only requires the estimated long-run covariance matrix $\hat{\Xi}_n$ satisfying $|\hat{\Xi}_n-\Xi|_\infty=o_{\rm p}(\delta_n)$ for some $\delta_n\rightarrow0$ as $n\rightarrow\infty$, where $\delta_n$ will be different for different selections of $\mathcal{A}$. There are various estimation methods for long-run covariance matrices, including the kernel-type estimators \citep{Andrews_1991} and utilizing moving block bootstraps \citep{Lahiri_2003}. See also \cite{DenHanLevin_1997} and \cite{Kieferetaj_2000}. Since the data sequence $\{X_t\}_{t=1}^n$ may be non-stationary, we suggest to adopt the kernel-type estimator for its long-run covariance matrix, that is
\begin{equation}
\hat{\Xi}_{n}=\sum_{j=-n+1}^{n-1}\mathcal{K}\bigg(\frac{j}{b_{n}}\bigg)\hat{H}_{j}\,,\label{eq:kernelest}
\end{equation}
where $\hat{H}_{j}=n^{-1}\sum_{t=j+1}^{n}(X_{t}-\bar{X})(X_{t-j}-\bar{X})^{\T}$ if $j\geq 0$ and $\hat{H}_{j}=n^{-1}\sum_{t=-j+1}^{n}(X_{t+j}-\bar{X})(X_{t}-\bar{X})^{\T}$ otherwise, with $\bar{X}=n^{-1}\sum_{t=1}^{n}X_{t}$. Here $\mathcal{K}(\cdot)$ is a symmetric kernel function that is continuous at 0 with $\mathcal{K}(0)=1$, and $b_{n}$ is the bandwidth diverging with $n$. Among a variety of kernel functions that guarantee the positive definiteness of the long-run covariance matrix estimators, \cite{Andrews_1991} derived an optimal kernel, i.e., the quadratic spectral kernel
\begin{equation}\label{eq:QSkernel}
\mathcal{K}_{\rm QS}(x)=\frac{25}{12\pi^{2}x^{2}}\bigg\{\frac{\sin(6\pi x/5)}{6\pi x/5}-\cos(6\pi x/5) \bigg\}\,,
\end{equation}
by minimizing the asymptotic truncated mean square error of the estimator. To study the property of $\hat{\Xi}_n$ given in \eqref{eq:kernelest}, we need the next condition imposed on the kernel $\mathcal{K}(\cdot)$.

\begin{ass}[Kernel regularity]
\label{as:kernel} The kernel function $\mathcal{K}(\cdot):\mathbb{R}\rightarrow [-1,1]$ is continuously differentiable with bounded derivatives on $\mathbb{R}$ and satisfies ${\rm(i)}\,\mathcal{K}(0)=1$, ${\rm(ii)}\,\mathcal{K}(x)=\mathcal{K}(-x)$ for any $x\in\mathbb{R}$, and ${\rm (iii)}\, |\mathcal{K}(x)|\lesssim |x|^{-\vartheta}$ as $|x|\to\infty$ for some constant $\vartheta>1$.
\end{ass}

\begin{theorem}[Bounds on $\Delta_{n,r}$]
	\label{tm:cov_compare}
	Assume $p\geq n^{\kappa}$ for some universal constant $\kappa>0$. Let Condition {\rm\ref{as:kernel}} hold and $b_n\asymp n^{\rho}$. 
	
	{\rm(i)} For $\alpha$-mixing sequence $\{X_t\}$ with Conditions {\rm\ref{as:tail}} and {\rm\ref{as:alpha-mixing}} being satisfied, if $0<\rho<(\vartheta-1)/(3\vartheta-2)$, there exist two constants $c_1>0$ depending only on $(\rho,\vartheta)$ and  $c_2>0$ depending only on $(\gamma_1,\gamma_2,\vartheta)$ such that $\Delta_{n,r}=O_{\rm p}\{B_n^2n^{-c_1}(\log p)^{c_2}\}+O(B_n^2n^{-\rho}).$
	
	{\rm(ii)} For $m$-dependent sequence $\{X_t\}$ with Condition {\rm\ref{as:tail}} being satisfied, if $0<\rho<(\vartheta-1)/(3\vartheta-2)$, there exist two constants $c_1>0$ depending only on $(\rho,\vartheta)$ and $c_2>0$ depending only on $(\gamma_1,\vartheta)$ such that $\Delta_{n,r}=O_{\rm p}\{B_n^2n^{-c_1}(\log p)^{c_2}\}+O(B_n^2m^2n^{-\rho}).$
	
	{\rm(iii)} For the sequence $\{X_t\}$ satisfying the model \eqref{eqn:weakly_dependent_ts_fdp},  assume Condition {\rm\ref{as:tail}} is satisfied and $\Phi_{\psi_\nu,0}<\infty$ and $\Psi_{2,\alpha}<\infty$ for some $\alpha,\nu\in(0,\infty)$, if $0<\rho<\min\{(2\alpha+2\vartheta-3)/(2\vartheta-2),1/(2\alpha)\}$ with $\alpha>(3-2\vartheta)/2$,  there exist two constants $c_1>0$ depending only on $(\rho,\vartheta,\alpha)$ and $c_2>0$ depending only on $(\alpha,\vartheta,\gamma_{1},\nu)$ such that $\Delta_{n,r}=O_{\rm p}\{(B_n\Phi_{\psi_{\nu,0}}+B_n^2+\Phi_{\psi_{\nu,0}}^2)n^{-c_1}(\log p)^{c_2}\}+O(n^{-\rho}\Psi_{2,0}\Psi_{2,\alpha}\varpi_{n})$ with $\varpi_n=(\log n)I(\alpha=1)+n^{1-\alpha}I(\alpha\neq 1)$.
\end{theorem}
The proof of Theorem \ref{tm:cov_compare} is given in Section \ref{sec:cov_compare} of the supplementary material.

\begin{remark}
Let  
$\Xi^* = \sum_{j=-n+1}^{n-1} \mathcal K(j/b_n)H_j$ with $H_{j}=n^{-1}\sum_{t=j+1}^{n}\mathbb{E}(X_{t}X_{t-j}^{\T})$ if $j\geq 0$ and $H_{j}=n^{-1}\sum_{t=-j+1}^{n}\mathbb{E}(X_{t+j}X_{t}^{\T})$ if $j<0$. The terms $O(B_n^2n^{-\rho})$, $O(B_n^2m^2n^{-\rho})$ and $O(n^{-\rho}\Psi_{2,0}\Psi_{2,\alpha}\varpi_n)$ in {\rm(i)}, {\rm(ii)} and {\rm(iii)} of Theorem \ref{tm:cov_compare} are, respectively, bounds on the bias $|\Xi^*-\Xi|_\infty$ in the three cases. Thus the bias terms do not depend on $p$ (at least directly). On the other hand, the terms $O_{\rm p}\{B_n^2n^{-c_1}(\log p)^{c_2}\}$, $O_{\rm p}\{B_n^2n^{-c_1}(\log p)^{c_2}\}$ and $O_{\rm p}\{(B_n\Phi_{\psi_{\nu,0}}+B_n^2+\Phi_{\psi_{\nu,0}}^2)n^{-c_1}(\log p)^{c_2}\}$ are bounds on $|\hat{\Xi}_n-\Xi^*|_\infty$ in the three cases, respectively. 
\end{remark}

Now, combining Remark \ref{remark:paraboot} and Theorem \ref{tm:cov_compare}, we see that our proposed parametric bootstrap procedure is asymptotically valid even if the dimension $p$ grows sub-exponentially fast in the sample size $n$. To implement the proposed parametric bootstrap, we need to solve two problems: (i) How to select bandwidth $b_n$ in practice? and (ii) How to generate $\hat{G}\sim N(0,\hat{\Xi}_n)$ efficiently when $p$ is large? For Problem (i), due to the positive definiteness of $\hat{\Xi}_n$ defined as \eqref{eq:kernelest} with the quadratic spectral kernel $\mathcal{K}_{{\rm QS}}(\cdot)$ defined as \eqref{eq:QSkernel}, we can use this kernel in our parametric bootstrap procedure. For $\mathcal{K}_{{\rm QS}}(\cdot)$,  \cite{Andrews_1991} suggested Algorithm \ref{alg:1} to select $b_n$. For Problem (ii), to generate a random vector $\hat{G}\sim N(0,\hat{\Xi}_n)$, the standard approach consists of three steps: (a) perform the Cholesky decomposition for the $p\times p$ matrix $\hat{\Xi}_n=L^\T L$; (b) generate independent standard normal random variables $Z_1,\ldots,Z_p$ and let $Z=(Z_1,\ldots,Z_p)^\T$; (c) perform the transformation $\hat{G}=L^\T Z$. However, the computation complexity of the standard approach is $O(np^2+p^3)$ and it also requires a large storage space for $\{X_t\}_{t=1}^n$ and the estimated matrix $\hat{\Xi}_n$. To circumvent the high computing cost with large $p$, we propose Algorithm \ref{alg:2} below which involves generating a random vector from an $n$-dimensional normal distribution instead. It is easy to check that such obtained $\hat{G}\sim N(0,\hat{\Xi}_n)$ conditionally on $\mathcal{X}_n$. Algorithm \ref{alg:2} was initially introduced in \cite{ChangYaoZhou_2017}.

\begin{algorithm}[H] \caption{Data-driven procedure of bandwidth selection for $\mathcal{K}_{\rm QS}(\cdot)$}
\begin{itemize}\label{alg:1}
\item[{\bf Step 1.}] For each $j\in[p]$, fit an AR(1) model to the $j$-th coordinate marginal sequence $\{X_{t,j}\}_{t=1}^n$. Denote by $\hat{\rho}_j$ and $\hat{\sigma}_j^2$, respectively, the estimated autoregressive coefficient and innovation variance. \\

\item[{\bf Step 2.}] Select $b_n=1.3221(\hat{a}n)^{1/5}$ with
$
\hat{a}=\{\sum_{j=1}^p4\hat{\rho}_j^2\hat{\sigma}_j^4(1-\hat{\rho}_j)^{-8}\}/\{\sum_{j=1}^p\hat{\sigma}_j^4(1-\hat{\rho}_j)^{-4}\}$.
\end{itemize}
\end{algorithm}

\begin{algorithm}[H] \caption{Generating $\hat{G}\sim N(0,\hat{\Xi}_n)$ for $\hat{\Xi}_n$ given in \eqref{eq:kernelest} with $\mathcal{K}_{{\rm QS}}(\cdot)$}
\begin{itemize}\label{alg:2}
\item[{\bf Step 1.}] Obtain the bandwidth $b_n$ by Algorithm \ref{alg:1}. Define $\Theta=(\theta_{i,j})_{n\times n}$ with $\theta_{i,j}=\mathcal{K}_{{\rm QS}}\{(i-j)/b_n\}$. \\

\item[{\bf Step 2.}] Generate $Z=(Z_1,\ldots,Z_n)^\T\sim N(0,\Theta)$ independent of $\{X_t\}_{t=1}^n$. Define
$
\hat{G}=n^{-1/2}\sum_{t=1}^nZ_tX_t$.
\end{itemize}
\end{algorithm}

\section{Applications}
\label{sec:apps}

 In this section, we discuss several statistical applications of the Gaussian and parametric bootstrap approximation results for high-dimensional dependent data developed in Sections \ref{sec:main_results} and \ref{sec:parametric_boostrap}.

\subsection{Testing high-dimensional mean vector}
\label{ex:test_mean}
Given data $\{X_t\}_{t=1}^{n_1}$ with $\mathbb{E}(X_t)= \theta_x\in\mathbb{R}^p$ for any $t\in[n_1]$, it is of general interest in testing the hypothesis
\begin{equation}\label{eq:one-sample}
H_0 : \theta_x = 0 \quad \mbox{versus} \quad H_1 : \theta_x \neq 0\,.
\end{equation}
If there is another group of data $\{Y_t\}_{t=1}^{n_2}$ with $\mathbb{E}(Y_t)=\theta_y\in\mathbb{R}^p$ for any $t\in[n_2]$, we are also interested in the hypothesis testing problem
\begin{equation}\label{eq:two-sample}
H_0 : \theta_x = \theta_y \quad \mbox{versus} \quad H_1 : \theta_x \neq \theta_y\,.
\end{equation}
Hypotheses \eqref{eq:one-sample} and \eqref{eq:two-sample} are called, respectively, one-sample and two-sample mean testing problems in the literature. Lots of statistical inference problems in practice can be formulated as \eqref{eq:one-sample} and \eqref{eq:two-sample}. Generally, the $\ell^2$-type and $\ell^\infty$-type statistics are used to test the hypotheses \eqref{eq:one-sample} and \eqref{eq:two-sample} in the high-dimensional settings. With independent data, we refer to \cite{ChenQin_2010} and \cite{CaiLiuXia_2014} for the uses of $\ell^2$-type statistic and $\ell^\infty$-type statistic in these testing problems, respectively. It has been well known that the $\ell^2$-type statistics are powerful for detecting relatively dense signals while the $\ell^\infty$-type statistics are preferable for detecting relatively sparse signals. In practice, we usually have less knowledge on whether the signals are dense or sparse. Let $\Gamma = \{\Cov(X_t)\}^{-1} = (\gamma_{i,j})_{p\times p}$. When $\Gamma$ is known, to combine the advantages of the $\ell^2$-type and $\ell^\infty$-type statistics, \cite{Zhang_2015} considered the following test statistic for \eqref{eq:one-sample} with independent data $\{X_t\}_{t=1}^{n_1}$:
\begin{equation}
\label{eqn:test_mean_vector_statistic}
T_n(s) =\max_{1 \leq j_1< \cdots < j_s \leq p} \sum_{k=1}^s {n_1\bar{Z}_{j_k}^2 \over \gamma_{j_k,j_k}}\,,
\end{equation}
where  $\bar{Z} = \Gamma \bar{X}:=(\bar{Z}_1,\ldots,\bar{Z}_p)^\T$ with $\bar{X} = n_1^{-1}\sum_{t=1}^{n_1} X_t$. When $\Gamma$ is unknown, \cite{Zhang_2015} proposed a feasible analogue for $T_n(s)$ by replacing $\Gamma$ by its estimator $\hat{\Gamma}$. To simplify our presentation, we assume $\Gamma$ is known in the rest of this subsection.

For dependent data, testing for white noise or serial correlation is a fundamental problem in statistical inference, as many testing problems in linear modelling can be transformed into a white noise test. Let $\{\varepsilon_t\}$ be a $d$-dimensional weakly stationary time series with mean zero. Denote by $\Sigma(k)=\cov(\varepsilon_{t+k},\varepsilon_t)$ the autocovariance of $\varepsilon_t$ at lag $k$. Given a prescribed integer $K$, the white noise hypothesis of $\{\varepsilon_t\}$ can be formulated as
\begin{align}\label{eq:whitenoise}
H_{0} : \Sigma(1) = \dots = \Sigma(K) = 0 \quad \mbox{versus} \quad H_1: H_0\,\,\mbox{is not true}\,.
\end{align}
Let $n_1=n-K$ and $X_t=\{{\rm vec}(\varepsilon_{t+1}\varepsilon_t^\T),\ldots,{\rm vec}(\varepsilon_{t+K}\varepsilon_t^\T)\}^\T$, where ${\rm vec}(A)$ denotes a row vector that collecting all the elements in $A$. Then the white noise hypothesis \eqref{eq:whitenoise} can be covered by the hypothesis \eqref{eq:one-sample} with $p=d^2K$. \cite{ChangYaoZhou_2017} proposed a bootstrap test based on the $\ell^{\infty}$-type statistic for the white noise hypothesis \eqref{eq:whitenoise} under the $\beta$-mixing assumption of $\{\varepsilon_t\}$. To enhance the power performance of \cite{ChangYaoZhou_2017}, we can use the test statistic $T_n(s)$ given in \eqref{eqn:test_mean_vector_statistic}.

Notice that the distribution function of $T_n(s)$ can be written in terms of probability of the random vector $n_1^{1/2}\bar{Z}$ over a class of convex subsets of the form $\{w \in \R^p : \sum_{j \in \Theta_s} \gamma_{j,j}^{-1}w_j^2 \leq t \}$ with $\Theta_s = \{w \in \R^p : |w|_0 = s\}$ which is a subset of the class of all $s$-sparsely convex sets in $\mathbb{R}^{p}$. Let $\hat{G}=(\hat{G}_1,\ldots,\hat{G}_p)^\T\sim N(0,\hat{\Xi}_n)$ with $\hat\Xi_n$ being the kernel-type estimator of the long-run covariance matrix $\cov(n_1^{1/2}\bar{Z})$. Using the results developed in Sections \ref{sec:main_results} and \ref{sec:parametric_boostrap}, the null-distribution of $T_n(s)$ can be approximated by that of
$$
\hat{T}_n=\max_{1 \leq j_1 < \cdots < j_s \leq p} \sum_{k=1}^s \frac{\hat{G}_{j_k}^2}{\gamma_{j_k,j_k}}$$
under both the $\alpha$-mixing assumption and physical dependency assumption, where $s$ can diverge with $n$ at some polynomial rate. For any $\delta\in(0,1)$, let $q_{\delta}$ be the upper $\delta$-quantile of the distribution of $\hat{T}_n$. Given the significant level $\delta$, we reject the null hypothesis of the white noise hypothesis \eqref{eq:whitenoise} if the test statistic $T_n(s)$ specified in \eqref{eqn:test_mean_vector_statistic} is larger than $q_\delta$. Our procedure allows arbitrary dependency among the components of $X_t$.

\subsection{Change point detection}
\label{ex:change_point_detection}
Consider the problem of change point detection for high-dimensional distributions in a location family
$
X_t= \theta \cdot I(t > m) + \xi_t$,
where $\theta \in \mathbb{R}^p$ is the location-shift parameter and $\{\xi_t\}$ is a sequence of stationary time series noise with mean zero. If $\theta=0$ or $m\geq n$, there is no change point in $\{X_t\}_{t=1}^n$. \cite{YuChen_2019} proposed a procedure to test whether there exists change point in the data based on the $U$-statistic
$$
U_n = (U_{n,1},\ldots,U_{n,p})^\T = {n \choose 2}^{-1} \sum_{1 \leq i < j \leq n} h(X_i, X_j)\,,$$
where $h : \mathbb{R}^p \times \mathbb{R}^p \to \mathbb{R}^p$ is an anti-symmetric kernel $h(x,y)=-h(y,x)$. The anti-symmetry
of the kernel $h$ plays a key role in testing for the change point in terms of noise cancellations so that after proper normalization the distribution of $U_n$ can be approximated by that of a Gaussian analogue. Specifically, under the null hypothesis that there is no change point and assuming independent and identically distributed  noise $\{\xi_t\}$ with distribution $F$, \cite{YuChen_2019} showed that $\mathbb{E}(U_n)= 0$ and the probability of $U_n$ on max-rectangles can be well approximated by that of $N(0,4\Gamma/3)$, where $\Gamma={\rm Cov}\{g(X_1)\}$ and $g(x) = \mathbb{E}\{h(x,X_2)\}$ is the $\mathcal{L}^{2}(F)$ projection of $h$ onto a linear subspace. On the other hand, under the alternative hypothesis when the change point location $m$ is known, $U_n$ is a two-sample Mann-Whitney test statistic (see e.g., Chapter 12 in \cite{vanderVaart_1998}), and  the signal distortion under certain nonlinear kernels can be controlled such that the between-sample change point signal is magnitude preserving. To practically calibrate the distribution of $\max_{j \in[p]} |n^{1/2} U_{n,j}|_\infty$, \cite{YuChen_2019} proposed a jackknife multiple bootstrap, which is powerful against alternatives with strong signals.

However, validity of the jackknife multiple bootstrap with a general nonlinear kernel heavily relies on the independent and identically distributed assumption of the noise sequence $\{\xi_t\}_{t=1}^n$ \citep{ChenKato_2020}. For time series data, with the linear kernel $h(x,y)=x-y$ we may write
$$
W_n=(W_{n,1},\ldots,W_{n,p})^\T = 2n^{-1/2}(n-1)^{-1} \sum_{t=1}^n (n-2t+1) X_t\,,$$
which can be viewed as one-pass CUSUM test statistic \citep{YuChen_2020}.  Thus we can enhance the power performance of the change point test of \cite{YuChen_2019} in the setting of linear kernel by using the test statistic
$$
T_n(s) = \max_{1 \leq j_1 < \dots < j_s \leq p} \sum_{k=1}^s W_{n,j_k}^2\,,$$
which allows $s$ to diverge with $n$ at some polynomial rate. Let $\hat{G} = (\hat{G}_1,\ldots,\hat{G}_p)^\T\sim N(0,\hat\Xi_n)$ with $\hat\Xi_n$ being the kernel-type estimator of the long-run covariance matrix ${\rm Cov}(W_n)$ of the weighted sequence $\{2(n-1)^{-1}(n-2t+1)X_t\}_{t=1}^n$. Based on the results developed in Sections \ref{sec:main_results} and  \ref{sec:parametric_boostrap}, the null-distribution of $T_n(s)$ can be calibrated by that of
$$
\hat{T}_n=\max_{1 \leq j_1 < \dots < j_s \leq p} \sum_{k=1}^s \hat{G}_{j_k}^2\,,$$
under both the $\alpha$-mixing assumption and physical dependence assumption. For any $\delta\in(0,1)$, let $q_{\delta}$ be the upper $\delta$-quantile of the distribution of $\hat{T}_n$. Given the significant level $\delta$, we reject the null hypothesis that there is no change point if the test statistic $T_n(s) = \max_{1 \leq j_1 < \dots < j_s \leq p} \sum_{k=1}^s W_{n,j_k}^2>q_\delta$. Our procedure does not need to impose any specific structure assumption on the dependency among different components of $X_t$.

\subsection{Confidence regions for the instantaneous covariance matrix and its inverse}
\label{ex:test_covariance_matrix}
Given $d$-dimensional dependent (and possibly non-stationary) data $\{Y_t\}_{t=1}^n$ with mean zero and instantaneous covariance $\Sigma$, i.e., $\mathbb{E}(Y_t)=0$ and ${\rm Cov}(Y_t)=\Sigma$ for any $t\in[n]$,  the instantaneous covariance matrix $\Sigma$ and the precision matrix $\Omega=\Sigma^{-1}=(\omega_{i,j})_{d\times d}$ quantify the dependence among the $d$ components of $Y_t$.
Confidence regions for $\Sigma$ and $\Omega$ can quantify the uncertainty in their estimates. For a given index set $\mathcal{S}\subset [d]^2$, denote by $\Sigma_{\mathcal{S}}$ and $\Omega_{\mathcal{S}}$ the vectors consisting, respectively, the entries of $\Sigma$ and $\Omega$ with their indices in $\mathcal{S}$. We are interested in constructing a class of confidence regions $\{\mathcal{C}_{\mathcal{S},\delta}\}_{0<\delta<1}$ for $\Sigma_{\mathcal{S}}$ such that
\begin{equation}\label{eq:confidenceregion}
\sup_{0<\delta<1}|\mathbb{P}(\Sigma_{\mathcal{S}}\in\mathcal{C}_{\mathcal{S},\delta})-\delta|\rightarrow0\quad\mbox{as $n,d\rightarrow\infty$}\,.
\end{equation}
We can also consider the confidence regions $\{\mathcal{C}_{\mathcal{S},\delta}\}_{0<\delta<1}$ for $\Omega_{\mathcal{S}}$ such that
$$
\sup_{0<\delta<1}|\mathbb{P}(\Omega_{\mathcal{S}}\in\mathcal{C}_{\mathcal{S},\delta})-\delta|\rightarrow0\quad\mbox{as $n,d\rightarrow\infty$}\,.$$

Given observations $\{Y_t\}_{t=1}^n$, we can estimate $\Sigma$ as $\hat{\Sigma}=n^{-1}\sum_{t=1}^nY_tY_t^\T$ and estimate $\Omega$ by fitting the node-wise regressions
$
Y_{j,t}=\sum_{k\neq j}\beta_{j,k}Y_{k,t}+\epsilon_{j,t}$ for each $j\in[d]$. For high-dimensional scenario, we need to use the regularization method to estimate the parameters in the node-wise regressions. More specifically, write $\beta_j=(\beta_{j,1},\ldots,\beta_{j,j-1},-1,\beta_{j,j+1},\ldots,\beta_{j,d})^\T$ which can be estimated as 
\[
\hat{\beta}_j\equiv(\hat{\beta}_{j,1},\ldots,\hat{\beta}_{j,j-1},-1,\hat{\beta}_{j,j+1},\hat{\beta}_{j,d})^\T=\arg\min_{\gamma\in\Theta_j}\bigg\{\frac{1}{n}\sum_{t=1}^n(\gamma^\T Y_t)^2+2\lambda_j|\gamma|_1\bigg\}\,,
\]
where $\Theta_j=\{\gamma=(\gamma_1,\ldots,\gamma_d)^\T\in\mathbb{R}^d:\gamma_j=-1\}$ and $\lambda_j>0$ is the tuning parameter. Let $V={\rm Cov}(\epsilon_t):=(v_{i,j})_{d\times d}$ with $\epsilon_t=(\epsilon_{1,t},\ldots,\epsilon_{d,t})^\T$, which can be estimated as $\hat{V}=(\hat{v}_{i,j})_{d\times d}$ with
\[
\hat{v}_{i,j}=
\bigg\{-\frac{1}{n}\sum_{t=1}^n(\hat{\epsilon}_{i,t}\hat{\epsilon}_{j,t}+\hat{\beta}_{i,j}\hat{\epsilon}_{j,t}^2+\hat{\beta}_{j,i}\hat{\epsilon}_{i,t}^2)\bigg\}I(i\neq j)+\bigg(\frac{1}{n}\sum_{t=1}^n\hat{\epsilon}_{i,t}\hat{\epsilon}_{j,t}\bigg)I(i=j)
\]
and $\hat{\epsilon}_{j,t}=-\hat{\beta}_j^\T Y_t$ for $j\in[d]$ and $t\in[n]$. Then $\hat{\Omega}=(\hat{\omega}_{i,j})_{d\times d}$ with $\hat{\omega}_{i,j}=\hat{v}_{i,j}/(\hat{v}_{i,i}\hat{v}_{j,j})$ provides an estimate of $\Omega$. Proposition 1 of \cite{ChangQiuYaoZou_2018} shows that $$\hat{\omega}_{i,j}-\omega_{i,j}=-\frac{1}{n}\sum_{t=1}^nv_{i,i}^{-1}v_{j,j}^{-1}(\epsilon_{i,t}\epsilon_{j,t}-v_{i,j})+o_{\rm p}\{(n\log d)^{-1/2}\}\,,$$ where the remainder term $o_{\rm p}\{(n\log d)^{-1/2}\}$  holds uniformly over $(i,j)\in[d]^2$. Let $p=|\mathcal{S}|$. We can see the leading terms of $\hat{\Sigma}_{\mathcal{S}}-\Sigma_{\mathcal{S}}$ and $\hat{\Omega}_{\mathcal{S}}-\Omega_{\mathcal{S}}$ can both be formulated as a general form $n^{-1}\sum_{t=1}^nX_t$ for some $p$-dimensional dependent sequence $\{X_t\}_{t=1}^n$. Write $\Xi={\rm Cov}(n^{-1/2}\sum_{t=1}^nX_t)$ and denote by $\hat{\Xi}_n$ the estimate of $\Xi$ given in Section  \ref{sec:parametric_boostrap}. Let $\hat{G}^{(\mathcal{S})}=\{\hat{G}_1^{(\mathcal{S})},\ldots,\hat{G}_p^{(\mathcal{S})}\}^\T\sim N(0,\hat{\Xi}_n)$ and define
$$
f_s\{\hat{G}^{(\mathcal{S})}\}=\max_{1\leq j_1<\cdots<j_s\leq p}\sum_{k=1}^sa_k\{\hat{G}_{j_k}^{(\mathcal{S})}\}^2\,,$$
where $a_1,\ldots,a_s>0$ denote the prescribed weights. For any $\delta\in(0,1)$, let $q_{\mathcal{S},\delta}$ be the upper $\delta$-quantile of the distribution of $f_s\{\hat{G}^{(\mathcal{S})}\}$, which can be determined by Monte Carlo simulation. Write $\hat{\Sigma}_{\mathcal{S}}=\{\hat{\sigma}_1^{(\mathcal{S})},\ldots,\hat{\sigma}_p^{(\mathcal{S})}\}^\T$. We can select the confidence region $\mathcal{C}_{\mathcal{S},\delta}$ for $\Sigma_{\mathcal{S}}$ as follows:
\begin{equation}\label{eq:conf}
\mathcal{C}_{\mathcal{S},\delta}=\bigg\{\xi=(\xi_1,\ldots,\xi_p)^\T:\max_{1\leq j_1<\cdots<j_s\leq p}\sum_{k=1}^sa_k\{\hat{\sigma}_{j_k}^{(\mathcal{S})}-\xi_{j_k}\}^2\leq q_{\mathcal{S},\delta}\bigg\}\,.
\end{equation}
Based on the results developed in Sections \ref{sec:main_results} and \ref{sec:parametric_boostrap}, we know such defined $\mathcal{C}_{\mathcal{S},\delta}$ satisfies \eqref{eq:confidenceregion}. Analogously, we can also obtain the confidence region $\mathcal{C}_{\mathcal{S},\delta}$ for $\Omega_{\mathcal{S}}$ in the same manner. \cite{ChangQiuYaoZou_2018} used this idea to construct the confidence region for $\Omega_{\mathcal{S}}$ with selecting $s=1$ in \eqref{eq:conf}, and established its validity under the $\beta$-mixing assumption. Using the results in Sections \ref{sec:main_results} and \ref{sec:parametric_boostrap}, we can show the validity of the confidence region defined as \eqref{eq:conf} in more general $\alpha$-mixing setting and physical dependence setting with diverging $s$. Such defined $\mathcal{C}_{\mathcal{S},\delta}$ can be applied for testing the structures and doing support recovering of $\Sigma$ and $\Omega$, respectively. See \cite{ChangQiuYaoZou_2018} for the usefulness of these confidence regions. As we have discussed in Section \ref{ex:test_mean}, involving a larger $s$ in \eqref{eq:conf} can enhance the power performance in finite-samples in comparison to the $\ell^\infty$-type statistic that with $s=1$.



%
%

\section{Proof of main results in Section \ref{sec:hyper}}\label{se:pftm1-3}

In this section, we provide proofs of the high-dimensional CLTs on hyper-rectangles in Section \ref{sec:hyper} under $\alpha$-mixing (Theorem {\rm\ref{tm:1}}), dependency graph (Theorem {\rm\ref{lem:GA_dependency_graph}}), and physical dependence (Theorem {\rm\ref{prop:GA_weakly-dependent}}) frameworks. Such quantitative CLTs on hyper-rectangles are the backbone for deriving CLTs on simple convex sets (Section \ref{sec:simple}) and sparsely convex sets (Section \ref{sec:sparse}). 

\subsection{Proof of Theorem {\rm\ref{tm:1}}}\label{sec:pfalphamix}

To prove Theorem \ref{tm:1}, we need the following lemma which is proved in Lemma C.5 by \cite{ChenKato_2019}. The proof is also implicit in the proof of Lemma C.1 in \cite{Chen_2018}, where a conditional version is given. 
\begin{lemma}
	\label{lem:gaussian_comparison}
	Let $Y$ and $W$ be centered Gaussian random vectors in $\R^{p}$ with covariance matrices $\Sigma_{y}= (\sigma_{j,k}^y)_{j,k\in[p]}$ and $\Sigma_{w} = (\sigma_{j,k}^{w})_{j,k\in[p]}$, respectively. If $\min_{
		j \in[p]} \sigma_{j,j}^{y} \vee \min_{j \in[p]} \sigma_{j,j}^{w} \geq c$ for some universal constant $c> 0$, it then holds that $\sup_{u\in\mathbb{R}^p}|\mathbb{P}(Y\leq u)-\mathbb{P}(W\leq u)|\leq C|\Sigma_{y} - \Sigma_{w} |_{\infty}^{1/3}(\log p)^{2/3}$ for some constant $C>0$ only depending on $c$.
\end{lemma}


Without loss of generality, we let $\log p=o(n^{2/21})$, $\log p=o\{n^{\gamma_{2}/(3+6\gamma_{2})}\}$ and $B_n^{2}(\log p)^{1/\gamma_{2}}=o(n^{1/3})$, since otherwise we can make the assertions trivial.
Let $Q=o(n)$ be a positive integer that will diverge with $n$. We first decompose the sequence $[n]$ to $L+1$ blocks with $L=\lfloor n/Q\rfloor$: $\mathcal{G}_\ell=\{(\ell-1)Q+1,\ldots,\ell Q\}$ for $\ell\in[L]$ and $\mathcal{G}_{L+1}=\{LQ+1,\ldots,n\}$. Additionally, let $b$ and $h$ be two nonnegative integers such that $Q=b+h$ and $h=o(b)$. We decompose each $\mathcal{G}_{\ell}$ for $\ell\in[L]$ to a ``large" block with length $b$ and a ``small" block with length $h$. Specifically, $\mathcal{I}_{\ell}=\{(\ell-1)Q+1,\ldots,(\ell-1)Q+b\}$ and $\mathcal{J}_{\ell}=\{(\ell-1)Q+b+1,\ldots,\ell Q\}$ for $\ell\in[L]$, and
$\mathcal{J}_{L+1}=\mathcal{G}_{L+1}$. Define $\tilde{X}_\ell=b^{-1/2}\sum_{t\in\mathcal{I}_\ell}X_t$ and $\check{X}_\ell=h^{-1/2}\sum_{t\in\mathcal{J}_\ell}X_t$ for $\ell\in[L]$, and $\check{X}_{L+1}=(n-LQ)^{-1/2}\sum_{t\in\mathcal{J}_{L+1}}X_t$. Let $\{Y_t\}_{t=1}^n$ be a sequence of independent normal random vectors with mean zero, where the covariance of $Y_t$ ($t\in\mathcal{I}_\ell$) is $\mathbb{E}(\tilde{X}_\ell\tilde{X}_\ell^\T)$ for each $\ell\in[L]$. Define $\tilde{Y}_\ell=b^{-1/2}\sum_{t\in\mathcal{I}_\ell}Y_t$ for $\ell\in[L]$. Let $S_{n,x}^{(1)}=L^{-1/2}\sum_{\ell=1}^L\tilde{X}_\ell$ and $S_{n,y}^{(1)}=L^{-1/2}\sum_{\ell=1}^L\tilde{Y}_\ell$. Write $\tilde{\Xi}=\cov\{S_{n,y}^{(1)}\}$. It holds that $\tilde{\Xi}=L^{-1}\sum_{\ell=	1}^{L}\mathbb{E}(\tilde{Y}_\ell \tilde{Y}_\ell ^{\T})=L^{-1}\sum_{\ell=	1}^{L}\mathbb{E}(\tilde{X}_\ell \tilde{X}_\ell ^{\T})$. Define
\begin{align}
\varrho_n^{(1)}:=\sup_{u\in\mathbb{R}^p,\nu\in[0,1]}|\mathbb{P}\{\sqrt{\nu}S_{n,x}^{(1)}+\sqrt{1-\nu}S_{n,y}^{(1)}\leq u\}-\mathbb{P}\{S_{n,y}^{(1)}\leq u\}|\,,\label{eq:alpharho1}\\
\varrho_n^{(2)}:=\sup_{u\in\mathbb{R}^p,\nu\in[0,1]}|\mathbb{P}\{\sqrt{\nu}S_{n,x}+\sqrt{1-\nu}S_{n,y}^{(1)}\leq u\}-\mathbb{P}\{S_{n,y}^{(1)}\leq u\}|\,.\label{eq:nnn}
\end{align}
Write $\tilde{X}_\ell=(\tilde{X}_{\ell,1},\ldots,\tilde{X}_{\ell,p})^\T$. We first present the following lemmas. The proof of Lemma \ref{lem:cov_infinity} is almost identical to the proof of Lemma L1 in \cite{ChangJiangShao_2021} with $m=1$ but using the condition $\mathbb{E}(|X_{t,j}|^4)\lesssim B_n^4$ in the steps based on Davydov's inequality \cite[Corollary 2]{Davydov_1968}.  We omit details here. The proofs of Lemmas \ref{pn:1} and \ref{pn:2} are given in Sections \ref{se:pfbd1alpha} and \ref{se:pfpn2} of the supplementary material, respectively. 

\begin{lemma}
	\label{lem:cov_infinity}
	Assume Conditions {\rm\ref{as:tail}--\ref{as:alpha-mixing}} hold. 
	Then $|\tilde{\Xi}-\Xi|_{\infty}\lesssim B_{n}^{2}(hb^{-1}+bn^{-1})$.
\end{lemma}

\begin{lemma}\label{pn:1}
	Assume Conditions {\rm\ref{as:tail}}--{\rm\ref{as:longrun}} hold. Let $\gamma=\gamma_{2}/(2\gamma_{2}+1)$ and $h=C\{\log(pn)\}^{1/\gamma_2}$ for some sufficiently large $C>0$. If $B_n^2h\ll b\ll nB_n^{-2}$, then $\varrho_n^{(1)}\lesssim B_{n}L^{-1/6}(\log p)^{7/6}$ provided that $(\log p)^{3+\gamma}=o(b^{3\gamma/2}L^{\gamma})$ and $\log p=o(L^{2/5})$.
\end{lemma}

\begin{lemma}\label{pn:2}
	Let $p\geq n^\kappa$ for some universal constant $\kappa>0$. Assume Conditions {\rm\ref{as:tail}}--{\rm\ref{as:longrun}} hold. Let $\gamma=\gamma_{2}/(2\gamma_{2}+1)$ and $h=C(\log p)^{1/\gamma_2}$ for some sufficiently large $C>0$. If $b$ satisfies $\min\{nB_n^{-2},n^{1/2}\}\gg b\gg \max\{n^{1/4}(\log p)^{(3-\gamma_2)/(4\gamma_2)},B_n^2h\}$, then $\varrho_n^{(2)}\lesssim B_{n}L^{-1/6}(\log p)^{7/6}$ provided that $(\log p)^{3+\gamma}=o(b^{3\gamma/2}L^{\gamma})$ and $\log p=o(L^{2/5})$.
\end{lemma}


Now we begin to prove Theorem \ref{tm:1}. Let $G\sim {N}(0,\Xi)$ be independent of $\mathcal{X}_n=\{X_1,\ldots,X_n\}$, where $\Xi={\rm Cov}(n^{-1/2}\sum_{t=1}^nX_t)$. Recall
$
\varrho_n=\sup_{u\in\mathbb{R}^p,\nu\in[0,1]}|\mathbb{P}(\sqrt{\nu}S_{n,x}+\sqrt{1-\nu}G\leq u)-\mathbb{P}(G\leq u)|$ and note $S_{n,x}$ is independent of $G$ and $S_{n,y}^{(1)}$. For any $u\in\mathbb{R}^p$, we have
\begin{align*}
&\,\,|\mathbb{P}(\sqrt{\nu}S_{n,x}+\sqrt{1-\nu}G\leq u)-\mathbb{P}\{\sqrt{\nu}S_{n,x}+\sqrt{1-\nu}S_{n,y}^{(1)}\leq u\}|\\
&~~=\bigg|\int\mathbb{P}(\sqrt{1-\nu}G\leq u-\sqrt{\nu}v)\,{\rm d}F_{S_{n,x}}(v)-\int\mathbb{P}\{\sqrt{1-\nu}S_{n,y}^{(1)}\leq u-\sqrt{\nu}v\}\,{\rm d}F_{S_{n,x}}(v)\bigg|\\
&~~\leq\int|\mathbb{P}(\sqrt{1-\nu}G\leq u-\sqrt{\nu}v)-\mathbb{P}\{\sqrt{1-\nu}S_{n,y}^{(1)}\leq u-\sqrt{\nu}v\}|\,{\rm d}F_{S_{n,x}}(v)\\
&~~\leq\sup_{u\in\mathbb{R}^p}|\mathbb{P}(G\leq u)-\mathbb{P}\{S_{n,y}^{(1)}\leq u\}|\,,
\end{align*}
where $F_{S_{n,x}}(\cdot)$ denotes the distribution function of $S_{n,x}$. By triangle inequality, it holds that
\begin{align}\label{eq:apprxga}
\varrho_n\leq&\,\,\sup_{u\in\mathbb{R}^p}|\mathbb{P}(\sqrt{\nu}S_{n,x}+\sqrt{1-\nu}G\leq u)-\mathbb{P}\{\sqrt{\nu}S_{n,x}+\sqrt{1-\nu}S_{n,y}^{(1)}\leq u\}|\notag\\
&+\sup_{u\in\mathbb{R}^p}|\mathbb{P}\{\sqrt{\nu}S_{n,x}+\sqrt{1-\nu}S_{n,y}^{(1)}\leq u\}-\mathbb{P}\{S_{n,y}^{(1)}\leq u\}|\\
&+\sup_{u\in\mathbb{R}^p}|\mathbb{P}\{S_{n,y}^{(1)}\leq u\}-\mathbb{P}(G\leq u)|\notag\\
\leq&\,\,\varrho_n^{(2)}+2\sup_{u\in\mathbb{R}^p}|\mathbb{P}(G\leq u)-\mathbb{P}\{S_{n,y}^{(1)}\leq u\}|\,.\notag
\end{align}
Notice that $S_{n,y}^{(1)}\sim {N}(0,\tilde{\Xi})$ and $G\sim {N}(0,\Xi)$. By Lemmas \ref{lem:gaussian_comparison} and \ref{lem:cov_infinity} with $b=o(n^{1/2})$, 
$$
\sup_{u\in\mathbb{R}^p}|\mathbb{P}\{S_{n,y}^{(1)}\leq u\}-\mathbb{P}(G\leq u)|\lesssim B_n^{2/3}h^{1/3}b^{-1/3}(\log p)^{2/3}\,.$$ Notice that $\log p=o(n^{2/21})$, $\log p=o\{n^{\gamma_{2}/(3+6\gamma_{2})}\}$ and $B_n^{2}(\log p)^{1/\gamma_{2}}=o(n^{1/3})$. Selecting $h\asymp (\log p)^{1/\gamma_{2}}$ and $b\asymp n^{1/3}$, if $(\log p)^{3-\gamma_2}=o(n^{\gamma_2/3})$, together with Lemma \ref{pn:2}, we have
$
\varrho_n\lesssim B_{n}^{2/3}n^{-1/9}(\log p)^{(1+2\gamma_2)/(3\gamma_2)}+B_{n}n^{-1/9}(\log p)^{7/6}$. We construct Theorem \ref{tm:1}. $\hfill\Box$

\subsection{Proof of Theorem \ref{lem:GA_dependency_graph}}\label{sec:5.2}
Without loss of generality, we assume $(D_{n}D_n^*)^2(\log p)^{7}=o(n)$, since otherwise we can make the assertions hold trivially.
Let $\{Y_{t}\}_{t=1}^n$, independent of $\mathcal{X}_n=\{X_1,\ldots,X_n\}$, be a centered Gaussian sequence such that $\Cov(Y_{t}, Y_{s}) = \Cov(X_{t}, X_{s})$ for all $t,s\in[n]$. Define $S_{n,y}=n^{-1/2}\sum_{t=1}^nY_t$. Then $S_{n,y}\overset{d}{=}G\sim {N}(0,\Xi)$. Let $\mathcal{W}_n=\{W_{1},\ldots,W_{n}\}$ be an independent copy of $\mathcal{Y}_n=\{Y_{1},\ldots,Y_{n}\}$ which is also independent of $\mathcal{X}_n$. Analogously, we can define $S_{n,w}$ based on $\mathcal{W}_n$. Write $X_t=(X_{t,1},\ldots,X_{t,p})^\T$, $Y_t=(Y_{t,1},\ldots,Y_{t,p})^\T$ and $W_t=(W_{t,1},\ldots,W_{t,p})^\T$. Recall $\mathcal{N}_{t} = \{s\in V_n: (t,s) \in E_{n}\}$ for any $t\in[n]$. For $\phi \geq 1$, define
\begin{align}\label{eq:Mphixm}
M_{n,x}(\phi) =&\,\,\max_{t\in[n]}\mathbb{E} \bigg\{ \max_{j\in[p],s \in \cup_{\ell\in\mathcal{N}_t}\mathcal{N}_{\ell}} |X_{s,j}|^{3}\notag\\
&\quad\quad\quad\quad\quad\quad\quad\quad \times I\bigg( \max_{j \in[p],s \in \cup_{\ell\in\mathcal{N}_t}\mathcal{N}_{\ell}} |X_{s,j}| > {\sqrt{n} \over 8(D_n^4D_{n}^*)^{1/3} \phi \log{p}} \bigg) \bigg\}\,,\\
\tilde{M}_{n,x}(\phi)=&\,\,\max_{t\in[n]}\max_{s\in(\cup_{\ell\in\mathcal{N}_t}\mathcal{N}_\ell)\backslash\mathcal{N}_t}\mathbb{E}\bigg\{\max_{j\in[p],s'\in\cup_{\ell\in\{s\}\cup\mathcal{N}_t}\mathcal{N}_\ell} |X_{s',j}|^{3}\notag\\
&\quad\quad\quad\quad\quad\quad\quad \times I\bigg(\max_{j\in[p],s'\in\cup_{\ell\in\{s\}\cup\mathcal{N}_t}\mathcal{N}_\ell} |X_{s',j}| > {\sqrt{n} \over 8(D_n^{4}D_{n}^*)^{1/3}\phi\log p} \bigg) \bigg\}\,.\notag
\end{align}
Similarly, we define $M_{n,y}(\phi)$ and $\tilde{M}_{n,y}(\phi)$ in the same manner with $X_{s}$ replaced by $Y_{s}$. Set $M_{n}(\phi) = M_{n,x}(\phi) + M_{n,y}(\phi)$ and $\tilde{M}_{n}(\phi) = \tilde{M}_{n,x}(\phi) + \tilde{M}_{n,y}(\phi)$.

Let $\beta=\phi\log p$. For a given $u=(u_1,\ldots,u_p)^\T\in\mathbb{R}^p$, define
$
F_\beta(v)=\beta^{-1}\log[\sum_{j=1}^p\exp\{\beta(v_j-u_j)\}]$
for any $v=(v_1,\ldots,v_p)^\T\in\mathbb{R}^p$. Such defined function $F_{\beta}(v)$ satisfies the property $0\leq F_{\beta}(v)-\max_{j\in[p]}(v_j-u_j)\leq \beta^{-1}\log p=\phi^{-1}$ for any $v\in\mathbb{R}^p$. Select a thrice continuously differentiable function $g_0:\mathbb{R}\rightarrow[0,1]$ whose derivatives up to the third order are all bounded such that $g_0(t)=1$ for $t\leq 0$ and $g_0(t)=0$ for $t\geq1$. Define $g(t):=g_0(\phi t)$ for any $t\in\mathbb{R}$, and $q(v):=g\{F_\beta(v)\}$ for any $v\in\mathbb{R}^p$. Define
$$
\mathcal{T}_{n} := q(\sqrt{\nu}S_{n,x} + \sqrt{1-\nu}S_{n,y}) - q(S_{n,w})\,.$$
Using the same arguments stated in Section \ref{se:pfbd1alpha} of the supplementary material, we have
$$
\varrho_n
\lesssim\phi^{-1}(\log p)^{1/2}+\sup_{u\in\mathbb{R}^p,\nu\in[0,1]}|\mathbb{E}(\mathcal{T}_n)|\,.$$ 
To specify the convergence rate of $\varrho_n$, we only need to bound $\sup_{u\in\mathbb{R}^p,\nu\in[0,1]}|\mathbb{E}(\mathcal{T}_n)|$. Define $Z(\omega) = \sum_{t=1}^{n}Z_{t}(\omega)$ for any $\omega\in[0,1]$, where
$$
Z_{t}(\omega) = n^{-1/2}\{ \sqrt{\omega} (\sqrt{\nu}X_{t} + \sqrt{1-\nu}Y_{t}) + \sqrt{1-\omega}W_t \}\,.$$
Then $Z(1) = \sqrt{\nu}S_{n,x} + \sqrt{1-\nu}S_{n,y}$ and $Z(0) =S_{n,w}$. For notational simplicity, we omit the dependence of $Z(\omega)$ and $Z_{t}(\omega)$ on $\omega$ in the rest of this proof. Let $\delta_{t} = \sum_{s\in\mathcal{N}_t}Z_s=(\delta_{t,1},\ldots,\delta_{t,p})^\T$, $Z^{(-t)} = Z - \delta_{t}$, and
$$\dot{Z}_{t} = n^{-1/2}\{ \omega^{-1/2}(\sqrt{\nu}X_{t} + \sqrt{1-\nu}Y_{t}) - (1-\omega)^{-1/2}W_{t}\}=(\dot{Z}_{t,1},\ldots,\dot{Z}_{t,p})^\T\,.$$
By Taylor expansion, we have 
$
2\mathbb{E}(\mathcal{T}_{n})={\rm I}+{\rm II}+{\rm III}$, where 
\begin{align*}
{\rm I}=&\,\,\sum_{j=1}^{p} \sum_{t=1}^{n} \int_{0}^{1} \mathbb{E}[\partial_{j}q\{Z^{(-t)}\} \dot{Z}_{t,j}]\,{\rm d}\omega\,,~{\rm II}=\sum_{j,k=1}^{p} \sum_{t=1}^{n} \int_{0}^{1} \mathbb{E}[\partial_{jk}q\{Z^{(-t)}\} \delta_{t,k} \dot{Z}_{t,j}]\,{\rm d}\omega\,,\\
{\rm III}=&\,\,\sum_{j,k,l=1}^{p} \sum_{t=1}^{n} \int_{0}^{1} \int_{0}^{1} (1-\tau) \mathbb{E}[\partial_{jkl}q\{Z^{(-t)} + \tau \delta_{t}\} \delta_{t,k} \delta_{t,l} \dot{Z}_{t,j}]\,{\rm d}\tau {\rm d}\omega\,.
\end{align*} 
Since $Z^{(-t)}$ and $\dot{Z}_{t,j}$ are independent, we have ${\rm I} = 0$. The following two lemmas give the upper bounds for ${\rm II}$ and ${\rm III}$, respectively, whose proofs are given in Sections \ref{sec:pfmIII} and \ref{se:pflemmll} of the supplementary material, respectively.

\begin{lemma}\label{lem:mIII}
	 $
	|{\rm III}| \lesssim D_{n}^{2} \phi^3 n^{-1/2}(\log p)^{2}\{ M_{n}(\phi) + B_{n}^{3}\phi^{-1}(\log{p})^{1/2} + B_{n}^{3}\varrho_{n}\}$.
\end{lemma}

\begin{lemma}\label{lem:mII}
	 $|{\rm II}|\lesssim {D_{n}D_n^* \phi^3n^{-1/2}}(\log p)^2\{ \tilde{M}_{n}(\phi) + B_{n}^{3}\phi^{-1}(\log{p})^{1/2} +B_{n}^{3} \varrho_{n}\}$.
\end{lemma}

Hence, by Lemmas \ref{lem:mIII} and \ref{lem:mII}, we have
\begin{align*}
\varrho_n\lesssim \phi^{-1}(\log p)^{1/2}+&D_{n}\phi^3n^{-1/2}(\log p)^{2}\{ D_nM_{n}(\phi)+ D_n^*\tilde{M}_n(\phi) \\
&\quad\quad\quad\quad\quad\quad\quad\quad\quad+ B_{n}^{3}D_n^*\phi^{-1}(\log{p})^{1/2}
+B_{n}^{3} D_n^*\varrho_{n}\}\,.
\end{align*}
Taking $\phi=C'n^{1/6}B_{n}^{-1}(D_nD_n^*)^{-1/3}(\log p)^{-2/3}$ for some sufficiently small $C'>0$, then 
\begin{align*}
\varrho_n\lesssim&\,\, B_{n}^{-3}M_n\{C'n^{1/6}B_{n}^{-1}(D_nD_n^*)^{-1/3}(\log p)^{-2/3}\}\\
&+B_{n}^{-3}\tilde{M}_n\{C'n^{1/6}B_{n}^{-1}(D_nD_n^*)^{-1/3}(\log p)^{-2/3}\}+n^{-1/6}B_{n}(D_nD_n^{*})^{1/3}(\log p)^{7/6}\,.
\end{align*}
Write $$X^{(t)}=\max_{j\in[p],s\in\cup_{\ell\in\mathcal{N}_t}\mathcal{N}_\ell} |X_{s,j}|~\textrm{and}~X^{(t),s}=\max_{j\in[p],s'\in\cup_{\ell\in\{s\}\cup\mathcal{N}_t}\mathcal{N}_\ell} |X_{s',j}|\,.$$ By Condition {\rm\ref{as:tail}}, 
$\mathbb{P}\{X^{(t)}>u\}\leq 2pD_n^*\exp(-u^{\gamma_1}B_{n}^{-\gamma_{1}})$ and $\mathbb{P}\{X^{(t),s}>u\}\leq 4pD_n^*\exp(-u^{\gamma_1}B_{n}^{-\gamma_{1}})$ for any $u>0$. 
Notice that $D_n\leq D_n^*\leq D_n^2$ and $\mathbb{E}\{|\xi|^3I(|\xi|>v)\}=v^3\mathbb{P}(|\xi|>v)+3\int_v^\infty u^2\mathbb{P}(|\xi|>u)\,{\rm d}u$. Due to $D_n^{3}(\log p)^{1+3/\gamma_1}=o(n)$ and $p\geq n^\kappa$ for some $\kappa>0$, we have
\begin{align*}
&\mathbb{E}[\{X^{(t)}\}^{3}I\{X^{(t)}> B_{n}n^{1/3}(8C')^{-1}\cdot D_{n}^{-1}(\log p)^{-1/3}\}]\\
&~~~~~~\lesssim B_{n}^{3}\exp\{-Cn^{\gamma_1/3}D_n^{-\gamma_1}(\log p)^{-\gamma_1/3}\}\lesssim B_{n}^{3}D_n^{1/2}n^{-1/6}(\log p)^{7/6}\,,
\end{align*}
 which implies   
$$
M_{n,x}\{C'n^{1/6}B_{n}^{-1}(D_nD_n^*)^{-1/3}(\log p)^{-2/3}\}
\lesssim n^{-1/6}B_{n}^{3}D_n^{1/2}(\log p)^{7/6}\,.$$ 
Recall $Y_t$ is a normal random vector. Since $\max_{t\in [n],j\in [p]}\mathbb{E}(|X_{t,j}|^{2})\lesssim B_{n}^{2}$ and $\Cov(Y_{t},Y_{s})=\Cov(X_{t},X_{s})$ for all $t,s\in[n]$, then $\max_{t\in[n],j\in [p]}\mathbb{P}(|Y_{t,j}|>u)\leq 2\exp(-Cu^2B_{n}^{-2})$ for any $u>0$.  Analogously,
$$
M_{n,y}\{C'n^{1/6}B_{n}^{-1}(D_nD_n^*)^{-1/3}(\log p)^{-2/3}\}\lesssim B_{n}^{3}D_n^{1/2}n^{-1/6}(\log p)^{7/6}$$
provided that $D_n^{3}(\log p)^{5/2}=o(n)$. Notice that $(D_{n}D_n^{*})^2(\log p)^{7}=o(n)$ and $\gamma_1\geq1$. Thus, $$
M_n\{C'n^{1/6}B_{n}^{-1}(D_nD_n^*)^{-1/3}(\log p)^{-2/3}\}\lesssim B_{n}^{3}D_n^{1/2}n^{-1/6}(\log p)^{7/6}\,.$$ 
By the same arguments, we can also show  
$$
\mathbb{E}[\{X^{t,(s)}\}^{3}I\{X^{t,(s)}> (8C')^{-1}B_{n}n^{1/3}D_{n}^{-1}(\log p)^{-1/3}\}]\lesssim B_{n}^{3}D_n^{1/2}n^{-1/6}(\log p)^{7/6}\,,$$ 
which implies $$
\tilde{M}_{n}\{C'n^{1/6}B_{n}^{-1}(D_nD_n^*)^{-1/3}(\log p)^{-2/3}\}\lesssim B_{n}^{3}D_n^{1/2}n^{-1/6}(\log p)^{7/6}\,.$$ Hence, 
$\varrho_n\lesssim n^{-1/6} B_n(D_nD_n^*)^{1/3}(\log p)^{7/6}$. We complete the proof of Theorem \ref{lem:GA_dependency_graph}. $\hfill\Box$

\subsection{Proof of Theorem {\rm\ref{prop:GA_weakly-dependent}}}\label{sec:pflma5}

Let $X_{t}^{(m)} = \mathbb{E}(X_{t}\,|\,\varepsilon_{t}, \dots, \varepsilon_{t-m})$ for any $m\geq1$. Then $\{X_{t}^{(m)}\}_{t=1}^{n}$ is an $m$-dependent sequence with mean zero. Let $\Xi^{(m)}=\cov\{S_{n,x}^{(m)}\}$ with $S_{n,x}^{(m)} = n^{-1/2} \sum_{t=1}^{n}X_{t}^{(m)}$. Recall $S_{n,x}=n^{-1/2}\sum_{t=1}^{n}X_{t}$. Write $S_{n,x}=(S_{n,x,1},\ldots,S_{n,x,p})^\T$ and $S_{n,x}^{(m)}=\{S_{n,x,1}^{(m)},\ldots,S_{n,x,p}^{(m)}\}^\T$.
\begin{lemma}
	\label{lem:tail_prob_bound_m-dep_approx}
	Let  $\{X_{t}\}$ be a sequence of centered random vectors generated from the model {\rm(\ref{eqn:weakly_dependent_ts_fdp})} such that $\Phi_{\psi_{\nu},\alpha} < \infty$ for some $\alpha, \nu \in (0, \infty)$. Then there exists a universal constant $C > 0$ depending only on $\nu$ such that
	$
	\max_{j\in[p]}\mathbb{P}\{|S_{n,x,j}^{(m)} - S_{n,x,j}| > u\} \leq C \exp\{-(4e)^{-1}(1+2\nu)(um^{\alpha}\Phi_{\psi_{\nu},\alpha}^{-1})^{2/(1+2\nu)}\}$
	for any $u>0$.
\end{lemma}

\begin{lemma}
	\label{lem:moment_bound_cov_mat_compare}
	Let $q \geq 2$. For each $j\in[p]$, it holds that $\|S_{n,x,j}\|_{q} \leq (q-1)^{1/2}\Theta_{0,q,j}$, $\|S_{n,x,j}^{(m)}\|_{q} \leq (q-1)^{1/2}\Theta_{0,q,j}$ and $\|S_{n,x,j}-S_{n,x,j}^{(m)}\|_{q} \leq (q-1)^{1/2}\Theta_{m+1,q,j}$.
\end{lemma}

The proof of Lemma \ref{lem:tail_prob_bound_m-dep_approx} essentially follows from the arguments in proving Lemma C.3 of \cite{ZhangWu_2017} with the necessary modification using the uniform functional dependence measure to non-stationarity of the sequence $\{X_{t}\}$. Details are omitted. The proof of Lemma \ref{lem:moment_bound_cov_mat_compare} is given in Section \ref{se:pflemmoment_bound_cov_mat_compare} of the supplementary material.

\subsubsection{Proof of Part {\rm(i)} of Theorem {\rm\ref{prop:GA_weakly-dependent}}}\label{sec:pfphysical_m}
Recall $\Xi=\cov(S_{n,x})$. We will apply the large-and-small-blocks technique stated in Appendix \ref{sec:pfalphamix} to derive the upper bound of $\varrho_n$. Without loss of generality, we assume $\Phi_{\psi_{\nu},0}=o\{n^{\alpha/(3+9\alpha)}\}$ and $\Psi_{2,\alpha}\Psi_{2,0}=o\{n^{\alpha/(1+3\alpha)}\}$, since otherwise the assertions hold trivially. Let $Q=o(n)$ be a positive integer that will diverge with $n$. We first decompose the sequence $\{X_{t}^{(m)}\}_{t=1}^{n}$ to $L+1$ blocks with $L=\lfloor n/Q\rfloor$: $\mathcal{G}_\ell=\{(\ell-1)Q+1,\ldots,\ell Q\}$ for $\ell\in[L]$ and $\mathcal{G}_{L+1}=\{LQ+1,\ldots,n\}$. Let $b\gg m$ be a nonnegative integer such that $Q=b+m$. We decompose each $\mathcal{G}_{\ell}$ for $\ell\in[L]$ to a ``large" block with length $b$ and a ``small" block with length $m$. Specifically, $\mathcal{I}_{\ell}=\{(\ell-1)Q+1,\ldots,(\ell-1)Q+b\}$ and $\mathcal{J}_{\ell}=\{(\ell-1)Q+b+1,\ldots,\ell Q\}$ for $\ell\in[L]$, and
$\mathcal{J}_{L+1}=\mathcal{G}_{L+1}$. Define $\tilde{X}_\ell^{(m)}=b^{-1/2}\sum_{t\in\mathcal{I}_\ell}X_t^{(m)}$ and $\check{X}_\ell^{(m)}=m^{-1/2}\sum_{t\in\mathcal{J}_\ell}X_t^{(m)}$ for $\ell\in[L]$, and $\check{X}_{L+1}^{(m)}=(n-LQ)^{-1/2}\sum_{t\in\mathcal{J}_{L+1}}X_t^{(m)}$. Since $\{X_t^{(m)}\}_{t=1}^n$ is an $m$-dependent sequence, we know $\{\tilde{X}_\ell^{(m)}\}_{\ell=1}^{L}$ is an independent sequence.
Let $\{Y_t^{(m)}\}_{t=1}^n$ be a sequence of independent normal random vectors with mean zero, where the covariance of $Y_t^{(m)}$ ($t\in\mathcal{I}_\ell$) is $\mathbb{E}[\tilde{X}_\ell^{(m)}\{\tilde{X}_\ell^{(m)}\}^\T]$ for each $\ell\in[L]$. Define $\tilde{Y}_\ell^{(m)}=b^{-1/2}\sum_{t\in\mathcal{I}_\ell}Y_t^{(m)}$ for $\ell\in[L]$. Let $\tilde S_{n,x}^{(m)}=L^{-1/2}\sum_{\ell=1}^L\tilde{X}_\ell^{(m)}$ and $\tilde S_{n,y}^{(m)}=L^{-1/2}\sum_{\ell=1}^L\tilde{Y}_\ell^{(m)}$. Write $\tilde{\Xi}=\cov\{\tilde S_{n,y}^{(m)}\}$. Then $\tilde{\Xi}=L^{-1}\sum_{\ell=	1}^{L}\mathbb{E}[\tilde{Y}_\ell^{(m)} \{\tilde{Y}_\ell^{(m)}\}^{\T}]=L^{-1}\sum_{\ell=	1}^{L}\mathbb{E}[\tilde{X}_\ell^{(m)} \{\tilde{X}_\ell ^{(m)}\}^{\T}]=\cov\{\tilde S_{n,x}^{(m)}\}$.
Define
\begin{align*}
\varrho_{n,1}^{(m)}&:=\sup_{u\in\mathbb{R}^p,\nu\in[0,1]}|\mathbb{P}\{\sqrt{\nu}\tilde S_{n,x}^{(m)}+\sqrt{1-\nu}\tilde S_{n,y}^{(m)}\leq u\}-\mathbb{P}\{\tilde S_{n,y}^{(m)}\leq u\}|\,,\\
\varrho_{n,2}^{(m)}&:=\sup_{u\in\mathbb{R}^p,\nu\in[0,1]}|\mathbb{P}\{\sqrt{\nu}S_{n,x}^{(m)}+\sqrt{1-\nu}\tilde S_{n,y}^{(m)}\leq u\}-\mathbb{P}\{\tilde S_{n,y}^{(m)}\leq u\}|\,.
\end{align*}
Write $\tilde{X}_\ell^{(m)}=\{\tilde{X}_{\ell,1}^{(m)},\ldots,\tilde{X}_{\ell,p}^{(m)}\}^\T$. We first present the following lemmas whose proofs are given in Sections \ref{sec:pfltpbm}--\ref{sec:pflpn2m} of the supplementary material, respectively.
\begin{lemma}
	\label{lem:tail_prob_bound_m}
	If $\Phi_{\psi_{\nu},0} < \infty$ for some $\nu\in(0,\infty)$, we have 
	$
	\max_{\ell\in [L]}\max_{j\in[p]}\mathbb{P}\{|\tilde X_{\ell,j}^{(m)}| > u\} \lesssim \exp\{-(4e)^{-1}(1+2\nu)(u\Phi_{\psi_{\nu},0}^{-1})^{2/(1+2\nu)}\}$
	for any $u>0$.
\end{lemma}

\begin{lemma}
	\label{lem:mom_m}
	If $\Phi_{\psi_{\nu},0} < \infty$ for some $\nu\in(0,\infty)$, we have  $\max_{\ell\in[L]}\max_{j\in [p]}\mathbb{E}\{|\tilde{X}_{\ell,j}^{(m)}|^{q}\} \lesssim \Phi_{\psi_{\nu,0}}^q$ for any positive integer $q\geq1$.
\end{lemma}

\begin{lemma}
	\label{lem:cov_infinity_m}
	If $\Phi_{\psi_{\nu},0}, \Psi_{2,\alpha} < \infty$ for some $\nu,\alpha\in(0,\infty)$, it holds that $|\tilde\Xi-\Xi|_\infty\lesssim \Phi_{\psi_{\nu},0}^2(mb^{-1}+bn^{-1})+m^{-\alpha} \Psi_{2,\alpha} \Psi_{2,0}$.
\end{lemma}

\begin{lemma}\label{pn:1_m}
Let $p\geq n^\kappa$ for some universal constant $\kappa>0$. 	Assume $\Phi_{\psi_{\nu},0},\Psi_{2,\alpha}  < \infty$ for some $\nu,\alpha\in(0,\infty)$, and $\min_{j\in[p]}V_{n,j}\geq C$ for some universal constant $C>0$, where $V_{n,j}$ is defined in \eqref{eq:Vtj}. If $\Phi_{\psi_{\nu},0}^2m\ll b\ll \Phi_{\psi_{\nu},0}^{-2}n$ and $m^{\alpha}\gg \Psi_{2,\alpha}\Psi_{2,0}$, then $\varrho_{n,1}^{(m)}\lesssim L^{-1/6}\Phi_{\psi_{\nu,0}}(\log p)^{7/6}$ provided that $(\log p)^{5+6\nu}=o(L^2)$.
\end{lemma}

\begin{lemma}\label{pn:2_m}
	Let $p\geq n^\kappa$ for some universal constant $\kappa>0$. Assume $\Phi_{\psi_{\nu},0},\Psi_{2,\alpha}  < \infty$ for some $\nu,\alpha\in(0,\infty)$, and $\min_{j\in[p]}V_{n,j}\geq C$ for some universal constant $C>0$, where $V_{n,j}$ is defined in \eqref{eq:Vtj}. If $\max\{m\Phi_{\psi_{\nu},0}^2,n^{1/4}m^{3/4}(\log p)^{(6\nu-1)/4}\}\ll b\ll \min\{n\Phi_{\psi_{\nu},0}^{-2},(mn)^{1/2}\}$ and $m^{\alpha}\gg \Psi_{2,\alpha}\Psi_{2,0}$, then $\varrho_{n,2}^{(m)}\lesssim L^{-1/6}\Phi_{\psi_{\nu,0}}(\log p)^{7/6}$ provided that $(\log p)^{5+6\nu}=o(L^2)$.
\end{lemma}

Recall $G\sim N(0,\Xi)$ and $\tilde{S}_{n,y}^{(m)}\sim N(0,\tilde{\Xi})$. By Lemma \ref{lem:gaussian_comparison}, $$\sup_{u\in\mathbb{R}^p}|\mathbb{P}(G\leq u)-\mathbb{P}\{\tilde S_{n,y}^{(m)}\leq u\}|\lesssim |\tilde\Xi-\Xi|_\infty^{1/3}(\log p)^{2/3}\,.$$ Let $D_n=C m^{-\alpha}\Phi_{\psi_{\nu},\alpha}(\log p)^{(1+2\nu)/2}$ for some sufficiently large constant $C>0$. Then $D_n(\log p)^{1/2}\lesssim m^{-\alpha}\Phi_{\psi_{\nu},\alpha}(\log p)^{1+\nu}$. Define the event $\mathcal{E}=\{|S_{n,x}-S_{n,x}^{(m)}|_\infty\leq D_n\}$. If $m^{\alpha}\gg \Psi_{2,\alpha}\Psi_{2,0}$ and $\max\{m\Phi_{\psi_{\nu},0}^2,n^{1/4}m^{3/4}(\log p)^{(6\nu-1)/4}\}\ll b\ll \min\{n\Phi_{\psi_{\nu},0}^{-2},(mn)^{1/2}\}$, then 
\begin{align}\label{eq:bdphys_m}
\varrho_n\leq&\,\sup_{u\in\mathbb{R}^p,\nu\in[0,1]}|\mathbb{P}\{\sqrt{\nu}S_{n,x}+\sqrt{1-\nu}\tilde S_{n,y}^{(m)}\leq u\}-\mathbb{P}\{\tilde S_{n,y}^{(m)}\leq u\}|\notag\\
&+2\sup_{u\in\mathbb{R}^p}|\mathbb{P}(G\leq u)-\mathbb{P}\{\tilde S_{n,y}^{(m)}\leq u\}|\notag\\
\lesssim&\,\,\varrho_{n,2}^{(m)}+\mathbb{P}(\mathcal{E}^c)+\sup_{u\in\mathbb{R}^p,\nu\in[0,1]}|\mathbb{P}\{\tilde S_{n,y}^{(m)}\leq u-\sqrt{\nu}D_n\}-\mathbb{P}\{\tilde S_{n,y}^{(m)}\leq u\}|\\
&+\sup_{u\in\mathbb{R}^p,\nu\in[0,1]}|\mathbb{P}\{\tilde S_{n,y}^{(m)}\leq u+\sqrt{\nu}D_n\}-\mathbb{P}\{\tilde S_{n,y}^{(m)}\leq u\}|+|\tilde\Xi-\Xi|_\infty^{1/3}(\log p)^{2/3}\notag\\
\lesssim&\,\, L^{-1/6}\Phi_{\psi_\nu,0}(\log p)^{7/6}+D_n(\log p)^{1/2}+\mathbb{P}(\mathcal{E}^c)+|\tilde\Xi-\Xi|_\infty^{1/3}(\log p)^{2/3}\,,\notag
\end{align}
provided that $(\log p)^{5+6\nu}=o(L^2)$, where the first step is identical to \eqref{eq:apprxga}, and the last step is based on Lemma \ref{pn:2_m} and Nazarov's inequality.
Lemma \ref{lem:tail_prob_bound_m-dep_approx} implies 
$$\mathbb{P}(\mathcal{E}^c)\lesssim p\exp\{-C(D_nm^\alpha\Phi_{\psi_{\nu},\alpha}^{-1})^{2/(1+2\nu)}\}\lesssim L^{-1/6}\Phi_{\psi_\nu,0}(\log p)^{7/6}\,.$$ Due to $b\ll(mn)^{1/2}$, by Lemma \ref{lem:cov_infinity_m}, $|\tilde\Xi - \Xi|_\infty \lesssim \Phi_{\psi_{\nu},0}^2mb^{-1}+m^{-\alpha} \Psi_{2,\alpha} \Psi_{2,0} $.  By \eqref{eq:bdphys_m} and recalling $L\asymp nb^{-1}$, we have 
\begin{align*}\varrho_n\lesssim&\,\, b^{1/6}n^{-1/6}\Phi_{\psi_\nu,0}(\log p)^{7/6}+m^{-\alpha}\Phi_{\psi_{\nu},\alpha}(\log p)^{1+\nu}\\
&~~~~+\Phi_{\psi_{\nu},0}^{2/3}m^{1/3}b^{-1/3}(\log p)^{2/3}+\Psi_{2,\alpha}^{1/3} \Psi_{2,0}^{1/3}m^{-\alpha/3}(\log p)^{2/3}\,.
\end{align*}
With selecting $b\asymp m^{2/3}n^{1/3}$, it then holds that 
\begin{align*}
\varrho_n
\lesssim m^{1/9}n^{-1/9}\Phi_{\psi_\nu,0}(\log p)^{7/6}+m^{-\alpha}\Phi_{\psi_{\nu},\alpha}(\log p)^{1+\nu}+\Psi_{2,\alpha}^{1/3} \Psi_{2,0}^{1/3}m^{-\alpha/3}(\log p)^{2/3}
\end{align*}
provided that $\Psi_{2,\alpha}^{1/\alpha}\Psi_{2,0}^{1/\alpha}\ll m\ll\min\{n\Phi_{\psi_\nu,0}^{-6},n(\log p)^{3(1-6\nu)},n(\log p)^{-3(5+6\nu)/4}\}$. Recall $\Phi_{\psi_{\nu},0}=o\{n^{\alpha/(3+9\alpha)}\}$ and $\Psi_{2,\alpha}\Psi_{2,0}=o\{n^{\alpha/(1+3\alpha)}\}$. Letting $m\asymp n^{1/(1+3\alpha)}$, we have $$\varrho_n\lesssim n^{-\alpha/(3+9\alpha)}(\log p)^{2/3}\{\Phi_{\psi_\nu,0}(\log p)^{1/2}+\Psi_{2,\alpha}^{1/3} \Psi_{2,0}^{1/3}\}+n^{-\alpha/(1+3\alpha)}\Phi_{\psi_\nu,\alpha}(\log p)^{1+\nu}$$ provided that $(\log p)^{6\nu-1}= o\{n^{\alpha/(1+3\alpha)}\}$ and $(\log p)^{(5+6\nu)/4}= o\{n^{\alpha/(1+3\alpha)}\}$. We have Part (i) of Theorem \ref{prop:GA_weakly-dependent}. $\hfill\Box$

\subsubsection{Proof of Part {\rm(ii)} of Theorem {\rm\ref{prop:GA_weakly-dependent}}}\label{sec:pfphysical}
Without loss of generality, we assume $\Phi_{\psi_{\nu},\alpha}^{4+2\alpha}=o(n^{\alpha})$, since otherwise the assertions hold trivially. Let $\{Y_{t}^{(m)}\}_{t=1}^n$ be a sequence of centered Gaussian random vectors in $\R^{p}$ such that $\Cov\{Y_{t}^{(m)}, Y_{s}^{(m)}\}  = \Cov\{X_{t}^{(m)}, X_{s}^{(m)}\}$ for all $t,s\in[n]$. Set $S_{n,y}^{(m)} = n^{-1/2} \sum_{t=1}^{n}Y_{t}^{(m)}$. Recall $S_{n,x}=n^{-1/2}\sum_{t=1}^nX_t$. Lemma \ref{lem:tail_prob_bound_m-dep_approx} implies 
$
\max_{j\in[p]}\mathbb{E}\{|S_{n,x,j}-S_{n,x,j}^{(m)}|^{2}\}
\lesssim \Phi_{\psi_{\nu},\alpha}^{2}{m}^{-2\alpha}
$. Since $\Var\{S_{n,x,j}^{(m)}\}\geq \Var(S_{n,x,j})+\mathbb{E}\{|S_{n,x,j}-S_{n,x,j}^{(m)}|^2\}-2\{\Var(S_{n,x,j})\}^{1/2}[\mathbb{E}\{|S_{n,x,j}-S_{n,x,j}^{(m)}|^2\}]^{1/2}$, if we select $m\geq C\Phi_{\psi_{\nu},\alpha}^{1/\alpha}$ for some sufficiently large $C>0$, we know that $\Var\{S_{n,x,j}^{(m)}\}$ is uniformly bounded away from zero. By H\"older's inequality, Jensen's inequality and Condition {\rm\ref{as:tail}}, $\|X_{t,j}^{(m)}\|_{\psi_{\gamma_{1}}}\leq \|X_{t,j}\|_{\psi_{\gamma_{1}}}\leq B_{n}$ for any  $t\in [n]$ and $j\in [p]$. Thus Corollary \ref{cor:GA_m-dependent_subexp} yields that \begin{align*}
\varrho_n\{S_{n,x}^{(m)},S_{n,y}^{(m)}\}:=&\,\,\sup_{u\in\mathbb{R}^p,\nu\in[0,1]}|\mathbb{P}\{\sqrt{\nu}S_{n,x}^{(m)}+\sqrt{1-\nu}S_{n,y}^{(m)}\leq u\}-\mathbb{P}\{S_{n,y}^{(m)}\leq u\}|\\
\lesssim&\,\, n^{-1/6}B_{n}m^{2/3}(\log p)^{7/6}\,.
\end{align*}
Recall $G\sim N(0,\Xi)$ with $\Xi=\cov(S_{n,x})$. Let $\tilde D_n=C m^{-\alpha}\Phi_{\psi_{\nu},\alpha}(\log p)^{(1+2\nu)/2}$ for some sufficiently large constant $C>0$. Then $\tilde D_n(\log p)^{1/2}\lesssim m^{-\alpha}\Phi_{\psi_{\nu},\alpha}(\log p)^{1+\nu}$. Consider the event $\mathcal{E}=\{|S_{n,x}-S_{n,x}^{(m)}|_\infty\leq \tilde{D}_n\}$. Write $\Xi^{(m)}=\cov\{S_{n,x}^{(m)}\}$. By Lemma \ref{lem:gaussian_comparison}, $\sup_{u\in\mathbb{R}^p}|\mathbb{P}(G\leq u)-\mathbb{P}\{S_{n,y}^{(m)}\leq u\}|\lesssim |\Xi^{(m)}-\Xi|_\infty^{1/3}(\log p)^{2/3}$. Then
\begin{align*}
\varrho_n\leq&\,\sup_{u\in\mathbb{R}^p,\nu\in[0,1]}|\mathbb{P}\{\sqrt{\nu}S_{n,x}+\sqrt{1-\nu}S_{n,y}^{(m)}\leq u\}-\mathbb{P}\{S_{n,y}^{(m)}\leq u\}|\notag\\
&+2\sup_{u\in\mathbb{R}^p}|\mathbb{P}(G\leq u)-\mathbb{P}\{S_{n,y}^{(m)}\leq u\}|\notag\\
\leq&\,\,\varrho_n\{S_{n,x}^{(m)},S_{n,y}^{(m)}\}+\mathbb{P}(\mathcal{E}^c)+\sup_{u\in\mathbb{R}^p,\nu\in[0,1]}|\mathbb{P}\{S_{n,y}^{(m)}\leq u-\sqrt{\nu}\tilde{D}_n\}-\mathbb{P}\{S_{n,y}^{(m)}\leq u\}|\\
&+\sup_{u\in\mathbb{R}^p,\nu\in[0,1]}|\mathbb{P}\{S_{n,y}^{(m)}\leq u+\sqrt{\nu}\tilde{D}_n\}-\mathbb{P}\{S_{n,y}^{(m)}\leq u\}|+|\Xi^{(m)}-\Xi|_\infty^{1/3}(\log p)^{2/3}\notag\\
\lesssim&\,\,n^{-1/6}B_{n}m^{2/3}(\log p)^{7/6}+\tilde{D}_n(\log p)^{1/2}+\mathbb{P}(\mathcal{E}^c)+|\Xi^{(m)}-\Xi|_\infty^{1/3}(\log p)^{2/3}\,,\notag
\end{align*}
where the first step is identical to \eqref{eq:apprxga}, and the last step is based on Nazarov's inequality.
By Lemma \ref{lem:tail_prob_bound_m-dep_approx}, 
$\mathbb{P}(\mathcal{E}^c)\lesssim p\exp\{-C(\tilde{D}_nm^\alpha\Phi_{\psi_{\nu},\alpha}^{-1})^{2/(1+2\nu)}\}\lesssim n^{-1/6}B_{n}m^{2/3}(\log p)^{7/6}$.
It follows from Cauchy-Schwarz inequality and Lemma \ref{lem:moment_bound_cov_mat_compare} that
\begin{align*}
&|\mathbb{E}(S_{n,x,j}S_{n,x,k})- \mathbb{E}\{S_{n,x,j}^{(m)}S_{n,x,k}^{(m)}\}|\\ &~~~\leq|\mathbb{E}[S_{n,x,j}\{S_{n,x,k}-S_{n,x,k}^{(m)}\}]| + |\mathbb{E}\{S_{n,x,k}^{(m)}\{S_{n,x,j}-S_{n,x,j}^{(m)}\}]|\\
&~~~\leq\|S_{n,x,j}\|_{2} \|S_{n,x,k}-S_{n,x,k}^{(m)}\|_{2} + \|S_{n,x,k}^{(m)}\|_{2}\|S_{n,x,j}-S_{n,x,j}^{(m)}\|_{2}\\
&~~~\leq \Theta_{0,2,j}\Theta_{m+1,2,k}+\Theta_{0,2,k}\Theta_{m+1,2,j}\,.
\end{align*}
Then  
$
|\Xi^{(m)} - \Xi|_\infty \leq2(\max_{j\in[p]} \Theta_{0,2,j})(\max_{j \in[p]} \Theta_{m+1,2,j})\lesssim m^{-\alpha}\Psi_{2,\alpha} \Psi_{2,0}$.  Thus, 
$$
\varrho_n
\lesssim n^{-1/6}B_{n}m^{2/3}(\log p)^{7/6}+m^{-\alpha}\Phi_{\psi_{\nu},\alpha}(\log p)^{1+\nu}+m^{-\alpha/3}\Psi_{2,\alpha}^{1/3} \Psi_{2,0}^{1/3}\cdot(\log p)^{2/3}\,.$$
Letting $m\asymp n^{1/(4+2\alpha)}$, we have 
$$
\varrho_n
\lesssim n^{-\alpha/(12+6\alpha)}\{B_{n}(\log p)^{7/6}+(\log p)^{2/3}\Psi_{2,\alpha}^{1/3} \Psi_{2,0}^{1/3}\}+n^{-\alpha/(4+2\alpha)}\Phi_{\psi_{\nu},\alpha}(\log p)^{1+\nu}\,.$$ We have Part (ii) of Theorem \ref{prop:GA_weakly-dependent}. $\hfill\Box$

\begin{acks}[Acknowledgments]
 The authors would like to thank the editor, the associate editor and two reviewers for their constructive suggestions which led to the improvements of the paper.
\end{acks}
\begin{funding}
Chang and Wu were supported in part by the National Natural Science Foundation of China (grant
nos. 71991472, 72125008 and 11871401). Chang was also supported by the Center of Statistical Research at Southwestern University of
Finance and Economics. Chen was supported by in part by the National Science Foundation (grant no. 1752614). 
%
\end{funding}

\begin{supplement}
\stitle{Additional technical proofs for "Central limit theorems for high dimensional dependent data"}
\sdescription{Technical proofs of Theorems 4--11, and Proposition 1.}
\end{supplement}




\newpage

\setcounter{equation}{0}
\setcounter{table}{0}
\setcounter{lemma}{0}
\setcounter{proposition}{0}
\setcounter{page}{1}
\setcounter{section}{0}
\setcounter{ass}{0}
\renewcommand{\thelemma}{L\arabic{lemma}}
\renewcommand{\theequation}{S.\arabic{equation}}
\renewcommand{\theproposition}{P\arabic{proposition}}
\renewcommand{\theass}{C\arabic{ass}}
\renewcommand{\thepage}{S\arabic{page}}
\renewcommand{\thetable}{S\arabic{table}}
\renewcommand{\thesection}{S\arabic{section}}
\begin{center}
{\bf \Large Supplementary material for "Central limit theorems for high dimensional dependent data" by Jinyuan Chang, Xiaohui Chen and Mingcong Wu}
\end{center}

In the sequel, we use $C$ to denote a generic positive constant that may have different values in different places.

\section{Proofs of the lemmas used in Section \ref{sec:pfalphamix}}
To prove Lemmas \ref{pn:1} and \ref{pn:2}, we need the following auxiliary lemmas.

\begin{lemma}\label{la:large_deviation}
	Let $\{Z_t\}_{t=1}^{\tilde{n}} $ be an $\alpha$-mixing sequence of centered random variables with $\alpha$-mixing coefficients $\{\tilde{\alpha}(k)\}_{k\geq1}$. Assume there exist some universal constants $a_1, a_2,b_1, b_2, r_1, r_2>0$ such that {\rm(i)} $
	\max_{t\in[\tilde{n}]}\mathbb{P}(|Z_t|>u)\leq b_1\exp(-b_2u^{r_1}\tilde{B}_{\tilde{n}}^{-r_1})$
	for any $u>0$, where $\tilde{B}_{\tilde{n}}>0$ may diverge with $\tilde{n}$; {\rm(ii)} $
	\tilde{\alpha}(k)\leq a_1\exp(-a_2\tilde{L}_{\tilde{n}}^{-r_2}|k-\tilde{j}_{\tilde{n}}|_+^{r_2})$
	for any integer $k\geq 1$, where $\tilde{L}_{\tilde{n}}\geq 1$ and $0\leq \tilde{j}_{\tilde{n}}=o(\tilde{n})$ may both diverge with $\tilde{n}$. Let $r=(2+|r_1^{-1}-1|_{+}+r_2^{-1})^{-1}$, where $|\cdot|_{+}=\max(\cdot,0)$. It then holds that
	$\mathbb{P}(|\sum_{t=1}^{\tilde{n}}Z_t|\geq x)\lesssim \exp\{-C(\tilde{L}_{\tilde{n}}+\tilde{j}_{\tilde{n}})^{-1}\tilde{B}_{\tilde{n}}^{-2}\tilde{n}^{-1}x^2\}+\exp\{-C(\tilde{L}_{\tilde{n}}+\tilde{j}_{\tilde{n}})^{-r}\tilde{B}_{\tilde{n}}^{-r}x^r\}$ for any $x>0$.
\end{lemma}

\begin{lemma}\label{la:mom}
	Under Conditions {\rm\ref{as:tail}} and {\rm\ref{as:alpha-mixing}}, we have $
	\max_{\ell\in[L]}\max_{j\in[p]}\mathbb{E}(|\tilde{X}_{\ell,j}|^3)\lesssim B_{n}^3$.
\end{lemma}

\subsection{Proof of Lemma \ref{la:large_deviation}}\label{sec:largedeviation}
In this proof, we let $\bar{C}$ be a universal constant depending only on $a_{1}$, $a_2$, $b_1$, $b_2$, $r_1$ and $r_2$. Write $S_{\tilde{n}}=\sum_{t=1}^{\tilde{n}}Z_t$.
Define $$\Lambda_{\tilde{n}}(u)=\max\bigg\{1,\max_{s\in[\tilde{n}]}\sum_{t=s}^{\tilde{n}}\tilde\alpha^{1/u}(t-s)\bigg\}\,.$$ For any integer $k\geq 2$, 
\begin{align*}
\Lambda_{\tilde{n}}\{2(k-1)\}\leq&\, 1+\max_{s\in[\tilde{n}]}\sum_{t=s}^{\tilde{n}}\{\tilde\alpha (t-s)\}^{1/(2k-2)}\leq 1+\sum_{m=0}^{\tilde{n}}\{\tilde\alpha (m)\}^{1/(2k-2)}\\
\leq&\, 1+(a_{1}\vee 1)\sum_{m=0}^{\tilde{n}}\exp\{-a_2(2k-2)^{-1}\tilde{L}_{\tilde{n}}^{-r_2}|m-\tilde{j}_{\tilde{n}}|_+^{r_2}\}\\
\leq&\, 1+(a_1\vee 1)(\tilde{j}_{\tilde{n}}+1)+(a_{1}\vee 1)\sum_{m=1}^{\tilde{n}}\exp\{-a_2(2k-2)^{-1}\tilde{L}_{\tilde{n}}^{-r_2}m^{r_2}\}\\
\leq&\, 1+(a_1\vee 1)(\tilde{j}_{\tilde{n}}+1)+(a_{1}\vee 1)\sum_{m=1}^{\tilde{n}}\int_{m-1}^{m}\exp\{-a_2(2k-2)^{-1}\tilde{L}_{\tilde{n}}^{-r_2}x^{r_2}\}\,{\rm d}x\\
\leq&\, 1+(a_1\vee 1)(\tilde{j}_{\tilde{n}}+1)+(a_{1}\vee 1)\int_{0}^{\infty}\exp\{-a_2(2k-2)^{-1}\tilde{L}_{\tilde{n}}^{-r_2}x^{r_2}\}\,{\rm d}x\\
\leq&\, \bar{C}(k-1)^{1/r_{2}}(\tilde{L}_{\tilde{n}}+\tilde{j}_{\tilde{n}})\,.
\end{align*}
Due to $k^k\leq k!e^k$ for any integer $k\geq 1$, it holds that 
\begin{align*}
\Lambda^{k-1}_{\tilde{n}}\{2(k-1)\}\leq&\,\, \bar{C}^{k-1}\{(k-1)^{k-1}\}^{1/r_{2}}(\tilde{L}_{\tilde{n}}+\tilde{j}_{\tilde{n}})^{k-1}\\
\leq&\,\, \bar{C}^{k-1}(k!)^{1/r_2}e^{(k-1)/r_2}(\tilde{L}_{\tilde{n}}+\tilde{j}_{\tilde{n}})^{k-1}\\
=&\,\,(k!)^{1/r_2}\{\bar{C}(\tilde{L}_{\tilde{n}}+\tilde{j}_{\tilde{n}})\}\{\bar{C}(\tilde{L}_{\tilde{n}}+\tilde{j}_{\tilde{n}})\}^{k-2}\,.
\end{align*}
Notice that $
\max_{t\in[\tilde{n}]}\mathbb{P}(|Z_t|>u)\leq b_1\exp(-b_2u^{r_1}\tilde{B}_{\tilde{n}}^{-r_1})$.
Denote by $\Gamma (\cdot)$ the gamma function. \\
{\indent\bf Case 1.} If $r_1 \geq 1$, for any integer $k\geq 2r_1$, we have by Stirling's formula that
\begin{align*}
\mathbb{E}(|Z_{t}|^k)\leq&\, \int_{0}^{\infty}b_1ku^{k-1}\exp(-b_2u^{r_1}\tilde{B}_{\tilde{n}}^{-r_1})\,{\rm d} u\\
\leq &\, k(\bar{C}\tilde{B}_{\tilde{n}})^k \Gamma(k/r_1)\leq k(\bar{C}\tilde{B}_{\tilde{n}})^k\sqrt{2\pi (kr_1^{-1}-1)}\{e^{-1}(kr_1^{-1}-1)\}^{kr_1^{-1}-1}\\
\leq&\, k(\bar{C}\tilde{B}_{\tilde{n}})^k\sqrt{2\pi (k-1)}\{e^{-1}(k-1)\}^{k-1}\leq k!(\bar{C}\tilde{B}_{\tilde{n}})^k\,,
\end{align*}
where the fourth step is due to the  function $x\rightarrow x^x$ is nondecreasing as soon as $x\geq e^{-1}$. Using the same arguments, we also have $\mathbb{E}(|Z_{t}|^k)\leq k(\bar{C}\tilde{B}_{\tilde{n}})^k\max_{1\leq k\leq \lfloor2r_1 \rfloor}\Gamma(k/r_1)\leq k!(\bar{C}\tilde{B}_{\tilde{n}})^k$ for any integer $k<2r_1$. Hence $\mathbb{E}(|Z_{t}|^k)\leq k!(\bar{C}\tilde{B}_{\tilde{n}})^k$ for any integer $k\geq 2$. Let $G_k(x)$ be the $k$-th order cumulant of the random variable $x$. By Theorem 4.17 of \cite{SS_1991} with $\delta=1$, we have $|G _k(S_{\tilde{n}})|\leq (k!)^{2+r_2^{-1}}\tilde{n}\{\bar{C}(\tilde{L}_{\tilde{n}}+\tilde{j}_{\tilde{n}})\tilde{B}_{\tilde{n}}^2\}\{\bar{C}(\tilde{L}_{\tilde{n}}+\tilde{j}_{\tilde{n}})\tilde{B}_{\tilde{n}}\}^{k-2}$. By Lemma 2.4 of \cite{SS_1991}, it holds that $$\mathbb{P}(|S_{\tilde{n}}|\geq x)\lesssim \exp\{-\bar{C}(\tilde{L}_{\tilde{n}}+\tilde{j}_{\tilde{n}})^{-1}\tilde{B}_{\tilde{n}}^{-2}\tilde{n}^{-1}x^2\}+\exp\{-\bar{C}(\tilde{L}_{\tilde{n}}+\tilde{j}_{\tilde{n}})^{-r}\tilde{B}_{\tilde{n}}^{-r}x^r\}$$ with $r=(2+r_2^{-1})^{-1}$. \\
{\indent\bf Case 2.} If $0<r_1<1$, for any integer $k\geq 2$,
\begin{align*}
\mathbb{E}(|Z_{t}|^k)\leq&\, k(\bar{C}\tilde{B}_{\tilde{n}})^k\sqrt{2\pi (kr_1^{-1}-1)}\{e^{-1}(kr_1^{-1}-1)\}^{kr_1^{-1}-1}\\
\leq&\, k(\bar{C}\tilde{B}_{\tilde{n}})^k\{\sqrt{2\pi k}(e^{-1}k)^{k}\}^{1/r_1}\leq (k!)^{r_1^{-1}}(\bar{C}\tilde{B}_{\tilde{n}})^k\,,
\end{align*}
where the last step is due to $k\leq 2^{k}$. Using the same arguments in Case 1, we also have $|G _k(S_{\tilde{n}})|\leq (k!)^{1+r_1^{-1}+r_2^{-1}}\tilde{n}\{\bar{C}(\tilde{L}_{\tilde{n}}+\tilde{j}_{\tilde{n}})\tilde{B}_{\tilde{n}}^2\}\{\bar{C}(\tilde{L}_{\tilde{n}}+\tilde{j}_{\tilde{n}})\tilde{B}_{\tilde{n}}\}^{k-2}$, and then $$\mathbb{P}(|S_{\tilde{n}}|\geq x)\lesssim \exp\{-\bar{C}(\tilde{L}_{\tilde{n}}+\tilde{j}_{\tilde{n}})^{-1}\tilde{B}_{\tilde{n}}^{-2}\tilde{n}^{-1}x^2\}+\exp\{-\bar{C}(\tilde{L}_{\tilde{n}}+\tilde{j}_{\tilde{n}})^{-r}\tilde{B}_{\tilde{n}}^{-r}x^r\}$$ with $r=(1+r_1^{-1}+r_2^{-1})^{-1}$. We complete the proof of Lemma \ref{la:large_deviation}. $\hfill\Box$

\subsection{Proof of Lemma \ref{la:mom}} Applying Lemma \ref{la:large_deviation} with $(\tilde{B}_{\tilde{n}},\tilde{L}_{\tilde{n}},\tilde{j}_{\tilde{n}},r_1,r_2)=(B_n,1,0,\gamma_{1},\gamma_{2})$, then
\begin{align}\label{eq:xtildetail}
	\max_{\ell\in[L]}\max_{j\in[p]}\mathbb{P}(|\tilde{X}_{\ell,j}|>\lambda)\lesssim \exp(-CB_n^{-2}\lambda^2)+\exp(-Cb^{\gamma/2}B_n^{-\gamma}\lambda^{\gamma})\,,
\end{align}
where $\gamma=\gamma_{2}/(2\gamma_{2}+1)$. Notice that $\mathbb{E}(|\tilde{X}_{\ell,j}|^3)=3\int_0^\infty \lambda^2\mathbb{P}(|\tilde{X}_{\ell,j}|>\lambda)\,{\rm d}\lambda$. By some basic calculations, we have $
\max_{\ell\in[L]}\max_{j\in[p]}\mathbb{E}(|\tilde{X}_{\ell,j}|^3)\lesssim B_{n}^3$. $\hfill\Box$

\subsection{Proof of Lemma \ref{pn:1}}\label{se:pfbd1alpha}
Let $\mathcal{W}_n=\{W_1,\ldots,W_{n}\}$ be a copy of $\mathcal{Y}_n=\{Y_1,\ldots,Y_{n}\}$. Write $\mathcal{X}_n=\{X_1,\ldots,X_n\}$. We assume $\mathcal{X}_n$, $\mathcal{Y}_n$ and $\mathcal{W}_n$ are independent. Let $S_{n,w}^{(1)}=L^{-1/2}\sum_{\ell=1}^L\tilde{W}_\ell$ with $\tilde{W}_\ell=b^{-1/2}\sum_{t\in\mathcal{I}_\ell}W_t$. Then $\varrho_n^{(1)}$ defined as \eqref{eq:alpharho1} can be reformulated as
$$
\varrho_n^{(1)}=\sup_{u\in\mathbb{R}^p,\nu\in[0,1]}|\mathbb{P}\{\sqrt{\nu}S_{n,x}^{(1)}+\sqrt{1-\nu}S_{n,y}^{(1)}\leq u\}-\mathbb{P}\{S_{n,w}^{(1)}\leq u\}|\,.$$
Recall $\tilde{X}_\ell=b^{-1/2}\sum_{t\in\mathcal{I}_\ell}X_t$ and $\tilde{Y}_\ell=b^{-1/2}\sum_{t\in\mathcal{I}_\ell}Y_t$. For any $\phi>0$, define
\begin{align}
M_{\tilde{x}}(\phi)=&\,\max_{\ell\in[L]}\mathbb{E}\bigg\{|\tilde{X}_\ell|_\infty^3I\bigg(|\tilde{X}_\ell|_\infty>\frac{\sqrt{L}}{4\phi\log p}\bigg)\bigg\}\,,\label{eq:mtildex}\\
M_{\tilde{y}}(\phi)=&\,\max_{\ell\in[L]}\mathbb{E}\bigg\{|\tilde{Y}_\ell|_\infty^3I\bigg(|\tilde{Y}_\ell|_\infty>\frac{\sqrt{L}}{4\phi\log p}\bigg)\bigg\}\,.\label{eq:mtildey}\end{align}

For a given $u=(u_1,\ldots,u_p)^\T\in\mathbb{R}^p$, define
$
F_\beta(v)=\beta^{-1}\log[\sum_{j=1}^p\exp\{\beta(v_j-u_j)\}]$ with $\beta:=\phi\log p$
for any $v=(v_1,\ldots,v_p)^\T\in\mathbb{R}^p$. Such defined function $F_{\beta}(v)$ satisfies the property $0\leq F_{\beta}(v)-\max_{j\in[p]}(v_j-u_j)\leq \beta^{-1}\log p=\phi^{-1}$ for any $v\in\mathbb{R}^p$. Select a thrice continuously differentiable function $g_0:\mathbb{R}\rightarrow[0,1]$ whose derivatives up to the third order are all bounded such that $g_0(t)=1$ for $t\leq 0$ and $g_0(t)=0$ for $t\geq1$. Define $g(t):=g_0(\phi t)$ for any $t\in\mathbb{R}$, and $q(v):=g\{F_\beta(v)\}$ for any $v\in\mathbb{R}^p$. Define
$$
\mathcal{T}_n:=q\{\sqrt{\nu}S_{n,x}^{(1)}+\sqrt{1-\nu}S_{n,y}^{(1)}\}-q\{S_{n,w}^{(1)}\}\,.$$
Write $v_n=\sqrt{\nu}S_{n,x}^{(1)}+\sqrt{1-\nu}S_{n,y}^{(1)}$. Notice that $\tilde\Xi=\cov\{S_{n,w}^{(1)}\}$. Due to $B_n^2h\ll b\ll nB_n^{-2}$, by Lemma \ref{lem:cov_infinity}, $|\tilde\Xi-\Xi|_{\infty}=o(1)$, which implies the all the elements in the main diagonal of $\tilde{\Xi}$ are uniformly bounded away from zero. Then
\begin{align*}
\mathbb{P}(v_n\leq u-\phi^{-1})\leq&\,\,\mathbb{P}\{F_\beta(v_n)\leq 0\}\leq\mathbb{E}\{q(v_n)\}\\
\leq&\,\,\mathbb{P}[F_\beta\{S_{n,w}^{(1)}\}\leq \phi^{-1}]+\mathbb{E}(\mathcal{T}_n)\leq\mathbb{P}\{S_{n,w}^{(1)}\leq u+\phi^{-1}\}+|\mathbb{E}(\mathcal{T}_n)|\\
\leq&\,\,\mathbb{P}\{S_{n,w}^{(1)}\leq u-\phi^{-1}\}+C\phi^{-1}(\log p)^{1/2}+|\mathbb{E}(\mathcal{T}_n)|\,,
\end{align*}
where the last step is based on Nazarov's inequality. Analogously, $\mathbb{P}(v_n\leq u-\phi^{-1})\geq\mathbb{P}\{S_{n,w}^{(1)}\leq u-\phi^{-1}\}-C\phi^{-1}(\log p)^{1/2}-|\mathbb{E}(\mathcal{I}_n)|$. Hence,
$$
\varrho_n^{(1)}\lesssim\phi^{-1}(\log p)^{1/2}+\sup_{u\in\mathbb{R}^p,\nu\in[0,1]}|\mathbb{E}(\mathcal{T}_n)|\,.$$
Let $\mathcal{F}_{-\ell}$ be the $\sigma$-filed generated by $\{\tilde{X}_s\}_{s\neq \ell}$. We will show in Section \ref{sec:ETnalpha} that
\begin{align}\label{eq:alphatoshow1}
\sup_{u\in\mathbb{R}^p,\nu\in[0,1]}|\mathbb{E}(\mathcal{T}_n)|\lesssim&\,\,\phi \sqrt{L}\max_{\ell\in[L]}\mathbb{E}\bigg\{\max_{j\in[p]}|\mathbb{E}(\tilde{X}_{\ell,j}\,|\,\mathcal{F}_{-\ell})|\bigg\}\notag\\
&+\phi^2(\log p)\max_{\ell\in[L]}\mathbb{E}\bigg[\max_{k,j\in[p]}|\mathbb{E}\{\tilde{X}_{\ell,k}\tilde{X}_{\ell,j}-\mathbb{E}(\tilde{X}_{\ell,k}\tilde{X}_{\ell,j})\,|\,\mathcal{F}_{-\ell}\}|\bigg]\\
&+\frac{\phi^3(\log p)^2}{\sqrt{L}} \max_{\ell\in[L]}\mathbb{E}\bigg[\max_{j\in[p]}|\mathbb{E}\{|\tilde{X}_{\ell,j}|^3-\mathbb{E}(|\tilde{X}_{\ell,j}|^3)\,|\,\mathcal{F}_{-\ell}\}|\bigg]\notag\\
&+\frac{\phi^3(\log p)^2}{\sqrt{L}}\big\{B_{n}^3\phi^{-1}(\log p)^{1/2}+B_{n}^3\varrho_n^{(1)}+M_{\tilde{x}}(\phi)+M_{\tilde{y}}(\phi)\big\}\,.\notag
\end{align}
With selecting $\phi=CB_{n}^{-1}L^{1/6}(\log p)^{-2/3}$ for some sufficiently small $C>0$, we have
\begin{align}\label{eq:rho1}
\varrho_n^{(1)}\lesssim&\,\,\frac{B_{n}(\log p)^{7/6}}{L^{1/6}}+\frac{1}{B_{n}^3}M_{\tilde{x}}\bigg\{\frac{CL^{1/6}}{B_{n}(\log p)^{2/3}}\bigg\}+\frac{1}{B_{n}^3}M_{\tilde{y}}\bigg\{\frac{CL^{1/6}}{B_{n}(\log p)^{2/3}}\bigg\}\notag\\
&+\frac{L^{2/3}}{B_{n}(\log p)^{2/3}}\max_{\ell\in[L]}\mathbb{E}\bigg\{\max_{j\in[p]}|\mathbb{E}(\tilde{X}_{\ell,j}\,|\,\mathcal{F}_{-\ell})|\bigg\}\\
&+\frac{L^{1/3}}{B_{n}^2(\log p)^{1/3}}\max_{\ell\in[L]}\mathbb{E}\bigg[\max_{k,j\in[p]}|\mathbb{E}\{\tilde{X}_{\ell,k}\tilde{X}_{\ell,j}-\mathbb{E}(\tilde{x}_{\ell,k}\tilde{x}_{\ell,j})\,|\,\mathcal{F}_{-\ell}\}|\bigg]\notag\\
&+\frac{1}{B_{n}^3}\max_{\ell\in[L]}\mathbb{E}\bigg[\max_{j\in[p]}|\mathbb{E}\{|\tilde{X}_{\ell,j}|^3-\mathbb{E}(|\tilde{X}_{\ell,j}|^3)\,|\,\mathcal{F}_{-\ell}\}|\bigg]\,.\notag
\end{align}
Let $\eta=B_{n}L^{1/3}(\log p)^{-1/3}$. It holds that
\begin{align}\label{eq:tailmom}
\mathbb{E}\bigg\{|\tilde{Y}_\ell|_\infty^3I\bigg(|\tilde{Y}_\ell|_\infty>\frac{\eta}{4C}\bigg)\bigg\}=\bigg(\frac{\eta}{4C}\bigg)^3\mathbb{P}\bigg(|\tilde{Y}_\ell|_\infty>\frac{\eta}{4C}\bigg)+3\int_{\frac{\eta}{4C}}^\infty u^2\mathbb{P}\big(|\tilde{Y}_\ell|_\infty>u\big)\,{\rm d}u\,.
\end{align}
Notice that $\tilde{Y}_\ell=(\tilde{Y}_{\ell,1},\ldots,\tilde{Y}_{\ell,p})^\T\sim {N}\{0,\mathbb{E}(\tilde{X}_{\ell}\tilde{X}_\ell^\T)\}$. Denote by $\sigma_{\ell,j,j}^2$ the  $j$-th element in the main diagonal of $\mathbb{E}(\tilde{X}_{\ell}\tilde{X}_\ell^\T)$. By Bonferroni inequality, we then have
\begin{align*}
\mathbb{P}\{|\tilde{Y}_\ell|_\infty>\eta(4C)^{-1}\}\leq&\,\, \sum_{j=1}^p\mathbb{P}\{|\tilde{Y}_{\ell,j}|>\eta(4C)^{-1}\}\lesssim \frac{1}{\eta}\sum_{j=1}^p\sigma_{\ell,j,j}\exp\{-\eta^2(32C^2\sigma_{\ell,j,j}^2)^{-1}\}
\end{align*}
and
\begin{align*}
\int_{\eta(4C)^{-1}}^\infty u^2\mathbb{P}(|\tilde{Y}_\ell|_\infty>u)\,{\rm d}u\leq&\,\,\sum_{j=1}^p\int_{\eta(4C)^{-1}}^\infty u^2\mathbb{P}\big(|\tilde{Y}_{\ell,j}|>u\big)\,{\rm d}u\\
\lesssim&\,\,\sum_{j=1}^p\sigma_{\ell,j,j}^3\exp\{-\eta^2(32C^2\sigma_{\ell,j,j}^2)^{-1}\}\,,
\end{align*}
which implies that
$$
\mathbb{E}[|\tilde{Y}_\ell|_\infty^3I\{|\tilde{Y}_\ell|_\infty>\eta(4C)^{-1}\}]\lesssim \sum_{j=1}^p(\eta^2\sigma_{\ell,j,j}+\sigma_{\ell,j,j}^3)\exp\{-\eta^2(32C^2\sigma_{\ell,j,j}^2)^{-1}\}\,.$$
By Lemma \ref{la:mom}, we know $\max_{\ell\in[L]}\max_{j\in[p]}\mathbb{E}(|\tilde{X}_{\ell,j}|^2)\lesssim B_{n}^2$. Then, if $\log p=o(L^{2/5})$, we have
\begin{align*}
M_{\tilde{y}}\bigg\{\frac{CL^{1/6}}{B_{n}(\log p)^{2/3}}\bigg\}=\max_{\ell\in[L]}\mathbb{E}\bigg\{|\tilde{Y}_\ell|_\infty^3I\bigg(|\tilde{Y}_\ell|_\infty>\frac{\eta}{4C}\bigg)\bigg\}\lesssim\frac{B_{n}^3(\log p)^{7/6}}{L^{1/6}}\,.
\end{align*}
On the other hand, by Bonferroni inequality and \eqref{eq:xtildetail}, 
$$
\max_{\ell\in[L]}\mathbb{P}(|\tilde{X}_{\ell}|_\infty>u)\lesssim p\exp(-CB_n^{-2}u^2)+p\exp(-Cb^{\gamma/2}B_n^{-\gamma}u^{\gamma})\,.$$
Same as (\ref{eq:tailmom}), if $(\log p)^{3+\gamma}=o(b^{3\gamma/2}L^{\gamma})$ and $\log p=o(L^{2/5})$, we have 
$
\mathbb{E}[|\tilde{X}_\ell|_\infty^3I\{|\tilde{X}_\ell|_\infty>\eta(4C)^{-1}\}]\lesssim B_{n}^3L^{-1/6}(\log p)^{7/6}$,
which implies 
\begin{align*}
M_{\tilde{x}}\bigg\{\frac{CL^{1/6}}{B_{n}(\log p)^{2/3}}\bigg\}=\max_{\ell\in[L]}\mathbb{E}\bigg\{|\tilde{X}_\ell|_\infty^3I\bigg(|\tilde{X}_\ell|_\infty>\frac{\eta}{4C}\bigg)\bigg\}\lesssim \frac{B_{n}^3(\log p)^{7/6}}{L^{1/6}}\,.
\end{align*}
It follows from (\ref{eq:rho1}) that
\begin{align}\label{eq:bou1}
\varrho_n^{(1)}\lesssim&\,\,\frac{B_{n}(\log p)^{7/6}}{L^{1/6}}+\frac{L^{2/3}}{B_{n}(\log p)^{2/3}}\max_{\ell\in[L]}\mathbb{E}\bigg\{\max_{j\in[p]}|\mathbb{E}(\tilde{X}_{\ell,j}\,|\,\mathcal{F}_{-\ell})|\bigg\}\notag\\
&+\frac{L^{1/3}}{B_{n}^2(\log p)^{1/3}}\max_{\ell\in[L]}\mathbb{E}\bigg[\max_{k,j\in[p]}|\mathbb{E}\{\tilde{X}_{\ell,k}\tilde{X}_{\ell,j}-\mathbb{E}(\tilde{X}_{\ell,k}\tilde{X}_{\ell,j})\,|\,\mathcal{F}_{-\ell}\}|\bigg]\\
&+\frac{1}{B_{n}^3}\max_{\ell\in[L]}\mathbb{E}\bigg[\max_{j\in[p]}|\mathbb{E}\{|\tilde{X}_{\ell,j}|^3-\mathbb{E}(|\tilde{X}_{\ell,j}|^3)\,|\,\mathcal{F}_{-\ell}\}|\bigg]\notag
\end{align}
provided that $(\log p)^{3+\gamma}=o(b^{3\gamma/2}L^{\gamma})$ and $\log p=o(L^{2/5})$.

Given $D_{1n}\rightarrow\infty$ as $n\rightarrow\infty$ such that $D_{1n}\gg B_n$, 
$
\mathbb{E}(\tilde{X}_{\ell,j}\,|\,\mathcal{F}_{-\ell})=\mathbb{E}\{\tilde{X}_{\ell,j}I(|\tilde{X}_{\ell,j}|\leq D_{1n})\,|\,\mathcal{F}_{-\ell}\}-\mathbb{E}\{\tilde{X}_{\ell,j}I(|\tilde{X}_{\ell,j}|\leq D_{1n})\}+\mathbb{E}\{\tilde{X}_{\ell,j}I(|\tilde{X}_{\ell,j}|> D_{1n})\,|\,\mathcal{F}_{-\ell}\}-\mathbb{E}\{\tilde{X}_{\ell,j}I(|\tilde{X}_{\ell,j}|> D_{1n})\}$,
which implies that
\begin{align}\label{eq:bound1}
\mathbb{E}\{|\mathbb{E}(\tilde{X}_{\ell,j}\,|\,\mathcal{F}_{-\ell})|\}\leq&\,\,\mathbb{E}[|\mathbb{E}\{\tilde{X}_{\ell,j}I(|\tilde{X}_{\ell,j}|\leq D_{1n})\,|\,\mathcal{F}_{-\ell}\}-\mathbb{E}\{\tilde{X}_{\ell,j}I(|\tilde{X}_{\ell,j}|\leq D_{1n})\}|]\notag\\
&+2\mathbb{E}\{|\tilde{X}_{\ell,j}|I(|\tilde{X}_{\ell,j}|> D_{1n})\}\,.
\end{align}
It follows from \eqref{eq:xtildetail} that
\begin{equation}\label{eq:bound11}
\begin{split}
\mathbb{E}\{|\tilde{X}_{\ell,j}|I(|\tilde{X}_{\ell,j}|> D_{1n})\}=&\,\,D_{1n}\mathbb{P}(|\tilde{X}_{\ell,j}|>D_{1n})+\int_{D_{1n}}^\infty\mathbb{P}(|\tilde{X}_{\ell,j}|>u)\,{\rm d}u\\\
\lesssim&\,\,B_{n}\exp(-Cb^{\gamma/2}D_{1n}^\gamma B_{n}^{-\gamma})+B_{n}\exp(-CD_{1n}^2B_{n}^{-2})\,.
\end{split}
\end{equation}
Note that
\begin{align*}
\mathbb{E}\{\tilde{X}_{\ell,j}I(|\tilde{X}_{\ell,j}|\leq D_{1n})\,|\,\mathcal{F}_{-\ell}\}=&\,D_{1n}\mathbb{P}(|\tilde{X}_{\ell,j}|>D_{1n}\,|\,\mathcal{F}_{-\ell})-2D_{1n}\mathbb{P}(\tilde{X}_{\ell,j}>D_{1n}\,|\,\mathcal{F}_{-\ell})\\
&+\int_{-D_{1n}}^{D_{1n}}\mathbb{P}(\tilde{X}_{\ell,j}>u\,|\,\mathcal{F}_{-\ell})\,{\rm d}u-D_{1n}
\end{align*}
and
\begin{align*}
\mathbb{E}\{\tilde{X}_{\ell,j}I(|\tilde{X}_{\ell,j}|\leq D_{1n})\}=&\,D_{1n}\mathbb{P}(|\tilde{X}_{\ell,j}|>D_{1n})-2D_{1n}\mathbb{P}(\tilde{X}_{\ell,j}>D_{1n})\\
&+\int_{-D_{1n}}^{D_{1n}}\mathbb{P}(\tilde{X}_{\ell,j}>u)\,{\rm d}u-D_{1n}\,.
\end{align*}
Equation (1.10c) of \cite{Rio_2013} implies that
$
\int_{-D_{1n}}^{D_{1n}}\mathbb{E}\{|\mathbb{P}(\tilde{X}_{\ell,j}>u\,|\,\mathcal{F}_{-\ell})-\mathbb{P}(\tilde{X}_{\ell,j}>u)|\}\,{\rm d}u\lesssim D_{1n}\alpha_n(h)$.
By triangle inequality and  \eqref{eq:xtildetail}, we have
\begin{align*}
&\mathbb{E}[|\mathbb{E}\{\tilde{X}_{\ell,j}I(|\tilde{X}_{\ell,j}|\leq D_{1n})\,|\,\mathcal{F}_{-\ell}\}-\mathbb{E}\{\tilde{X}_{\ell,j}I(|\tilde{X}_{\ell,j}|\leq D_{1n})\}|]\\
&~~~\leq6D_{1n}\mathbb{P}(|\tilde{X}_{\ell,j}|>D_{1n})+\int_{-D_{1n}}^{D_{1n}}\mathbb{E}\{|\mathbb{P}(\tilde{X}_{\ell,j}>u\,|\,\mathcal{F}_{-\ell})-\mathbb{P}(\tilde{X}_{\ell,j}>u)|\}\,{\rm d}u\\
&~~~\lesssim B_{n} \exp(-Cb^{\gamma/2}D_{1n}^\gamma B_{n}^{-\gamma})+B_{n}\exp(-CD_{1n}^2B_{n}^{-2})+D_{1n}\alpha_n(h)
\end{align*}
holds uniformly over $\ell\in[L]$ and $j\in[p]$. Together with (\ref{eq:bound11}), (\ref{eq:bound1}) implies that
\begin{align}\label{eq:bound12}
&\max_{\ell\in[L]}\mathbb{E}\bigg\{\max_{j\in[p]}|\mathbb{E}(\tilde{X}_{\ell,j}\,|\,\mathcal{F}_{-\ell})|\bigg\}\leq\max_{\ell\in[L]}\sum_{j=1}^p\mathbb{E}\{|\mathbb{E}(\tilde{X}_{\ell,j}\,|\,\mathcal{F}_{-\ell})|\}\notag\\
&~~\lesssim pB_{n}\exp(-Cb^{\gamma/2}D_{1n}^\gamma B_{n}^{-\gamma})+pB_{n}\exp(-CD_{1n}^2B_{n}^{-2})+pD_{1,n}\alpha_n(h)\,.
\end{align}
Recall $\eta=B_{n}L^{1/3}(\log p)^{-1/3}$. Selecting $h=C'\{\log(pn)\}^{1/\gamma_2}$ and $D_{1n}=C'\eta$ for some sufficiently large $C'>0$, due to $(\log p)^{3+\gamma}=o(b^{3\gamma/2}L^{\gamma})$ and $\log p=o(L^{2/5})$, 
$$
\max_{\ell\in[L]}\mathbb{E}\bigg\{\max_{j\in[p]}|\mathbb{E}(\tilde{X}_{\ell,j}\,|\,\mathcal{F}_{-\ell})|\bigg\}\lesssim B_{n}^{7/2}\eta^{-5/2}
\,,$$
which implies
\begin{align*}
\frac{L^{2/3}}{B_{n}(\log p)^{2/3}}\max_{\ell\in[L]}\mathbb{E}\bigg\{\max_{j\in[p]}|\mathbb{E}(\tilde{X}_{\ell,j}\,|\,\mathcal{F}_{-\ell})|\bigg\}\lesssim B_{n}^{1/2}\eta^{-1/2}\lesssim \frac{B_{n}(\log p)^{7/6}}{L^{1/6}}\,.
\end{align*}
By Lemma 2 of \cite{CTW_2013} and \eqref{eq:xtildetail}, 
$
\max_{\ell\in[L]}\max_{k,j\in[p]}\mathbb{P}\{|\tilde{X}_{\ell,k}\tilde{X}_{\ell,j}-\mathbb{E}(\tilde{X}_{\ell,k}\tilde{X}_{\ell,j})|>\lambda\}\lesssim \exp(-CB_n^{-2}\lambda)+\exp(-Cb^{\gamma/2}B_n^{-\gamma}\lambda^{\gamma/2})$
for any $\lambda\geq C''B_{n}^2$, where $C''>0$ is a sufficiently large constant. Similarly, 
\begin{align*}
\max_{\ell\in[L]}\max_{j\in[p]}\mathbb{P}\{||\tilde{X}_{\ell,j}|^3-\mathbb{E}(|\tilde{X}_{\ell,j}|^3)|>\lambda\}\lesssim \exp(-CB_n^{-2}\lambda^{2/3})+\exp(-Cb^{\gamma/2}B_n^{-\gamma}\lambda^{\gamma/3})
\end{align*}
for any $\lambda\geq C''B_{n}^3$. Let $D_{2n}=C''\eta^2$ and $D_{3n}=C''\eta^3$. Due to $\log p=o(L^{2/5})$, we know $D_{2n}\geq C''B_{n}^2$ and $D_{3n}\geq C''B_{n}^3$ for sufficiently large $n$.
Applying the identical arguments for (\ref{eq:bound12}), we have
\begin{align}\label{eq:conditional_bound_2}
&\max_{\ell\in[L]}\mathbb{E}\bigg[\max_{k,j\in[p]}|\mathbb{E}\{\tilde{X}_{\ell,k}\tilde{X}_{\ell,j}-\mathbb{E}(\tilde{X}_{\ell,k}\tilde{X}_{\ell,j})\,|\,\mathcal{F}_{-\ell}\}|\bigg]\notag\\
&~~~\lesssim p^2B_{n}^2\exp(-Cb^{\gamma/2}\eta^{\gamma}B_{n}^{-\gamma})+p^2B_{n}^2\exp(-C\eta^2B_{n}^{-2})+p^2\eta^2\alpha_n(h)\,,\\
&\max_{\ell\in[L]}\mathbb{E}\bigg[\max_{j\in[p]}|\mathbb{E}\{|\tilde{X}_{\ell,j}|^3-\mathbb{E}(|\tilde{X}_{\ell,j}|^3)|\mathcal{F}_{-\ell}\}|\bigg]\notag\\
&~~~\lesssim pB_{n}^3\exp(-Cb^{\gamma/2}\eta^{\gamma}B_{n}^{-\gamma})+pB_{n}^{3}\exp(-C\eta^{2}B_{n}^{-2})+p\eta^3\alpha_n(h)\notag\,.
\end{align}
Recall $h=C'\{\log(pn)\}^{1/\gamma_2}$. If $(\log p)^{3+\gamma}=o(b^{3\gamma/2}L^{\gamma})$ and $\log p=o(L^{2/5})$, it holds that
\begin{align*}
&\frac{L^{1/3}}{B_{n}^2(\log p)^{1/3}}\max_{\ell\in[L]}\mathbb{E}\bigg[\max_{k,j\in[p]}|\mathbb{E}\{\tilde{X}_{\ell,k}\tilde{X}_{\ell,j}-\mathbb{E}(\tilde{X}_{\ell,k}\tilde{X}_{\ell,j})\,|\,\mathcal{F}_{-\ell}\}|\bigg]\lesssim B_{n}^{1/2}\eta^{-1/2}\lesssim\frac{B_{n}(\log p)^{7/6}}{L^{1/6}}\,,\\
&\,\,\,\,\,\quad\,\quad\frac{1}{B_{n}^3}\max_{\ell\in[L]}\mathbb{E}\bigg[\max_{j\in[p]}|\mathbb{E}\{|\tilde{X}_{\ell,j}|^3-\mathbb{E}(|\tilde{X}_{\ell,j}|^3)|\mathcal{F}_{-\ell}\}|\bigg]\lesssim B_{n}^{1/2}\eta^{-1/2}\lesssim\frac{B_{n}(\log p)^{7/6}}{L^{1/6}}\,.
\end{align*}
By (\ref{eq:bou1}), we have $\varrho_n^{(1)}\lesssim B_{n}L^{-1/6}(\log p)^{7/6}$ provided that $(\log p)^{3+\gamma}=o(b^{3\gamma/2}L^{\gamma})$ and $\log p=o(L^{2/5})$. $\hfill\Box$

\subsection{Proof of Lemma \ref{pn:2}}\label{se:pfpn2} Recall $
S_{n,x}=S_{n,x}^{(1)}+\delta_n$ with $
\delta_n=n^{-1/2}\sum_{\ell=1}^{L+1}\sum_{t\in\mathcal{J}_\ell}{X}_t+\{n^{-1/2}-(Lb)^{-1/2}\}\sum_{\ell=1}^L\sum_{t\in\mathcal{I}_\ell}{X}_t$.
For some ${D}_{n}>0$, consider the event $\mathcal{E}=\{|\delta_n|_\infty\leq {D}_{n}\}$. Then
\begin{align}\label{eq:bb2}
\varrho_n^{(2)}\leq&\,\,\varrho_n^{(1)}+\sup_{u\in\mathbb{R}^p}|\mathbb{P}\{S_{n,y}^{(1)}\leq u-\sqrt{\nu}D_n\}-\mathbb{P}\{S_{n,y}^{(1)}\leq u\}|\notag\\
&+\sup_{u\in\mathbb{R}^p}|\mathbb{P}\{S_{n,y}^{(1)}\leq u+\sqrt{\nu}D_n\}-\mathbb{P}\{S_{n,y}^{(1)}\leq u\}|+\mathbb{P}(\mathcal{E}^c)\\
\lesssim&\,\,B_{n}L^{-1/6}(\log p)^{7/6}+{D}_{n}(\log p)^{1/2}+\mathbb{P}(\mathcal{E}^c)\,,\notag
\end{align}
where the last step is based on Lemma \ref{pn:1} and Nazarov's inequality. In the sequel, we will bound $\mathbb{P}(\mathcal{E}^c)$. Due to $b\ll n^{1/2}$ and $L\asymp nb^{-1}$, we have $n-Lb\leq 2Lh$ which implies
$
D_n\sqrt{nLb}(\sqrt{n}-\sqrt{Lb})^{-1}\geq D_n\sqrt{n}bh^{-1}$. By Bonferroni inequality, it holds that
\begin{align*}
\mathbb{P}(\mathcal{E}^c)\leq&\,\,\mathbb{P}\bigg(\bigg|\sum_{\ell=1}^{L+1}\sum_{t\in\mathcal{J}_\ell}{X}_t\bigg|_\infty>\frac{D_n\sqrt{n}}{2}\bigg)+\mathbb{P}\bigg\{\bigg|\sum_{\ell=1}^L\sum_{t\in\mathcal{I}_\ell}{X}_t\bigg|_\infty>\frac{D_n\sqrt{nLb}}{2(\sqrt{n}-\sqrt{Lb})}\bigg\}\\
\leq&\,\,\sum_{j=1}^p\mathbb{P}\bigg(\bigg|\sum_{\ell=1}^{L+1}\sum_{t\in\mathcal{J}_\ell}X_{t,j}\bigg|>\frac{D_n\sqrt{n}}{2}\bigg)+\sum_{j=1}^p\mathbb{P}\bigg\{\bigg|\sum_{\ell=1}^L\sum_{t\in\mathcal{I}_\ell}X_{t,j}\bigg|>\frac{D_n\sqrt{n}b}{2h}\bigg\}\,.
\end{align*}
For each $j\in[p]$, $\{X_{t,j}\}_{t\in\mathcal{I}_\ell,\ell\in[L]}$ and $\{X_{t,j}\}_{t\in\mathcal{J}_\ell,\ell\in[L+1]}$ are also $\alpha$-mixing processes. Note that  $h=o(b)$. If $D_{n}\sqrt{n}\rightarrow\infty$ as $n\rightarrow\infty$,
by Lemma \ref{la:large_deviation}, we have
\begin{align*}
&\max_{j\in[p]}\mathbb{P}\bigg(\bigg|\sum_{\ell=1}^L\sum_{t\in\mathcal{I}_\ell}X_{t,j}\bigg|>\frac{D_n\sqrt{n}b}{2h}\bigg)+\max_{j\in[p]}\mathbb{P}\bigg(\bigg|\sum_{\ell=1}^{L+1}\sum_{t\in\mathcal{J}_\ell}X_{t,j}\bigg|>\frac{D_n\sqrt{n}}{2}\bigg)\\
&~~~\lesssim\exp(-CD_n^{\gamma}n^{\gamma/2}B_{n}^{-\gamma})+\exp(-CD_n^2bh^{-1}B_{n}^{-2})\,.
\end{align*}
With selecting $D_n=C'B_{n}L^{-1/6}(\log p)^{2/3}$ for some sufficiently large constant $C'>0$, it holds that $D_n(\log p)^{1/2}\lesssim B_{n}L^{-1/6}(\log p)^{7/6}$ and
$$
\mathbb{P}(\mathcal{E}^c)\lesssim p\exp\{-Cn^{\gamma/3}b^{\gamma/6}(\log p)^{2\gamma/3}\}+p\exp\{-Cn^{-1/3}b^{4/3}h^{-1}(\log p)^{4/3}\}\,.$$
Note that $p\geq n^\kappa$ for some $\kappa>0$ and $h\asymp\{\log(pn)\}^{1/\gamma_2}\asymp(\log p)^{1/\gamma_2}$. To make $\mathbb{P}(\mathcal{E}^c)\rightarrow0$, we require two restrictions: (i)
$b\geq C''n^{1/4}h^{3/4}(\log p)^{-1/4}$ for some sufficiently large constant $C''>0$, and (ii) $
(\log p)^{3-2\gamma}=O(n^{\gamma}b^{\gamma/2})$.
Under such restrictions, 
\begin{align*}
p\exp\{-Cn^{-1/3}b^{4/3}h^{-1}(\log p)^{4/3}\}\lesssim&\,\, \exp(-C\log p)\lesssim B_nL^{-1/6}(\log p)^{7/6}\,,\\ p\exp\{-Cn^{\gamma/3}b^{\gamma/6}(\log p)^{2\gamma/3}\} \lesssim&\,\, B_{n}L^{-1/6}(\log p)^{7/6}\,.
\end{align*}
Hence, $
\mathbb{P}(\mathcal{E}^c)\lesssim B_{n}L^{-1/6}(\log p)^{7/6}$. From (\ref{eq:bb2}), we complete the proof. $\hfill\Box$

\subsection{Proof of \eqref{eq:alphatoshow1}}\label{sec:ETnalpha}
Let $Z(\omega)=\sum_{\ell=1}^LZ_\ell(\omega)$, where
$
Z_\ell(\omega)=L^{-1/2}\{\sqrt{\omega}(\sqrt{\nu}\tilde{X}_\ell+\sqrt{1-\nu}\tilde{Y}_{\ell})+\sqrt{1-\omega}\tilde{W}_\ell\}$.
Then $Z(1)=\sqrt{\nu}S_{n,x}^{(1)}+\sqrt{1-\nu}S_{n,y}^{(1)}$ and $Z(0)=S_{n,w}^{(1)}$. Let $Z^{(-\ell)}(\omega)=Z(\omega)-Z_{\ell}(\omega)$ and
\begin{align*}
\dot{Z}_\ell(\omega)&=L^{-1/2}\{\omega^{-1/2}(\sqrt{\nu}\tilde{X}_\ell+\sqrt{1-\nu}\tilde{Y}_\ell)-(1-\omega)^{-1/2}\tilde{W}_\ell\}=:\{\dot{Z}_{\ell,1}(\omega),\ldots,\dot{Z}_{\ell,p}(\omega)\}^\T\,.
\end{align*}
Notice that  
\begin{align*}
\mathcal{T}_n=q\{\sqrt{\nu}S_{n,x}^{(1)}+\sqrt{1-\nu}S_{n,y}^{(1)}\}-q\{S_{n,w}^{(1)}\}=\frac{1}{2}\sum_{j=1}^p\sum_{\ell=1}^L\int_0^1\partial_{j}q\{Z(\omega)\}\dot{Z}_{\ell,j}(\omega)\,{\rm d}\omega\,.
\end{align*}

Write $Z_\ell(\omega)=\{Z_{\ell,1}(\omega),\ldots,Z_{\ell,p}(\omega)\}^\T$. By Taylor expansion,
\begin{align*}
2\mathbb{E}(\mathcal{T}_n)=&\,\,\underbrace{\sum_{j=1}^p\sum_{\ell=1}^L\int_0^1\mathbb{E}[\partial_{j}q\{Z^{(-\ell)}(\omega)\}\dot{Z}_{\ell,j}(\omega)]\,{\rm d}\omega}_{\textrm{I}}\\
&+\underbrace{\sum_{j,k=1}^p\sum_{\ell=1}^L\int_0^1\mathbb{E}[\partial_{jk}q\{Z^{(-\ell)}(\omega)\}Z_{\ell,k}(\omega)\dot{Z}_{\ell,j}(\omega)]\,{\rm d}\omega}_{\textrm{II}}\\
&+\underbrace{\sum_{j,k,l=1}^p\sum_{\ell=1}^L\int_0^1\int_0^1(1-\tau)\mathbb{E}[\partial_{jkl}q\{Z^{(-\ell)}(\omega)+\tau Z_\ell(\omega)\} Z_{\ell,k}(\omega)Z_{\ell,l}(\omega)\dot{Z}_{\ell,j}(\omega)]\,{\rm d}\tau{\rm d}\omega}_{\textrm{III}}\,.
\end{align*}
In the sequel, we will bound the three terms $\textrm{I}$, $\textrm{II}$ and $\textrm{III}$, respectively. To simplify the notation and without causing much confusion, we write $Z_\ell(\omega)$, $Z^{(-\ell)}(\omega)$, $Z_{\ell,j}(\omega)$ and $\dot{Z}_{\ell,j}(\omega)$ as $Z_\ell$, $Z^{(-\ell)}$, $Z_{\ell,j}$ and $\dot{Z}_{\ell,j}$, respectively.

Note that
$
\sqrt{L}\partial_{j}q\{Z^{(-\ell)}\}\dot{Z}_{\ell,j}=\sqrt{\nu}\omega^{-1/2}\tilde{X}_{\ell,j}\partial_{j}q\{Z^{(-\ell)}\}+\sqrt{1-\nu}\omega^{-1/2}\tilde{Y}_{\ell,j}\partial_{j}q\{Z^{(-\ell)}\}-(1-\omega)^{-1/2}\tilde{W}_{\ell,j}\partial_{j}q\{Z^{(-\ell)}\}$.
Since $\tilde{Y}_\ell$ is independent of $\{\tilde{X}_\ell\}_{\ell=1}^L$, $\{\tilde{W}_\ell\}_{\ell=1}^L$ and $\{\tilde{Y}_s\}_{s\neq \ell}$, we know $\mathbb{E}[\partial_{j}q\{Z^{(-\ell)}\}\tilde{Y}_{\ell,j}]=0$, which implies
$
\sum_{j=1}^p\sum_{\ell=1}^L\mathbb{E}[\partial_{j}q\{Z^{(-\ell)}\}\tilde{Y}_{\ell,j}]=0$.
Analogously, we also have $
\sum_{j=1}^p\sum_{\ell=1}^L\mathbb{E}[\partial_{j}q\{Z^{(-\ell)}\}\tilde{W}_{\ell,j}]=0$. Thus,
\begin{align*}
{\rm I}=\frac{\sqrt{\nu}}{\sqrt{L}}\sum_{j=1}^p\sum_{\ell=1}^L\int_0^1\omega^{-1/2}\mathbb{E}[\tilde{X}_{\ell,j}\partial_{j}q\{Z^{(-\ell)}\}]\,{\rm d}\omega\,.
\end{align*}
Let $\mathcal{F}_{-\ell}^*$ be the $\sigma$-filed generated by $\{\tilde{X}_s,\tilde{Y}_s,\tilde{W}_s\}_{s\neq \ell}$. Note that $\sum_{j=1}^p|\partial_{j}q(v)|\leq C\phi$ for any $v\in\mathbb{R}^p$, and $\mathbb{E}[\tilde{X}_{\ell,j}\partial_{j}q\{Z^{(-\ell)}\}]=\mathbb{E}[\partial_{j}q\{Z^{(-\ell)}\}\mathbb{E}(\tilde{X}_{\ell,j}\,|\,\mathcal{F}_{-\ell}^*)]$. Recall $\mathcal{F}_{-\ell}$ is the $\sigma$-filed generated by $\{\tilde{X}_s\}_{s\neq \ell}$. Since $\tilde{X}_\ell$ is independent of $\{\tilde{Y}_s,\tilde{W}_s\}_{s\neq \ell}$, we have
\begin{align}\label{eq:I}
|\textrm{I}|
\lesssim&\,\,\frac{1}{\sqrt{L}}\sum_{\ell=1}^L\sum_{j=1}^p\int_0^1\frac{1}{\sqrt{\omega}}\mathbb{E}\bigg[|\partial_{j}q\{Z^{(-\ell)}\}|\max_{j\in[p]}|\mathbb{E}(\tilde{X}_{\ell,j}\,|\,\mathcal{F}_{-\ell}^*)|\bigg]\,{\rm d}\omega\notag\\
\lesssim&\,\,\phi \sqrt{L}\max_{\ell\in[L]}\mathbb{E}\bigg\{\max_{j\in[p]}|\mathbb{E}(\tilde{X}_{\ell,j}\,|\,\mathcal{F}_{-\ell}^*)|\bigg\}=\phi \sqrt{L}\max_{\ell\in[L]}\mathbb{E}\bigg\{\max_{j\in[p]}|\mathbb{E}(\tilde{X}_{\ell,j}\,|\,\mathcal{F}_{-\ell})|\bigg\}\,.
\end{align}
For $\textrm{II}$, 
since $\tilde{Y}_\ell$ is independent of $\tilde{W}_\ell$, and $(\tilde{Y}_\ell,\tilde{W}_\ell)$ is independent of $\{\tilde{Y}_s,\tilde{W}_s\}_{s\neq\ell}$ and  $\{\tilde{X}_\ell\}_{\ell=1}^L$,  
\begin{align*}
&\mathbb{E}[\partial_{jk}q\{Z^{(-\ell)}\}Z_{\ell,k}\dot{Z}_{\ell,j}]\\
&~~~=\nu L^{-1}\mathbb{E}[\partial_{jk}q\{Z^{(-\ell)}\}\tilde{X}_{\ell,k}\tilde{X}_{\ell,j}]+(1-\nu)L^{-1}\mathbb{E}[\partial_{jk}q\{Z^{(-\ell)}\}]\mathbb{E}(\tilde{Y}_{\ell,j}\tilde{Y}_{\ell,k})\\
&~~~~-L^{-1}\mathbb{E}[\partial_{jk}q\{Z^{(-\ell)}\}]\mathbb{E}(\tilde{W}_{\ell,j}\tilde{W}_{\ell,k})\,.
\end{align*}
Due to $\mathbb{E}(\tilde{X}_{\ell,j}\tilde{X}_{\ell,k})=\mathbb{E}(\tilde{Y}_{\ell,j}\tilde{Y}_{\ell,k})=\mathbb{E}(\tilde{W}_{\ell,j}\tilde{W}_{\ell,k})$, it holds that
$$
\mathbb{E}[\partial_{jk}q\{Z^{(-\ell)}\}Z_{\ell,k}\dot{Z}_{\ell,j}]=\nu L^{-1}\mathbb{E}[\partial_{jk}q\{Z^{(-\ell)}\}\mathbb{E}\{\tilde{X}_{\ell,k}\tilde{X}_{\ell,j}-\mathbb{E}(\tilde{X}_{\ell,k}\tilde{X}_{\ell,j})\,|\,\mathcal{F}_{-\ell}^*\}]\,.$$
Lemmas A.5 and A.6 of \cite{CCK_2013} show that $\sum_{j,k=1}^p|\partial_{jk}q(v)|\lesssim\phi^2+\phi\beta$ for any $v\in\mathbb{R}^p$. Then
\begin{align}\label{eq:II}
|\textrm{II}|\lesssim&\,\,(\phi^2+\phi\beta)\max_{\ell\in[L]}\mathbb{E}\bigg[\max_{k,j\in[p]}|\mathbb{E}\{\tilde{X}_{\ell,k}\tilde{X}_{\ell,j}-\mathbb{E}(\tilde{X}_{\ell,k}\tilde{X}_{\ell,j})\,|\,\mathcal{F}_{-\ell}^*\}|\bigg]\notag\\
\lesssim&\,\,\phi^2(\log p)\max_{\ell\in[L]}\mathbb{E}\bigg[\max_{k,j\in[p]}|\mathbb{E}\{\tilde{X}_{\ell,k}\tilde{X}_{\ell,j}-\mathbb{E}(\tilde{X}_{\ell,k}\tilde{X}_{\ell,j})\,|\,\mathcal{F}_{-\ell}\}|\bigg]\,.
\end{align}

To bound $\textrm{III}$, we first define
$
\chi_\ell=I\{\max_{j\in[p]}(|\tilde{X}_{\ell,j}|\vee|\tilde{Y}_{\ell,j}|\vee|\tilde{W}_{\ell,j}|)\leq \sqrt{L}/(4\beta)\}$. Then
\begin{align*}
\textrm{III}=&\,\underbrace{\sum_{j,k,l=1}^p\sum_{\ell=1}^L\int_0^1\int_0^1(1-\tau)\mathbb{E}[\chi_\ell \partial_{jkl}q\{Z^{(-\ell)}+\tau Z_\ell\}Z_{\ell,k}Z_{\ell,l}\dot{Z}_{\ell,j}]\,{\rm d}\tau{\rm d}\omega}_{\textrm{III}_1}\\
&+\underbrace{\sum_{j,k,l=1}^p\sum_{\ell=1}^L\int_0^1\int_0^1(1-\tau)\mathbb{E}[(1-\chi_\ell) \partial_{jkl}q\{Z^{(-\ell)}+\tau Z_\ell\}Z_{\ell,k}Z_{\ell,l}\dot{Z}_{\ell,j}]\,{\rm d}\tau{\rm d}\omega}_{\textrm{III}_2}\,.
\end{align*}
Define
$
\kappa(\omega)=(\sqrt{\omega}\wedge\sqrt{1-\omega})^{-1}$ and $
h(v,\omega)=I\{-\phi^{-1}-\omega/\beta<\max_{j\in[p]}(v_j-u_j)\leq \phi^{-1}+\omega/\beta\}$
for any $v=(v_1,\ldots,v_p)^\T\in\mathbb{R}^p$ and $\omega>0$. Lemmas A.5 and A.6 of \cite{CCK_2013} show that $\sum_{j,k,l=1}^p|\partial_{jkl}q(v)|\lesssim\phi\beta^2$ for any $v\in\mathbb{R}^p$. Then
\begin{align}\label{eq:III2alpha}
|{\rm III}_2|\lesssim \phi\beta^2\sum_{\ell=1}^L\int_0^1\mathbb{E}\bigg\{(1-\chi_\ell)\max_{j,k,l\in[p]}|Z_{\ell,k}Z_{\ell,l}\dot{Z}_{\ell,j}|\bigg\}\,{\rm d}\omega\,.
\end{align}
Observe that
$
\max_{j,k,l\in[p]}|Z_{\ell,k}Z_{\ell,l}\dot{Z}_{\ell,j}|\lesssim \kappa(\omega)L^{-3/2}(|\tilde{X}_{\ell}|_\infty^3\vee|\tilde{Y}_{\ell}|_\infty^3\vee|\tilde{W}_{\ell}|_\infty^3)$
and
$
1-\chi_\ell\leq I\{|\tilde{X}_{\ell}|_\infty>\sqrt{L}/(4\beta)\}+I\{|\tilde{Y}_{\ell}|_\infty>\sqrt{L}/(4\beta)\}+I\{|\tilde{W}_{\ell}|_\infty>\sqrt{L}/(4\beta)\}$.
Thus, by Chebyshev's association inequality \citep[Lemma B.1]{CCK_2017}, 
\begin{align}\label{eq:bdident}
\mathbb{E}\bigg\{(1-\chi_\ell)\max_{j,k,l\in[p]}|Z_{\ell,k}Z_{\ell,l}\dot{Z}_{\ell,j}|\bigg\}\lesssim \frac{\kappa(\omega)}{L^{3/2}}\{M_{\tilde{x}}(\phi)+M_{\tilde{y}}(\phi)\}
\end{align}
with $M_{\tilde{x}}(\cdot)$ and $M_{\tilde{y}}(\cdot)$ defined in \eqref{eq:mtildex} and \eqref{eq:mtildey}, respectively,
which implies
\begin{align}\label{eq:III2}
|{\rm III}_2|\lesssim \frac{\phi\beta^2}{\sqrt{L}}\{M_{\tilde{x}}(\phi)+M_{\tilde{y}}(\phi)\}=\frac{\phi^3(\log p)^2}{\sqrt{L}}\{M_{\tilde{x}}(\phi)+M_{\tilde{y}}(\phi)\}\,.
\end{align}
If $\chi_\ell=1$, then
$
\max_{j\in[p]}|Z_{\ell,j}|
\leq(4\beta)^{-1}\{\sqrt{\omega}(\sqrt{\nu}+\sqrt{1-\nu})+\sqrt{1-\omega}\}\leq(4\beta)^{-1}\sqrt{3}$
for any $\omega\in[0,1]$.
If $h\{Z^{(-\ell)},1\}=0$, then $\max_{j\in[p]}\{Z_{j}^{(-\ell)}-u_j\}\leq -\phi^{-1}-\beta^{-1}$ or $\max_{j\in[p]}\{Z_{j}^{(-\ell)}-u_j\}>\phi^{-1}+\beta^{-1}$. When $\max_{j\in[p]}\{Z_{j}^{(-\ell)}-u_j\}\leq -\phi^{-1}-\beta^{-1}$ and $\chi_\ell=1$, 
$
F_\beta\{Z^{(-\ell)}+\tau Z_\ell\}\leq -\beta^{-1}(1-\sqrt{3}/4)
$ for any $\tau\in[0,1]$, which implies that $q\{Z^{(-\ell)}+\tau Z_{\ell}\}=1$ for any $\omega\in[0,1]$ and $\tau\in[0,1]$. When $\max_{j\in[p]}\{Z_{j}^{(-\ell)}-u_j\}> \phi^{-1}+\beta^{-1}$ and $\chi_\ell=1$, $F_\beta\{Z^{(-\ell)}+\tau Z_{\ell}\}>\phi^{-1}+\beta^{-1}(1-\sqrt{3}/4)$ for any $\tau\in[0,1]$, which implies $q\{Z^{(-\ell)}+\tau Z_{\ell}\}=0$ for any $\omega\in[0,1]$ and $\tau\in[0,1]$. Therefore, if $\chi_\ell=1$ and $h\{Z^{(-\ell)},1\}=0$, we have $\partial_{jkl}q\{Z^{(-\ell)}+\tau Z_{\ell}\}=0$ for any $\omega\in[0,1], \tau\in[0,1]$ and $j,k,l\in[p]$. Lemmas A.5 and A.6 of \cite{CCK_2013} indicate that there exist $U_{jkl}(v)$ such that $|\partial_{jkl}q(v)|\leq U_{jkl}(v)$ and $\sum_{j,k,l=1}^pU_{jkl}(v)\lesssim\phi\beta^2$ for any $v\in\mathbb{R}^p$. Then we have 
\begin{align*}
\chi_\ell |\partial_{jkl}q\{Z^{(-\ell)}+\tau Z_{\ell}\}|=&\,\,\chi_\ell h\{Z^{(-\ell)},1\} |\partial_{jkl}q\{Z^{(-\ell)}+\tau Z_{\ell}\}|\\
\lesssim&\,\,\chi_\ell h\{Z^{(-\ell)},1\} U_{jkl}\{Z^{(-\ell)}\}\leq h\{Z^{(-\ell)},1\} U_{jkl}\{Z^{(-\ell)}\}\,,
\end{align*}
which implies
\begin{align}\label{eq:III1}
|\textrm{III}_1|
\lesssim&\sum_{j,k,l=1}^p\sum_{\ell=1}^L\int_0^1\mathbb{E}[h\{Z^{(-\ell)},1\} U_{jkl}\{Z^{(-\ell)}\}|Z_{\ell,k}Z_{\ell,l}\dot{Z}_{\ell,j}|]\,{\rm d}\omega\notag\\
\lesssim&\sum_{j,k,l=1}^p\sum_{\ell=1}^L\int_0^1\frac{\kappa(\omega)}{L^{3/2}}\mathbb{E}[h\{Z^{(-\ell)},1\} U_{jkl}\{Z^{(-\ell)}\}(|\tilde{X}_{\ell,k}|^3+|\tilde{Y}_{\ell,k}|^3+|\tilde{W}_{\ell,k}|^3)]\,{\rm d}\omega\,.
\end{align}
In the sequel, we will show
\begin{align}\label{eq:alphatoprove1}
R:=&\,\sum_{j,k,l=1}^p\sum_{\ell=1}^L\int_0^1\frac{\kappa(\omega)}{L^{3/2}}\mathbb{E}[h\{Z^{(-\ell)},1\} U_{jkl}\{Z^{(-\ell)}\}(|\tilde{X}_{\ell,k}|^3+|\tilde{Y}_{\ell,k}|^3+|\tilde{W}_{\ell,k}|^3)]\,{\rm d}\omega\notag\\
\lesssim&\,\,\frac{\phi^3(\log p)^2}{\sqrt{L}}\{B_{n}^3\phi^{-1}(\log p)^{1/2}+B_{n}^3\varrho_n^{(1)}+M_{\tilde{x}}(\phi)+M_{\tilde{y}}(\phi)\}\\
&~~~~~~+\frac{\phi^3(\log p)^2}{\sqrt{L}} \max_{\ell\in[L]}\mathbb{E}\bigg[\max_{k\in[p]}|\mathbb{E}\{|\tilde{X}_{\ell,k}|^3-\mathbb{E}(|\tilde{X}_{\ell,k}|^3)\,|\,\mathcal{F}_{-\ell}\}|\bigg]\,.\notag
\end{align}
Recall that ${\rm III}={\rm III}_1+{\rm III}_2$. Together with (\ref{eq:III2}) and \eqref{eq:III1}, we have
\begin{align*}
|{\rm III}|\lesssim&\,\,\frac{\phi^3(\log p)^2}{\sqrt{L}}\{B_{n}^3\phi^{-1}(\log p)^{1/2}+B_{n}^3\varrho_n^{(1)}+M_{\tilde{x}}(\phi)+M_{\tilde{y}}(\phi)\}\\
&+\frac{\phi^3(\log p)^2}{\sqrt{L}} \max_{\ell\in[L]}\mathbb{E}\bigg[\max_{j\in[p]}|\mathbb{E}\{|\tilde{X}_{\ell,j}|^3-\mathbb{E}(|\tilde{X}_{\ell,j}|^3)\,|\,\mathcal{F}_{-\ell}\}|\bigg]\,.
\end{align*}
Due to $\mathbb{E}(\mathcal{T}_n)=2^{-1}({\rm I}+{\rm II}+{\rm III})$, together with (\ref{eq:I}) and (\ref{eq:II}), we have \eqref{eq:alphatoshow1}.

Now we begin to prove \eqref{eq:alphatoprove1}. Write $u_{\ell,k}=|\tilde{X}_{\ell,k}|^3+|\tilde{Y}_{\ell,k}|^3+|\tilde{W}_{\ell,k}|^3-\mathbb{E}(|\tilde{X}_{\ell,k}|^3)-\mathbb{E}(|\tilde{Y}_{\ell,k}|^3)-\mathbb{E}(|\tilde{W}_{\ell,k}|^3)$. Then
\begin{align}\label{eq:sd}
R=&\underbrace{\sum_{j,k,l=1}^p\sum_{\ell=1}^L\int_0^1\frac{\kappa(\omega)}{L^{3/2}}\mathbb{E}[\chi_\ell h\{Z^{(-\ell)},1\} U_{jkl}\{Z^{(-\ell)}\}]\mathbb{E}(|\tilde{X}_{\ell,k}|^3+|\tilde{Y}_{\ell,k}|^3+|\tilde{W}_{\ell,k}|^3)\,{\rm d}\omega}_{R_1}\notag\\
&+\underbrace{\sum_{j,k,l=1}^p\sum_{\ell=1}^L\int_0^1\frac{\kappa(\omega)}{L^{3/2}}\mathbb{E}[(1-\chi_\ell)h\{Z^{(-\ell)},1\}\notag U_{jkl}\{Z^{(-\ell)}\}]\mathbb{E}(|\tilde{X}_{\ell,k}|^3+|\tilde{Y}_{\ell,k}|^3+|\tilde{W}_{\ell,k}|^3)\,{\rm d}\omega}_{R_2}\notag\\
&+\underbrace{\sum_{j,k,l=1}^p\sum_{\ell=1}^L\int_0^1\frac{\kappa(\omega)}{L^{3/2}}\mathbb{E}[h\{Z^{(-\ell)},1\} U_{jkl}\{Z^{(-\ell)}\}u_{\ell,k}]\,{\rm d}\omega}_{R_3}\,.
\end{align}
Observe that if $\chi_\ell=1$ and $h(Z,2)=0$, then $h\{Z^{(-\ell)},1\}=0$. Due to $\sum_{j,k,l=1}^pU_{jkl}\{Z^{(-\ell)}\}\lesssim\phi\beta^2$ for any $\omega\in[0,1]$ and $\ell\in[L]$, by Lemma \ref{la:mom} and the fact that $\mathbb{E}(|\tilde Y_{\ell,j}|^3)\lesssim\{\mathbb{E}( |\tilde Y_{\ell,j}|^2)\}^{3/2}=\{\mathbb{E}(|\tilde X_{\ell,j}|^2)\}^{3/2}\leq\mathbb{E}(|\tilde X_{\ell,j}|^3)$, we have
\begin{align}\label{eq:III111}
R_1\leq&\,\sum_{j,k,l=1}^p\sum_{\ell=1}^L\int_0^1\frac{\kappa(\omega)}{L^{3/2}}\mathbb{E}[\chi_\ell h(Z,2)U_{jkl}\{Z^{(-\ell)}\}]\mathbb{E}(|\tilde{X}_{\ell,k}|^3+2|\tilde{Y}_{\ell,k}|^3)\,{\rm d}\omega\notag\\
\lesssim&\,\,\frac{\phi\beta^2}{L^{3/2}}\sum_{\ell=1}^L\int_0^1\kappa(\omega)\mathbb{E}\{h(Z,2)\}\max_{k\in[p]}\mathbb{E}(|\tilde{X}_{\ell,k}|^3+|\tilde{Y}_{\ell,k}|^3)\,{\rm d}\omega\\
\lesssim&\,\,\frac{\phi\beta^2B_{n}^3}{\sqrt{L}}\int_0^1\kappa(\omega)\mathbb{E}\{h(Z,2)\}\,{\rm d}\omega\,.\notag
\end{align}
Recall that $\varrho_n^{(1)}=\sup_{u\in\mathbb{R}^p,\nu\in[0,1]}|\mathbb{P}\{\sqrt{\nu}S_{n,x}^{(1)}+\sqrt{1-\nu}S_{n,y}^{(1)}\leq u\}-\mathbb{P}\{S_{n,y}^{(1)}\leq u\}|$ with $S_{n,x}^{(1)}=L^{-1/2}\sum_{\ell=1}^L\tilde{X}_\ell$ and $S_{n,y}^{(1)}=L^{-1/2}\sum_{\ell=1}^L\tilde{Y}_\ell$. Since
$
Z=L^{-1/2}\sum_{\ell=1}^L\{\sqrt{\omega\nu}\tilde{X}_\ell+\sqrt{\omega(1-\nu)}\tilde{Y}_\ell+\sqrt{1-\omega}\tilde{W}_\ell\}\overset{d}{=}\sqrt{\omega\nu}S_{n,x}^{(1)}+\sqrt{1-\omega\nu}{S}_{n,y}^{(1)}$,
we have $\mathbb{P}(Z\leq u+\phi^{-1}+2\beta^{-1})\leq \mathbb{P}\{S_{n,y}^{(1)}\leq u+\phi^{-1}+2\beta^{-1}\}+\varrho_n^{(1)}$ and $\mathbb{P}(Z\leq u-\phi^{-1}-2\beta^{-1})\geq \mathbb{P}\{S_{n,y}^{(1)}\leq u-\phi^{-1}-2\beta^{-1}\}-\varrho_n^{(1)}$. Therefore,
\begin{align*}
\mathbb{E}\{h(Z,2)\}=&\,\,\mathbb{P}(Z\leq u+\phi^{-1}+2\beta^{-1})-\mathbb{P}(Z\leq u-\phi^{-1}-2\beta^{-1})\\
\leq&\,\,\mathbb{P}\{S_{n,y}^{(1)}\leq u+\phi^{-1}+2\beta^{-1}\}-\mathbb{P}\{S_{n,y}^{(1)}\leq u-\phi^{-1}-2\beta^{-1}\}+2\varrho_n^{(1)}\\
\lesssim&\,\,\phi^{-1}(\log p)^{1/2}+\varrho_n^{(1)}\,,
\end{align*}
where the last step is based on Nazarov's inequality. By \eqref{eq:III111}, we have
\begin{align}\label{eq:III11}
R_1
\lesssim\frac{\phi^3B_{n}^3(\log p)^2}{\sqrt{L}}\{\phi^{-1}(\log p)^{1/2}+\varrho_n^{(1)}\}\,.
\end{align}
On the other hand, it holds that
\begin{align}\label{eq:bdterm1}
R_2\lesssim&\,\,\frac{\phi\beta^2}{L^{3/2}}\sum_{\ell=1}^L\int_0^1\kappa(\omega)\mathbb{E}[(1-\chi_\ell) h\{Z^{(-\ell)},1\}]\mathbb{E}\bigg\{\max_{k\in[p]}(|\tilde{X}_{\ell,k}|^3\vee|\tilde{Y}_{\ell,k}|^3\vee|\tilde{W}_{\ell,k}|^3)\bigg\}\,{\rm d}\omega\notag\\
\leq&\,\,\frac{\phi\beta^2}{L^{3/2}}\sum_{\ell=1}^L\int_0^1\kappa(\omega)\mathbb{E}(1-\chi_\ell)\mathbb{E}\bigg\{\max_{k\in[p]}(|\tilde{X}_{\ell,k}|^3\vee|\tilde{Y}_{\ell,k}|^3\vee|\tilde{W}_{\ell,k}|^3)\bigg\}\,{\rm d}\omega\\
\leq&\,\,\frac{\phi\beta^2}{L^{3/2}}\sum_{\ell=1}^L\int_0^1\kappa(\omega)\mathbb{E}\bigg\{(1-\chi_\ell)\max_{k\in[p]}(|\tilde{X}_{\ell,k}|^3\vee|\tilde{Y}_{\ell,k}|^3\vee|\tilde{W}_{\ell,k}|^3)\bigg\}\,{\rm d}\omega\notag\\
\lesssim&\,\,\frac{\phi^3(\log p)^2}{\sqrt{L}}\{M_{\tilde{x}}(\phi)+M_{\tilde{y}}(\phi)\}\,,\notag
\end{align}
where the third step and last step are based on Chebyshev's association inequality \citep[Lemma B.1]{CCK_2017}.
Recall $\mathcal{F}_{-\ell}^*$ is the $\sigma$-filed generated by $\{\tilde{X}_s,\tilde{Y}_s,\tilde{W}_s\}_{s\neq \ell}$. Since $(\tilde{Y}_\ell,\tilde{W}_\ell)$ is independent of $\{\tilde{X}_s,\tilde{Y}_s,\tilde{W}_s\}_{s\neq \ell}$, and $\{\tilde{X}_\ell\}_{\ell=1}^L$ is independent of $\{\tilde{Y}_s,\tilde{W}_s\}_{s\neq \ell}$, we have
\begin{align*}
R_3=&\sum_{j,k,l=1}^p\sum_{\ell=1}^L\int_0^1\frac{\kappa(\omega)}{L^{3/2}}\mathbb{E}[h\{Z^{(-\ell)},1\} U_{jkl}\{Z^{(-\ell)}\}\{|\tilde{X}_{\ell,k}|^3-\mathbb{E}(|\tilde{X}_{\ell,k}|^3)\}]\,{\rm d}\omega\\
=&\sum_{j,k,l=1}^p\sum_{\ell=1}^L\int_0^1\frac{\kappa(\omega)}{L^{3/2}}\mathbb{E}[h\{Z^{(-\ell)},1\} U_{jkl}\{Z^{(-\ell)}\}\mathbb{E}\{|\tilde{X}_{\ell,k}|^3-\mathbb{E}(|\tilde{X}_{\ell,k}|^3)\,|\,\mathcal{F}_{-\ell}^*\}]\,{\rm d}\omega\\
\lesssim&\,\,\frac{\phi^3(\log p)^2}{\sqrt{L}}\max_{\ell\in[L]}\mathbb{E}\bigg[\max_{k\in[p]}|\mathbb{E}\{|\tilde{X}_{\ell,k}|^3-\mathbb{E}(|\tilde{X}_{\ell,k}|^3)\,|\,\mathcal{F}_{-\ell}^*\}|\bigg]\\
=&\,\,\frac{\phi^3(\log p)^2}{\sqrt{L}}\max_{\ell\in[L]}\mathbb{E}\bigg[\max_{k\in[p]}|\mathbb{E}\{|\tilde{X}_{\ell,k}|^3-\mathbb{E}(|\tilde{X}_{\ell,k}|^3)\,|\,\mathcal{F}_{-\ell}\}|\bigg]\,,
\end{align*}
where $\mathcal{F}_{-\ell}$ is the $\sigma$-field generated by $\{\tilde{X}_s\}_{s\neq\ell}$. Together with \eqref{eq:III11} and \eqref{eq:bdterm1},  (\ref{eq:sd}) implies \eqref{eq:alphatoprove1} holds. $\hfill\Box$

\section{Proofs of the lemmas used in Section \ref{sec:5.2}}

\subsection{Proof of Lemma \ref{lem:mIII}} \label{sec:pfmIII}
Let
\begin{align*}
\chi_t=I\bigg\{\max_{j\in[p],s\in\cup_{\ell\in\mathcal{N}_t}\mathcal{N}_\ell}(|X_{s,j}|\vee|Y_{s,j}|\vee|W_{s,j}|)\leq \frac{\sqrt{n}}{8\beta (D_n^4D_{n}^*)^{1/3}}\bigg\}
\end{align*}
and define
\begin{align*}
{\rm III}_{1} &= \sum_{j,k,l=1}^{p} \sum_{t=1}^{n} \int_{0}^{1} \int_{0}^{1} (1-\tau) \mathbb{E}[ \chi_{t} \partial_{jkl}q\{Z^{(-t)} + \tau \delta_t\} \delta_{t,k} \delta_{t,l} \dot{Z}_{t,j}]\,{\rm d}\tau {\rm d}\omega\,,\\
{\rm III}_2 &=\sum_{j,k,l=1}^{p} \sum_{t=1}^{n} \int_{0}^{1} \int_{0}^{1} (1-\tau) \mathbb{E}[(1-\chi_{t}) \partial_{jkl}q\{Z^{(-t)} + \tau \delta_t\} \delta_{t,k} \delta_{t,l} \dot{Z}_{t,j}]\,{\rm d}\tau {\rm d}\omega\,.\end{align*}
Identical to \eqref{eq:III2alpha} in Section \ref{sec:ETnalpha},
$$
|{\rm III}_{2}|\lesssim \phi \beta^{2} \sum_{t=1}^{n} \int_{0}^{1} \mathbb{E}\bigg\{(1-\chi_{t}) \max_{j,k,l \in[p]} |\delta_{t,k}\delta_{t,l} \dot{Z}_{t,j}|\bigg\}\,{\rm d}\omega\,.$$
Define $\kappa(\omega)=(\sqrt{\omega}\wedge\sqrt{1-\omega})^{-1}$. Observe that
$
|\delta_{t,k}\delta_{t,l}\dot{Z}_{t,j}| \leq \kappa(\omega) n^{-3/2}\sum_{s,s' \in \mathcal{N}_{t}} (|X_{s,k}| + |Y_{s,k}| + |W_{s,k}|)(|X_{s',l}| + |Y_{s',l}| + |W_{s',l}|) (|X_{t,j}| + |Y_{t,j}| + |W_{t,j}|)$ so that
\begin{align*}
\max_{j,k,l\in[p]} |\delta_{t,k}\delta_{t,l} \dot{Z}_{t,j}| \lesssim&\,\,{D_{n}^{2} \kappa(\omega) \over n^{3/2}} \bigg\{ \max_{j\in[p],s\in\cup_{\ell\in\mathcal{N}_t}\mathcal{N}_\ell} |X_{s,j}|^{3} \\
&\quad\quad\quad\quad\quad+ \max_{j\in[p],s\in\cup_{\ell\in\mathcal{N}_t}\mathcal{N}_\ell} |Y_{s,j}|^{3} + \max_{j\in[p],s\in\cup_{\ell\in\mathcal{N}_t}\mathcal{N}_\ell} |W_{s,j}|^{3} \bigg\}\,.
\end{align*}
Identical to \eqref{eq:bdident} in Section \ref{sec:ETnalpha},
$\mathbb{E}\{(1-\chi_{t}) \max_{j,k,l \in [p]} |\delta_{t,k}\delta_{t,l} \dot{Z}_{t,j}|\}\lesssim D_{n}^{2}n^{-3/2} \kappa(\omega)M_{n}(\phi)$
with $M_n(\phi)$ defined in \eqref{eq:Mphixm}. Since $\int_{0}^{1} \kappa(\omega)\,{\rm d}\omega \lesssim 1$, then
$
|{\rm III}_{2}| \lesssim D_{n}^{2}n^{-1/2} M_{n}(\phi) \phi \beta^{2} $.

Next, we deal with ${\rm III}_{1}$. Define $\delta^{(\mathcal{N}_t)}=\sum_{s\in\cup_{\ell\in\mathcal{N}_t}\mathcal{N}_\ell}Z_s$ and $$Z^{(-\mathcal{N}_t)}=Z-\delta^{(\mathcal{N}_t)}=Z^{(-t)}+\delta_t-\delta^{(\mathcal{N}_t)}\,.$$ If $\chi_t=1$, we have
$
|(\tau-1)\delta_t+\delta^{(\mathcal{N}_t)}|_\infty\leq\sqrt{3}n^{-1/2}\sum_{s\in\cup_{\ell\in\mathcal{N}_t}\mathcal{N}_\ell}\max_{j\in[p]}(|X_{s,j}|\vee|Y_{s,j}|\vee|W_{s,j}|)\leq\sqrt{3}/(8\beta)<\beta^{-1}$.
Recall $Z^{(-t)}+\tau\delta_t=Z^{(-\mathcal{N}_t)}+(\tau-1)\delta_t+\delta^{(\mathcal{N}_t)}$.
For $h(v,\omega)$ specified in Section \ref{sec:ETnalpha}, 
same as the arguments for the first step of \eqref{eq:III1} in Section \ref{sec:ETnalpha}, 
\begin{align}\label{eq:bdIII1m}
|{\rm III}_{1}| \lesssim&\,\sum_{j,k,l=1}^{p} \sum_{t=1}^{n} \int_{0}^{1} \mathbb{E}[h\{Z^{(-\mathcal{N}_t)},1\}U_{jkl}\{Z^{(-\mathcal{N}_t)}\}|\delta_{t,k}\delta_{t,l} \dot{Z}_{t,j}|]\,{\rm d}\omega \\
=& \sum_{j,k,l=1}^{p} \sum_{t=1}^{n} \int_{0}^{1}\mathbb{E}[ h\{Z^{(-\mathcal{N}_t)},1\}U_{jkl}\{Z^{(-\mathcal{N}_t)}\}] \mathbb{E}(|\delta_{t,k}\delta_{t,l} \dot{Z}_{t,j}|)\,{\rm d}\omega\,,\notag
\end{align}
where $U_{jkl}(v)$ is specified in Section \ref{sec:ETnalpha}, and the second step follows from the independence between $Z^{(-\mathcal{N}_t)}$ and $\delta_{t,k}\delta_{t,l}\dot{Z}_{t,j}$. It holds that
\begin{align*}
&\sum_{j,k,l=1}^{p} \sum_{t=1}^{n} \int_{0}^{1}\mathbb{E}[ (1-\chi_{t}) h\{Z^{(-\mathcal{N}_t)},1\}U_{jkl}\{Z^{(-\mathcal{N}_t)}\}]\mathbb{E}(|\delta_{t,k}\delta_{t,l} \dot{Z}_{t,j}|)\,{\rm d}\omega \\
&~~~\lesssim \phi\beta^2\sum_{t=1}^{n} \int_{0}^{1}\mathbb{E}(1-\chi_{t}) \mathbb{E}\bigg(\max_{j,k,l\in[p]}|\delta_{t,k}\delta_{t,l} \dot{Z}_{t,j}|\bigg)\,{\rm d}\omega\lesssim {D_{n}^{2} M_{n}(\phi) \phi \beta^{2} \over n^{1/2}}\,,\end{align*}
where the last step is obtained by Chebyshev's association inequality.
Observe that if $\chi_{t}=1$ and $h(Z,2)=0$, then $h\{Z^{(-\mathcal{N}_t)},1\}=0$. Thus we have
\begin{align*}
& \sum_{j,k,l=1}^{p} \sum_{t=1}^{n}\int_0^1 \mathbb{E}[ \chi_{t} h\{Z^{(-\mathcal{N}_t)},1\}U_{jkl}\{Z^{(-\mathcal{N}_t)}\}]\mathbb{E}(|\delta_{t,k}\delta_{t,l}\dot{Z}_{t,j}|)\,{\rm d}\omega \\
&~~~\leq\sum_{j,k,l=1}^{p} \sum_{t=1}^{n}\int_0^1 \mathbb{E}[ \chi_{t} h(Z,2) U_{jkl}\{Z^{(-\mathcal{N}_t)}\}]\mathbb{E}(|\delta_{t,k}\delta_{t,l}\dot{Z}_{t,j}|)\,{\rm d}\omega \\
&~~~\lesssim\phi \beta^{2} \sum_{t=1}^n\int_0^1 \mathbb{E}\{h(Z,2)\}\max_{j,k,l\in[p]}\mathbb{E}(|\delta_{t,k}\delta_{t,l} \dot{Z}_{t,j}|)\,{\rm d}\omega\,.
\end{align*}
Notice that $\max_{t\in [n],j\in [p]}\mathbb{E}(|X_{t,j}|^{3})\lesssim B_{n}^{3}$. Hence, $\mathbb{E}(|W_{t,j}|^3)=\mathbb{E}(|Y_{t,j}|^3)\lesssim\{\mathbb{E}( |Y_{t,j}|^2)\}^{3/2}=\{\mathbb{E}(|X_{t,j}|^2)\}^{3/2}\leq\mathbb{E}(|X_{t,j}|^3)\lesssim B_{n}^{3}$ holds uniformly over $t\in [n]$ and $j\in [p]$,
which implies 
$$
\max_{j,k,l \in[p], t\in[n]}\mathbb{E}(|\delta_{t,k}\delta_{t,l} \dot{Z}_{t,j}|)
\lesssim\kappa(\omega) B_{n}^{3}n^{-3/2}D_{n}^{2}\,.$$
Notice that $Z\overset{d}{=}n^{-1/2} \sum_{t=1}^{n} (\sqrt{\omega v} X_{t} + \sqrt{1-\omega v} Y_{t})$ and
$
\mathbb{E}\{h(Z,2)\} = \mathbb{P}(Z\leq u + \phi^{-1} + 2 \beta^{-1}) - \mathbb{P}(Z\leq u- \phi^{-1} - 2 \beta^{-1})$. Then 
$
\mathbb{P}(Z\leq u + \phi^{-1} + 2 \beta^{-1}) \leq \mathbb{P}(S_{n,y} \leq u + \phi^{-1} + 2 \beta^{-1}) + \varrho_{n}$ and $\mathbb{P}(Z\leq u- \phi^{-1} - 2 \beta^{-1}) \geq \mathbb{P}(S_{n,y} \leq u- \phi^{-1} - 2 \beta^{-1}) - \varrho_{n}$.
By Nazarov's inequality, we have
$
\mathbb{P}(S_{n,y} \leq u + \phi^{-1} + 2 \beta^{-1}) - \mathbb{P}(S_{n,y} \leq u- \phi^{-1} - 2 \beta^{-1}) \lesssim\phi^{-1}(\log{p})^{1/2}$. Hence, we have $
\mathbb{E}\{h(Z,2)\} \lesssim \phi^{-1}(\log{p})^{1/2} + \varrho_{n}$, which implies 
\begin{align*}
&\sum_{j,k,l=1}^{p} \sum_{t=1}^{n}\int_0^1 \mathbb{E}[ \chi_{t} h\{Z^{(-\mathcal{N}_t)},1\}U_{jkl}\{Z^{(-\mathcal{N}_t)}\}]\mathbb{E}(|\delta_{t,k}\delta_{t,l}\dot{Z}_{t,j}|)\,{\rm d}\omega \\
&~~~~~\lesssim B_{n}^{3}D_n^2n^{-1/2}\phi\beta^2\{\varrho_n+\phi^{-1}(\log p)^{1/2}\}\,.
\end{align*}
Then 
$
|{\rm III}_{1}| \lesssim D_{n}^{2} n^{-1/2}\phi \beta^{2}\{ M_{n}(\phi) + B_{n}^{3}\phi^{-1}(\log{p})^{1/2} + B_{n}^{3}\varrho_{n}\}$.
Together with the upper bound for $|{\rm III}_2|$, we complete the proof of Lemma \ref{lem:mIII}. $\hfill\Box$

\subsection{Proof of Lemma \ref{lem:mII}} \label{se:pflemmll}

Write $Z_s=(Z_{s,1},\ldots,Z_{s,p})^\T$. Due to $\mathbb{E}(X_{s,k}X_{t,j}) = \mathbb{E}(Y_{s,k}Y_{t,j}) = \mathbb{E}(W_{s,k}W_{t,j})$, then
$
\mathbb{E}(\delta_{t,k} \dot{Z}_{t,j}) = \sum_{s \in \mathcal{N}_{t}} \mathbb{E}(Z_{s,k} \dot{Z}_{t,j}) = n^{-1} \sum_{s \in \mathcal{N}_{t}} \mathbb{E}\{ \nu X_{s,k} X_{t,j} + (1-\nu) Y_{s,k} Y_{t,j} - W_{s,k}W_{t,
	j} \}= 0$. Write $\Delta_t=\delta^{(\mathcal{N}_t)}-\delta_t=(\Delta_{t,1},\ldots,\Delta_{t,p})^\T$. Then $Z^{(-t)}-Z^{(-\mathcal{N}_t)}=\Delta_t$. Since $Z^{(-\mathcal{N}_t)}$ is independent of $\delta_{t,k}\dot{Z}_{t,j}$, we know $\mathbb{E}[\partial_{jk}q\{Z^{(-\mathcal{N}_t)}\}\delta_{t,k}\dot{Z}_{t,j}]=0$. By Taylor expansion,
\begin{align*}
{\rm II}=&\,\sum_{j,k,l=1}^p\sum_{t=1}^n\int_0^1\int_0^1\mathbb{E}[\partial_{jkl}q\{Z^{(-\mathcal{N}_t)}+\tau\Delta_t\}\delta_{t,k}\Delta_{t,l}\dot{Z}_{t,j}]\,{\rm d}\tau{\rm d}\omega\\
=&\,\sum_{j,k,l=1}^p\sum_{t=1}^n\sum_{s\in(\cup_{\ell\in\mathcal{N}_t}\mathcal{N}_\ell)\backslash\mathcal{N}_t}\int_0^1\int_0^1\mathbb{E}[\partial_{jkl}q\{Z^{(-\mathcal{N}_t)}+\tau\Delta_t\}\delta_{t,k}Z_{s,l}\dot{Z}_{t,j}]\,{\rm d}\tau{\rm d}\omega\,.
\end{align*}
For given $(t,s)$, let
\begin{align*}
\chi_{t,s}=I\bigg\{\max_{j\in[p],s'\in\cup_{\ell\in\{s\}\cup\mathcal{N}_t}\mathcal{N}_\ell}(|X_{s',j}|\vee|Y_{s',j}|\vee|W_{s',j}|)\leq\frac{\sqrt{n}}{8\beta(D_n^4D_n^*)^{1/3}}\bigg\}\,.
\end{align*}
Then we have
\begin{align*}
&\sum_{j,k,l=1}^p\int_0^1\int_0^1\mathbb{E}[\partial_{jkl}q\{Z^{(-\mathcal{N}_t)}+\tau\Delta_t\}\delta_{t,k}Z_{s,l}\dot{Z}_{t,j}]\,{\rm d}\tau{\rm d}\omega\\
&~~~~=\underbrace{\sum_{j,k,l=1}^p\int_0^1\int_0^1\mathbb{E}[\chi_{t,s}\partial_{jkl}q\{Z^{(-\mathcal{N}_t)}+\tau\Delta_t\}\delta_{t,k}Z_{s,l}\dot{Z}_{t,j}]\,{\rm d}\tau{\rm d}\omega}_{{\rm II}_1(t,s)}\\
&~~~~~+\underbrace{\sum_{j,k,l=1}^p\int_0^1\int_0^1\mathbb{E}[(1-\chi_{t,s})\partial_{jkl}q\{Z^{(-\mathcal{N}_t)}+\tau\Delta_t\}\delta_{t,k}Z_{s,l}\dot{Z}_{t,j}]\,{\rm d}\tau{\rm d}\omega}_{{\rm II}_2(t,s)}\,.
\end{align*}
Identical to \eqref{eq:III2alpha} in Section \ref{sec:ETnalpha},
$$
|{\rm II}_2(t,s)|\lesssim \phi\beta^2\int_0^1\mathbb{E}\{(1-\chi_{t,s})\max_{j,k,l\in[p]}|\delta_{t,k}Z_{s,l}\dot{Z}_{t,j}|\}\,{\rm d}\omega\,.$$
Notice that
\begin{align*}
|\delta_{t,k}Z_{s,l}\dot{Z}_{t,j}| \leq&\,\frac{\kappa(\omega)}{n^{3/2}}\sum_{s' \in \mathcal{N}_{t}} (|X_{s',k}| + |Y_{s',k}| + |W_{s',k}|)(|X_{s,l}| + |Y_{s,l}| + |W_{s,l}|)\\
&~~~~~~~~~\times(|X_{t,j}| + |Y_{t,j}| + |W_{t,j}|)
\end{align*}
with $\kappa(\omega)=(\sqrt{\omega}\wedge\sqrt{1-\omega})^{-1}$ so that
\begin{align*}
\max_{j,k,l\in[p]} |\delta_{t,k}Z_{s,l} \dot{Z}_{t,j}| \lesssim&\,{D_{n}\kappa(\omega) \over n^{3/2}} \bigg\{ \max_{j\in[p],s'\in\cup_{\ell\in\{s\}\cup\mathcal{N}_t}\mathcal{N}_\ell} |X_{s',j}|^{3} \\
&\quad\quad+ \max_{j\in[p],s'\in\cup_{\ell\in\{s\}\cup\mathcal{N}_t}\mathcal{N}_\ell} |Y_{s',j}|^{3} + \max_{j\in[p],s'\in\cup_{\ell\in\{s\}\cup\mathcal{N}_t}\mathcal{N}_\ell} |W_{s',j}|^{3} \bigg\}\,.
\end{align*}
Identical to \eqref{eq:bdident} in Section \ref{sec:ETnalpha},
$$\mathbb{E}\bigg\{(1-\chi_{t,s}) \max_{j,k,l \in [p]} |\delta_{t,k}Z_{s,l} \dot{Z}_{t,j}|\bigg\}\lesssim D_{n} \kappa(\omega)n^{-3/2} \tilde{M}_n(\phi)\,,$$
which implies
$
\sum_{t=1}^n\sum_{s\in(\cup_{\ell\in\mathcal{N}_t}\mathcal{N}_\ell)\backslash\mathcal{N}_t}|{\rm II}_2(t,s)|\lesssim D_nD_{n}^*n^{-1/2}\tilde{M}_n(\phi)\phi\beta^2$.

Now we deal with ${\rm II}_1(t,s)$. Define $Z^{(-\mathcal{N}_t\cup\{s\})}=Z-\sum_{s'\in\cup_{\ell\in\{s\}\cup\mathcal{N}_t}\mathcal{N}_\ell}Z_{s'}$. It then holds that
$
Z^{(-\mathcal{N}_t)}+\tau\Delta_t
=Z^{(-\mathcal{N}_t\cup\{s\})}+\sum_{s'\in\mathcal{N}_s\cap(\cup_{\ell\in\mathcal{N}_t}\mathcal{N}_\ell)^c}Z_{s'}+\tau\sum_{s'\in(\cup_{\ell\in\mathcal{N}_t}\mathcal{N}_\ell)\backslash\mathcal{N}_t}Z_{s'}$.
If $\chi_{t,s}=1$, we have
$|\sum_{s'\in\mathcal{N}_s\cap(\cup_{\ell\in\mathcal{N}_t}\mathcal{N}_\ell)^c}Z_{s'}+\tau\sum_{s'\in(\cup_{\ell\in\mathcal{N}_t}\mathcal{N}_\ell)\backslash\mathcal{N}_t}Z_{s'}|_\infty
\leq\sqrt{3}/(4\beta)<\beta^{-1}$.
Same as \eqref{eq:bdIII1m}, 
$$
|{\rm II}_1(t,s)|\lesssim
\sum_{j,k,l=1}^p\int_0^1\mathbb{E}[h\{Z^{(-\mathcal{N}_t\cup\{s\})},1\}U_{jkl}\{Z^{(-\mathcal{N}_t\cup\{s\})}\}]\mathbb{E}(|\delta_{t,k}Z_{s,l}\dot{Z}_{t,j}|)\,{\rm d}\omega\,.$$
The left steps to derive the upper bound of $|{\rm II}_1(t,s)|$ are fully identical to that for bounding $|{\rm III}_1|$ stated in Section \ref{sec:pfmIII}. We then have
$$
\sum_{t=1}^n\sum_{s\in(\cup_{\ell\in\mathcal{N}_t}\mathcal{N}_\ell)\backslash\mathcal{N}_t}|{\rm II}_1(t,s)|\lesssim {D_{n}D_{n}^* \phi \beta^{2} n^{-1/2}}\{ \tilde{M}_{n}(\phi)+ B_{n}^{3}\phi^{-1}(\log{p})^{1/2} +B_{n}^{3} \varrho_{n}\}\,.$$
Since 
$|{\rm II}|\leq\sum_{t=1}^n\sum_{s\in(\cup_{\ell\in\mathcal{N}_t}\mathcal{N}_\ell)\backslash\mathcal{N}_t}|{\rm II}_1(t,s)|+\sum_{t=1}^n\sum_{s\in(\cup_{\ell\in\mathcal{N}_t}\mathcal{N}_\ell)\backslash\mathcal{N}_t}|{\rm II}_2(t,s)|$, we complete the proof of Lemma \ref{lem:mII}. $\hfill\Box$

\section{Proofs of the lemmas used in Section \ref{sec:pflma5}}

\subsection{Proof of Lemma \ref{lem:moment_bound_cov_mat_compare}}\label{se:pflemmoment_bound_cov_mat_compare}
Write $\calF_{t} = \sigma(\varepsilon_{t},\varepsilon_{t-1},\dots)$ and let $\calP_{t}(\cdot) = \mathbb{E}(\cdot\,|\,\calF_{t}) - \mathbb{E}(\cdot\,|\,\calF_{t-1})$ be a projection operator. In this proof, we shall slightly abuse the notation by using $\calF_{t}$ to also denote the corresponding input innovation sequence $\{\varepsilon_{t},\varepsilon_{t-1},\dots\}$. Since $X_{t} = \sum_{l=0}^\infty \calP_{t-l}(X_{t})$, then $S_{n,x} =n^{-1/2} \sum_{l=0}^\infty \sum_{t=1}^{n} \calP_{t-l}(X_{t})$. Let $\calF_{t,l}' = \sigma(\varepsilon_{t},\dots,\varepsilon_{t-l+1}, \varepsilon'_{t-l}, \varepsilon_{t-l-1},\dots)$. By the independence of the innovations, we have
\begin{align*}
\calP_{t-l}(X_{t}) =&\,\,\mathbb{E}\{f_{t}(\calF_{t})\,|\,\calF_{t-l}\} - \mathbb{E}\{f_{t}(\calF_{t})\,|\,\calF_{t-l-1}\}\\
=&\,\,\mathbb{E}\{f_{t}(\calF_{t})\,|\,\calF_{t-l}\} - \mathbb{E}\{f_{t}(\calF_{t,l}')\,|\,\calF_{t-l-1}\} =\mathbb{E}[X_t - X'_{t,\{l\}}\,|\,\calF_{t-l}]\,.
\end{align*}
Due to $q\geq2$, by Jensen's inequality, we have $\|\calP_{t-l}(X_{t,j})\|_{q} \leq \|X_{t,j} - X'_{t,j,\{l\}}\|_{q} \leq \theta_{l,q,j}$ for any $j\in[p]$. In addition, observe that $\{\calP_{t-l}(X_{t})\}_{t=1}^{n}$ is a martingale difference sequence with respect to the filtration $\{\calF_{t-l}\}_{t=1}^{n}$. By triangle inequality, Burkholder's inequality \citep[Theorem 2.1]{Rio_2009} and Cauchy-Schwarz inequality, we have that
\begin{align*}
\|S_{n,x,j}\|_{q} \leq&\,\,\frac{1}{n^{1/2}}\sum_{l=0}^\infty \bigg\|\sum_{t=1}^{n} \calP_{t-l}(X_{t,j})\bigg\|_{q} \leq \frac{1}{n^{1/2}}\sum_{l=0}^\infty (q-1)^{1/2} \bigg\{\sum_{t=1}^{n} \|\calP_{t-l}(X_{t,j})\|_{q}^{2}\bigg\}^{1/2} \\
\leq&\,\,\frac{(q-1)^{1/2}}{n^{1/2}} \sum_{l=0}^\infty \bigg(\sum_{t=1}^{n} \theta_{l,q,j}^{2}\bigg)^{1/2} = (q-1)^{1/2}\sum_{l=0}^\infty \theta_{l,q,j} = (q-1)^{1/2}\Theta_{0,q,j}\,.
\end{align*}

Next, we deal with $\|S_{n,x,j}^{(m)}\|_{q}$. Recall $X_{t}^{(m)} = \mathbb{E}(X_{t}\,|\,\varepsilon_{t}, \dots, \varepsilon_{t-m})$. Note that if $l < 0$ or $l > m$, then $\calP_{t-l}\{X_{t}^{(m)}\} = 0$. For $0 \leq l \leq m$, we have
\begin{align*}
\calP_{t-l}\{X_{t}^{(m)}\} =&\,\,\mathbb{E}\{X_{t}^{(m)}\,|\,\calF_{t-l}\} - \mathbb{E}\{\mathbb{E}(X_{t}\,|\,\varepsilon_{t},\dots,\varepsilon_{t-m})\,|\,\calF_{t-l-1}\} =\mathbb{E}[X_{t}^{(m)} - X'^{,(m)}_{t,\{l\}}\,|\,\calF_{t-l} ]\,,
\end{align*}
where $X'^{,(m)}_{t,\{l\}} = \mathbb{E}[X'_{t,\{l\}}\,|\,\varepsilon_{t},\dots,\varepsilon_{t-l+1}, \varepsilon'_{t-l}, \varepsilon_{t-l-1},\dots,\varepsilon_{t-m}]$. Write $$X_t^{(m)}=\{X_{t,1}^{(m)},\ldots,X_{t,p}^{(m)}\}^\T~\textrm{and}~X'^{,(m)}_{t,\{l\}}=[X'^{,(m)}_{t,1,\{l\}},\ldots,X'^{,(m)}_{t,p,\{l\}}]^\T\,.$$ By Jensen's inequality twice, we have
\begin{align*}
\|\calP_{t-l}\{X_{t,j}^{(m)}\}\|_{q} \leq\|X_{t,j}^{(m)} - X'^{,(m)}_{t,j,\{l\}}\|_{q}=&\,\,\| \mathbb{E}[X_{t,j} - X'_{t,j,\{l\}}\,|\,\varepsilon_{t},\dots, \varepsilon_{t-m}, \varepsilon'_{t-l} ] \|_{q}\\
\leq&\,\,\|X_{t,j} - X'_{t,j,\{l\}}\|_{q} \leq \theta_{l,q,j}\,.
\end{align*}
Notice that $S_{n,x}^{(m)} =n^{-1/2} \sum_{l=0}^\infty \sum_{t=1}^{n} \calP_{t-l}\{X_{t}^{(m)}\}$. By the same arguments for bounding $\|S_{n,x,j}\|_{q}$, we get $\|S_{n,x,j}^{(m)}\|_{q} \leq (q-1)^{1/2}\Theta_{0,q,j}$. Due to $S_{n,x,j}-S_{n,x,j}^{(m)} = \sum_{l=m+1}^\infty\{S_{n,x,j}^{(l)}-S_{n,x,j}^{(l-1)}\}$, then
\begin{align*}
&n^{1/2}\|S_{n,x,j}-S_{n,x,j}^{(m)}\|_q\leq n^{1/2}\sum_{l=m+1}^\infty\|S_{n,x,j}^{(l)}-S_{n,x,j}^{(l-1)}\|_q\\
&~~~~~~=\sum_{l=m+1}^\infty\bigg\|\sum_{t=1}^n\{X_{t,j}^{(l)}-X_{t,j}^{(l-1)}\}\bigg\|_q=\sum_{l=m+1}^\infty\bigg\|\sum_{t=0}^{n-1}\{X_{n-t,j}^{(l)}-X_{n-t,j}^{(l-1)}\}\bigg\|_q\,.
\end{align*}
Let $\mathcal{F}_{t,l}^*=\sigma(\varepsilon_{n-t-l},\varepsilon_{n-t-l+1},\ldots)$. Then $\{X_{n-t,j}^{(l)}-X_{n-t,j}^{(l-1)}\}_{t=0}^{n-1}$ is a martingale difference sequence with respect to $\{\mathcal{F}_{t,l}^*\}_{t=0}^{n-1}$. Since
\begin{align*}
X_{n-t}^{(l)}-X_{n-t}^{(l-1)}=&\,\,\mathbb{E}\{f_{n-t}(\mathcal{F}_{n-t})\,|\,\varepsilon_{n-t},\ldots,\varepsilon_{n-t-l}\}-\mathbb{E}\{f_{n-t}(\mathcal{F}_{n-t})\,|\,\varepsilon_{n-t},\ldots,\varepsilon_{n-t-l+1}\}\\
=&\,\,\mathbb{E}\{f_{n-t}(\mathcal{F}_{n-t})\,|\,\varepsilon_{n-t},\ldots,\varepsilon_{n-t-l}\}-\mathbb{E}\{f_{n-t}(\mathcal{F}_{n-t,l}')\,|\,\varepsilon_{n-t},\ldots,\varepsilon_{n-t-l}\}\\
=&\,\,\mathbb{E}[X_{n-t}-X_{n-t,\{l\}}'\,|\,\varepsilon_{n-t},\ldots,\varepsilon_{n-t-l}]\,,
\end{align*}
it follows from Jensen's inequality that $\|X_{n-t,j}^{(l)}-X_{n-t,j}^{(l-1)}\|_q\leq\theta_{l,q,j}$.
By the same arguments for bounding $\|S_{n,x,j}\|_q$, we get $\|S_{n,x,j}-S_{n,x,j}^{(m)}\|_q\leq (q-1)^{1/2}\Theta_{m+1,q,j}$. $\hfill\Box$

\subsection{Proof of Lemma \ref{lem:tail_prob_bound_m}}\label{sec:pfltpbm}
Denote $X_{t}^{(m)}=\mathbb{E}(X_t\,|\,\varepsilon_{t},\ldots,\varepsilon_{t-m})=g_t(\varepsilon_t,\varepsilon_{t-1},\ldots)$. Let $\varepsilon'_i$ be an independent copy of $\varepsilon_i$ and $X_{t,\{i\}}^{(m)'}=g_t(\varepsilon_t,\ldots,\varepsilon_{t-i+1},\varepsilon_{t-i}',\varepsilon_{t-i-1},\ldots)$ be the coupled version of $X_t^{(m)}$ at the lag $i$ with $\varepsilon_{t-i}$ replacing by $\varepsilon_{t-i}'$. Recall $
X'_{t,\{i\}} =f_{t}(\varepsilon_{t},\dots,\varepsilon_{t-i+1},\varepsilon'_{t-i},\varepsilon_{t-i-1},\dots)$
be the coupled version of $X_{t}$ at the time lag $i$ with $\varepsilon_{t-i}$ replacing by $\varepsilon'_{t-i}$. Then $X_{t,\{i\}}^{(m)'}=X_{t}^{(m)}$ for $i>m$ and $X_{t,\{i\}}^{(m)'}=\mathbb{E}(X_t\,|\,\varepsilon_{t},\ldots,\varepsilon_{t-i}',\ldots,\varepsilon_{t-m})=\mathbb{E}[X_{t,\{i\}}'\,|\,\varepsilon_{t},\ldots,\varepsilon_{t-m}]$ for $0\leq i\leq m$. For $q\geq1$, by Jensen's inequality, $$\|X_{t,j}^{(m)}-X_{t,j,\{i\}}^{(m)'}\|_{q}\leq \|\mathbb{E}[X_{t,j}-X_{t,j,\{i\}}'\,|\,\varepsilon_{t},\ldots,\varepsilon_{t-m}]\|_{q}\leq \|X_{t,j}-X_{t,j,\{i\}}'\|_{q}\leq\theta_{i,q,j}\,.$$ Define $\Phi_{\psi_{\nu},\alpha}^{(m)}$ in the same manner as  $\Phi_{\psi_{\nu},\alpha}$ in \eqref{eq:aggnorm} by replacing $\{X_{t}\}$ with $\{X_t^{(m)}\}$. Then $\Phi_{\psi_{\nu},0}^{(m)}\leq \Phi_{\psi_{\nu},0}<\infty$. By Lemma \ref{lem:tail_prob_bound_m-dep_approx_2}, there exists a universal constant $C>0$ depending only on $\nu$ such that $
\max_{j\in[p]}\mathbb{P}\{|\tilde X_{\ell,j}^{(m)}| > u\} \leq C \exp\{-(4e)^{-1}(1+2\nu)(u\Phi_{\psi_{\nu,0}}^{-1})^{2/(1+2\nu)}\}$
for any $u>0$.$\hfill\Box$

\subsection{Proof of Lemma \ref{lem:mom_m}}

Note that $\mathbb{E}\{|\tilde{X}_{\ell,j}^{(m)}|^q\}=q\int_0^\infty u^{q-1}\mathbb{P}\{|\tilde{X}_{\ell,j}^{(m)}|>u\}\,{\rm d}u$. By Lemma \ref{lem:tail_prob_bound_m}, we have Lemma \ref{lem:mom_m}. $\hfill\Box$

\subsection{Proof of Lemma \ref{lem:cov_infinity_m}}\label{sec:pflcim}
By Cauchy-Schwarz inequality and Lemma \ref{lem:moment_bound_cov_mat_compare}, we know
\begin{align*}
&|\mathbb{E}(S_{n,x,j}S_{n,x,k})- \mathbb{E}\{S_{n,x,j}^{(m)}S_{n,x,k}^{(m)}\}|\\
&~~~~~\leq|\mathbb{E}[S_{n,x,j}\{S_{n,x,k}-S_{n,x,k}^{(m)}\}]| + |\mathbb{E}\{S_{n,x,k}^{(m)}\{S_{n,x,j}-S_{n,x,j}^{(m)}\}]|\\
&~~~~~\leq\|S_{n,x,j}\|_{2} \|S_{n,x,k}-S_{n,x,k}^{(m)}\|_{2} + \|S_{n,x,k}^{(m)}\|_{2} \|S_{n,x,j}-S_{n,x,j}^{(m)}\|_{2}\\
&~~~~~\leq\Theta_{0,2,j}\Theta_{m+1,2,k}+\Theta_{0,2,k}\Theta_{m+1,2,j}\,.
\end{align*}
Therefore,
$
|\Xi^{(m)} - \Xi|_\infty \leq2(\max_{j\in[p]} \Theta_{0,2,j})(\max_{j \in[p]} \Theta_{m+1,2,j})\leq 2m^{-\alpha} \Psi_{2,\alpha} \Psi_{2,0}$.   By triangle inequality, $|\tilde\Xi - \Xi|_\infty\leq|\tilde\Xi - \Xi^{(m)}|_\infty+2m^{-\alpha} \Psi_{2,\alpha} \Psi_{2,0}$.
We adopt the convention $\cI_{L+1}=\emptyset$ and set $\tilde{X}_{L+1}^{(m)}=0$. Define $\tilde{a}_\ell=b^{1/2}$ for any $\ell\in[L+1]$, $\check{a}_\ell=m^{1/2}$ for any $\ell\in[L]$ and $\check{a}_{L+1}=(n-LQ)^{1/2}$. Recall $\tilde\Xi=L^{-1}\sum_{\ell=	1}^{L}\mathbb{E}[\tilde X_\ell^{(m)} \{\tilde X^{(m)}_\ell \}^{\T}]$ and $\Xi^{(m)}=\cov\{n^{-1/2}\sum_{t=1}^{n}X_{t}^{(m)}\}$. By triangle inequality,
\begin{align*}
|\tilde\Xi-\Xi^{(m)}|_\infty =&\,\max_{j_1,j_2\in[p]}\bigg| \frac{1}{Lb}\sum_{\ell=1}^L\tilde{a}_\ell^2\mathbb{E}\{\tilde{X}_{\ell,j_1}^{(m)}\tilde{X}_{\ell,j_2}^{(m)}\} - \frac{1}{n}\mathbb{E}\bigg[\bigg\{\sum_{t=1}^n X_{t,j_1}^{(m)}\bigg\}\bigg\{\sum_{t=1}^n X_{t,j_2}^{(m)}\bigg\} \bigg] \bigg| \\
\leq&\,\max_{j_1,j_2\in[p]}\underbrace{\bigg| \frac{1}{Lb}\sum_{\ell=1}^L  \tilde{a}_{\ell}^2\mathbb{E}\{\tilde{X}_{\ell,j_1}^{(m)}\tilde{X}_{\ell,j_2}^{(m)}\}- \frac{1}{n}\sum_{\ell=1}^L \tilde{a}_{\ell}^2\mathbb{E}\{\tilde{X}_{\ell,j_1}^{(m)}\tilde{X}_{\ell,j_2}^{(m)}\} \bigg|}_{R_1(j_1,j_2)}\\
&\,+\max_{j_1,j_2\in[p]}\underbrace{\bigg| \frac{1}{n}\sum_{\ell=1}^L\tilde{a}_\ell^2\mathbb{E}\{\tilde{X}_{\ell,j_1}^{(m)}\tilde{X}_{\ell,j_2}^{(m)}\} - \frac{1}{n}\mathbb{E}\bigg[\bigg\{\sum_{t=1}^n X_{t,j_1}^{(m)}\bigg\}\bigg\{\sum_{t=1}^n X_{t,j_2}^{(m)}\bigg\} \bigg] \bigg|}_{R_2(j_1,j_2)}\,.
\end{align*}
As we will show in Sections \ref{sec:R1j1j2} and \ref{sec:R2j1j2} that $\max_{j_1,j_2\in[p]}R_1(j_1,j_2)\lesssim \Phi_{\psi_{\nu,0}}^2(mb^{-1}+bn^{-1})$ and $\max_{j_1,j_2\in[p]}R_2(j_1,j_2)\lesssim \Phi_{\psi_{\nu,0}}^2(mb^{-1}+bn^{-1})$, it holds that $$|\tilde\Xi-\Xi^{(m)}|_\infty\lesssim \Phi_{\psi_{v},0}^2(mb^{-1}+bn^{-1})\,.$$ Recall $|\tilde\Xi - \Xi|_\infty\leq|\tilde\Xi - \Xi^{(m)}|_\infty+2m^{-\alpha} \Psi_{2,\alpha} \Psi_{2,0}$.
We complete the proof of Lemma \ref{lem:cov_infinity_m}. $\hfill\Box$

\subsubsection{Convergence rate of $\max_{j_1,j_2\in[p]}R_1(j_1,j_2)$}\label{sec:R1j1j2}

By Lemma \ref{lem:mom_m}, $$\max_{\ell\in[L]}\max_{j\in [p]}\mathbb{E}\{|\tilde{X}_{\ell,j}^{(m)}|^{2}\}\lesssim \Phi_{\psi_{\nu,0}}^2\,.$$ For any $j_1,j_2\in[p]$, by triangle inequality and Cauchy-Schwarz inequality, we have
\begin{align}\label{eq:R1}
R_1(j_1,j_2)\leq&\, \dfrac{n-Lb}{nLb}\sum_{\ell=1}^{L}\tilde a_{\ell}^2|\mathbb{E}\{\tilde X_{\ell,j_{1}}^{(m)}\tilde X_{\ell,j_{2}}^{(m)}\}|\notag\\
\leq&\, \dfrac{n-Lb}{nLb}\sum_{\ell=1}^{L}\tilde a_{\ell}^2[\mathbb{E}\{|\tilde X_{\ell,j_{1}}^{(m)}|^2\}]^{1/2}[\mathbb{E}\{|\tilde X_{\ell,j_{2}}^{(m)}|^2\}]^{1/2}\\
\lesssim&\,\, \Phi_{\psi_{\nu,0}}^2(mb^{-1}+bn^{-1})\,,\notag
\end{align}
where the last step is due to the facts $n-Lb\leq Lm+Q$ and $L\asymp nb^{-1}$. $\hfill\Box$

\subsubsection{Convergence rate of $\max_{j_1,j_2\in[p]}R_2(j_1,j_2)$}\label{sec:R2j1j2}

For any given $j_1,j_2\in[p]$, it holds that
\begin{align*}
R_2(j_1,j_2)\leq&\,\underbrace{\bigg|\frac{1}{n}\sum_{\ell=1}^{L+1}\tilde{a}_\ell\check{a}_\ell\mathbb{E}\{\tilde{X}_{\ell,j_1}^{(m)}\check{X}_{\ell,j_2}^{(m)}\}\bigg|}_{R_{2,1}(j_1,j_2)}
+\underbrace{\bigg|\frac{1}{n}\sum_{\ell=1}^{L+1}\check{a}_\ell\tilde{a}_\ell\mathbb{E}\{\check{X}_{\ell,j_1}^{(m)}\tilde{X}_{\ell,j_2}^{(m)}\}\bigg|}_{R_{2,2}(j_1,j_2)}\\
&+ \underbrace{\bigg|\frac{1}{n}\sum_{\ell=1}^{L+1}\check{a}_\ell^2\mathbb{E}\{\check{X}_{\ell,j_1}^{(m)}\check{X}_{\ell,j_2}^{(m)}\}\bigg|}_{R_{2,3}(j_1,j_2)} + \underbrace{\bigg|\frac{1}{n}\sum_{\ell_1\neq\ell_2}\tilde{a}_{\ell_1}\tilde{a}_{\ell_2}\mathbb{E}\{\tilde{X}_{\ell_1,j_1}^{(m)}\tilde{X}_{\ell_2,j_2}^{(m)}\}\bigg|}_{R_{2,4}(j_1,j_2)}\\
&+ \underbrace{\bigg|\frac{1}{n}\sum_{\ell_1\neq \ell_2}\tilde{a}_{\ell_1}\check{a}_{\ell_2}\mathbb{E}\{\tilde{X}_{\ell_1,j_1}^{(m)}\check{X}_{\ell_2,j_2}^{(m)}\}\bigg|}_{R_{2,5}(j_1,j_2)}+\underbrace{\bigg|\frac{1}{n}\sum_{\ell_1\neq \ell_2}\check{a}_{\ell_1}\tilde{a}_{\ell_2}\mathbb{E}\{\check{X}_{\ell_1,j_1}^{(m)}\tilde{X}_{\ell_2,j_2}^{(m)}\}\bigg|}_{R_{2,6}(j_1,j_2)}\\
&+ \underbrace{\bigg|\frac{1}{n}\sum_{\ell_1\neq \ell_2}\check{a}_{\ell_1}\check{a}_{\ell_2}\mathbb{E}\{\check{X}_{\ell_1,j_1}^{(m)}\check{X}_{\ell_2,j_2}^{(m)}\}\bigg|}_{R_{2,7}(j_1,j_2)} \,.   \nonumber
\end{align*}
Since $\tilde X_{\ell_1}^{(m)}$ is independent of $\tilde X_{\ell_2}^{(m)}$ for any $\ell_1\neq \ell_2$, then $R_{2,4}(j_1,j_2)=0$.

\underline{\it Convergence rates of $\max_{j_1,j_2\in[p]}R_{2,1}(j_1,j_2)$ and $\max_{j_1,j_2\in[p]}R_{2,2}(j_1,j_2)$.}  Notice that $X_{t}^{(m)}=\sum_{l=0}^\infty \calP_{t-l}\{X_{t}^{(m)}\}$, where $\calP_{t}(\cdot) = \mathbb{E}(\cdot\,|\,\calF_{t}) - \mathbb{E}(\cdot\,|\,\calF_{t-1})$ with $\calF_{t} = \sigma(\varepsilon_{t},\varepsilon_{t-1},\dots)$. For $X_{t,\{l\}}^{(m)'}$ defined in Section \ref{sec:pfltpbm}, we have $$
\calP_{t-l}\{X_{t}^{(m)}\} =\mathbb{E}\{X_{t}^{(m)}\,|\,\calF_{t-l}\} - \mathbb{E}\{X_{t,\{l\}}^{(m)'}\,|\,\calF_{t-l-1}\} =\mathbb{E}[X_{t}^{(m)} - X_{t,\{l\}}^{(m)'}\,|\,\calF_{t-l}]\,.$$
For $q\geq2$, by Jensen's inequality, $\|\calP_{t-l}\{X_{t,j}^{(m)}\}\|_{q} \leq \|X_{t,j}^{(m)} - X_{t,j,\{l\}}^{(m)'}\|_{q}\leq \theta_{l,q,j}$ for any $j\in[p]$, where the last step is due to the result shown in Section \ref{sec:pfltpbm}. By triangle inequality and Cauchy-Schwarz inequality, for any $t\in[n]$ and $k\geq 1$, we have
\begin{align}\label{eq:phycovbd}
|\mathbb{E}\{X_{t,j_{1}}^{(m)}X_{t-k,j_{2}}^{(m)}\}|=&\,\bigg|\sum_{l=0}^{\infty}\mathbb{E}[\calP_{t-l}\{X_{t,j_1}^{(m)}\}\calP_{t-l}\{X_{t-k,j_2}^{(m)}\}]\bigg|\notag\\
\leq&\, \sum_{l=0}^{\infty}|\mathbb{E}[\calP_{t-l}\{X_{t,j_1}^{(m)}\}\calP_{t-l}\{X_{t-k,j_2}^{(m)}\}]|\\
\leq&\,\sum_{l=0}^{\infty}\|\calP_{t-l}\{X_{t,j_1}^{(m)}\}\|_2\|\calP_{t-l}\{X_{t-k,j_2}^{(m)}\}\|_2\leq\sum_{l=0}^{\infty}\theta_{l,2,j_1}\theta_{l-k,2,j_2}\,,\notag
\end{align}
which implies
\begin{align*}
R_{2,1}(j_1,j_2)
\leq &\,\,\frac{1}{n}\sum_{\ell=1}^{L+1}\sum_{t_1\in\cI_\ell}\sum_{t_2\in\cJ_\ell} |\mathbb{E}\{X_{t_1,j_1}^{(m)}X_{t_2,j_2}^{(m)}\}|\\
\leq&\,\,\frac{1}{n}\sum_{\ell=1}^{L}\sum_{t\in\cJ_\ell}\sum_{k=1}^{b+m} |\mathbb{E}\{X_{t,j_2}^{(m)}X_{t-k,j_1}^{(m)}\}|\leq \frac{mL}{n}\sum_{k=1}^{b+m}\sum_{l=0}^{\infty}\theta_{l,2,j_2}\theta_{l-k,2,j_1}\,.
\end{align*}
Due to $\theta_{l,2,j}=0$ for any $l<0$, then 
$
\sum_{k=1}^{b+m}\sum_{l=0}^{\infty}\theta_{l,2,j_2}\theta_{l-k,2,j_1}\leq\Theta_{0,2,j_1}\Theta_{0,2,j_2}\lesssim\Phi_{\psi_{v},0}^2$.
Since $L\asymp nb^{-1}$, we have $\max_{j_1,j_2\in[p]}R_{2,1}(j_1,j_2)\lesssim \Phi_{\psi_{v},0}^2mb^{-1}$.
Analogously, $\max_{j_1,j_2\in[p]}R_{2,2}(j_1,j_2)\lesssim \Phi_{\psi_{v},0}^2mb^{-1}$.

\underline{\it Convergence rate of $\max_{j_1,j_2\in[p]}R_{2,3}(j_1,j_2)$.} Using the same arguments for deriving Lemma \ref{lem:mom_m}, 
$\max_{\ell\in[L+1]}\max_{j\in [p]}\mathbb{E}\{|\check{X}_{\ell,j}^{(m)}|^{2}\}\lesssim \Phi_{\psi_{\nu,0}}^2$. By triangle inequality and Cauchy-Schwarz inequality, $R_{2,3}(j_1,j_2)\leq n^{-1}\sum_{\ell=1}^{L+1}\check{a}_\ell^2|\mathbb{E}\{\check{X}_{\ell,j_1}^{(m)}\check{X}_{\ell,j_2}^{(m)}\}|\lesssim  \Phi_{\psi_{v},0}^2(mb^{-1}+bn^{-1})$. 

\underline{\it Convergence rates of $\max_{j_1,j_2\in[p]}R_{2,5}(j_1,j_2)$ and $\max_{j_1,j_2\in[p]}R_{2,6}(j_1,j_2)$.} Notice that $\tilde X_{\ell_1}^{(m)}$ is independent of $\check X_{\ell_2}^{(m)}$ for any $(\ell_1,\ell_2)\in\{(\ell_1,\ell_2):\ell_1\neq\ell_2,\ell_{2}+1\}$. By triangle inequality,
\begin{align*}
R_{2,5}(j_1,j_2)=&\,\bigg|\frac{1}{n}\sum_{\ell= 1}^{L}\tilde{a}_{\ell+1}\check{a}_{\ell}\mathbb{E}\{\tilde{X}_{\ell+1,j_1}^{(m)}\check{X}_{\ell,j_2}^{(m)}\}\bigg|\leq \frac{1}{n}\sum_{\ell=1}^{L}\sum_{t_1\in\mathcal{I}_{\ell+1},t_2\in\mathcal{J}_{\ell}}|\mathbb{E}\{X_{t_1,j_1}^{(m)}X_{t_2,j_2}^{(m)}\}|\\
\leq&\,\frac{1}{n}\sum_{\ell=1}^{L}\sum_{t\in\mathcal{J}_{\ell}}\sum_{k=1}^{b+m}|\mathbb{E}\{X_{t+k,j_1}^{(m)}X_{t,j_2}^{(m)}\}|\lesssim \Phi_{\psi_{v},0}^2mb^{-1}\,,
\end{align*}
where the last step is due to \eqref{eq:phycovbd} and the fact $
\sum_{k=1}^{b+m}\sum_{l=0}^{\infty}\theta_{l,2,j_2}\theta_{l-k,2,j_1}\lesssim\Phi_{\psi_{v},0}^2$. Then $$\max_{j_1,j_2\in[p]}R_{2,5}(j_1,j_2)\lesssim \Phi_{\psi_{v},0}^2mb^{-1}\,.$$ Analogously, $\max_{j_1,j_2\in[p]}R_{2,6}(j_1,j_2)\lesssim \Phi_{\psi_{v},0}^2mb^{-1}$.

\underline{\it Convergence rate of $\max_{j_1,j_2\in[p]}R_{2,7}(j_1,j_2)$.} Notice that $\check X_{\ell_1}^{(m)}$ is independent of $\check X_{\ell_2}^{(m)}$ for any $(\ell_1,\ell_2)\notin\{(L,L+1),(L+1,L)\}$. Hence, by triangle inequality,
\begin{align*}
R_{2,7}(j_1,j_2)\lesssim \frac{1}{n}\sum_{t_1\in\mathcal{J}_{L},t_2\in\mathcal{J}_{L+1}}|\mathbb{E}\{X_{t_1,j_1}^{(m)}X_{t_2,j_2}^{(m)}\}|+\frac{1}{n}\sum_{t_1\in\mathcal{J}_{L+1},t_2\in\mathcal{J}_{L}}|\mathbb{E}\{X_{t_1,j_1}^{(m)}X_{t_2,j_2}^{(m)}\}|\,.
\end{align*}
Note that $$\frac{1}{n}\sum_{t_1\in\mathcal{J}_{L},t_2\in\mathcal{J}_{L+1}}|\mathbb{E}\{X_{t_1,j_1}^{(m)}X_{t_2,j_2}^{(m)}\}|\leq \frac{1}{n}\sum_{t\in\mathcal{J}_{L+1}}\sum_{k=1}^{n-LQ+m}|\mathbb{E}\{X_{t,j_2}^{(m)}X_{t-k,j_1}^{(m)}\}|\lesssim \Phi_{\psi_{v},0}^2bn^{-1}\,,$$ where the last step is due to \eqref{eq:phycovbd} and the fact $
\sum_{k=1}^{n-LQ+m}\sum_{l=0}^{\infty}\theta_{l,2,j_2}\theta_{l-k,2,j_1}\lesssim\Phi_{\psi_{v},0}^2$. Analogously, we also have $n^{-1}\sum_{t_1\in\mathcal{J}_{L+1},t_2\in\mathcal{J}_{L}}|\mathbb{E}\{X_{t_1,j_1}^{(m)}X_{t_2,j_2}^{(m)}\}|\lesssim \Phi_{\psi_{v},0}^2bn^{-1}$, which implies $$\max_{j_1,j_2\in[p]}R_{2,7}(j_1,j_2)\lesssim \Phi_{\psi_{v},0}^2bn^{-1}\,.$$ Thus, $\max_{j_1,j_2\in[p]}R_{2}(j_1,j_2)\lesssim \Phi_{\psi_{v},0}^2(mb^{-1}+bn^{-1})$. $\hfill\Box$

\subsection{Proof of Lemma \ref{pn:1_m}}\label{sec:pfpn1m}
By Lemma \ref{lem:cov_infinity_m}, we have $$|\tilde\Xi - \Xi|_\infty\lesssim \Phi_{\psi_{v},0}^2(mb^{-1}+bn^{-1})+m^{-\alpha} \Psi_{2,\alpha} \Psi_{2,0} =o(1)$$ provided that $\Phi_{\psi_{v},0}^2m\ll b\ll \Phi_{\psi_{v},0}^{-2}n$ and $m^{\alpha}\gg \Psi_{2,\alpha}\Psi_{2,0}$. Since $\min_{j\in[p]}V_{n,j}\geq C$, it holds that $\min_{j\in [p]}L^{-1}\sum_{\ell=1}^{L}\mathbb{E}\{|\tilde X_{\ell,j}^{(m)}|^2\}\geq C/2$.
By Lemma \ref{lem:mom_m}, $$\max_{\ell\in[L]}\max_{j\in [p]}\mathbb{E}\{|\tilde{X}_{\ell,j}^{(m)}|^{3}\}\lesssim \Phi_{\psi_{\nu,0}}^3\,,$$ which implies $\max_{j\in [p]}L^{-1}\sum_{\ell=1}^{L}\mathbb{E}\{|\tilde X_{\ell,j}^{(m)}|^3\}\lesssim \Phi_{\psi_{\nu,0}}^3$.
Recall $\{\tilde{Y}_{\ell}^{(m)}\}_{\ell=1}^{L}$ is a sequence of independent centered Gaussian random vectors with $\mathbb{E}[\tilde{Y}_{\ell}^{(m)}\{\tilde{Y}_{\ell}^{(m)}\}^\T]=\mathbb{E}[\tilde{X}_{\ell}^{(m)}\{\tilde{X}_{\ell}^{(m)}\}^\T]$ for $\ell\in [L]$. By Theorem 2.1 of \cite{CCK_2017}, it holds that
\begin{align}
\varrho_{n,1}^{(m)}\lesssim&\, \frac{\Phi_{\psi_{\nu,0}}(\log p)^{7/6}}{L^{1/6}}+\frac{1}{\Phi_{\psi_{\nu,0}}^3}\max_{\ell\in [L]}\mathbb{E}\bigg[|\tilde X_{\ell}^{(m)}|_{\infty}^3I\bigg\{|\tilde X_{\ell}^{(m)}|_{\infty}>\frac{CL^{1/3}\Phi_{\psi_{\nu,0}}}{(\log p)^{1/3}}\bigg\}\bigg]\notag\\
&\,+\frac{1}{\Phi_{\psi_{\nu,0}}^3}\max_{\ell\in [L]}\mathbb{E}\bigg[|\tilde Y_{\ell}^{(m)}|_{\infty}^3I\bigg\{|\tilde Y_{\ell}^{(m)}|_{\infty}>\frac{CL^{1/3}\Phi_{\psi_{\nu,0}}}{(\log p)^{1/3}}\bigg\}\bigg]\,.\label{eq:rhon1m}
\end{align}
Lemma \ref{lem:tail_prob_bound_m} implies that $\max_{\ell\in[L]}\mathbb{P}\{|\tilde{X}_{\ell}^{(m)}|_{\infty}> u\}\lesssim p\exp\{-C(u\Phi_{\psi_{\nu},0}^{-1})^{2/(1+2\nu)}\}$ for any $u>0$. Since $\mathbb{E}\{|\xi|^3I(|\xi|>v)\}=v^3\mathbb{P}(|\xi|>v)+3\int_v^\infty u^2\mathbb{P}(|\xi|>u)\,{\rm d}u$, it holds that
\begin{align*}
&\frac{1}{\Phi_{\psi_{\nu,0}}^3}\mathbb{E}\bigg[|\tilde X_{\ell}^{(m)}|_{\infty}^3I\bigg\{|\tilde X_{\ell}^{(m)}|_{\infty}>\frac{CL^{1/3}\Phi_{\psi_{\nu,0}}}{(\log p)^{1/3}}\bigg\}\bigg]\\
&~~\lesssim\frac{L}{\log p}\mathbb{P}\bigg\{|\tilde X_{\ell}^{(m)}|_{\infty}>\frac{CL^{1/3}\Phi_{\psi_{\nu,0}}}{(\log p)^{1/3}}\bigg\}+\frac{1}{\Phi_{\psi_{\nu,0}}^3}\int_{\frac{CL^{1/3}\Phi_{\psi_{\nu,0}}}{(\log p)^{1/3}}}^\infty u^2\mathbb{P}\{|\tilde X_{\ell}^{(m)}|_{\infty}>u\}\,{\rm d}u\,,
\end{align*}
where
\begin{align*}
&\mathbb{P}\bigg\{|\tilde X_{\ell}^{(m)}|_{\infty}>\frac{CL^{1/3}\Phi_{\psi_{\nu,0}}}{(\log p)^{1/3}}\bigg\}\lesssim p\exp\bigg\{-
\frac{CL^{2/(3+6\nu)}}{(\log p)^{2/(3+6\nu)}}\bigg\}\lesssim \frac{(\log p)^{7/6}}{L^{7/6}}
\end{align*}
and
\begin{align*}
&\int_{\frac{CL^{1/3}\Phi_{\psi_{\nu,0}}}{(\log p)^{1/3}}}^\infty u^2\mathbb{P}\{|\tilde X_{\ell}^{(m)}|_{\infty}>u\}\,{\rm d}u\lesssim p\int_{\frac{CL^{1/3}\Phi_{\psi_{\nu,0}}}{(\log p)^{1/3}}}^\infty u^2\exp\{-C(u\Phi_{\psi_{\nu},0}^{-1})^{2/(1+2\nu)}\}\,{\rm d}u\\
&~~~~~~~\lesssim\dfrac{L\Phi_{\psi_{\nu,0}}^3p}{\log p}\exp\bigg\{-
\frac{CL^{2/(3+6\nu)}}{(\log p)^{2/(3+6\nu)}}\bigg\}\lesssim \frac{\Phi_{\psi_{\nu,0}}^3(\log p)^{1/6}}{L^{1/6}}
\end{align*}
provided that $(\log p)^{5+6\nu}=o(L^2)$. Thus, $$\Phi_{\psi_{\nu,0}}^{-3}\mathbb{E}[|\tilde X_{\ell}^{(m)}|_{\infty}^3I\{|\tilde X_{\ell}^{(m)}|_{\infty}>CL^{1/3}\Phi_{\psi_{\nu,0}}(\log p)^{-1/3}\}]\\\lesssim L^{-1/6}(\log p)^{1/6}\,.$$ Recall $\tilde{Y}_{\ell,j}^{(m)}$ is a normal random variable with zero mean. By Lemma \ref{lem:mom_m}, we have 
$$\max_{\ell\in [L]}\max_{j\in [p]}\mathbb{E}[\{\tilde Y_{\ell,j}^{(m)}\}^2]=\max_{\ell\in [L]}\max_{j\in [p]}\mathbb{E}[\{\tilde X_{\ell,j}^{(m)}\}^2]\lesssim\Phi_{\psi_{\nu,0}}^2 \,,$$ which implies that $\max_{\ell\in [L]}\max_{j\in [p]}\mathbb{P}\{|\tilde Y_{\ell,j}^{(m)}|>u\}\lesssim\exp(-Cu^2\Phi_{\psi_{\nu,0}}^{-2})$. By Bonferroni inequality, it holds that $\max_{\ell\in [L]}\mathbb{P}\{|\tilde Y_{\ell}^{(m)}|_{\infty}>u\}\lesssim p\exp(-Cu^2\Phi_{\psi_{\nu,0}}^{-2})$ for any $u>0$. Then
\begin{align*}
&\mathbb{P}\bigg\{|\tilde Y_{\ell}^{(m)}|_{\infty}>\frac{CL^{1/3}\Phi_{\psi_{\nu,0}}}{(\log p)^{1/3}}\bigg\}\lesssim p\exp\bigg\{-
\frac{CL^{2/3}}{(\log p)^{2/3}}\bigg\}\lesssim \frac{(\log p)^{7/6}}{L^{7/6}}
\end{align*}
and
\begin{align*}
\int_{\frac{CL^{1/3}\Phi_{\psi_{\nu,0}}}{(\log p)^{1/3}}}^\infty u^2\mathbb{P}\{|\tilde Y_{\ell}^{(m)}|_{\infty}>u\}\,{\rm d}u\lesssim&\,\,p\int_{\frac{CL^{1/3}\Phi_{\psi_{\nu,0}}}{(\log p)^{1/3}}}^\infty u^2\exp(-Cu^2\Phi_{\psi_{\nu},0}^{-2})\,{\rm d}u\\
\lesssim&\,\,\dfrac{L\Phi_{\psi_{\nu,0}}^3p}{\log p}\exp\bigg\{-\frac{CL^{2/3}}{(\log p)^{2/3}}\bigg\}\lesssim \frac{\Phi_{\psi_{\nu,0}}^3(\log p)^{1/6}}{L^{1/6}}
\end{align*}
provided that $(\log p)^5=o(L^2)$. Thus, $$\Phi_{\psi_{\nu,0}}^{-3}\mathbb{E}[|\tilde Y_{\ell}^{(m)}|_{\infty}^3I\{|\tilde Y_{\ell}^{(m)}|_{\infty}>CL^{1/3}\Phi_{\psi_{\nu,0}}(\log p)^{-1/3}\}]\lesssim L^{-1/6}(\log p)^{1/6}\,.$$ Therefore, by \eqref{eq:rhon1m}, we have $\varrho_{n,1}^{(m)}\lesssim L^{-1/6}\Phi_{\psi_{\nu,0}}(\log p)^{7/6}$.$\hfill\Box$

\subsection{Proof of Lemma \ref{pn:2_m}}\label{sec:pflpn2m}
Note that $
S_{n,x}^{(m)}=\tilde S_{n,x}^{(m)}+\delta_n$ with $$
\delta_n=\frac{1}{\sqrt{n}}\sum_{\ell=1}^{L+1}\sum_{t\in\mathcal{J}_\ell}{X}_t^{(m)}+\bigg(\frac{1}{\sqrt{n}}-\frac{1}{\sqrt{Lb}}\bigg)\sum_{\ell=1}^L\sum_{t\in\mathcal{I}_\ell}{X}_t^{(m)}\,.$$
For some ${D}_{n}>0$, consider the event $\mathcal{E}=\{|\delta_n|_\infty\leq {D}_{n}\}$. Identical to \eqref{eq:bb2}, we have
\begin{align}\label{eq:bb2_m}
\varrho_{n,2}^{(m)}\leq&\,\,\varrho_{n,1}^{(m)}+\sup_{u\in\mathbb{R}^p,\nu\in[0,1]}|\mathbb{P}\{\tilde S_{n,y}^{(m)}\leq u-\sqrt{\nu}D_n\}-\mathbb{P}\{\tilde S_{n,y}^{(m)}\leq u\}|\notag\\
&+\sup_{u\in\mathbb{R}^p,\nu\in[0,1]}|\mathbb{P}\{\tilde S_{n,y}^{(m)}\leq u+\sqrt{\nu}D_n\}-\mathbb{P}\{\tilde S_{n,y}^{(m)}\leq u\}|+\mathbb{P}(\mathcal{E}^c)\\
\lesssim&\,\,L^{-1/6}\Phi_{\psi_{\nu,0}}(\log p)^{7/6}+{D}_{n}(\log p)^{1/2}+\mathbb{P}(\mathcal{E}^c)\,,\notag
\end{align}
where the last step is based on Lemma \ref{pn:1_m} and Nazarov's inequality. In the sequel, we will bound $\mathbb{P}(\mathcal{E}^c)$. Due to $b\ll (mn)^{1/2}$ and $L\asymp nb^{-1}$, we have $n-Lb\lesssim Lm$ which implies
$
D_n\sqrt{nLb}(\sqrt{n}-\sqrt{Lb})^{-1}\gtrsim D_n\sqrt{n}bm^{-1}$. By Bonferroni inequality, it holds that
\begin{align*}
\mathbb{P}(\mathcal{E}^c)\leq&\,\,\mathbb{P}\bigg\{\bigg|\sum_{\ell=1}^{L+1}\sum_{t\in\mathcal{J}_\ell}{X}_t^{(m)}\bigg|_\infty>\frac{D_n\sqrt{n}}{2}\bigg\}+\mathbb{P}\bigg\{\bigg|\sum_{\ell=1}^L\sum_{t\in\mathcal{I}_\ell}{X}_t^{(m)}\bigg|_\infty>\frac{D_n\sqrt{nLb}}{2(\sqrt{n}-\sqrt{Lb})}\bigg\}\\
\leq&\,\,\sum_{j=1}^p\mathbb{P}\bigg\{\bigg|\sum_{\ell=1}^{L+1}\sum_{t\in\mathcal{J}_\ell}X_{t,j}^{(m)}\bigg|>\frac{D_n\sqrt{n}}{2}\bigg\}+\sum_{j=1}^p\mathbb{P}\bigg\{\bigg|\sum_{\ell=1}^L\sum_{t\in\mathcal{I}_\ell}X_{t,j}^{(m)}\bigg|>\frac{CD_n\sqrt{n}b}{m}\bigg\}\,.
\end{align*}
For each $j\in[p]$, $\{X_t^{(m)}\}_{t\in\mathcal{I}_{\ell},\ell\in[L]}$ and $\{X_t^{(m)}\}_{t\in\mathcal{J}_{\ell},\ell\in[L+1]}$ are also physical sequences. Similar to the proof of Lemma \ref{lem:tail_prob_bound_m}, due to $b\gg m$, we have
\begin{align*}
&\max_{j\in[p]}\mathbb{P}\bigg\{\bigg|\sum_{\ell=1}^L\sum_{t\in\mathcal{I}_\ell}X_{t,j}^{(m)}\bigg|>\frac{CD_n\sqrt{n}b}{m}\bigg\}+\max_{j\in[p]}\mathbb{P}\bigg\{\bigg|\sum_{\ell=1}^{L+1}\sum_{t\in\mathcal{J}_\ell}X_{t,j}^{(m)}\bigg|>\frac{D_n\sqrt{n}}{2}\bigg\}\\
&~~~~~\lesssim\exp\bigg\{-C\bigg(\frac{b^{1/2}D_n}{m^{1/2}\Phi_{\psi_\nu,0}}\bigg)^{2/(1+2\nu)}\bigg\}\,.
\end{align*}
With selecting $D_n=C'L^{-1/6}\Phi_{\psi_\nu,0}(\log p)^{2/3}$ for some sufficiently large constant $C'>0$, then $D_n(\log p)^{1/2}\lesssim L^{-1/6}\Phi_{\psi_\nu,0}(\log p)^{7/6}$ and
$$
\mathbb{P}(\mathcal{E}^c)\lesssim p\exp[-C\{b^{2/3}n^{-1/6}(\log p)^{2/3}m^{-1/2}\}^{2/(1+2\nu)}]\,.$$
Due to $b\gg n^{1/4}m^{3/4}(\log p)^{(6\nu-1)/4}
$, we have $$p\exp[-C\{b^{2/3}n^{-1/6}(\log p)^{2/3}m^{-1/2}\}^{2/(1+2\nu)}]\lesssim L^{-1/6}\Phi_{\psi_\nu,0}(\log p)^{7/6}\,,$$ which implies $\mathbb{P}(\mathcal{E}^c)\lesssim L^{-1/6}\Phi_{\psi_\nu,0}(\log p)^{7/6}$. From (\ref{eq:bb2_m}), we complete the proof. $\hfill\Box$

\section{Proof of Proposition \ref{pn:zhangwu2017}}\label{sec:xic1c4}
Without loss of generality, we assume $\log p\leq n^c$ for some sufficiently small $c>0$. Otherwise, the assertion holds trivially. Let $\lambda = n^{-c_1}, \eta = n^{-c_2}$, $w \asymp n^{c_3}$ and $m\asymp n^{c_4}$ for some constants $c_1,c_2,c_3, c_4>0$ that can be optimized according to the right-hand side of \eqref{eqn:zhang_wu_rate_phystical_dependence}. Recall $p\geq n^\kappa$ for some $\kappa>0$. Since $m=o(M)$ and $w\asymp n/(M+m)$, we have $M\asymp n^{1-c_3}$. As we will show below in Cases 1--3, with the optimal selection of $(c_1,c_2,c_3,c_4)$, it holds that $\omega_n\lesssim n^{-\alpha/(11\alpha+3)}\cdot \mbox{polylog}(p)$ for $\alpha>0$.

{\indent\bf Case 1.} If $\alpha > 1$, by \eqref{eqn:zhang_wu_rate_phystical_dependence}, we have
\begin{align*}
\omega_n \lesssim &\,\, p \exp\{-C n^{\beta(\alpha c_4-c_2)}\}+p\exp\{-C n^{\beta(1-2c_2-c_3-c_4)/2}\}\\
&+ \{n^{-c_1} + n^{-c_2}+ n^{-c_3/8}  + n^{-\alpha c_4/3}+n^{-(1-c_3-c_4)/3}\}\cdot \mbox{polylog}(p)\,.
\end{align*}
To make the upper bound in above inequality converges to zero, we need to require $\alpha^{-1}c_2<c_4<1-2c_2-c_3$ and $c_2<\alpha(1-c_3)/(1+2\alpha)$. Define $$\xi(c_1,c_2,c_3,c_4)=\inf\{c_1,c_2,c_3/8,\alpha c_4/3,(1-c_3-c_4)/3\}$$ and $\mathcal{C}=\{(c_1,c_2,c_3,c_4):c_1>0,0<c_3<1, 0<c_2<\alpha(1-c_3)/(1+2\alpha), \alpha^{-1}c_2<c_4<1-2c_2-c_3\}$. Write $\tilde{\mathcal{C}}=\{(c_2,c_3,c_4):0<c_3<1, 0<c_2<\alpha(1-c_3)/(1+2\alpha), \alpha^{-1}c_2<c_4<1-2c_2-c_3\}$ and define $$\tilde{\xi}(c_2,c_3,c_4)=\inf\{c_2,c_3/8,\alpha c_4/3,(1-c_3-c_4)/3\}\,.$$ Notice that $\sup_{(c_1,c_2,c_3,c_4)\in\mathcal{C}}\xi(c_1,c_2,c_3,c_4)=\sup_{(c_2,c_3,c_4)\in\tilde{\mathcal{C}}}\tilde{\xi}(c_2,c_3,c_4)$.

(i) If $0<c_2<\alpha(1-c_3)/(2+2\alpha)$, we have $1-2c_2-c_3> (1-c_3)/(1+\alpha)>\alpha^{-1}c_2$, which implies
\begin{align*}
&\sup_{0<c_3<1}\sup_{0<c_2<\frac{\alpha(1-c_3)}{2+2\alpha}}\sup_{\alpha^{-1}c_2<c_4<1-2c_2-c_3}\tilde{\xi}(c_2,c_3,c_4)\\
&~~~=\sup_{0<c_3<1}\sup_{0<c_2<\frac{\alpha(1-c_3)}{2+2\alpha}}\inf\bigg\{c_2,\frac{c_3}{8},\frac{\alpha(1-c_3)}{3\alpha+3}\bigg\}\\
&~~~=\sup_{0<c_3<1}\inf\bigg\{\frac{c_3}{8},\frac{\alpha(1-c_3)}{3\alpha+3}\bigg\}=\frac{\alpha}{11\alpha+3}\,.
\end{align*}

(ii) If $\alpha(1-c_3)/(2+2\alpha)\leq c_2<\alpha(1-c_3)/(1+2\alpha)$, we have $1-2c_2-c_3\leq (1-c_3)/(1+\alpha)$, which implies
\begin{align*}
&\sup_{0<c_3<1}\sup_{\frac{\alpha(1-c_3)}{2+2\alpha}\leq c_2<\frac{\alpha(1-c_3)}{1+2\alpha}}\sup_{\alpha^{-1}c_2<c_4<1-2c_2-c_3}\tilde{\xi}(c_2,c_3,c_4)\\
&~~~=\sup_{0<c_3<1}\sup_{\frac{\alpha(1-c_3}{2+2\alpha}\leq c_2<\frac{\alpha(1-c_3)}{1+2\alpha}}\inf\bigg\{c_2,\frac{c_3}{8},\frac{\alpha(1-2c_2-c_3)}{3}\bigg\}\\
&~~~=\sup_{0<c_3<1}\sup_{\frac{\alpha(1-c_3)}{2+2\alpha}\leq c_2<\frac{\alpha(1-c_3)}{1+2\alpha}}\inf\bigg\{\frac{c_3}{8},\frac{\alpha(1-2c_2-c_3)}{3}\bigg\}\\
&~~~=\sup_{0<c_3<1}\inf\bigg\{\frac{c_3}{8},\frac{\alpha(1-c_3)}{3\alpha+3}\bigg\}=\frac{\alpha}{11\alpha+3}\,.
\end{align*}
Hence, we have $$\sup_{(c_1,c_2,c_3,c_4)\in\mathcal{C}}\xi(c_1,c_2,c_3,c_4)=\frac{\alpha}{11\alpha+3}$$ and it can be achieved by selecting $c_3=8\alpha/(11\alpha+3)$, $c_1=c_2=\alpha(1-c_3)/(3\alpha+3)$ and $c_4=(1-c_3)/(\alpha+1)$. With such selected $(c_1,c_2,c_3,c_4)$, we know $\beta(\alpha c_4-c_2)>0$ and $\beta(1-2c_2-c_3-c_4)/2>0$, which implies that $p \exp\{-C n^{\beta(\alpha c_4-c_2)}\}+p\exp\{-C n^{\beta(1-2c_2-c_3-c_4)/2}\} \lesssim n^{-\alpha/(11\alpha+3)}$ if $\log p\leq n^c$ for some sufficiently small $c>0$. Therefore, for $\alpha>1$, we have $\omega_n\lesssim n^{-\alpha/(11\alpha+3)}\cdot \mbox{polylog}(p)$.

{\indent\bf Case 2.} If $\alpha\in(0,1)$, by \eqref{eqn:zhang_wu_rate_phystical_dependence}, we have
\begin{align}
\omega_n \lesssim &\,\,p \exp\{-C n^{\beta(\alpha c_4-c_2)}\}+p\exp\{-C n^{\beta(1-2c_2-c_3-c_4)/2}\}\notag \\
&~~+ \{n^{-c_1} + n^{-c_2}+ n^{-c_3/8}  + n^{-\alpha c_4/3}\label{eq:wnbd1}\\
&~~~~~+n^{-\alpha(1-c_3)/3}+n^{-(1-c_3-c_4)/3}\}\cdot \mbox{polylog}(p)\,.\notag
\end{align}
We still need to require $\alpha^{-1}c_2<c_4<1-2c_2-c_3$ and $c_2<\alpha(1-c_3)/(1+2\alpha)$. Following the same arguments in Case 1, we can select $c_1=c_2$. Define $$\xi(c_2,c_3,c_4)=\inf\{c_2,c_3/8,\alpha c_4/3,(1-c_3-c_4)/3,\alpha(1-c_3)/3\}$$ and $\mathcal{C}=\{(c_2,c_3,c_4):0<c_3<1, 0<c_2<\alpha(1-c_3)/(1+2\alpha), \alpha^{-1}c_2<c_4<1-2c_2-c_3\}$.

(i) If $0<c_2<\alpha(1-c_3)/(2+2\alpha)$, we have $1-2c_2-c_3> (1-c_3)/(1+\alpha)>\alpha^{-1}c_2$, which implies
\begin{align*}
&\sup_{0<c_3<1}\sup_{0<c_2<\frac{\alpha(1-c_3)}{2+2\alpha}}\sup_{\alpha^{-1}c_2<c_4<1-2c_2-c_3}\xi(c_2,c_3,c_4)\\
&~~~=\sup_{0<c_3<1}\sup_{0<c_2<\frac{\alpha(1-c_3)}{2+2\alpha}}\inf\bigg\{c_2,\frac{c_3}{8},\frac{\alpha(1-c_3)}{3\alpha+3},\frac{\alpha(1-c_3)}{3}\bigg\}\\
&~~~=\sup_{0<c_3<1}\sup_{0<c_2<\frac{\alpha(1-c_3)}{2+2\alpha}}\inf\bigg\{c_2,\frac{c_3}{8},\frac{\alpha(1-c_3)}{3\alpha+3}\bigg\}\\
&~~~=\sup_{0<c_3<1}\inf\bigg\{\frac{c_3}{8},\frac{\alpha(1-c_3)}{3\alpha+3}\bigg\}=\frac{\alpha}{11\alpha+3}\,.
\end{align*}

(ii) If $\alpha(1-c_3)/(2+2\alpha)\leq c_2<\alpha(1-c_3)/(1+2\alpha)$, we have $1-2c_2-c_3\leq(1-c_3)/(1+\alpha)$, which implies
\begin{align*}
&\sup_{0<c_3<1}\sup_{\frac{\alpha(1-c_3)}{2+2\alpha}\leq c_2<\frac{\alpha(1-c_3)}{1+2\alpha}}\sup_{\alpha^{-1}c_2<c_4<1-2c_2-c_3}\xi(c_2,c_3,c_4)\\
&~~~=\sup_{0<c_3<1}\sup_{\frac{\alpha(1-c_3)}{2+2\alpha}\leq c_2<\frac{\alpha(1-c_3)}{1+2\alpha}}\inf\bigg\{c_2,\frac{c_3}{8},\frac{\alpha(1-2c_2-c_3)}{3},\frac{\alpha(1-c_3)}{3}\bigg\}\\
&~~~=\sup_{0<c_3<1}\sup_{\frac{\alpha(1-c_3)}{2+2\alpha}\leq c_2<\frac{\alpha(1-c_3)}{1+2\alpha}}\inf\bigg\{\frac{c_3}{8},\frac{\alpha(1-2c_2-c_3)}{3}\bigg\}\\
&~~~=\sup_{0<c_3<1}\inf\bigg\{\frac{c_3}{8},\frac{\alpha(1-c_3)}{3\alpha+3}\bigg\}=\frac{\alpha}{11\alpha+3}\,.
\end{align*}
Hence, we have $$\sup_{(c_2,c_3,c_4)\in\mathcal{C}}\xi(c_2,c_3,c_4)=\frac{\alpha}{11\alpha+3}$$ and it can be achieved by selecting $c_3=8\alpha/(11\alpha+3)$, $c_2=\alpha(1-c_3)/(3\alpha+3)$ and $c_4=(1-c_3)/(\alpha+1)$. With such selected $(c_2,c_3,c_4)$, we know $\beta(\alpha c_4-c_2)>0$ and $\beta(1-2c_2-c_3-c_4)/2>0$, which implies that $p \exp\{-C n^{\beta(\alpha c_4-c_2)}\}+p\exp\{-C n^{\beta(1-2c_2-c_3-c_4)/2}\} \lesssim n^{-\alpha/(11\alpha+3)}$ if $\log p\leq n^c$ for some sufficiently small $c>0$. Therefore, for $\alpha\in(0,1)$, we also have $\omega_n\lesssim n^{-\alpha/(11\alpha+3)}\cdot \mbox{polylog}(p)$.

{\indent\bf Case 3.} If $\alpha = 1$, by \eqref{eqn:zhang_wu_rate_phystical_dependence}, we have
\begin{align*}
\omega_n \lesssim&\,\, p \exp\{-C n^{\beta(c_4-c_2)}\}+p\exp\{-C n^{\beta(1-2c_2-c_3-c_4)/2}\}\\
&+ \{n^{-c_1} + n^{-c_2}+ n^{-c_3/8}  + n^{-c_4/3}+n^{-(1-c_3-c_4)/3}\}\cdot \mbox{polylog}(p)\,,
\end{align*}
which is identical to the upper bound specified in \eqref{eq:wnbd1} with $\alpha=1$. Hence,
following the arguments for Case 1, we have $
\omega_n\lesssim n^{-1/14}\cdot \mbox{polylog}(p)$ for $\alpha=1$. $\hfill\Box$

\section{Proof of Theorem \ref{pn:GA_simple_convex}}\label{sec:pf_simple_convex}
For any $A\in\mathcal{A}$, there exists  a $K$-generated $A^{K}$ such that $A^{K}\subset A\subset A^{K,\epsilon}$ with $K\leq (pn)^d$ and $\epsilon=a/n$. Define
\begin{align*}
\rho^{(K)}=&\,\,|\mathbb{P}(S_{n,x}\in A^{K})-\mathbb{P}(G\in A^{K})|\,,\\ \rho^{(K),\epsilon}=&\,\,|\mathbb{P}(S_{n,x}\in A^{K,\epsilon})-\mathbb{P}(G\in A^{K,\epsilon})|\,.
\end{align*}
Recall $|\mathcal{V}(A^K)|=K$. Without loss of generality, we write $\mathcal{V}(A^K)=\{v_1,\ldots,v_K\}$. Define $G^{\rm si}=(G^{\rm si}_{1},\ldots,G^{\rm si}_{K})^{\T}=(v_1^{\T}G,\ldots,v_K^\T G)^\T$. By Condition \ref{as:longrun_simple_convex}, we have $$\min_{i\in [K]}\var(G^{\rm si}_{i})=\min_{v\in\mathcal{V}(A^{K})}\var\bigg(\frac{1}{\sqrt{n}}\sum_{t=1}^{n}v^{\T}X_{t}\bigg)\gtrsim K_{4}\,.$$ Write $\mathcal{S}_{A^{K}}=\{\mathcal{S}_{A^{K}}(v_1),\ldots,\mathcal{S}_{A^{K}}(v_K)\}^\T\in\mathbb{R}^{K}$. By Nazarov's inequality,
\begin{align*}
\mathbb{P}(G\in A^{K,\epsilon})=&\,\,\mathbb{P}(G^{\rm si}\leq\mathcal{S}_{A^{K}}+\epsilon)\\
\leq &\,\,\mathbb{P}(G^{\rm si}\leq\mathcal{S}_{A^{K}})+C\epsilon(\log K)^{1/2}=\mathbb{P}(G\in A^{K})+C\epsilon(\log K)^{1/2}\,.
\end{align*}
Notice that
$\mathbb{P}(S_{n,x}\in A)\leq\mathbb{P}(S_{n,x}\in A^{K,\epsilon})\leq \mathbb{P}(G\in A^{K,\epsilon})+\rho^{(K),\epsilon}$. Then
\begin{align*}
\mathbb{P}(S_{n,x}\in A)\leq&\,\,\mathbb{P}(G\in A^{K})+C\epsilon(\log K)^{1/2}+\rho^{(K),\epsilon}\\
\leq&\,\, \mathbb{P}(G\in A)+C\epsilon(\log K)^{1/2}+\rho^{(K),\epsilon}\,.
\end{align*}
Analogously, we also have
$
\mathbb{P}(S_{n,x}\in A)\geq \mathbb{P}(G\in A)-C\epsilon(\log K)^{1/2}-\rho^{(K)}$.
Hence, 
\begin{align}\label{eq:simple_upper_bound}
	|\mathbb{P}(S_{n,x}\in A)-\mathbb{P}(G\in A)|\lesssim \epsilon(\log K)^{1/2}+\rho^{(K)}+\rho^{(K),\epsilon}\,.
\end{align}
Let $X^{\rm si}_{t}=(X^{\rm si}_{t,1},\ldots,X^{\rm si}_{t,K})^{\T}=(v_1^{\T}X_{t},\ldots,v_K^\T X_t)^\T$ for each $t\in [n]$. Notice that $\cov(G^{\rm si})=\cov{(n^{-1/2}\sum_{i=1}^{n}X^{\rm si}_{t})}$. Applying Theorem \ref{tm:1} to $\{X^{\rm si}_{t}\}_{t=1}^{n}$, we have
\begin{align*}
\rho^{(K)}	=&\,\,\bigg|\mathbb{P}\bigg(\frac{1}{\sqrt{n}}\sum_{i=1}^{n}X^{\rm si}_{t}\leq\mathcal{S}_{A^{K}}\bigg)-\mathbb{P}(G^{\rm si}\leq\mathcal{S}_{A^{K}})\bigg|\\
\lesssim&\,\,  n^{-1/9}\{B_{n}^{2/3}(\log K)^{(1+2\gamma_2)/(3\gamma_2)}+B_{n}(\log K)^{7/6}\}
\end{align*}
provided that $(\log K)^{3-\gamma_2}=o(n^{\gamma_2/3})$.
Analogously, we also have $$\rho^{(K),\epsilon}\lesssim n^{-1/9}B_{n}(\log K)^{7/6}+n^{-1/9}B_{n}^{2/3}(\log K)^{(1+2\gamma_2)/(3\gamma_2)}\,.$$ Due to $K\leq (pn)^d$, $\epsilon=a/n$ and $p\geq n^{\kappa}$, by \eqref{eq:simple_upper_bound}, we have
$
|\mathbb{P}(S_{n,x}\in A)-\mathbb{P}(G\in A)|\lesssim   n^{-1}a(d\log p)^{1/2} +n^{-1/9}\{B_{n}^{2/3}(d\log p)^{(1+2\gamma_2)/(3\gamma_2)}+B_{n}(d\log p)^{7/6}\}$
for any $A\in\mathcal{A}$. $\hfill\Box$ 

\section{Proof of Theorem \ref{pn:GA_dependency_graph_simple_convex}}\label{sec:gadpsc}
The proof of Theorem \ref{pn:GA_dependency_graph_simple_convex} is almost identical to that of Theorem \ref{pn:GA_simple_convex} except the upper bounds for $\rho^{(K)}$ and $\rho^{(K),\epsilon}$. By Theorem \ref{lem:GA_dependency_graph}, we have $\rho^{(K)}+\rho^{(K),\epsilon}\lesssim n^{-1/6}B_{n}(D_nD_n^*)^{1/3}(d\log p)^{7/6}$. $\hfill\Box$

\section{Proof of Theorem \ref{pn:GA_weakly-dependent_simple_convex}}\label{sec:gawdsc}
The proof of Theorem \ref{pn:GA_weakly-dependent_simple_convex} is almost identical to that of Theorem \ref{pn:GA_simple_convex} except the upper bounds for $\rho^{(K)}$ and $\rho^{(K),\epsilon}$. For part {\rm(i)} of Theorem \ref{pn:GA_weakly-dependent_simple_convex}, by Theorem \ref{prop:GA_weakly-dependent}{\rm(i)}, we conclude that 
\begin{align*}
\rho^{(K)}+\rho^{(K),\epsilon}\lesssim&\,\,  n^{-\alpha/(3+9\alpha)}(d\log p)^{2/3}\{\Phi_{\psi_\nu,0,\mathcal{A}}(d\log p)^{1/2}+\Psi_{2,\alpha,\mathcal{A}}^{1/2} \Psi_{2,0,\mathcal{A}}^{1/3}\}\\
&+n^{-\alpha/(1+3\alpha)}\Phi_{\psi_\nu,\alpha,\mathcal{A}}(d\log p)^{1+\nu}
\end{align*}
provided that $(d\log p)^{\max\{6\nu-1,(5+6\nu)/4\}}= o\{n^{\alpha/(1+3\alpha)}\}$. For part {\rm(ii)} of Theorem \ref{pn:GA_weakly-dependent_simple_convex}, by Theorem \ref{prop:GA_weakly-dependent}{\rm(ii)}, we have 
\begin{align*}
\rho^{(K)}+\rho^{(K),\epsilon}\lesssim &\,\, n^{-\alpha/(12+6\alpha)}\{B_{n}(d\log p)^{7/6}+(\Psi_{2,\alpha,\mathcal{A}}\Psi_{2,0,\mathcal{A}})^{1/3}(d\log p)^{2/3}\}\\
&+n^{-\alpha/(4+2\alpha)}\Phi_{\psi_{\nu},\alpha,\mathcal{A}}(d\log p)^{1+\nu}\,.
\end{align*}
We complete the proof of Theorem \ref{pn:GA_weakly-dependent_simple_convex}. 
$\hfill\Box$

\section{Proof of Theorem \ref{pn:GA_s-sparsely}}\label{sec:pf_s-sparsely}  Denote by $\mathcal{U}$ the class of all convex sets in $\mathbb{R}^{p}$ and define
$
\zeta_{n}:=\sup_{U\in\mathcal{U}}|\mathbb{P}(S_{n,x}\in U)-\mathbb{P}(G\in U)|$.
Recall $\Xi={\rm Cov}(S_{n,x})$. To construct the upper bound for $\zeta_n$, we need the following condition.
\begin{ass}\label{as:eigenvalue1}
	There exists a universal constant $K_{6}>0$ such that $\lambda_{\min}(\Xi)\geq K_{6}$.
\end{ass}
Proposition \ref{pn:BE_bound_alphamixing} gives an upper bound of $\zeta_n$ when $\{X_t\}$ is an $\alpha$-mixing sequence, whose proof is given in Section \ref{sec:pfpropP1}.
\begin{proposition}\label{pn:BE_bound_alphamixing}
	Assume $\{X_t\}$ is an $\alpha$-mixing sequence. Under Conditions {\rm\ref{as:tail}}--{\rm\ref{as:alpha-mixing}} and {\rm\ref{as:eigenvalue1}}, it holds that
	$
	\zeta_{n}\lesssim B_{n}^{4/3}p^{17/6}n^{-1/6}\{\log (pn)\}^{1/(6\gamma_{2})}(\log n)^{1/6}$.
\end{proposition}

Let $B=\{w=(w_1,\ldots,w_p)^\T\in\mathbb{R}^{p}:\max_{j\in[p]}|w_{j}|\leq pn^{5/2}\}$. Without loss of generality, we assume $B_{n}^{2}(\log p)^{(1+2\gamma_{2})/\gamma_{2}}=o(n^{1/3})$ since otherwise the assertions are trivial. By Markov inequality, it holds that 
\begin{align*}
\mathbb{P}(S_{n,x}\in A\cap B^{c})\leq&\,\,\mathbb{P}(S_{n,x}\in B^{c})=\mathbb{P}\bigg(\max_{j\in [p]}\bigg|\frac{1}{\sqrt{n}}\sum_{t=1}^{n}X_{t,j}\bigg|> pn^{5/2}\bigg)\\
\leq&\,\,\mathbb{P}\bigg(\max_{t\in[n],j\in[p]}|X_{t,j}|>pn^2\bigg)\leq p^{-1}n^{-2}\mathbb{E}\bigg(\max_{t\in[n],j\in[p]}|X_{t,j}|\bigg)\\
\leq&\,\, n^{-1}\max_{t\in[n],j\in[p]}\mathbb{E}(|X_{t,j}|)\lesssim B_{n}n^{-1}\lesssim n^{-1/2}\,.
\end{align*}
Recall $G\sim N(0,\Xi)$ with $\Xi=\Cov(S_{n,x})$.   Let $\{Y_{t}\}_{t=1}^n$, independent of $\mathcal{X}_n=\{X_1,\ldots,X_n\}$, be a centered Gaussian sequence such that $\Cov(Y_{t}, Y_{s}) = \Cov(X_{t}, X_{s})$ for all $t,s\in[n]$. Then $n^{-1/2}\sum_{t=1}^{n}Y_{t}\overset{d}{=}G$. Write $Y_{t}=(Y_{t,1},\ldots,Y_{t,p})^{\T}$. By Markov inequality again,
\begin{align*}
\mathbb{P}(G\in A\cap B^{c})\leq&\,\,\mathbb{P}(G\in B^{c})=\mathbb{P}\bigg(\max_{j\in [p]}\bigg|\frac{1}{\sqrt{n}}\sum_{t=1}^{n}Y_{t,j}\bigg|\geq pn^{5/2}\bigg)\\
\leq&\,\, n^{-1}\max_{t\in[n],j\in[p]}\mathbb{E}(|Y_{t,j}|)\leq n^{-1}\max_{t\in[n],j\in[p]}\{\mathbb{E}(|X_{t,j}|^2)\}^{1/2}
\lesssim 
n^{-1/2}\,.
\end{align*}
Thus 
$
|\mathbb{P}(S_{n,x}\in A)-\mathbb{P}(G\in A)|\leq |\mathbb{P}(S_{n,x}\in A\cap B)-\mathbb{P}(G\in A\cap B)|+Cn^{-1/2}$. 
Let $A^*=A\cap B$. Since $A$ is an $s$-sparsely convex set, $A^*$ has a sparse representation $A^{*}=\cap_{q=1}^{K_*}A_{q}^{*}$ such that $\max_{q\in [K^*]}\sup_{w\in A^{*}_{q}}\max_{j\in J(A_{q}^{*})}|w_{j}|\leq pn^{5/2}$, where $K_*\lesssim p^s$ and $J(A_{q}^{*})$ is the set containing the indexes of the main components of $A_{q}^{*}$. Write $s_{q}=|J(A_{q}^{*})|$. For any convex set $U\in\mathbb{R}^{p}$ and $\varepsilon>0$, let  $U^{\varepsilon}=\{w\in\mathbb{R}^{p}:\inf_{v\in U}|v-w|_{2}\leq \varepsilon\}$ and $U^{-\varepsilon}=\{w\in\mathbb{R}^{p}:B(w,\varepsilon)\subset U\}$ with $B(w,\varepsilon)=\{v\in\mathbb{R}^{p}:|v-w|_{2}\leq\varepsilon\}$. There are three cases for $A^{*}$.

{\bf Case 1.} ({\it $A^{*}$ contains a ball with radius $\epsilon=n^{-1}$ and center at some $w^{*}\in A^{*}$}) Letting $D:=A^{*}-w^{*}=\{w\in\mathbb{R}^{p}:w+w^{*}\in A^{*}\}$, $D$ must contain a ball with radius $\epsilon$ and center at the origin. Then $D$ admits a sparse representation $D=\cap_{q=1}^{K_*}D_{q}$ with $D_{q}=A_{q}^{*}-w^{*}$ also containing a ball with radius $\epsilon$ and center at the origin. Next, we will show $-D_{q}\subset\mu D_{q}$ for some $\mu\geq 1$. Since $I(w\in D_q)$ depends on at most $s$ components of the vector $w\in\mathbb{R}^p$, it suffices to show $-\mathring{D}_{q}\subset\mu\mathring{D}_{q}$, where $\mathring{D}_{q}\subset \mathbb{R}^{s_q}$ is the set generated only by the main components of $D_{q}$. Clearly, $\mathring{D}_{q}$ also contains a ball in $\mathbb{R}^{s_q}$ with radius $\epsilon$ and center at the origin. Hence, $2p^{3/2}n^{5/2}\epsilon^{-1} \mathring{D}_{q}$ must contain a ball in $\mathbb{R}^{s_q}$ with radius $2p^{3/2}n^{5/2}$ and center at the origin. Note that $\max_{q\in [K^*]}\sup_{w\in \mathring{D}_{q}}|w|_{2}\leq 2p^{3/2}n^{5/2}$. Then each $\mathring{D}_{q}$ is contained in a ball in $\mathbb{R}^{s_q}$ with radius $2p^{3/2}n^{5/2}$ and center at the origin, so is $-\mathring{D}_{q}$. Thus, $-\mathring{D}_{q}\subset 2p^{3/2}n^{5/2}\epsilon^{-1} \mathring{D}_{q}=2p^{3/2}n^{7/2}\mathring{D}_{q}$, which implies $-D_{q}\subset 2p^{3/2}n^{7/2}D_{q}$. Notice that $\sup_{w\in D}|w|_2\leq 2p^{3/2}n^{5/2}$. By applying Lemma D.1 of Chernozhukov, Chetverikov and Kato (2017) to the set $D$ with $R=2p^{3/2}n^{5/2}$ and $\mu=2p^{3/2}n^{7/2}$, and recalling $A^{*}=D+w^{*}$, we know $A^{*}\in\mathcal{A}^{\rm si}(a,d)$ is a simple convex set with $a=1$ and $d=Cs^{2}$ for some positive constant $C$, and the corresponding $K$-generated convex set $A^{K}:=A^{K}(A^*)$ satisfies $\max_{v\in\mathcal{V}(A^{K})}|v|_{0}\leq s$. We will use Theorem \ref{pn:GA_simple_convex} to $A^{*}$ for deriving an upper bound of $|\mathbb{P}(S_{n,x}\in A^{*})-\mathbb{P}(G\in A^{*})|$. For any $v=(v_{1},\ldots,v_{p})^{\T}\in\mathcal{V}(A^{K})$, by Condition \ref{as:tail}, we have $
\max_{v\in\mathcal{V}(A^{K})}\|v^{\T}X_{t}\|_{\psi_{\gamma_{1}}}\leq sB_{n}$. Applying Theorem \ref{pn:GA_simple_convex} with replacing $B_n$ by $sB_n$, 
\begin{align*}
|\mathbb{P}(S_{n,x}\in A^{*})-\mathbb{P}(G\in A^{*})|\lesssim&\,\, n^{-1/9}B_{n}^{2/3}s^{(2+6\gamma_{2})/(3\gamma_{2})}(\log p)^{(1+2\gamma_2)/(3\gamma_2)}\\
&+n^{-1/9}B_{n}s^{10/3}(\log p)^{7/6}
\end{align*}
provided that $(s^2\log p)^{3-\gamma_2}=o(n^{\gamma_2/3})$.

{\bf Case 2.} ({\it $A^{*}=\cap_{q=1}^{K_{*}}A_{q}^*$ doesn't contain a ball with radius $\epsilon=n^{-1}$, and there exists one $A_{q}^{*}$ containing no ball with radius $\epsilon$}) For this given $q$, we define $\tau_{q}(w)=(w_{j})_{j\in J(A_{q}^{*})}\in\mathbb{R}^{s_{q}}$ for any $w=(w_{1},\ldots,\,w_{p})^{\T}\in\mathbb{R}^{p}$. Let $\mathring{A}_{q}^{*}=\{\tau_q(w):w\in A_q^*\}$.   Clearly, $\mathring{A}_{q}^{*}$ doesn't contain a ball in $\mathbb{R}^{s_q}$ with radius $\epsilon$. Then $(\mathring{A}_{q}^{*})^{-\epsilon}=\emptyset$.
Applying Lemma A.2 of Chernozhukov, Chetverikov and Kato (2017) to $\tau_{q}(G)$ with $Q=\mathring{A}_{q}^{*}$ and $h_{1}=h_{2}=\epsilon$, we have
$
\mathbb{P}(G\in A_{q}^*)=\mathbb{P}\{\tau_{q}(G)\in \mathring{A}_{q}^*\}\lesssim \epsilon\|\Xi_{s_{q}}^{-1}\|_{\rm HS}^{1/2}$,
where $\Xi_{s_{q}}=\Cov\{\tau_{q}(G)\}\in\mathbb{R}^{s_{q}\times s_{q}}$ and $\|\cdot\|_{\rm HS}$ is the Hilbert-Schmidt norm. Let $\lambda_{1}\geq\cdots\geq\lambda_{s_{q}}>0$ be the $s_{q}$ eigenvalues of $\Xi_{s_{q}}$. Since
$
\|\Xi_{s_{q}}^{-1}\|_{\rm HS}=(\sum_{i=1}^{s_{q}}\lambda_{i}^{-2})^{1/2}\leq s_{q}^{1/2}\lambda_{s_{q}}^{-1}\lesssim s^{1/2}$,
we have $\mathbb{P}(G\in A^*)\leq \mathbb{P}(G\in A_{q}^*)\lesssim s^{1/4} n^{-1}$. Notice that $\tau_{q}(S_{n,x})=n^{-1/2}\sum_{t=1}^{n}\tau_{q}(X_{t})$ is a partial sum of an $s_q$-dimensional sequence. Applying the Berry-Essen bound stated in Proposition \ref{pn:BE_bound_alphamixing} for such $s_q$-dimensional $\tau_q(S_{n,x})$ by setting the general parameter $p$ (that is different from the notation $p$ used in current proof) appeared in Proposition \ref{pn:BE_bound_alphamixing} equals to $s_q$, we have
\begin{align*}
|\mathbb{P}(S_{n,x}\in A^{*}_{q})-\mathbb{P}(G\in A^{*}_{q})|=&\,\,|\mathbb{P}\{\tau_q(S_{n,x})\in \mathring{A}_{q}^{*}\}-\mathbb{P}\{\tau_q(G)\in \mathring{A}_{q}^{*}\}|\\
\lesssim&\,\,    n^{-1/6}B_{n}^{4/3}s^{17/6}(\log p)^{(1+\gamma_{2})/(6\gamma_{2})}\,,
\end{align*}
which implies 
\begin{align*}
\mathbb{P}(S_{n,x}\in A^{*})\leq \mathbb{P}(S_{n,x}\in A^{*}_{q})\leq&\,\,|\mathbb{P}(S_{n,x}\in A^{*}_{q})-\mathbb{P}(G\in A^{*}_{q})|+\mathbb{P}(G\in A^{*}_{q}) \\
\lesssim&\,\,   n^{-1/6}B_{n}^{4/3}s^{17/6}(\log p)^{(1+\gamma_{2})/(6\gamma_{2})}\,.
\end{align*}
Due to $B_{n}^{2}(\log p)^{(1+2\gamma_{2})/\gamma_{2}}=o(n^{1/3})$, it holds that $$n^{-1/6}B_{n}^{4/3}s^{17/6}(\log p)^{(1+\gamma_{2})/(6\gamma_{2})}\lesssim n^{-1/9}B_{n}s^{10/3}(\log p)^{7/6}\,.$$ By triangle inequality, we have
$$
|\mathbb{P}(S_{n,x}\in A^{*})-\mathbb{P}(G\in A^{*})|\leq\mathbb{P}(S_{n,x}\in A^{*})+\mathbb{P}(G\in A^{*})\lesssim  n^{-1/9}B_{n}s^{10/3}(\log p)^{7/6}\,.$$

{\bf Case 3.} ({\it each $A_{q}^{*}$ contains a ball with radius $\epsilon=n^{-1}$, but they are not all the same}) In this case, each $\mathring{A}_{q}^{*}$ must contain a ball in $\mathbb{R}^{s_q}$ with radius $\epsilon$ and $\sup_{w\in \mathring{A}^{*}_{q}}|w|_{\infty}\leq pn^{5/2}$. Using the same arguments in Case 1, for each $\mathring{A}_{q}^{*}$, there exists a $K_q$-generated convex set $\mathring{A}_{q}^{K_q}$ with $K_q\leq(pn)^{Cs^{2}}$ such that $\mathring{A}_{q}^{K_q}\subset \mathring{A}_{q}^{*} \subset \mathring{A}_{q}^{K_q,\epsilon}$ and $\max_{v\in\mathcal{V}(\mathring{A}^{K_q}_{q})}|v|_{0}\leq s$. The definition of $\mathring{A}_q^{K_q,\epsilon}$ is given in \eqref{eq:inclusion}. Hence, for each $A_{q}^{*}$, there also exists a $K_q$-generated convex set $A_{q}^{K_q}$ with $K_q\leq(pn)^{Cs^{2}}$ such that $A_{q}^{K_q}\subset A_{q}^{*} \subset A_{q}^{K_q,\epsilon}$ and $\max_{v\in\mathcal{V}(A^{K_q}_{q})}|v|_{0}\leq s$. Let  $A^{0}:=\cap_{q=1}^{K_{*}}A_{q}^{K_q,\epsilon}$. We know $A^0\in\mathcal{A}^{\rm si}(a,d)$ with $a=0$ and $d= Cs^{2}$. Notice that $\max_{v\in\mathcal{V}(A^{0})}|v|_{0}\leq s$. Applying Theorem \ref{pn:GA_simple_convex} with replacing $B_n$ by $sB_n$, 
$
|\mathbb{P}(S_{n,x}\in A^{0})-\mathbb{P}(G\in A^{0})|\lesssim n^{-1/9}\{B_{n}^{2/3}s^{(2+6\gamma_{2})/(3\gamma_{2})}(\log p)^{(1+2\gamma_2)/(3\gamma_2)}+B_{n}s^{10/3}(\log p)^{7/6}\}$ 
provided that $(s^2\log p)^{3-\gamma_2}=o(n^{\gamma_2/3})$.
For any $x\in\cap_{q=1}^{K_{*}}A_{q}^{K_q,-\epsilon}$ where $A_{q}^{K_q,-\epsilon}$ follows the definition in \eqref{eq:inclusion} with replacing $\epsilon$ by $-\epsilon$, then $v^{\T}x\leq\mathcal{S}_{A_{q}^{K_q}}(v)-\epsilon$ for all $v\in\mathcal{V}(A_{q}^{K_q})$ and $q\in [K_{*}]$. For any $y\in B(x,\epsilon)$, we have $v^{\T}y=v^{\T}(y-x)+v^{\T}x\leq\mathcal{S}_{A^{K_q}_{q}}(v)$ for all $v\in\mathcal{V}(A_{q}^{K_q})$ and $q\in [K_{*}]$, which implies $B(x,\epsilon)\subset\cap_{q=1}^{K_{*}}A_{q}^{K_q}$ and $x\in(\cap_{q=1}^{K_{*}}A_{q}^{K_q})^{-\epsilon}$. Hence
$
\cap_{q=1}^{K_{*}}A_{q}^{K_q,-\epsilon}\subset(\cap_{q=1}^{K_{*}}A_{q}^{K_q})^{-\epsilon}\subset (A^*)^{-\epsilon}=\emptyset$,
where the second step is due to $\cap_{q=1}^{K_{*}}A_{q}^{K_q}\subset\cap_{q=1}^{K_{*}}A^*_{q}=A^*$ and the third step is due to our assumption that $A^{*}$ contains no ball with radius $\epsilon$. Due to $p\geq n^\kappa$ for some $\kappa>0$, by Nazarov's inequality,
$
\mathbb{P}(G\in A^{0})=\mathbb{P}\{v^{\T}G\leq\mathcal{S}_{A_{q}^{K_q}}(v)+\epsilon\,\text{ for all }v\in\mathcal{V}(A_{q}^{K_q})\,,q\in[K_{*}]\}-\mathbb{P}\{v^{\T}G\leq\mathcal{S}_{A_{q}^{K_q}}(v)-\epsilon\,\text{ for all }v\in\mathcal{V}(A_{q}^{K_q})\,,q\in[K_{*}]\}\lesssim  n^{-1}(s^2\log p)^{1/2}$,
where the last step is due to $K_{*}\lesssim p^s$ and $\epsilon=n^{-1}$. Then we have 
\begin{align*}
\mathbb{P}(S_{n,x}\in A^{0})\lesssim n^{-1/9}B_{n}^{2/3}s^{(2+6\gamma_{2})/(3\gamma_{2})}(\log p)^{(1+2\gamma_2)/(3\gamma_2)}+n^{-1/9}B_{n}s^{10/3}(\log p)^{7/6}\,.
\end{align*}
Since $A^*=\cap_{q=1}^{K_{*}}A_{q}^{*}\subset\cap_{q=1}^{K_{*}}A_{q}^{K,\epsilon}=A^{0}$, $\mathbb{P}(G\in A^{*})\leq \mathbb{P}(G\in A^{0})\lesssim n^{-1}(s^2\log p)^{1/2}$ and 
\begin{align*}
\mathbb{P}(S_{n,x}\in A^{*})\leq&\,\,\mathbb{P}(S_{n,x}\in A^{0})\\
\lesssim&\,\, n^{-1/9}B_{n}^{2/3}s^{(2+6\gamma_{2})/(3\gamma_{2})}(\log p)^{(1+2\gamma_2)/(3\gamma_2)}+n^{-1/9}B_{n}s^{10/3}(\log p)^{7/6}\,,
\end{align*}
by triangle inequality, we have that 
\begin{align*}
&|\mathbb{P}(S_{n,x}\in A^{*})-\mathbb{P}(G\in A^{*})|\leq \mathbb{P}(S_{n,x}\in A^{*})+\mathbb{P}(G\in A^{*})\\
&~~~~~\lesssim  n^{-1/9}\{B_{n}s^{10/3}(\log p)^{7/6}+B_{n}^{2/3}s^{(2+6\gamma_{2})/(3\gamma_{2})}(\log p)^{(1+2\gamma_2)/(3\gamma_2)}\}\,.
\end{align*}

Combining the results in the three cases, we complete the proof of Theorem \ref{pn:GA_s-sparsely}.$\hfill\Box$

\section{Proof of Proposition \ref{pn:BE_bound_alphamixing}}\label{sec:pfpropP1}
We will still use the data blocking technique used in the proof of Theorem \ref{tm:1} in Section \ref{sec:pfalphamix}. Following the notation used in Section \ref{sec:pfalphamix}, define
$
\zeta_{n}^{(1)}:=\sup_{U\in\mathcal{U}}|\mathbb{P}\{S_{n,x}^{(1)}\in U\}-\mathbb{P}\{G^{(1)}\in U\}|$ and $
\zeta_{n}^{(2)}:=\sup_{U\in\mathcal{U}}|\mathbb{P}(S_{n,x}\in U)-\mathbb{P}\{G^{(1)}\in U\}|$,
where $S_{n,x} = n^{-1/2} \sum_{t=1}^n X_t$, $S_{n,x}^{(1)}=L^{-1/2}\sum_{\ell=1}^{L}\tilde{X}_{\ell}$ with $\tilde{X}_\ell=b^{-1/2}\sum_{t\in\mathcal{I}_\ell}X_t$, and $G^{(1)}\sim N\{0,\Xi^{(1)}\}$ with $\Xi^{(1)}=\cov\{S_{n,x}^{(1)}\}$. Without loss of generality, we assume $B_{n}^{8}p^{17}\{\log(pn)\}^{1/\gamma_{2}}\log n=o(n)$. Otherwise the assertions hold trivially. We first present next lemmas whose proofs are given in Sections \ref{sec:pflem13}--\ref{sec:pflem15}, respectively.

\begin{lemma}\label{la:BE_bound_Block}
	Assume Conditions {\rm\ref{as:tail}}--{\rm\ref{as:alpha-mixing}} and {\rm\ref{as:eigenvalue1}} hold.
	Let $h=C\{\log(pn)\}^{1/\gamma_{2}}$ for some sufficiently large constant $C>0$.
	If $B_{n}^{2}p\{\log (pn)\}^{1/\gamma_{2}}\ll b\ll B_{n}^{-2}p^{-1}n$,
	then $ \zeta_{n}^{(1)}\lesssim L^{-1/4}B_{n}^{3/2}p^{13/4}$.
\end{lemma}
\begin{lemma}\label{la:BE_bound_non_Block}
	Assume Conditions  {\rm\ref{as:tail}}--{\rm\ref{as:alpha-mixing}} and {\rm\ref{as:eigenvalue1}} hold. Let $h=C\{\log (pn)\}^{1/\gamma_{2}}$ for some sufficiently large constant $C>0$. If $b$ satisfies 
	$
	\max[B_{n}^{-2/3}n^{1/3}p^{-10/3}\{\log(pn) \}^{2(\gamma_{2}+1)/(3\gamma_{2})}, B_{n}^{2}p\{\log (pn)\}^{1/\gamma_{2}}]\ll b\ll \min(n^{1/2},B_{n}^{-2}p^{-1}n)$,
	then $\zeta_{n}^{(2)}\lesssim L^{-1/4}B_{n}^{3/2}p^{13/4}$.
\end{lemma}
\begin{lemma}\label{la:comparison_convex}
	Let $X\sim N(0,\Gamma_{1})$ and $Y\sim N(0,\Gamma_{2})$ be two p-dimensional Gaussian random vectors. If $\lambda_{\min}(\Gamma_{2})\geq c$ for some constant $c>0$, then there exists a universal constant $C>0$ depending only on $c$ such that $\sup_{U\in\mathcal{U}}|\mathbb{P}(X\in U)-\mathbb{P}(Y\in U)|\leq Cp^{2}|\Gamma_{1}-\Gamma_{2}|_{\infty}^{1/2}$, where $\mathcal{U}$ is the class of all convex sets in $\mathbb{R}^{p}$.
\end{lemma}

Recall $G\sim N(0,\Xi)$. By Condition \ref{as:eigenvalue1} and Lemma \ref{la:comparison_convex}, it holds that
$
\sup_{U\in\mathcal{U}}|\mathbb{P}\{G^{(1)}\in U\}-\mathbb{P}(G\in U)|\lesssim p^{2}|\Xi^{(1)}-\Xi|_{\infty}^{1/2}$.
Select $h=C\{\log (pn)\}^{1/\gamma_{2}}$ for some sufficiently large $C>0$. As shown in the proof of Lemma \ref{la:xtildestar_tail} in Section \ref{sec:pflem15}, 
$|\Xi^{(1)}-\Xi|_{\infty}\leq B_{n}^{2}(hb^{-1}+bn^{-1})$. Since $L\asymp nb^{-1}$, by Lemma \ref{la:BE_bound_non_Block}, we have
\begin{align*}
\zeta_{n}\leq&\,\, \zeta_{n}^{(2)}+\sup_{U\in\mathcal{U}}|\mathbb{P}\{G^{(1)}\in U\}-\mathbb{P}(G\in U)|\\
\lesssim&\,\, B_{n}^{3/2}p^{13/4}n^{-1/4}b^{1/4}+B_{n}p^{2}\{\log (pn)\}^{1/(2\gamma_{2})}b^{-1/2}
\end{align*}
provided that $b$ satisfying $\max[B_{n}^{-2/3}n^{1/3}p^{-10/3}\{\log(pn) \}^{2(\gamma_{2}+1)/(3\gamma_{2})},B_{n}^{2}p\{\log (pn)\}^{1/\gamma_{2}}] \ll b\ll  \min(n^{1/2}, B_{n}^{-2}p^{-1}n)$.
Selecting $b\asymp B_{n}^{-2/3}p^{-5/3}n^{1/3}\{\log (pn)\}^{2/(3\gamma_{2})}(\log n)^{2/3}$, we conclude that
$
\zeta_{n}\lesssim B_{n}^{4/3}n^{-1/6}p^{17/6} \{\log (pn)\}^{1/(6\gamma_{2})}(\log n)^{1/6}$.
We complete the proof of Proposition \ref{pn:BE_bound_alphamixing}. $\hfill\Box$

\subsection{Two auxiliary lemmas}\label{sec:pflem15}
For any $\ell\in[L]$, define $\tilde{X}_\ell^*=b^{-1/2}\sum_{t\in\mathcal{I}_\ell}[\{\Xi^{(1)}\}^{-1/2}X_{t}]=(\tilde{X}_{\ell,1}^*,\ldots,\tilde{X}_{\ell,p}^*)^\T$ with $\Xi^{(1)}={\rm Cov}\{S_{n,x}^{(1)}\}$ and $S_{n,x}^{(1)}=L^{-1/2}\sum_{\ell=1}^{L}\tilde{X}_{\ell}$. Let
$
S_{n,x}^{(1)*}=L^{-1/2}\sum_{\ell=1}^{L}\tilde{X}^{*}_{\ell}$. To prove Lemma \ref{la:BE_bound_Block}, we need the following auxiliary lemmas.
\begin{lemma}\label{la:xtildestar_tail}
	Assume Conditions {\rm\ref{as:tail}--\ref{as:alpha-mixing}} and {\rm\ref{as:eigenvalue1}} hold. Let $\gamma=\gamma_{2}/(2\gamma_{2}+1)$ and $h=C\{\log(pn)\}^{1/\gamma_{2}}$ for some sufficiently large $C>0$. Then there exists a universal constant $C>0$ such that $\max_{\ell\in[L]}\max_{j\in [p]}\mathbb{P}(|\tilde{X}_{\ell,j}^{*}|>\lambda)\lesssim p\exp(-Cb^{\gamma/2}p^{-\gamma/2}\lambda^{\gamma}B_{n}^{-\gamma})+p\exp(-Cp^{-1}\lambda^{2}B_{n}^{-2})$
	for any $\lambda>0$, provided that $B_{n}^{2}p\{\log (pn)\}^{1/\gamma_{2}}\ll b\ll B_{n}^{-2}p^{-1}n$.
\end{lemma}
\noindent{\it Proof.} Write $A=\{\Xi^{(1)}\}^{-1/2}=(a_1,\ldots,a_p)$ where $a_j\in\mathbb{R}^p$ is the $j$-th column. Recall that $\tilde{X}_\ell=b^{-1/2}\sum_{t\in\mathcal{I}_\ell}X_t$. For any $\ell\in [L]$ and $j\in [p]$, we have
\begin{align}\label{eq:tuta_x}
|\tilde{X}_{\ell,j}^*|=|a_j^{\T}\tilde{X}_{\ell}|\leq |a_{j}|_1|\tilde{X}_{\ell}|_{\infty}\leq\sqrt{p}|a_{j}|_2 |\tilde{X}_{\ell}|_{\infty}\leq\sqrt{p}\lambda_{\min}^{-1/2}\{\Xi^{(1)}\}|\tilde{X}_{\ell}|_{\infty}\,.
\end{align}
Recall $\tilde{\Xi}=L^{-1}\sum_{\ell=1}^{L}\mathbb{E}(\tilde X_\ell \tilde X_\ell ^{\T})$. Then
\begin{align}\label{eq:cov_infty_x}
|\Xi^{(1)}-\tilde{\Xi}|_{\infty}&=\max_{j_{1},j_{2}\in [p]}\bigg|\frac{1}{L}\sum_{\ell_{1},\ell_{2}=1}^{L}\mathbb{E}(\tilde{X}_{\ell_{1},j_{1}}\tilde{X}_{\ell_{2},j_{2}})-\frac{1}{L}\sum_{\ell=1}^{L}\mathbb{E}(\tilde{X}_{\ell,j_{1}}\tilde{X}_{\ell,j_{2}})\bigg|\notag\\
&=\max_{j_{1},j_{2}\in [p]}\bigg|\frac{1}{L}\sum_{\ell_{1}\neq\ell_{2}}\mathbb{E}(\tilde{X}_{\ell_{1},j_{1}}\tilde{X}_{\ell_{2},j_{2}})\bigg|\,.
\end{align}
For any $\ell_{1}>\ell_{2}+1$ and $j_{1},j_{2}\in [p]$, by Davydov's inequality,
\begin{align*}
|\mathbb{E}(\tilde{X}_{\ell_{1},j_{1}}\tilde{X}_{\ell_{2},j_{2}})|\lesssim&\, B_{n}^2\exp\{-C(|\ell_{1}-\ell_{2}|Q-b)^{\gamma_{2}}\}\\
\lesssim&\, B_{n}^2\exp\{-C(|\ell_{1}-\ell_{2}|-1)^{\gamma_{2}}b^{\gamma_{2}}\}\\
\lesssim&\, B_{n}^{2}\exp\{-C(|\ell_{1}-\ell_{2}|-1)^{\gamma_{2}}\}\exp(-Cb^{\gamma_{2}})
\end{align*}
holds uniformly over $j_{1},j_{2}\in [p]$, where the last step is due to $(|\ell_{1}-\ell_{2}|-1)^{\gamma_{2}}+b^{\gamma_{2}}\lesssim (|\ell_{1}-\ell_{2}|-1)^{\gamma_{2}}b^{\gamma_{2}}$ for $\ell_{1}>\ell_{2}+1$. Then
$$\bigg|\frac{1}{L}\sum_{\ell_{1}>\ell_{2}+1}\mathbb{E}(\tilde{X}_{\ell_{1},j_{1}}\tilde{X}_{\ell_{2},j_{2}})\bigg|\lesssim B_{n}^{2}\exp(-Cb^{\gamma_{2}})\lesssim B_{n}^{2}hb^{-1}\,.$$
Due to $h=C\{\log(pn)\}^{1/\gamma_{2}}$, applying Davydov's inequality again,
$$
\bigg|\frac{1}{L}\sum_{\ell=1}^{L-1}\mathbb{E}(\tilde{X}_{\ell+1,j_{1}}\tilde{X}_{\ell,j_{2}})\bigg|\lesssim B_{n}^{2}\exp(-Ch^{\gamma_{2}})\lesssim B_{n}^{2}bn^{-1}\,.$$
Thus $\max_{j_1,j_2\in[p]}|L^{-1}\sum_{\ell_{1}>\ell_{2}}\mathbb{E}(\tilde{X}_{\ell_{1},j_{1}}\tilde{X}_{\ell_{2},j_{2}})|\lesssim B_{n}^{2}(hb^{-1}+bn^{-1})$. Analogously, we also have $\max_{j_1,j_2\in[p]}|L^{-1}\sum_{\ell_{1}<\ell_{2}}\mathbb{E}(\tilde{X}_{\ell_{1},j_{1}}\tilde{X}_{\ell_{2},j_{2}})|\lesssim B_{n}^{2}(hb^{-1}+bn^{-1})$. It then follows from \eqref{eq:cov_infty_x} that $|\Xi^{(1)}-\tilde{\Xi}|_{\infty}\lesssim B_{n}^{2}(hb^{-1}+bn^{-1})$.
Together with Lemma \ref{lem:cov_infinity}, we have
$$|\Xi^{(1)}-\Xi|_{\infty}\leq|\Xi^{(1)}-\tilde{\Xi}|_{\infty}+|\tilde{\Xi}-\Xi|_{\infty}\lesssim B_{n}^{2}(hb^{-1}+bn^{-1})\,.$$
For any $v\in\mathbb{S}^{p-1}$,
$
|v^{\T}\{\Xi^{(1)}-\Xi\}v|\leq |\Xi^{(1)}-\Xi|_{\infty}|v|_{1}^2\leq B_{n}^{2}p(hb^{-1}+bn^{-1})=o(1)$
provided that $B_{n}^{2}p\{\log (pn)\}^{1/\gamma_{2}}\ll b\ll B_{n}^{-2}p^{-1}n$. By Condition \ref{as:eigenvalue1}, we have $\lambda_{\min}(\Xi)\geq K_6$. Therefore,
\begin{align}\label{eq:eigen_min}
\lambda_{\min}\{\Xi^{(1)}\}=\min_{v\in\mathbb{S}^{p-1}}v^{\T}\Xi^{(1)}v\geq \min_{v\in\mathbb{S}^{p-1}}v^{\T}\{\Xi^{(1)}-\Xi\}v+\lambda_{\min}(\Xi)\geq K_6/2
\end{align}
for any sufficiently large $n$. Then \eqref{eq:tuta_x} implies $\max_{\ell\in[L]}\max_{j\in [p]}|\tilde{X}_{\ell,j}^{*}|\lesssim \sqrt{p}|\tilde{X}_{\ell}|_{\infty}$. By \eqref{eq:xtildetail}, it holds that
\begin{align*}
\max_{\ell\in[L]}\max_{j\in [p]}\mathbb{P}(|\tilde{X}_{\ell,j}^{*}|>\lambda)&\leq p\max_{\ell\in[L]}\max_{j\in [p]}\mathbb{P}(|\tilde{X}_{\ell,j}|>Cp^{-1/2}\lambda)\\
&\lesssim p\exp(-Cp^{-1}\lambda^{2}B_{n}^{-2})+p\exp(-Cb^{\gamma/2}p^{-\gamma/2}\lambda^{\gamma}B_{n}^{-\gamma})
\end{align*}
for any $\lambda>0$. We complete the proof of Lemma \ref{la:xtildestar_tail}.$\hfill\Box$

\bigskip

Notice that $\mathbb{E}(|\tilde{X}^*_{\ell,j}|^3)=3\int^{\infty}_{0}\lambda^2\mathbb{P}(|\tilde{X}^*_{\ell,j}|>\lambda)\,{\rm d}\lambda$. By Lemma \ref{la:xtildestar_tail}, we have the following lemma.
\begin{lemma}\label{la:xtildestar_mom}
	Assume Conditions {\rm\ref{as:tail}--\ref{as:alpha-mixing}} and {\rm\ref{as:eigenvalue1}} hold. Let $h=C\{\log(pn)\}^{1/\gamma_{2}}$ for some sufficiently large constant $C>0$. Then $\max_{\ell\in[L]}\max_{j\in[p]}\mathbb{E}(|\tilde{X}_{\ell,j}^*|^3)\lesssim B_{n}^{3}p^{5/2}$
	provided that $B_{n}^{2}p\{\log (pn)\}^{1/\gamma_{2}}\ll b\ll B_{n}^{-2}p^{-1}n$.
\end{lemma}

\subsection{Proof of Lemma \ref{la:BE_bound_Block}}\label{sec:pflem13}
Without loss of generality, we assume $G^{(1)}$ is independent of $S_{n,x}^{(1)}$. Let $A=\{\Xi^{(1)}\}^{-1/2}$. Then $S_{n,x}^{(1)*}=AS_{n,x}^{(1)}$. By the definition of $\zeta_n^{(1)}$, we have
\begin{align}\label{eq:transform_BE_block}
\zeta_{n}^{(1)}&=\sup_{U\in\mathcal{U}}|\mathbb{P}\{S_{n,x}^{(1)*}\in A\circ U\}-\mathbb{P}(Z\in A \circ U)|=\sup_{U\in\mathcal{U}}|\mathbb{P}\{S_{n,x}^{(1)*}\in U\}-\mathbb{P}(Z\in U)|\,,
\end{align}
where $A\circ U:=\{Au:u\in U\}$, and $Z$ is a $p$-dimensional normal random vector with mean zero and covariance ${\rm I}_p$ that is independent of $S_{n,x}^{(1)*}$.  Given $U\in\mathcal{U}$, we define $h(x)=I(x\in U)$ and $\tilde{h}(x)=h(x)-\int_{\mathbb{R}^p} h(u)\phi(u)\,{\rm d}u$, where $\phi(\cdot)$ is the density function of $p$-dimensional standard normal distribution. Then $\zeta_{n}^{(1)}=\sup_{U\in\mathcal{U}}|\mathbb{E}[\tilde{h}\{S_{n,x}^{(1)*}\}]|$. For $a>0$, let $\mathscr{T}_{a}\tilde{h}(x)=\int_{\mathbb{R}^{p}} \tilde{h}(e^{-a}x+\sqrt{1-e^{-2a}}u)\phi(u)\,{\rm d}u$ and $\psi_{a}(x)=-\int_{a}^{\infty}\int_{\mathbb{R}^p}\tilde{h}(e^{-s}x+\sqrt{1-e^{-2s}}u)\phi(u)\,{\rm d}u{\rm d}s$. As in (1.14) and (3.1) of Bhattacharya and Holmes (2010), we have
\begin{align}\label{eq:expression}
\mathbb{E}[\mathscr{T}_{a}\tilde{h}\{S_{n,x}^{(1)*}\}]=\mathbb{E}[\Delta\psi_{a}\{S_{n,x}^{(1)*}\}-S_{n,x}^{(1)*}\cdot\nabla\psi_{a}\{S_{n,x}^{(1)*}\}]\,,
\end{align}
where $\Delta$ and $\nabla$ denote, respectively, the Laplace operator and gradient operator, and $\cdot$ is the inner product of two $p$-dimensional vectors. Our proof includes two steps: (i) specify an upper bound of $\sup_{U\in\mathcal{U}}|\mathbb{E}[\mathscr{T}_{a}\tilde{h}\{S_{n,x}^{(1)*}\}]|$ via the formula  \eqref{eq:expression}, and (ii) derive the upper bound of $\zeta_{n}^{(1)}$ based on the smoothing inequality and the upper bound given in Step (i).

Based on the definition of $\psi_a(\cdot)$, we can reformulate it as follows:
\begin{align*}
\psi_a(x)&=-\int_a^\infty\int_{\mathbb{R}^p}\frac{\tilde{h}(y)}{(1-e^{-2s})^{p/2}}\phi\bigg(\frac{y-e^{-s}x}{\sqrt{1-e^{-2s}}}\bigg)\,{\rm d}y{\rm d}s\,,
\end{align*}
which implies that
\begin{align}
\partial_{j}\psi_{a}(x)&=\int_a^\infty\int_{\mathbb{R}^p}\frac{\tilde{h}(y)}{(1-e^{-2s})^{p/2}}\partial_j\phi\bigg(\frac{y-e^{-s}x}{\sqrt{1-e^{-2s}}}\bigg)\frac{e^{-s}}{\sqrt{1-e^{-2s}}}\,{\rm d}y{\rm d}s\notag\\
&=\int_{a}^{\infty}\frac{e^{-s}}{\sqrt{1-e^{-2s}}}\bigg\{\int_{\mathbb{R}^p}\tilde{h}(e^{-s}x+\sqrt{1-e^{-2s}}u)\partial_j\phi(u)\,{\rm d}u\bigg\}\,{\rm d}s\,,\label{eq:psi_first_partial}\\
\partial_{jk}\psi_{a}(x)&=-\int_a^\infty\int_{\mathbb{R}^p}\frac{\tilde{h}(y)}{(1-e^{-2s})^{p/2}}\partial_{jk}\phi\bigg(\frac{y-e^{-s}x}{\sqrt{1-e^{-2s}}}\bigg)
\bigg(\frac{e^{-s}}{\sqrt{1-e^{-2s}}}\bigg)^2\,{\rm d}y{\rm d}s\notag\\
&=-\int_{a}^{\infty}\bigg(\frac{e^{-s}}{\sqrt{1-e^{-2s}}}\bigg)^2\bigg\{\int_{\mathbb{R}^p}\tilde{h}(e^{-s}x+\sqrt{1-e^{-2s}}u)\partial_{jk}\phi(u)\,{\rm d}u\bigg\}\,{\rm d}s\,,\label{eq:psi_second_partial}\\
\partial_{jkl}\psi_{a}(x)&=\int_a^\infty\int_{\mathbb{R}^p}\frac{\tilde{h}(y)}{(1-e^{-2s})^{p/2}}\partial_{jkl}\phi\bigg(\frac{y-e^{-s}x}{\sqrt{1-e^{-2s}}}\bigg)
\bigg(\frac{e^{-s}}{\sqrt{1-e^{-2s}}}\bigg)^3\,{\rm d}y{\rm d}s\notag\\
&=\int_{a}^{\infty}\bigg(\frac{e^{-s}}{\sqrt{1-e^{-2s}}}\bigg)^3\bigg\{\int_{\mathbb{R}^p}\tilde{h}(e^{-s}x+\sqrt{1-e^{-2s}}u)\partial_{jkl}\phi(u)\,{\rm d}u\bigg\}\,{\rm d}s\,.\label{eq:psi_third_partial}
\end{align}
Recall $S_{n,x}^{(1)*}=L^{-1/2}\sum_{\ell=1}^L\tilde{X}_{\ell}^*$ with ${\rm Cov}\{S_{n,x}^{(1)*}\}={\rm I}_{p}$. Let $\{\tilde{X}_{\ell}^{*\dag}\}$ be an independent copy of $\{\tilde{X}_{\ell}^{*}\}$. Then $\mathbb{E}\{(L^{-1/2}\sum_{\ell=1}^{L}\tilde{X}_{\ell,j}^{*\dag})(L^{-1/2}\sum_{\ell=1}^{L}\tilde{X}_{\ell,k}^{*\dag})\}=I(j=k)$.
By Taylor expansion,
\begin{align*}
&\mathbb{E}[\Delta\psi_{a}\{S_{n,x}^{(1)*}\}]=\sum_{j,k=1}^{p}\mathbb{E}[\partial_{jk}\psi_{a}\{S_{n,x}^{(1)*}\}]\mathbb{E}\bigg\{\bigg(\frac{1}{\sqrt{L}}\sum_{\ell=1}^{L}\tilde{X}_{\ell,j}^{*\dag}\bigg)\bigg(\frac{1}{\sqrt{L}}\sum_{\ell=1}^{L}\tilde{X}_{\ell,k}^{*\dag}\bigg)\bigg\}\\
&~~=\frac{1}{L}\sum_{\ell=1}^{L}\sum_{j,k=1}^{p}\mathbb{E}[\tilde{X}_{\ell,j}^{*\dag}\tilde{X}_{\ell,k}^{*\dag}\partial_{jk}\psi_{a}\{S_{n,x}^{(1)*}\}]+\frac{1}{L}\sum_{\ell_{1}\neq \ell_{2}}\sum_{j,k=1}^{p}\mathbb{E}[\tilde{X}_{\ell_{1},j}^{*\dag}\tilde{X}_{\ell_2,k}^{*\dag}\partial_{jk}\psi_{a}\{S_{n,x}^{(1)*}\}]\\
&~~=\frac{1}{L}\sum_{\ell=1}^{L}\sum_{j,k=1}^{p}\mathbb{E}[\tilde{X}_{\ell,j}^{*\dag}\tilde{X}_{\ell,k}^{*\dag}\partial_{jk}\psi_{a}\{S_{n,x}^{(1)*}-L^{-1/2}\tilde{X}_{\ell}^*\}]\\
&~~~+\frac{1}{L^{3/2}}\sum_{\ell=1}^{L}\sum_{j,k,l=1}^{p}\mathbb{E}\bigg[\tilde{X}_{\ell,j}^{*\dag}\tilde{X}_{\ell,k}^{*\dag}\tilde{X}_{\ell,l}^{*}\int_{0}^{1}\partial_{jkl}\psi_{a}\{S_{n,x}^{(1)*}-(1-v)L^{-1/2}\tilde{X}_{\ell}^*\}\,{\rm d }v\bigg]\\
&~~~+\frac{1}{L}\sum_{\ell_{1}\neq \ell_{2}}\sum_{j,k=1}^{p}\mathbb{E}[\tilde{X}_{\ell_{1},j}^{*\dag}\tilde{X}_{\ell_2,k}^{*\dag}\partial_{jk}\psi_{a}\{S_{n,x}^{(1)*}\}]\,,\\
&\mathbb{E}[S_{n,x}^{(1)*}\cdot\nabla\psi_{a}\{S_{n,x}^{(1)*}\}]=\frac{1}{\sqrt{L}}\sum_{\ell=1}^{L}\sum_{j=1}^{p}\mathbb{E}[\tilde{X}_{\ell,j}^*\partial_{j}\psi_{a}\{S_{n,x}^{(1)*}\}]\\
&~~=\frac{1}{\sqrt{L}}\sum_{\ell=1}^{L}\sum_{j=1}^{p}\mathbb{E}[\tilde{X}_{\ell,j}^*\partial_{j}\psi_{a}\{S_{n,x}^{(1)*}-L^{-1/2}\tilde{X}_{\ell}^*\}]\\
&~~~+\frac{1}{L}\sum_{\ell=	1}^{L}\sum_{j,k=1}^{p}\mathbb{E}[\tilde{X}_{\ell,j}^*\tilde{X}_{\ell,k}^*\partial_{jk}\psi_{a}\{S_{n,x}^{(1)*}-L^{-1/2}\tilde{X}_{\ell}^*\}]\\
&~~~+\frac{1}{L^{3/2}}\sum_{\ell=1}^{L}\sum_{j,k,l=1}^{p}\mathbb{E}\bigg[\tilde{X}_{\ell,j}^*\tilde{X}_{\ell,k}^*\tilde{X}_{\ell,l}^*\int_{0}^{1}(1-v)\partial_{jkl}\psi_{a}\{S_{n,x}^{(1)*}-(1-v)L^{-1/2}\tilde{X}^*_{\ell}\}\,{\rm d}v\bigg]\,.
\end{align*}
Recall $\mathcal{F}_{-\ell}$ is the $\sigma$-filed generated by $\{\tilde{X}_s\}_{s\neq \ell}$ and $\tilde{X}_\ell^*=A\tilde{X}_{\ell}$. Note that 
\begin{align*}
\mathbb{E}[\tilde{X}_{\ell,j}^*\tilde{X}_{\ell,k}^*\partial_{jk}\psi_{a}\{S_{n,x}^{(1)*}-L^{-1/2}\tilde{X}_{\ell}^*\}]=&\,\,\mathbb{E}[\partial_{jk}\psi_{a}\{S_{n,x}^{(1)*}-L^{-1/2}\tilde{X}_{\ell}^*\}\mathbb{E}(\tilde{X}_{\ell,j}^*\tilde{X}_{\ell,k}^*\,|\,\mathcal{F}_{-\ell})]\,,\\ \mathbb{E}[\tilde{X}_{\ell,j}^{*\dag}\tilde{X}_{\ell,k}^{*\dag}\partial_{jk}\psi_{a}\{S_{n,x}^{(1)*}-L^{-1/2}\tilde{X}_{\ell}^*\}]=&\,\,\mathbb{E}[\partial_{jk}\psi_{a}\{S_{n,x}^{(1)*}-L^{-1/2}\tilde{X}_{\ell}^*\}]\mathbb{E}(\tilde{X}_{\ell,j}^{*\dag}\tilde{X}_{\ell,k}^{*\dag})\\
=&\,\,\mathbb{E}[\partial_{jk}\psi_{a}\{S_{n,x}^{(1)*}-L^{-1/2}\tilde{X}_{\ell}^*\}]\mathbb{E}(\tilde{X}_{\ell,j}^*\tilde{X}_{\ell,k}^*)\,.
\end{align*}
Then it holds that
\begin{align*}
&\frac{1}{L}\sum_{\ell=1}^{L}\sum_{j,k=1}^{p}(\mathbb{E}[\tilde{X}_{\ell,j}^{*\dag}\tilde{X}_{\ell,k}^{*\dag}\partial_{jk}\psi_{a}\{S_{n,x}^{(1)*}-L^{-1/2}\tilde{X}_{\ell}^*\}]-\mathbb{E}[\tilde{X}_{\ell,j}^*\tilde{X}_{\ell,k}^*\partial_{jk}\psi_{a}\{S_{n,x}^{(1)*}-L^{-1/2}\tilde{X}_{\ell}^*\}])\\
=&-\frac{1}{L}\sum_{\ell=1}^{L}\sum_{j,k=1}^{p}\mathbb{E}[\partial_{jk}\psi_{a}\{S_{n,x}^{(1)*}-L^{-1/2}\tilde{X}_{\ell}^*\}\mathbb{E}\{\tilde{X}_{\ell,j}^*\tilde{X}_{\ell,k}^*-\mathbb{E}(\tilde{X}_{\ell,j}^*\tilde{X}_{\ell,k}^*)\,|\,\mathcal{F}_{-\ell}\}]\,.
\end{align*}
By \eqref{eq:expression}, we have
\begin{align}\label{eq:five_terms}
\mathbb{E}[\mathscr{T}_{a}\tilde{h}\{S_{n,x}^{(1)*}\}]&=\underbrace{\frac{1}{L}\sum_{\ell=1}^{L}\sum_{j,k=1}^{p}\mathbb{E}[\partial_{jk}\psi_{a}\{S_{n,x}^{(1)*}-L^{-1/2}\tilde{X}_{\ell}^*\}\mathbb{E}\{\tilde{X}_{\ell,j}^*\tilde{X}_{\ell,k}^*-\mathbb{E}(\tilde{X}_{\ell,j}^*\tilde{X}_{\ell,k}^*)\,|\,\mathcal{F}_{-\ell}\}]}_{I_{1}}\notag\\
&~-\underbrace{\frac{1}{\sqrt{L}}\sum_{\ell=1}^{L}\sum_{j=1}^{p}\mathbb{E}[\tilde{X}_{\ell,j}^*\partial_{j}\psi_{a}\{S_{n,\bx}^{(1)*}-L^{-1/2}\tilde{X}_{\ell}^*\}]}_{I_{2}}\notag\\
&~+\underbrace{\frac{1}{L}\sum_{\ell_{1}\neq \ell_{2}}\sum_{j,k=1}^{p}\mathbb{E}[\tilde{X}_{\ell_{1},j}^{*\dag}\tilde{X}_{\ell_{2},k}^{*\dag}\partial_{jk}\psi_{a}\{S_{n,x}^{(1)*}\}]}_{I_{3}}\notag\\
&~+\frac{1}{L^{3/2}}\sum_{\ell=1}^{L}\sum_{j,k,l=1}^{p}\mathbb{E}\bigg[\tilde{X}_{\ell,j}^{*\dag}\tilde{X}_{\ell,k}^{*\dag}\tilde{X}_{\ell,l}^*\\
&~~\underbrace{~~~~~~\times\int_{0}^{1}\partial_{jkl}\psi_{a}\{S_{n,x}^{(1)*}-(1-v)L^{-1/2}\tilde{X}_{\ell}^*\}\,{\rm d }v\bigg]}_{I_{4}}\notag\\
&~-\frac{1}{L^{3/2}}\sum_{\ell=1}^{L}\sum_{j,k,l=1}^{p}\mathbb{E}\bigg[\tilde{X}_{\ell,j}^*\tilde{X}_{\ell,k}^*\tilde{X}_{\ell,l}^*\notag\\
&~~\underbrace{~~~~~~\times\int_{0}^{1}(1-v)\partial_{jkl}\psi_{a}\{S_{n,x}^{(1)*}-(1-v)L^{-1/2}\tilde{X}_{\ell}^*\}\,{\rm d}v\bigg]}_{I_{5}}\notag\,.
\end{align}

Our first step is to bound $I_1, I_2$, $I_3$, $I_4$ and $I_5$ in \eqref{eq:five_terms}, respectively.
Notice that $|\tilde{h}(\cdot)|=|h(\cdot)-\int_{\mathbb{R}^p} h(u)\phi(u)\,{\rm d}u|\leq 1$. For any $j,k\in [p]$,
by \eqref{eq:psi_second_partial}, we have
$
|\partial_{jk}\psi_{a}(x)|\leq\int_{a}^{\infty}e^{-2s}(1-e^{-2s})^{-1}\,{\rm d}s\cdot\int_{\mathbb{R}^{p}}|\partial_{jk}\phi(u)|\,{\rm d}u$ for any $x\in\mathbb{R}^p$.
Since $\partial_{j}\phi(u)=-u_{j}\phi(u)$, we know
$
\partial_{jk}\phi(u)=
\phi(u)\{(u_{j}^{2}-1)I(j= k)+u_{j}u_{k}I(j\neq k)\}$,
which yields $$\max_{j,k\in[p]}\int_{\mathbb{R}^{p}}|\partial_{jk}\phi(u)|{\rm d}u\leq C\,.$$ Moreover,
$
\int_{a}^{\infty}e^{-2s}(1-e^{-2s})^{-1}\,{\rm d}s\leq e^{-a}(1-e^{-2a})^{-1/2}\int_{0}^{\infty}e^{-s}(1-e^{-2s})^{-1/2}\,{\rm d}s\lesssim a^{-1/2}$,
where the last step is based on the facts $e^{-a}(1-e^{-2a})^{-1/2}\leq (2a)^{-1/2}$ for $a>0$ and $\int_{0}^{\infty}e^{-s}(1-e^{-2s})^{-1/2}\,{\rm d}s=\pi/2$. Hence,
\begin{align}\label{eq:second_partial_bound}
\sup_{x\in\mathbb{R}^p}\max_{j,k\in [p]}|\partial_{jk}\psi_{a}(x)|\lesssim a^{-1/2}\,.
\end{align}
By H\"older's inequality and Lemma \ref{la:xtildestar_mom}, we have $\max_{\ell\in[L]}\max_{j,k\in[p]}|\mathbb{E}(\tilde{X}_{\ell,j}^*\tilde{X}_{\ell,k}^*)|\lesssim B_{n}^{2}p^{5/3}$. Lemma \ref{la:xtildestar_tail} yields
\begin{align*}
    &\max_{\ell\in[L]}\max_{j,k\in[p]}\mathbb{P}\{|\tilde{X}_{\ell,j}^*\tilde{X}_{\ell,k}^*-\mathbb{E}(\tilde{X}_{\ell,j}^*\tilde{X}_{\ell,k}^*)|>\lambda\}\\
    &~~~~~\lesssim p\exp(-Cp^{-1}\lambda B_{n}^{-2})+p\exp(-Cb^{\gamma/2}p^{-\gamma/2}\lambda^{\gamma/2}B_{n}^{-\gamma})
    \end{align*}
for any $\lambda\geq C'B_{n}^{2}p^{5/3}$, where $C'>0$ is a sufficiently large constant. Letting $D_{4n}=C'B_{n}^{2}(pn)^{5/3}$, using the same technique for deriving \eqref{eq:bound12}, it holds that \begin{align*}
&\max_{\ell\in[L]}\max_{j,k\in[p]}\mathbb{E}[|\mathbb{E}\{\tilde{X}_{\ell,j}^{*}\tilde{X}_{\ell,k}^{*}-\mathbb{E}(\tilde{X}_{\ell,j}^{*}\tilde{X}_{\ell,k}^{*})\,|\,\mathcal{F}_{-\ell}\}|]\\
&~~~\lesssim B_n^2 p^2\{\exp(-Cp^{-1}D_{4n}B_{n}^{-2})+\exp(-Cb^{\gamma/2}p^{-\gamma/2}D_{4n}^{\gamma/2}B_{n}^{-\gamma})\}+D_{4n}\alpha_{n}(h)\,.
\end{align*}
Select $h=C''\{\log(pn)\}^{1/\gamma_{2}}$ for a sufficiently large constant $C''>0$. Notice that $B_n\leq n$. By \eqref{eq:second_partial_bound},
\begin{align*}
|I_{1}|&\leq a^{-1/2}p^{2}\max_{\ell\in[L]}\max_{j,k\in [p]}\mathbb{E}[|\mathbb{E}\{\tilde{X}_{\ell,j}^*\tilde{X}_{\ell,k}^*-\mathbb{E}(\tilde{X}_{\ell,j}^*\tilde{X}_{\ell,k}^*)\,|\,\mathcal{F}_{-\ell}\}|]\lesssim a^{-1/2}(pn)^{-C'''}
\end{align*}
for some sufficiently large constant $C'''>0$. By \eqref{eq:psi_first_partial}, we have
$
\max_{j\in [p]}|\partial_{j}\psi_{a}(x)|\leq C$.
Analogously, with selecting $D_{5n}=C'B_{n}pn$ for some sufficiently large constant $C'>0$, by Lemma \ref{la:xtildestar_tail},
\begin{align*}
&\max_{\ell\in [L]}\max_{j\in [p]}|\mathbb{E}[\tilde{X}_{\ell,j}^*\partial_{j}\psi_{a}\{S_{n,x}^{(1)*}-L^{-1/2}\tilde{X}^*_{\ell}\}]|
\lesssim \max_{\ell\in [L]}\max_{j\in [p]}|\mathbb{E}\{|\mathbb{E}(\tilde{X}_{\ell,j}^*\,|\,\mathcal{F}_{-\ell})|\}\\
&~~~\lesssim  B_{n}p^{3/2}\{\exp(-Cp^{-1}D_{5n}^{2}B_{n}^{-2})+\exp(-Cb^{\gamma/2}p^{-\gamma/2}D_{5n}^{\gamma}B_{n}^{-\gamma})\}+D_{5n}\alpha_{n}(h)\,.
\end{align*}
Due to $h=C''\{\log(pn)\}^{1/\gamma_{2}}$ for a sufficiently large constant $C''>0$, we have
\begin{align*}
|I_{2}|\leq \sqrt{L}p\max_{\ell\in[L]}\max_{j\in[p]}|\mathbb{E}[\tilde{X}_{\ell,j}^*\partial_{j}\psi_{a}\{S_{n,x}^{(1)*}-L^{-1/2}\tilde{X}_{\ell}^*\}]|\lesssim (pn)^{-C'''}
\end{align*}
for some sufficiently large constant $C'''>0$.  By Davydov's inequality and Lemma \ref{la:xtildestar_mom}, we have $
|L^{-1}\sum_{\ell=1}^{L-1}\mathbb{E}(\tilde{X}_{\ell+1,j}^{*\dag}\tilde{X}_{\ell,k}^{*\dag})|=|L^{-1}\sum_{\ell=1}^{L-1}\mathbb{E}(\tilde{X}_{\ell+1,j}^{*}\tilde{X}_{\ell,k}^{*})|\lesssim B_{n}^{2}p^{5/3}\exp(-Ch^{\gamma_{2}})$ and
$
|\mathbb{E}(\tilde{X}_{\ell_{1},j}^{*\dag}\tilde{X}_{\ell_{2},k}^{*\dag})|=|\mathbb{E}(\tilde{X}_{\ell_{1},j}^{*}\tilde{X}_{\ell_{2},k}^{*})|
\lesssim B_{n}^{2}p^{5/3}\exp\{-C(|\ell_{1}-\ell_{2}|-1)^{\gamma_{2}}\}\exp(-Cb^{\gamma_{2}})$ for any $(\ell_1,\ell_2)$ such that $\ell_{1}> \ell_{2}+1$, which yields $\max_{j,k\in [p]}|L^{-1}\sum_{\ell_{1}> \ell_{2}}\mathbb{E}(\tilde{X}_{\ell_{1},j}^{*\dag}\tilde{X}_{\ell_{2},k}^{*\dag})|\lesssim B_{n}^{2}p^{5/3}\exp(-Ch^{\gamma_{2}})$. Analogously, we also have $\max_{j,k\in [p]}|L^{-1}\sum_{\ell_{1}< \ell_{2}}\mathbb{E}(\tilde{X}_{\ell_{1},j}^{*\dag}\tilde{X}_{\ell_{2},k}^{*\dag})|\lesssim B_{n}^{2}p^{5/3}\exp(-Ch^{\gamma_{2}})$. Hence, 
$$
\max_{j,k\in [p]}\bigg|\frac{1}{L}\sum_{\ell_{1}\neq \ell_{2}}\mathbb{E}(\tilde{X}_{\ell_{1},j}^{*\dag}\tilde{X}_{\ell_{2},k}^{*\dag})\bigg|\lesssim B_{n}^{2}p^{5/3}\exp(-Ch^{\gamma_{2}})\,.$$
Together with \eqref{eq:second_partial_bound}, we have
\begin{align*}
|I_{3}|
&\lesssim \sum_{j,k=1}^{p}\mathbb{E}[|\partial_{jk}\psi_{a}\{S_{n,x}^{(1)*}\}|]\cdot\max_{j,k\in[p]}\bigg|\frac{1}{L}\sum_{\ell_{1}\neq \ell_{2}}\mathbb{E}(\tilde{X}_{\ell_{1},j}^{*\dag}\tilde{X}_{\ell_{2},k}^{*\dag})\bigg|\lesssim a^{-1/2}(pn)^{-C'''}
\end{align*}
for some sufficiently large $C'''>0$. Notice that $\max_{j,k,l\in [p]}\int_{\mathbb{R}^{p}}|\partial_{jkl}\phi(u)|{\rm d}u\leq C$. By \eqref{eq:psi_third_partial}, $
\sup_{x\in\mathbb{R}^p}\max_{j,k,l\in [p]}|\partial_{jkl}\psi_{a}(x)|\lesssim\int^{\infty}_{a}e^{-3s}(1-e^{-2s})^{-3}\,{\rm d}s\lesssim \int_{a}^{\infty}s^{-3/2}\,{\rm d}s\lesssim a^{-1/2}$. By H\"older's inequality and Lemma \ref{la:xtildestar_mom}, $$\max_{\ell\in[L]}\max_{j,k,l\in [p]}\mathbb{E}(|\tilde{X}_{\ell,j}^{*\dag}\tilde{X}_{\ell,k}^{*\dag}\tilde{X}_{\ell,l}^*|)\leq \max_{\ell\in[L]}\max_{j\in [p]}\mathbb{E}(|\tilde{X}_{\ell,j}^*|^3)\lesssim B_{n}^{3}p^{5/2}\,,$$ which implies
\begin{align*}
|I_{4}|\lesssim\frac{1}{L^{3/2}a^{1/2}}\sum_{\ell=1}^{L}\sum_{j,k,l=1}^{p}\mathbb{E}(|\tilde{X}_{\ell,j}^{*\dag}\tilde{X}_{\ell,k}^{*\dag}\tilde{X}_{\ell,l}^*|)\lesssim a^{-1/2}L^{-1/2}B_{n}^3p^{11/2}\,.
\end{align*}
Analogously, we also have
$
|I_{5}|\lesssim a^{-1/2}L^{-1/2}B_{n}^3p^{11/2}$.
Therefore, by \eqref{eq:five_terms}, it holds that
\begin{align}\label{eq:smoothing_bound}
|\mathbb{E}[\mathscr{T}_{a}\tilde{h}\{S_{n,x}^{(1)*}\}]|\lesssim a^{-1/2}L^{-1/2}B_{n}^3p^{11/2}+(pn)^{-C'''}
\end{align}
for some sufficiently large constant $C'''>0$.

Our second step is to bound $\zeta_{n}^{(1)}$ by the smoothing inequality. Let $\mathcal{H}=\{I_U(\cdot):U\in\mathcal{U}\}$ with $I_U(\cdot)=I(\cdot\in U)$. For some constant $c>0$ and any distribution function $G(\cdot)$, we write $G_{c}(\cdot)$ as the distribution function of $cW$ with $W\sim G(\cdot)$. Denote by $F(\cdot)$ the distribution function of $S_{n,x}^{(1)*}$. Since $S_{n,x}^{(1)*}$ is independent of $Z$, the distribution function of $e^{-a}S_{n,x}^{(1)*}+\sqrt{1-e^{-2a}}Z$ is $F_{e^{-a}}\star\Phi_{\sqrt{1-e^{-2a}}}(\cdot)$, where $\star$ is the convolution operator. Since $\Phi_{e^{-a}}\star\Phi_{\sqrt{1-e^{-2a}}}(\cdot)=\Phi(\cdot)$ and $\tilde{h}(\cdot)=h(\cdot)-\int h\,{\rm d}\Phi$, we have
\begin{align}\label{eq:convolution}
\mathbb{E}[\mathscr{T}_{a}\tilde{h}\{S_{n,x}^{(1)*}\}]&=\mathbb{E}[\tilde{h}\{e^{-a}S_{n,x}^{(1)*}+\sqrt{1-e^{-2a}}Z\}]=\int_{\mathbb{R}^p}\tilde{h}\,{\rm d}F_{e^{-a}}\star\Phi_{\sqrt{1-e^{-2a}}}\notag\\
&=\int_{\mathbb{R}^p}\tilde{h}\,{\rm d}\{F_{e^{-a}}-\Phi_{e^{-a}}\}\star\Phi_{\sqrt{1-e^{-2a}}}+\int_{\mathbb{R}^p}\tilde{h}\,{\rm d}\Phi\\
&=\int_{\mathbb{R}^p}\tilde{h}\,{\rm d}\{F_{e^{-a}}-\Phi_{e^{-a}}\}\star\Phi_{\sqrt{1-e^{-2a}}}=\int_{\mathbb{R}^p}h\,{\rm d}\{F_{e^{-a}}-\Phi_{e^{-a}}\}\star\Phi_{\sqrt{1-e^{-2a}}}\,.\notag
\end{align}
On the other hand, since $\mathcal{U}$ is invariant under  affine transformations, it holds that
\begin{align}\label{eq:rhobestar}
\zeta_{n}^{(1)}&=\sup_{h\in\mathcal{H}}|\mathbb{E}[\tilde{h}\{S_{n,x}^{(1)*}\}]|=\sup_{h\in\mathcal{H}}\bigg|\int_{\mathbb{R}^p}h\,{\rm d}(F-\Phi)\bigg|=\sup_{h\in\mathcal{H}}\bigg|\int_{\mathbb{R}^p}h\,{\rm d}\{F_{e^{-a}}-\Phi_{e^{-a}}\}\bigg|\,.
\end{align}
Notice that $\{F_{e^{-a}}-\Phi_{e^{-a}}\}\star\Phi_{\sqrt{1-e^{-2a}}}$ is a perturbed of $F_{e^{-a}}-\Phi_{e^{-a}}$ by convolution with $\Phi_{\sqrt{1-e^{-2a}}}$. Let $\epsilon>0$ satisfy
$
\int_{|u|_2<\epsilon}{\rm d}\Phi_{\sqrt{1-e^{-2a}}}(u)=7/8$.
Then $\epsilon\lesssim (pa)^{1/2}$. Write $M_{h}(x;\epsilon)=\sup\{h(y):y\in\mathbb{R}^p, |x-y|_{2}\leq\epsilon\}$, $m_{h}(x;\epsilon)=\inf\{h(y):y\in\mathbb{R}^p, |x-y|_{2}\leq\epsilon\}$ and $w_{h}(x;\epsilon)=M_{h}(x;\epsilon)-m_{h}(x;\epsilon)$ for any $h\in\mathcal{H}$. Given $h\in\mathcal{H}$, there exists $U\in\mathcal{U}$ such that $h(\cdot)=I(\cdot\in U)$. Applying the smoothing inequality \citep[Corollary 11.5]{BhattacharyaRao_2010} with $\mu=F_{e^{-a}}$, $\nu=\Phi_{e^{-a}}$, $K_{\epsilon}=\Phi_{\sqrt{1-e^{-2a}}}$, $f(\cdot)=h(\cdot)$, $\alpha=7/8$ and $\epsilon$ specified above, we have
\begin{align}\label{eq:smoothing}
\bigg|\int_{\mathbb{R}^p}h\,{\rm d}\{F_{e^{-a}}-\Phi_{e^{-a}}\}\bigg|\lesssim \gamma^*(h;\epsilon)+w_{h}^*(2\epsilon;\Phi_{e^{-a}})\,,
\end{align}
where $\gamma^*(h;\epsilon)=\sup_{y\in\mathbb{R}^p}\gamma(h_{y};\epsilon)$ and $w_{h}^*(2\epsilon;\Phi_{e^{-a}})=\sup_{y\in\mathbb{R}^p}\bar{w}_{h_{y}}(2\epsilon;\Phi_{e^{-a}})$ with $h_{y}(x)=h(x+y)$, $\bar{w}_{h_y}(\epsilon;\Phi_{e^{-a}})=\int w_{h_y}(x;\epsilon)\,{\rm d}\Phi_{e^{-a}}(x)$ and
\begin{align*}
&\gamma(h_y;\epsilon)=\max\bigg[\int_{\mathbb{R}^p} M_{h_y}(x;\epsilon)\,{\rm d}\{F_{e^{-a}}-\Phi_{e^{-a}}\}\star\Phi_{\sqrt{1-e^{-2a}}}(x)\,,\\
&~~~~~~~~~~~~~-\int_{\mathbb{R}^p} m_{h_y}(x;\epsilon)\,{\rm d}\{F_{e^{-a}}-\Phi_{e^{-a}}\}\star\Phi_{\sqrt{1-e^{-2a}}}(x)\bigg]\,.
\end{align*}
Next, we bound $\gamma^*(h;\epsilon)$ and $w_{h}^*(2\epsilon;\Phi_{e^{-a}})$, respectively. For any $\delta>0$ and $U_*\in\mathcal{U}$, we write $U_*^{\delta}=\{w\in\mathbb{R}^{p}:\rho(w,U_*)\leq \delta\}$ with $\rho(w,U_*)=\inf_{v\in U_*}|v-w|_{2}$, and $U_*^{-\delta}=\{w\in\mathbb{R}^{p}:B(w,\delta)\subset U_*\}$ with $B(w,\delta)=\{v\in\mathbb{R}^{p}:|v-w|_{2}\leq\delta\}$.
Due to $M_{h_{y}}(x;\epsilon)=I\{x\in (U-y)^{\epsilon}\}\in\mathcal{H}$ and $m_{h_{y}}(x;\epsilon)=I\{x\in (U-y)^{-\epsilon}\}\in\mathcal{H}$, by \eqref{eq:convolution}, we have
\begin{align*}
\bigg|\int_{\mathbb{R}^p} M_{h_{y}}(x;\epsilon)\,{\rm d}\{F_{e^{-a}}-\Phi_{e^{-a}}\}\star\Phi_{\sqrt{1-e^{-2a}}}(x)\bigg|\leq&\,\sup_{h\in\mathcal{H}}|\mathbb{E}[\mathscr{T}_{a}\tilde{h}\{S_{n,x}^{(1)*}\}]|\,,\\
\bigg|\int_{\mathbb{R}^p} m_{h_{y}}(x;\epsilon)\,{\rm d}\{F_{e^{-a}}-\Phi_{e^{-a}}\}\star\Phi_{\sqrt{1-e^{-2a}}}(x)\bigg|\leq&\,\sup_{h\in\mathcal{H}}|\mathbb{E}[\mathscr{T}_{a}\tilde{h}\{S_{n,x}^{(1)*}\}]|\,,
\end{align*}
which implies
$
|\gamma(h_{y};\epsilon)|\leq \sup_{h\in\mathcal{H}}|\mathbb{E}[\mathscr{T}_{a}\tilde{h}\{S_{n,x}^{(1)*}\}]$
for any $y\in\mathbb{R}^p$. Therefore,
\begin{align}\label{eq:gamma_star}
|\gamma^*(h;\epsilon)|\leq\sup_{y\in\mathbb{R}^p}|\gamma(h_{y};\epsilon)|\leq \sup_{h\in\mathcal{H}}|\mathbb{E}[\mathscr{T}_{a}\tilde{h}\{S_{n,x}^{(1)*}\}]|\,.
\end{align}
Let $\tilde{U}=e^a(U-y)$. Recall $Z\sim N(0,{\rm I}_p)$. By Corollary 3.2 of \cite{BhattacharyaRao_2010}, 
\begin{align*}
\bar{w}_{h_{y}}(2\epsilon;\Phi_{e^{-a}})&=\int_{\mathbb{R}^p} I\{x\in(U-y)^{2\epsilon}\backslash( U-y)^{-2\epsilon}\}\,{\rm d}\Phi_{e^{-a}}(x)\\
&=\mathbb{P}\{e^{-a}Z\in(U-y)^{2\epsilon}\backslash (U-y)^{-2\epsilon}\}=\mathbb{P}(Z\in \tilde{U}^{2e^a\epsilon}\backslash \tilde{U}^{-2e^a\epsilon})\\
&\lesssim \frac{\Gamma\{(p-1)/2\}}{\Gamma(p/2)}pe^a\epsilon\lesssim p^{1/2}e^a\epsilon
\end{align*}
holds uniformly over $y\in\mathbb{R}^p$, where $\Gamma(\cdot)$ is the gamma function. Thus, we have
$
w_{h}^*(2\epsilon;\Phi_{e^{-a}})=\sup_{y\in\mathbb{R}^p}\bar{w}_{h_{y}}(2\epsilon;\Phi_{e^{-a}})\lesssim p^{1/2}e^a\epsilon$.
Since $\epsilon\lesssim (pa)^{1/2}$, we have
$
w_{h}^*(2\epsilon;\Phi_{e^{-a}})\lesssim pa^{1/2}e^a$.
Together with \eqref{eq:smoothing_bound}, \eqref{eq:rhobestar}, \eqref{eq:smoothing} and \eqref{eq:gamma_star}, we have
\begin{align}\label{eq:smoothinginequalitybound}
\zeta_{n}^{(1)}&\lesssim \sup_{h\in\mathcal{H}}|\mathbb{E}[\mathscr{T}_{a}\tilde{h}\{S_{n,x}^{(1)*}\}]|+pa^{1/2}e^{a}\\
&\lesssim a^{-1/2}L^{-1/2}B_{n}^3p^{11/2}+(pn)^{-C'''}+pa^{1/2}e^{a}\notag
\end{align}
for any $a> 0$. Considering $a\leq 1$, we have
$
\zeta_{n}^{(1)}\lesssim a^{-1/2}L^{-1/2}B_{n}^3p^{11/2}+(pn)^{-C'''}+pa^{1/2}$.
Select $a_*=L^{-1/2}B_{n}^3p^{9/2}$. If $a_*\leq 1$, we have $
\zeta_{n}^{(1)}\lesssim L^{-1/4}B_{n}^{3/2}p^{13/4}$. If $a_*>1$, such upper bound holds automatically.
We complete the proof of Lemma \ref{la:BE_bound_Block}. $\hfill\Box$

\subsection{Proof of Lemma \ref{la:BE_bound_non_Block}}\label{sec:f.5}
Recall $
S_{n,x}=S_{n,x}^{(1)}+\delta_n$ with $
\delta_n=n^{-1/2}\sum_{\ell=1}^{L+1}\sum_{t\in\mathcal{J}_\ell}{X}_t+\{n^{-1/2}-(Lb)^{-1/2}\}\sum_{\ell=1}^L\sum_{t\in\mathcal{I}_\ell}{X}_t$. For some $D_{n}> 0$, consider the event $\mathcal{E}:=\{|\delta_{n}|_{\infty}\leq D_{n}\}$. Restricted on $\mathcal{E}$, we have
$
|S_{n,x}-S_{n,x}^{(1)}|_{2}=|\delta_{n}|_{2}\leq p^{1/2}D_{n}=:R$.
For any convex set $U\in\mathcal{U}$ and $\epsilon>0$, denote $U^{\epsilon}=\{w\in\mathbb{R}^{p}:\rho(w,U)\leq \epsilon\}$ with $\rho(w,U)=\inf_{v\in U}|v-w|_{2}$, and $U^{-\epsilon}=\{w\in\mathbb{R}^{p}:B(w,\epsilon)\subset U\}$ with $B(w,\epsilon)=\{v\in\mathbb{R}^{p}:|v-w|_{2}\leq\epsilon\}$. Recall $G^{(1)}\sim N\{0,\Xi^{(1)}\}$ with $\Xi^{(1)}=\cov\{S_{n,x}^{(1)}\}$.
Then we have
\begin{align*}
&\mathbb{P}(S_{n,x}\in U)-\mathbb{P}\{G^{(1)}\in U\}=\mathbb{P}(S_{n,x}\in U,\mathcal{E})+\mathbb{P}(S_{n,x}\in U,\mathcal{E}^{c})-\mathbb{P}\{G^{(1)}\in U\}\\
&~~~~\leq\mathbb{P}\{S_{n,x}^{(1)}\in U^{R}\}+\mathbb{P}(\mathcal{E}^c)-[\mathbb{P}\{G^{(1)}\in U^{R}\}-\mathbb{P}\{G^{(1)}\in U^{R}\backslash U\}]\\
&~~~~\leq\zeta_{n}^{(1)}+\mathbb{P}(\mathcal{E}^c)+\mathbb{P}\{G^{(1)}\in U^{R}\backslash U\}\,.
\end{align*}
Due to $\{(U^c)^{R}\}^c=U^{-R}$, we have
\begin{align*}
&\mathbb{P}(S_{n,x}\in U)-\mathbb{P}\{G^{(1)}\in U\}\geq\mathbb{P}(\mathcal{E})- \mathbb{P}(S_{n,x}\in U^c,\mathcal{E})-\mathbb{P}\{G^{(1)}\in U\}\\
&~~~~\geq\mathbb{P}(\mathcal{E})- \mathbb{P}\{S_{n,x}^{(1)}\in (U^c)^{R}\}-\mathbb{P}\{G^{(1)}\in U\}\\
&~~~~=-\mathbb{P}(\mathcal{E}^c)+ \mathbb{P}\{S_{n,x}^{(1)}\in U^{-R}\}-[\mathbb{P}\{G^{(1)}\in U^{-R}\}+\mathbb{P}\{G^{(1)}\in U\backslash U^{-R}\}]\\
&~~~~\geq-\zeta_{n}^{(1)}-\mathbb{P}(\mathcal{E}^c)-\mathbb{P}\{G^{(1)}\in U\backslash U^{-R}\}\,.
\end{align*}
Therefore, it holds that
\begin{align}\label{eq:substituteblock}
&|\mathbb{P}(S_{n,x}\in U)-\mathbb{P}\{G^{(1)}\in U\}|\leq \zeta_{n}^{(1)}+\mathbb{P}(\mathcal{E}^c)+\mathbb{P}\{G^{(1)}\in U^{R}\backslash U\}+\mathbb{P}\{G^{(1)}\in U\backslash U^{-R}\}\,.
\end{align}
By Lemma A.2 of \cite{CCK_2017}, we have
$
\sup_{U\in\mathcal{U}}\mathbb{P}\{G^{(1)}\in U^{R}\backslash U\}\lesssim R\|\{\Xi^{(1)}\}^{-1}\|_{\rm HS}^{1/2}$, where $\|\cdot\|_{\rm HS}$ is the Hilbert-Schmidt norm. Let $\lambda_{1}\geq\lambda_{2}\geq\cdots\geq\lambda_{p}>0$ be the eigenvalues of $\Xi^{(1)}$. By \eqref{eq:eigen_min}, 
$
\|\{\Xi^{(1)}\}^{-1}\|_{\rm HS}=(\sum_{i=1}^{p}\lambda_{i}^{-2})^{1/2}\leq p^{1/2}\lambda_{p}^{-1}\leq Cp^{1/2}$,
which implies that $
\mathbb{P}\{G^{(1)}\in U^{R}\backslash U\}\lesssim Rp^{1/4}=p^{3/4}D_{n}$.
Analogously, we also have $\mathbb{P}\{G^{(1)}\in U\backslash U^{-R}\}\lesssim p^{3/4}D_{n}$. Together with Lemma \ref{la:BE_bound_Block}, \eqref{eq:substituteblock} yields
$$
\zeta_{n}^{(2)}\lesssim B_{n}^{3/2}p^{13/4}L^{-1/4}+\mathbb{P}(\mathcal{E}^c)+p^{3/4}D_{n}$$ 
provided that $B_{n}^{2}p\{\log (pn)\}^{1/\gamma_{2}}\ll b\ll B_{n}^{-2}p^{-1}n$. As shown in the proof of Lemma \ref{pn:2}, if $b=o(n^{1/2})$, then 
$$
\mathbb{P}(\mathcal{E}^{c})\lesssim p\exp(-CD_n^{\gamma}n^{\gamma/2}B_{n}^{-\gamma})+p\exp(-CD_n^2bh^{-1}B_{n}^{-2})\,.$$ 
Select $D_{n}= C'B_{n}^{3/2}p^{5/2}L^{-1/4}$ for some sufficiently large constant $C'>0$. Then $p^{3/4}D_{n}\asymp B_{n}^{3/2}p^{13/4}L^{-1/4}$. Since $L\asymp nb^{-1}$ and $D_{n}\sqrt{n}\asymp B_{n}^{3/2}p^{5/2}n^{1/4}b^{1/4}$, it holds that $$p\exp(-CD_{n}^{\gamma}n^{\gamma/2}B_{n}^{-\gamma})\lesssim B_{n}^{3/2}p^{13/4}L^{-1/4}\,.$$ Due to $h\asymp \{\log(pn)\}^{1/\gamma_2}$, if $b\gg B_{n}^{-2/3}n^{1/3}p^{-10/3}\{\log(pn) \}^{2(\gamma_{2}+1)/(3\gamma_{2})}$, it then holds that $p\exp(-CD_{n}^2 bh^{-1}B_{n}^{-2})\lesssim B_{n}^{3/2}p^{13/4}L^{-1/4}$. Hence, $\zeta_{n}^{(2)}\lesssim B_{n}^{3/2}p^{13/4}L^{-1/4}$.  $\hfill\Box$

\subsection{Proof of Lemma \ref{la:comparison_convex}} \label{sec:pflem15}

Without loss of generality, we assume $p^{2}|\Gamma_{1}-\Gamma_{2}|_{\infty}< 1$, since otherwise the assertions are trivial.
Note that
$
\sup_{U\in\mathcal{U}}|\mathbb{P}(X\in U)-\mathbb{P}(Y\in U)|
\leq \sup_{U\in\mathcal{U}}|\mathbb{P}(\Gamma_{2}^{-1/2}X\in U)-\mathbb{P}(Z\in U)|$ with a $p$-dimensional standard normal random vector $Z$. We assume $Z$ is independent of $X$. Let $\mathring{X}=\Gamma_{2}^{-1/2}X=(\mathring{X}_{1},\ldots,\mathring{X}_{p})^{\T}$. For $\mathscr{T}_a\tilde{h}(\cdot)$ and $\psi_a(\cdot)$ defined in the proof of Lemma \ref{la:BE_bound_Block}, identical to \eqref{eq:expression}, it holds that
$$
\mathbb{E}\{\mathscr{T}_{a}\tilde{h}(\mathring{X})\}=\sum_{j=1}^{p}\mathbb{E}\{\partial_{jj}\psi_{a}(\mathring{X})\}-\sum_{j=1}^{p}\mathbb{E}\{\mathring{X}_{j}\partial_{j}\psi_{a}(\mathring{X})\}\,.$$ 
By Stein's identity \citep[Lemma 2]{CCK_2015}, we have $\mathbb{E}\{\mathring{X}_{j}\partial_{j}\psi_{a}(\mathring{X})\}\\=\sum_{k=1}^{p}\mathbb{E}(\mathring{X}_{j}\mathring{X}_{k})\mathbb{E}\{\partial_{jk}\psi_{a}(\mathring{X})\}$, which implies
\begin{align}\label{eq:stein_identity}
\mathbb{E}\{\mathscr{T}_{a}\tilde{h}(\mathring{X})\}
&=\sum_{j,k=1}^{p}\{I(j=k)-\mathbb{E}(\mathring{X}_{j}\mathring{X}_{k})\}\mathbb{E}\{\partial_{jk}\psi_{a}(\mathring{X})\}\,.
\end{align}
Write $\Gamma_{2}^{-1/2}=(\gamma_1,\ldots,\gamma_p)$. Due to $\cov(\mathring{X})=\Gamma_{2}^{-1/2}\Gamma_{1}\Gamma_{2}^{-1/2}$, we have
$
|I(j=k)-\mathbb{E}(\mathring{X}_{j}\mathring{X}_{k})|=|\gamma_j^{\T}(\Gamma_{2}-\Gamma_{1})\gamma_k|\leq p|\Gamma_{2}-\Gamma_{1}|_{\infty}|\gamma_j|_{2}|\gamma_k|_{2}\lesssim p|\Gamma_{2}-\Gamma_{1}|_{\infty}$
holds uniformly over $j,k\in [p]$, where the last step is based on the same arguments for the last step in \eqref{eq:tuta_x}. Together with \eqref{eq:stein_identity} and \eqref{eq:second_partial_bound}, we have
$
|\mathbb{E}\{\mathscr{T}_{a}\tilde{h}(\mathring{X})\}|\lesssim a^{-1/2}p^{3}|\Gamma_{1}-\Gamma_{2}|_{\infty}$.
For $\mathcal{H}$ specified in the proof of Lemma \ref{la:BE_bound_Block}, identical to the first step of  \eqref{eq:smoothinginequalitybound}, we have
\begin{align*}
\sup_{U\in\mathcal{U}}|\mathbb{P}(\mathring{X}\in U)-\mathbb{P}(Z\in U)|\lesssim&\, \sup_{h\in\mathcal{H}}|\mathbb{E}\{\mathscr{T}_{a}\tilde{h}(\mathring{X})\}|+pa^{1/2}e^{a}\\
\lesssim&\, a^{-1/2}p^{3}|\Gamma_{1}-\Gamma_{2}|_{\infty}+p a^{1/2}
\end{align*}
for any $a\in(0,1)$. We then obtain Lemma \ref{la:comparison_convex} by selecting $a=p^{2}|\Gamma_{1}-\Gamma_{2}|_{\infty}$. $\hfill\Box$

\section{Proof of Theorem \ref{pn:GA_dependency_graph_s_sparsely}}\label{sec:gas}
Recall $\mathcal{U}$ is the class of all convex sets in $\mathbb{R}^{p}$ and
$
\zeta_{n}:=\sup_{U\in\mathcal{U}}|\mathbb{P}(S_{n,x}\in U)-\mathbb{P}(G\in U)|$.
Proposition \ref{pn:BE_bound_graph} gives an upper bound of $\zeta_n$ based on the maximum degrees $D_n$ and $D_n^*$ of the dependency graph determined by the underlying sequence $\{X_t\}_{t=1}^n$, whose proof is given in Section \ref{sec:pfpropP2}.

\begin{proposition}\label{pn:BE_bound_graph}
	Under Conditions {\rm\ref{as:tail}} and {\rm\ref{as:eigenvalue1}}, then $\zeta_{n}\lesssim B_{n}^{3/2}(D_nD_n^*)^{1/2}p^{13/4}n^{-1/4}$,
	where $D_n$ and $D^*_n$ are the maximum degrees of the first-degree and second-degree of connections in the dependency graph generated by the sequence $\{X_t\}_{t=1}^n$, respectively.
\end{proposition}

Without loss of generality, we assume $B_n^6(D_nD_n^*)^2\leq n$, since otherwise the assertions are trivial. Theorem \ref{pn:GA_dependency_graph_s_sparsely} can be derived following the same strategy stated in Section \ref{sec:pf_s-sparsely} for the proof of Theorem \ref{pn:GA_s-sparsely}. We only point out the differences here. In Case 1,  applying Theorem \ref{pn:GA_dependency_graph_simple_convex} with replacing $B_n$ by $sB_n$, it holds that $$|\mathbb{P}(S_{n,x}\in A^{*})-\mathbb{P}(G\in A^{*})|\lesssim n^{-1/6}s^{10/3}B_{n}(D_nD_n^*)^{1/3}(\log p)^{7/6}\,.$$ In Case 2, by Proposition \ref{pn:BE_bound_graph} with setting the general parameter $p$ (that is different from the notation $p$ here) in Proposition \ref{pn:BE_bound_graph} as $s$, we then have 
\begin{align*}
|\mathbb{P}(S_{n,x}\in A^{*}_{q})-\mathbb{P}(G\in A^{*}_{q})|\lesssim&\,\, n^{-1/4}s^{13/4}B_{n}^{3/2}(D_nD_n^*)^{1/2}\\
\lesssim&\,\, n^{-1/6}s^{10/3}B_{n}(D_nD_n^*)^{1/3}(\log p)^{7/6}\,.
\end{align*}
Furthermore, following the arguments for Case 2 in Section \ref{sec:pf_s-sparsely}, we have $|\mathbb{P}(S_{n,x}\in A^{*})-\mathbb{P}(G\in A^{*})|\lesssim n^{-1/6}s^{10/3}B_{n}(D_nD_n^*)^{1/3}(\log p)^{7/6}$. In Case 3, applying Theorem \ref{pn:GA_dependency_graph_simple_convex} again with replacing $B_n$ by $sB_n$, we have $$|\mathbb{P}(S_{n,x}\in A^{0})-\mathbb{P}(G\in A^{0})|\lesssim n^{-1/6}s^{10/3}B_{n}(D_nD_n^*)^{1/3}(\log p)^{7/6}\,.$$ By the arguments for Case 3 in Section \ref{sec:pf_s-sparsely}, we have $$
|\mathbb{P}(S_{n,x}\in A^*)-\mathbb{P}(G\in A^*)|\lesssim n^{-1/6}s^{10/3}B_{n}(D_nD_n^*)^{1/3}(\log p)^{7/6}\,.$$
We complete the proof. 
$\hfill\Box$

\subsection{Proof of Proposition \ref{pn:BE_bound_graph}}\label{sec:pfpropP2}
Without loss of generality, we assume $p^{9/2}D_nD_n^*B_{n}^{3}\leq n^{1/2}$, since otherwise the assertion holds trivially. For any $t\in [n]$, let  $X^*_t=(X_{t,1}^*,\ldots,X_{t,p}^*)^\T :=\Xi^{-1/2}X_{t}$
with $\Xi=\cov(S_{n,x})$ and $S_{n,x}^{*}=n^{-1/2}\sum_{t=1}^n X_{t}^*$. Recall $S_{n,x} = n^{-1/2} \sum_{t=1}^n X_t$. Then $S_{n,x}=\Xi^{1/2}S_{n,x}^*$. Similar to \eqref{eq:transform_BE_block}, we have $$\zeta_{n}=\sup_{U\in\mathcal{U}}|\mathbb{P}(S_{n,x}^{*}\in U)-\mathbb{P}(Z\in U)|$$
with $Z\sim N(0,{\rm I}_{p})$.
For $a>0$, $\mathscr{T}_{a}\tilde{h}(\cdot)$ and $\psi_{a}(\cdot)$ defined in the proof of Lemma \ref{la:BE_bound_Block}, identical to \eqref{eq:expression}, we have
\begin{align}\label{eq:Tt_graph}
\mathbb{E}[\mathscr{T}_{a}\tilde{h}(S_{n,x}^*)]=\mathbb{E}[\Delta\psi_{a}(S_{n,x}^{*})-S_{n,x}^{*}\cdot\nabla\psi_{a} (S_{n,x}^{*})]\,.
\end{align}
Recall $\mathcal{N}_{t} = \{s\in V_n: (t,s) \in E_{n}\}$ for any $t\in[n]$. Let $\{X_{t}^{*\dag}\}$ be an independent copy of $\{X_{t}^{*}\}$. Write $X_{t}^{*\dag}=(X_{t,1}^{*\dag},\ldots,X_{t,p}^{*\dag})^\T$. Define $\delta_{t}^*=\sum_{s\in\mathcal{N}_{t}}X_{s}^*=(\delta_{t,1}^*,\ldots,\delta_{t,p}^*)^{\T}$, $\delta_{t}^{*\dag}=\sum_{s\in\mathcal{N}_{t}}X_{s}^{*\dag}=(\delta_{t,1}^{*\dag},\ldots,\delta_{t,p}^{*\dag})^{\T}$  and $S_{n,x}^{*(-t)}=S_{n,x}^{*}-n^{-1/2}\delta_{t}^*$. Since $X_t^{*\dag}$ is independent of $\{X_{s}^{*\dag}\}_{s\notin\mathcal{N}_t}$, by $\cov(n^{-1/2}\sum_{t=1}^{n}X_{t}^{*\dag})={\rm I}_{p}$, we have $n^{-1}\sum_{t=1}^{n}\mathbb{E}(X_{t,j}^{*\dag}\delta_{t,k}^{*\dag})=I(j=k)$.
By Taylor expansion, it holds that
\begin{align*}
&\mathbb{E}\{\Delta\psi_{a}(S_{n,x}^{*})\}=\sum_{j,k=1}^{p}\mathbb{E}\{\partial_{jk}\psi_{a}(S_{n,x}^{*})\}\bigg\{\frac{1}{n}\sum_{t=1}^{n}\mathbb{E}(X_{t,j}^{*\dag}\delta_{t,k}^{*\dag})\bigg\}\\
&~~~=\frac{1}{n}\sum_{t=1}^{n}\sum_{j,k=1}^{p}\mathbb{E}\big\{X_{t,j}^{*\dag}\delta_{t,k}^{*\dag}\partial_{jk}\psi_{a}(S_{n,x}^{*})\big\}\\
&~~~=\frac{1}{n}\sum_{t=1}^{n}\sum_{j,k=1}^{p}\mathbb{E}\big[X_{t,j}^{*\dag}\delta_{t,k}^{*\dag}\partial_{jk}\psi_{a}\{S_{n,x}^{*(-t)}\}\big]\\
&~~~~+\frac{1}{n^{3/2}}\sum_{t=1}^{n}\sum_{j,k,l=1}^{p}\mathbb{E}\bigg[X_{t,j}^{*\dag}\delta_{t,k}^{*\dag}\delta_{t,l}^{*}\int_{0}^{1}\partial_{jkl}\psi_{a}\{S_{n,x}^{*(-t)}+vn^{-1/2}\delta_{t}^*\}{\rm d }v\bigg]\,,\\
&\mathbb{E}\{S_{n,x}^{*}\cdot\nabla\psi_{a}(S_{n,x}^{*})\}=\frac{1}{n^{1/2}}\sum_{t=1}^{n}\sum_{j=1}^{p}\mathbb{E}\big\{X_{t,j}^*\partial_{j}\psi_{a}(S_{n,x}^{*})\big\}\\
&~~~=\frac{1}{n^{1/2}}\sum_{t=1}^{n}\sum_{j=1}^{p}\mathbb{E}\big[X_{t,j}^*\partial_{j}\psi_{a}\{S_{n,x}^{*(-t)}\}\big]+\frac{1}{n}\sum_{t=	1}^{n}\sum_{j,k=1}^{p}\mathbb{E}\big[X_{t,j}^*\delta_{t,k}^*\partial_{jk}\psi_{a}\{S_{n,x}^{*(-t)}\}\big]\\
&~~~~+\frac{1}{n^{3/2}}\sum_{t=1}^{n}\sum_{j,k,l=1}^{p}\mathbb{E}\bigg[X_{t,j}^*\delta_{t,k}^*\delta_{t,l}^*\int_{0}^{1}(1-v)\partial_{jkl}\psi_{a}\{S_{n,x}^{*(-t)}+vn^{-1/2}\delta_{t}^*\}{\rm d}v\bigg]\,.
\end{align*}
Notice that $\mathbb{E}[X_{t,j}^*\partial_{j}\psi_{a}\{S_{n,x}^{*(-t)}\}]=\mathbb{E}(X_{t,j}^*)\mathbb{E}[\partial_{j}\psi_{a}\{S_{n,x}^{*(-t)}\}]=0$.
By \eqref{eq:Tt_graph}, we have
\begin{align}\label{eq:graphI123}
\mathbb{E}\{\mathscr{T}_{a}\tilde{h}(S_{n,x}^{*})\}
=&\,\underbrace{\frac{1}{n}\sum_{t=1}^{n}\sum_{j,k=1}^{p}\mathbb{E}\big[(X_{t,j}^{*\dag}\delta_{t,k}^{*\dag}-X_{t,j}^*\delta_{t,k}^*)\partial_{jk}\psi_{a}\{S_{n,x}^{*(-t)}\}\big]}_{J_{1}}\notag\\
&+\underbrace{\frac{1}{n^{3/2}}\sum_{t=1}^{n}\sum_{j,k,l=1}^{p}\mathbb{E}\bigg[X_{t,j}^{*\dag}\delta_{t,k}^{*\dag}\delta_{t,l}^{*}\int_{0}^{1}\partial_{jkl}\psi_{a}\{S_{n,x}^{*(-t)}+vn^{-1/2}\delta_{t}^*\}{\rm d }v\bigg]}_{J_{2}}\\
&-\underbrace{\frac{1}{n^{3/2}}\sum_{t=1}^{n}\sum_{j,k,l=1}^{p}\mathbb{E}\bigg[X_{t,j}^*\delta_{t,k}^*\delta_{t,l}^*\int_{0}^{1}(1-v)\partial_{jkl}\psi_{a}\{S_{n,x}^{*(-t)}+vn^{-1/2}\delta_{t}^*\}\,{\rm d}v\bigg]}_{J_{3}}\,.\notag
\end{align}

In the sequel, we will first bound $J_{1}$, $J_{2}$ and $J_{3}$, respectively. Define $$\delta^{*}_{\mathcal{N}_{t}}=(\delta_{\mathcal{N}_t,1}^*,\ldots,\delta_{\mathcal{N}_t,p}^*)^\T:=\sum_{s\in\cup_{\ell\in\mathcal{N}_{t}}\mathcal{N}_{\ell}}X_{s}^*$$ and $S_{n,x}^{*(-\mathcal{N}_{t})}=S_{n,x}^*-n^{-1/2}\delta_{\mathcal{N}_{t}}^*$. Notice that $\partial_{jk}\psi_{a}\{S_{n,x}^{*(-\mathcal{N}_{t})}\}$ is independent of $X_{t,j}^{*\dag}\delta_{t,k}^{*\dag}-X_{t,j}^*\delta_{t,k}^*$ and $\mathbb{E}(X_{t,j}^{*\dag}\delta_{t,k}^{*\dag}-X_{t,j}^*\delta_{t,k}^*)=0$. Due to $S_{n,x}^{*(-t)}-S_{n,x}^{*(-\mathcal{N}_{t})}=n^{-1/2}(\delta_{\mathcal{N}_{t}}^*-\delta_{t}^*)$,
by Taylor expansion, 
\begin{align*}
J_{1}=&\,n^{-3/2}\sum_{t=1}^{n}\sum_{j,k,l=1}^{p}\mathbb{E}\bigg[(X_{t,j}^{*\dag}\delta_{t,k}^{*\dag}-X_{t,j}^*\delta_{t,k}^*)(\delta_{\mathcal{N}_{t},l}^*-\delta_{t,l}^*)\\
&~~~~~~~~~\times\int_{0}^{1}\partial_{jkl}\psi_{a}\{S_{n,x}^{*(-\mathcal{N}_{t})}+vn^{-1/2}(\delta_{\mathcal{N}_{t}}^*-\delta_{t}^*)\}{\rm d}v\bigg]\,.
\end{align*}
As shown in the proof of Lemma \ref{la:BE_bound_Block}, $
\sup_{x\in\mathbb{R}^p}\max_{j,k,l\in [p]}|\partial_{jkl}\psi_{a}(x)|\lesssim a^{-1/2}$. Then
\begin{align}\label{eq:J1graph}
|J_{1}|
\lesssim&\, \frac{1}{a^{1/2}n^{3/2}}\sum_{t=1}^{n}\sum_{j,k,l=1}^{p}\big\{\mathbb{E}\big(|X_{t,j}^{*\dag}\delta_{t,k}^{*\dag}||\delta_{\mathcal{N}_{t},l}^*-\delta_{t,l}^*|\big)+\mathbb{E}\big(|X_{t,j}^{*}\delta_{t,k}^{*}||\delta_{\mathcal{N}_{t},l}^*-\delta_{t,l}^*|\big)\big\}\,.
\end{align}
Due to $\delta_{\mathcal{N}_{t}}^*-\delta_{t}^*=\sum_{s\in\{\cup_{\ell\in\mathcal{N}_{t}}\mathcal{N}_{\ell}\}\backslash\mathcal{N}_{t}}X_{s}^*$, by triangle inequality, we have
\begin{align*}
&\frac{1}{a^{1/2}n^{3/2}}\sum_{t=1}^{n}\sum_{j,k,l=1}^{p}\mathbb{E}\big(|X_{t,j}^{*}\delta_{t,k}^{*}||\delta_{\mathcal{N}_{t},l}^*-\delta_{t,l}^*|\big)\\
&~~~\leq\frac{1}{a^{1/2}n^{3/2}}\sum_{t=1}^{n}\sum_{j,k,l=1}^{p}\mathbb{E}\bigg\{|X_{t,j}^{*}|\bigg(\sum_{s\in\mathcal{N}_{t}}|X_{s,k}^*|\bigg)\bigg(\sum_{s\in\{\cup_{\ell\in\mathcal{N}_{t}}\mathcal{N}_{\ell}\}\backslash\mathcal{N}_{t}}|X_{s,l}^*|\bigg)\bigg\}\\
&~~~=\frac{1}{a^{1/2}n^{3/2}}\sum_{t=1}^{n}\sum_{j,k,l=1}^{p}\sum_{s_{1}\in\mathcal{N}_{t}}\sum_{s_{2}\in\{\cup_{\ell\in\mathcal{N}_{t}}\mathcal{N}_{\ell}\}\backslash\mathcal{N}_{t}}\mathbb{E}(|X_{t,j}^{*}X_{s_{1},k}^*X_{s_{2},l}^*|)\,.
\end{align*}
Identical to \eqref{eq:tuta_x}, by Condition \ref{as:eigenvalue1}, $
\max_{j\in[p]}|X_{t,j}^*|\leq\sqrt{p}\lambda_{\min}^{-1/2}(\Xi)|X_{t}|_{\infty}\lesssim \sqrt{p}|X_{t}|_{\infty}$ for any $t\in [n]$. Together with Condition \ref{as:tail},  $$
\max_{t\in [n]}\max_{j\in [p]}\mathbb{P}(|X_{t,j}^*|>u)\leq 2p\exp(-Cp^{-\gamma_{1}/2}u^{\gamma_{1}}B_{n}^{-\gamma_1})$$ for any $u>0$,
which implies $\max_{t\in [n]}\max_{j\in [p]}\mathbb{E}(|X_{t,j}^*|^3)\lesssim p^{5/2}B_{n}^3$. By H\"older's inequality, $a^{-1/2}n^{-3/2}\sum_{t=1}^{n}\sum_{j,k,l=1}^{p}\mathbb{E}(|X_{t,j}^{*}\delta_{t,k}^{*}||\delta_{\mathcal{N}_{t},l}^*-\delta_{t,l}^*|)\lesssim a^{-1/2}n^{-1/2}p^{11/2}D_{n}D_n^*B_{n}^3$. Analogously, $a^{-1/2}n^{-3/2}\sum_{t=1}^{n}\sum_{j,k,l=1}^{p}\mathbb{E}(|X_{t,j}^{*\dag}\delta_{t,k}^{*\dag}||\delta_{\mathcal{N}_{t},l}^*-\delta_{t,l}^*|)\lesssim a^{-1/2}n^{-1/2}p^{11/2}D_{n}D_n^*B_{n}^3$. By \eqref{eq:J1graph},
	$
	|J_{1}| \lesssim a^{-1/2}n^{-1/2}p^{11/2}D_nD_n^*B_{n}^3$.
Using the same arguments, we have
$
|J_{2}|
\lesssim a^{-1/2}n^{-1/2}p^{11/2}D_{n}^{2}B_{n}^3$ and $|J_{3 }|\lesssim a^{-1/2}n^{-1/2}p^{11/2}D_{n}^{2}B_{n}^3$.
Thus, by \eqref{eq:graphI123}, $$
	|\mathbb{E}\{\mathscr{T}_{a}\tilde{h}(S_{n,x}^{*})\}|\lesssim  a^{-1/2}n^{-1/2}p^{11/2}D_{n}D_n^*B_{n}^3\,.$$
For $\mathcal{H}$ specified in the proof of Lemma \ref{la:BE_bound_Block}, identical to the first step of  \eqref{eq:smoothinginequalitybound}, it holds that $\zeta_{n}\lesssim\sup_{h\in\mathcal{H}}|\mathbb{E}\{\mathscr{T}_{a}\tilde{h}(S_{n,x}^{*})\}|+pa^{1/2}e^{a}\lesssim a^{-1/2}n^{-1/2}p^{11/2}D_nD_n^*B_{n}^{3}+pa^{1/2}$
for any $a\in (0,1)$. By choosing $a=n^{-1/2}p^{9/2}D_nD_n^*B_{n}^{3}\leq 1$, $\zeta_n\lesssim n^{-1/4}p^{13/4}(D_nD_n^*)^{1/2}B_n^{3/2}$. We complete the proof of Proposition \ref{pn:BE_bound_graph}.$\hfill\Box$

\section{Proof of Theorem \ref{pn:GA_weakly-dependent_s_sparsely}}\label{sec:gawds}

Recall $\mathcal{U}$ is the class of all convex sets in $\mathbb{R}^{p}$ and
$
\zeta_{n}:=\sup_{U\in\mathcal{U}}|\mathbb{P}(S_{n,x}\in U)-\mathbb{P}(G\in U)|$. 
Proposition \ref{pn:BE_bound_physical} gives an upper bound of $\zeta_n$ when $\{X_t\}$ is a sequence with physical dependence, whose proof is given in Section \ref{sec:pfpropP3}.
\begin{proposition}\label{pn:BE_bound_physical}
	Assume $\Phi_{\psi_{\nu},\alpha} < \infty$ for some $\alpha, \nu \in (0,\infty)$.
	
	{\rm (i)} Under Condition {\rm\ref{as:eigenvalue1}}, then 
	\begin{align*}
	\zeta_{n}\lesssim&\, n^{-\alpha/(2+4\alpha)}p^{3/4}\Phi_{\psi_\nu,0}\{p^{7/2}\Phi_{\psi_\nu,0}^2+(\log n)^{(1+2\nu)/2}\}\\
	&+n^{-\alpha/(1+2\alpha)}p^{3/4}\Phi_{\psi_{\nu},\alpha}(\log n)^{(1+2\nu)/2}+n^{-\alpha/(2+4\alpha)}p^2\Psi_{2,\alpha}^{1/2} \Psi_{2,0}^{1/2}\,.
	\end{align*}
	
	{\rm (ii)} Under Conditions {\rm\ref{as:tail}} and {\rm\ref{as:eigenvalue1}}, then 
	\begin{align*}\zeta_{n}\lesssim&\, n^{-\alpha/(8+4\alpha)}(B_{n}^{3/2}p^{13/4}+p^{2}\Psi_{2,\alpha}^{1/2} \Psi_{2,0}^{1/2})\\
	&+n^{-\alpha/(4+2\alpha)}p^{3/4}\Phi_{\psi_{\nu},\alpha}(\log n)^{(1+2\nu)/2}\,.
	\end{align*}
\end{proposition}

Theorem \ref{pn:GA_weakly-dependent_s_sparsely} can be derived following the same strategy stated in Section \ref{sec:pf_s-sparsely} for the proof of Theorem \ref{pn:GA_s-sparsely}. We only point out the differences here.
\subsection{Proof of Part {\rm(i)} of Theorem {\rm\ref{pn:GA_weakly-dependent_s_sparsely}}}
Assume $s^{4}\Psi_{2,\alpha,\Omega_{s,c}}\Psi_{2,0,\Omega_{s,c}}=o\{n^{\alpha/(1+3\alpha)}\}$. Otherwise the assertions hold trivially. In Case 1, by Theorem \ref{pn:GA_weakly-dependent_simple_convex}{\rm (i)}, 
\begin{align*}
&|\mathbb{P}(S_{n,x}\in A^{*})-\mathbb{P}(G\in A^{*})|\\
&~~~~\lesssim n^{-\alpha/(3+9\alpha)}(s^2\log p)^{2/3}\{\Phi_{\psi_{\nu},0,\Omega_{s,c}}(s^2\log p)^{1/2}+\Psi_{2,\alpha,\Omega_{s,c}}^{1/3}\Psi_{2,0,\Omega_{s,c}}^{1/3}\}\\
&~~~~~+n^{-\alpha/(1+3\alpha)}(s^2\log p)^{1+\nu}\Phi_{\psi_{\nu},\alpha,\Omega_{s,c}}
\end{align*}
provided that $(s^2\log p)^{\max\{6\nu-1,(5+6\nu)/4\}}=o\{n^{\alpha/(1+3\alpha)}\}$. In Case 2, using Proposition \ref{pn:BE_bound_physical}{\rm (i)} with setting $p$ as $s$, 
\begin{align*}
&|\mathbb{P}(S_{n,x}\in A^{*}_{q})-\mathbb{P}(G\in A^{*}_{q})|\\
&~~~~\lesssim 
n^{-\alpha/(3+9\alpha)}(s^2\log p)^{2/3}\{\Phi_{\psi_{\nu},0,\Omega_{s,c}}(s^2\log p)^{1/2}+\Psi_{2,\alpha,\Omega_{s,c}}^{1/3}\Psi_{2,0,\Omega_{s,c}}^{1/3}\}\\
&~~~~~+n^{-\alpha/(1+3\alpha)}\Phi_{\psi_{\nu},\alpha,\Omega_{s,c}}(s^2\log p)^{1+\nu}+n^{-\alpha/(2+4\alpha)}s^{17/4}\Phi_{\psi_\nu,0}^3\,.
\end{align*}
By the arguments for Case 2 in Section \ref{sec:pf_s-sparsely}, 
\begin{align*}
&|\mathbb{P}(S_{n,x}\in A^{*})-\mathbb{P}(G\in A^{*})|\\
&~~~~\lesssim n^{-\alpha/(3+9\alpha)}(s^2\log p)^{2/3}\{\Phi_{\psi_{\nu},0,\Omega_{s,c}}s(\log p)^{1/2}+\Psi_{2,\alpha,\Omega_{s,c}}^{1/3} \Psi_{2,0,\Omega_{s,c}}^{1/3}\}\\
&~~~~~+n^{-\alpha/(1+3\alpha)}\Phi_{\psi_{\nu},\alpha,\Omega_{s,c}}(s^2\log p)^{1+\nu}+n^{-\alpha/(2+4\alpha)} s^{17/4}\Phi_{\psi_\nu,0}^3\,.
\end{align*}
In Case 3, by Theorem \ref{pn:GA_weakly-dependent_simple_convex}{\rm (i)},
\begin{align*}
&|\mathbb{P}(S_{n,x}\in A^{0})-\mathbb{P}(G\in A^{0})|\\
&~~~~\lesssim n^{-\alpha/(3+9\alpha)}(s^2\log p)^{2/3}\{\Phi_{\psi_{\nu},0,\Omega_{s,c}}s(\log p)^{1/2}+\Psi_{2,\alpha,\Omega_{s,c}}^{1/3}\Psi_{2,0,\Omega_{s,c}}^{1/3}\}\\
&~~~~~+n^{-\alpha/(1+3\alpha)}\Phi_{\psi_{\nu},\alpha,\Omega_{s,c}}(s^2\log p)^{1+\nu}
\end{align*}
provided that $(s^2\log p)^{\max\{6\nu-1,(5+6\nu)/4\}}=o\{n^{\alpha/(1+3\alpha)}\}$. By the arguments for Case 3 in Section \ref{sec:pf_s-sparsely}, 
\begin{align*}
&|\mathbb{P}(S_{n,x}\in A^*)-\mathbb{P}(G\in A^*)|\\
&~~~~\lesssim n^{-\alpha/(3+9\alpha)}(s^2\log p)^{2/3}\{\Phi_{\psi_{\nu},0,\Omega_{s,c}}s(\log p)^{1/2}+\Psi_{2,\alpha,\Omega_{s,c}}^{1/3}\Psi_{2,0,\Omega_{s,c}}^{1/3}\}\\
&~~~~~+n^{-\alpha/(1+3\alpha)}\Phi_{\psi_{\nu},\alpha,\Omega_{s,c}}(s^2\log p)^{1+\nu}\,.
\end{align*}
We complete the proof of Part {\rm(i)} of Theorem {\rm\ref{pn:GA_weakly-dependent_s_sparsely}}. 
$\hfill\Box$

\subsection{Proof of Part {\rm(ii)} of Theorem {\rm\ref{pn:GA_weakly-dependent_s_sparsely}}}

Let $B_{n}=o\{n^{\alpha/(12+6\alpha)}\}$ and $s^{4}\Psi_{2,\alpha,\Omega_{s,c}}\Psi_{2,0,\Omega_{s,c}}=o\{n^{\alpha/(4+2\alpha)}\}$. Otherwise the assertions hold trivially. In Case 1, by Theorem \ref{pn:GA_weakly-dependent_simple_convex}{\rm (ii)} with replacing $B_n$ by $sB_n$, we have
\begin{align*}
|\mathbb{P}(S_{n,x}\in A^{*})-\mathbb{P}(G\in A^{*})|\lesssim&\, n^{-\alpha/(12+6\alpha)}(s^2\log p)^{2/3}\{s^{2}B_{n}(\log p)^{1/2}+\Psi_{2,\alpha,\Omega_{s,c}}^{1/3}\Psi_{2,0,\Omega_{s,c}}^{1/3}\}\\
&+n^{-\alpha/(4+2\alpha)}\Phi_{\psi_{\nu},\alpha,\Omega_{s,c}}(s^2\log p)^{1+\nu}\,.
\end{align*}
In Case 2, by Proposition \ref{pn:BE_bound_physical}{\rm (ii)} with setting $p$ as $s$, we have 
\begin{align*}
|\mathbb{P}(S_{n,x}\in A^{*}_{q})-\mathbb{P}(G\in A^{*}_{q})|
\lesssim&\, n^{-\alpha/(12+6\alpha)}(s^2\log p)^{2/3}\{s^{2}B_{n}(\log p)^{1/2}+\Psi_{2,\alpha,\Omega_{s,c}}^{1/3}\Psi_{2,0,\Omega_{s,c}}^{1/3}\}\\
&+n^{-\alpha/(4+2\alpha)}\Phi_{\psi_{\nu},\alpha,\Omega_{s,c}}(s^2\log p)^{1+\nu}\,.
\end{align*}
Following the arguments for Case 2 in Section \ref{sec:pf_s-sparsely}, we have 
\begin{align*}
    |\mathbb{P}(S_{n,x}\in A^{*})-\mathbb{P}(G\in A^{*})|\lesssim&\, n^{-\alpha/(12+6\alpha)}(s^2\log p)^{2/3}\{s^{2}B_{n}(\log p)^{1/2}+\Psi_{2,\alpha,\Omega_{s,c}}^{1/3}\Psi_{2,0,\Omega_{s,c}}^{1/3}\}\\
    &+n^{-\alpha/(4+2\alpha)}\Phi_{\psi_{\nu},\alpha,\Omega_{s,c}}(s^2\log p)^{1+\nu}\,.
    \end{align*}
    In Case 3, by Theorem \ref{pn:GA_weakly-dependent_simple_convex}{\rm (ii)} with replacing $B_n$ by $sB_n$, we have 
\begin{align*}
|\mathbb{P}(S_{n,x}\in A^{0})-\mathbb{P}(G\in A^{0})|\lesssim&\, n^{-\alpha/(12+6\alpha)}(s^2\log p)^{2/3}\{s^{2}B_{n}(\log p)^{1/2}+\Psi_{2,\alpha,\Omega_{s,c}}^{1/3}\Psi_{2,0,\Omega_{s,c}}^{1/3}\}\\
&+n^{-\alpha/(4+2\alpha)}\Phi_{\psi_{\nu},\alpha,\Omega_{s,c}}(s^2\log p)^{1+\nu}\,.
\end{align*}
Following the arguments for Case 3 in Section \ref{sec:pf_s-sparsely}, it holds that  
\begin{align*}
|\mathbb{P}(S_{n,x}\in A^*)-\mathbb{P}(G\in A^*)|\lesssim&\, n^{-\alpha/(12+6\alpha)}(s^2\log p)^{2/3}\{s^{2}B_{n}(\log p)^{1/2}+\Psi_{2,\alpha,\Omega_{s,c}}^{1/3}\Psi_{2,0,\Omega_{s,c}}^{1/3}\}\\
&+n^{-\alpha/(4+2\alpha)}\Phi_{\psi_{\nu},\alpha,\Omega_{s,c}}(s^2\log p)^{1+\nu}\,.
\end{align*}
We complete the proof of Part {\rm(ii)} of Theorem {\rm\ref{pn:GA_weakly-dependent_s_sparsely}}. 
$\hfill\Box$

\section{Proof of Proposition {\rm\ref{pn:BE_bound_physical}}}\label{sec:pfpropP3}

\subsection{Proof of Part {\rm(i)} of Proposition {\rm\ref{pn:BE_bound_physical}}} Without loss of generality, we assume $p^{17/2}\Phi_{\psi_{\nu},0}^6+p^4\Psi_{2,\alpha} \Psi_{2,0}=o\{n^{\alpha/(1+2\alpha)}\}$, since otherwise the assertions hold trivially. In the proof of Part (i), we still follow the notation used in Section \ref{sec:pfphysical_m}. Define
\begin{align*}
\zeta_{n,1}^{(m)}=&\,\sup_{U\in\mathcal{U}}|\mathbb{P}\{\tilde S_{n,x}^{(m)}\in U\}-\mathbb{P}\{\tilde S_{n,y}^{(m)}\in U\}|\,, \\
\zeta_{n,2}^{(m)}=&\,\sup_{U\in\mathcal{U}}|\mathbb{P}\{S_{n,x}^{(m)}\in U\}-\mathbb{P}\{\tilde S_{n,y}^{(m)}\in U\}|\,.
\end{align*}
Recall that $\Xi=\cov(S_{n,x})$ and $\tilde{\Xi}=\cov\{\tilde S_{n,x}^{(m)}\}$. Write $\tilde X_\ell^{(m)*}=\tilde{\Xi}^{-1/2}\tilde X_\ell^{(m)}=\{\tilde X_{\ell,1}^{(m)*},\ldots, \tilde X_{\ell,p}^{(m)*}\}^{\T}$. Similar to \eqref{eq:transform_BE_block}, we have $$\zeta_{n,1}^{(m)}=\sup_{U\in\mathcal{U}}|\mathbb{P}\{\tilde S_{n,x}^{(m)*}\in U\}-\mathbb{P}(Z\in U)|\,,$$
where $\tilde S_{n,x}^{(m)*}=L^{-1/2}\sum_{\ell=1}^L \tilde X_\ell^{(m)*}$ with $\cov\{S_{n,x}^{(m)*}\}={\rm I}_p$ and $Z\sim N(0,{\rm I}_{p})$. By Lemma \ref{lem:cov_infinity_m}, $$\|\tilde\Xi-\Xi\|_{2}\leq p|\tilde\Xi-\Xi|_\infty\lesssim p\Phi_{\psi_{\nu},0}^2(mb^{-1}+bn^{-1})+pm^{-\alpha} \Psi_{2,\alpha} \Psi_{2,0}\,.$$ Thus, Weyl's inequality yields
$
|\lambda_{\min}(\tilde\Xi)-\lambda_{\min}(\Xi)|\leq\|\tilde\Xi-\Xi\|_{2}=o(1)$ provided that $p\Phi_{\psi_{\nu},0}^2m \ll b\ll p^{-1}\Phi_{\psi_{\nu},0}^{-2}n$ and $m^{\alpha}\gg p\Psi_{2,\alpha} \Psi_{2,0}$. By Condition \ref{as:eigenvalue1}, we have $\lambda_{\min}(\tilde\Xi)$ is bounded away from zero. Same as \eqref{eq:tuta_x}, $
|\tilde X_{\ell}^{(m)*}|_\infty\lesssim \sqrt{p}|\tilde X_{\ell}^{(m)}|_{\infty}$ for any $\ell\in[L]$, which yields $L^{-3/2}\sum_{\ell=1}^{L}\mathbb{E}\{|\tilde X_{\ell}^{(m)*}|_2^3\}\lesssim p^3L^{-3/2}\sum_{\ell=1}^{L}\mathbb{E}\{|\tilde X_{\ell}^{(m)}|_{\infty}^3\}$ for any $\ell\in [L]$. By Lemma \ref{lem:tail_prob_bound_m}, we have  $$\max_{\ell\in[L]}\mathbb{P}\{|\tilde{X}_{\ell}^{(m)}|_{\infty}>u\}\lesssim p\exp\{-C(u\Phi_{\psi_{\nu},0}^{-1})^{2/(1+2\nu)}\}$$ for any $u>0$, which implies $\max_{\ell\in [L]}\mathbb{E}\{|\tilde X_{\ell}^{(m)}|_{\infty}^3\}\lesssim p\Phi_{\psi_{\nu},0}^3$.
Thus, $$\frac{1}{L^{3/2}}\sum_{\ell=1}^{L}\mathbb{E}\{|\tilde X_{\ell}^{(m)*}|_2^3\}\lesssim L^{-1/2}p^4\Phi_{\psi_{\nu},0}^3\,.$$  Notice that $\{\tilde{X}_\ell^{(m)*}\}_{\ell=1}^L$ is an independent sequence. By Theorem 1.1 of Rai\v{c} (2019), we have $\zeta_{n,1}^{(m)}\lesssim L^{-1/2}p^{17/4}\Phi_{\psi_{\nu},0}^3$. Notice that $
S_{n,x}^{(m)}=\tilde S_{n,x}^{(m)}+\delta_n$ with $$
\delta_n=n^{-1/2}\sum_{\ell=1}^{L+1}\sum_{t\in\mathcal{J}_\ell}{X}_t^{(m)}\\ +\{n^{-1/2}-(Lb)^{-1/2}\}\sum_{\ell=1}^L\sum_{t\in\mathcal{I}_\ell}{X}_t^{(m)}\,.$$
Let $D_{n1}=C b^{-1/2}m^{1/2}\Phi_{\psi_\nu,0}(\log n)^{(1+2\nu)/2}$ for some sufficiently large constant $C>0$. Consider the event $\mathcal{E}_1=\{|\delta_n|_\infty\leq {D}_{n1}\}$.
Restricted on $\mathcal{E}_1$, we have $|S_{n,x}^{(m)}-\tilde S_{n,x}^{(m)}|_{2}=|\delta_n|_{2}\leq \sqrt{p}D_{n1}=:R_1$. Identical to \eqref{eq:substituteblock},  
\begin{align*}
|\mathbb{P}\{S_{n,x}^{(m)}\in U\}-\mathbb{P}\{\tilde S_{n,y}^{(m)}\in U\}|\leq&\,\, \zeta_{n,1}^{(m)}+\mathbb{P}(\mathcal{E}_1^c)+\mathbb{P}\{\tilde S_{n,y}^{(m)}\in U^{R_1}\backslash U\}+\mathbb{P}\{\tilde S_{n,y}^{(m)}\in U\backslash U^{-R_1}\}
\end{align*}
for any $U\in\mathcal{U}$. By Lemma A.2 of Chernozhukov, Chetverikov and Kato (2017), 
\begin{align*}
&\sup_{U\in\mathcal{U}}\mathbb{P}\{\tilde S_{n,y}^{(m)}\in U^{R_1}\backslash U\}+\sup_{U\in\mathcal{U}}\mathbb{P}\{\tilde S_{n,y}^{(m)}\in U\backslash U^{-R_1}\}\\
&~~~~\lesssim p^{3/4}D_{n1}\lesssim p^{3/4}b^{-1/2}m^{1/2}\Phi_{\psi_\nu,0}(\log n)^{(1+2\nu)/2}\,.
\end{align*}
Following the same arguments for deriving the upper bound of $\mathbb{P}(\mathcal{E}^c)$ in Section \ref{sec:pflpn2m}, if $b\lesssim (mn)^{1/2}$, it then holds that
$$\mathbb{P}(\mathcal{E}_1^c)\lesssim p\exp\{-C(b^{1/2}D_{n1}m^{-1/2}\Phi_{\psi_\nu,0}^{-1})^{2/(1+2\nu)}\}\lesssim L^{-1/2}p^{17/4}\Phi_{\psi_\nu,0}^3\,.$$
Hence, we have $\zeta_{n,2}^{(m)}\lesssim L^{-1/2}p^{17/4}\Phi_{\psi_{\nu},0}^3+p^{3/4}b^{-1/2}m^{1/2}\Phi_{\psi_\nu,0}(\log n)^{(1+2\nu)/2}$ provided that $p\Phi_{\psi_{\nu},0}^2m \ll b\ll p^{-1}\Phi_{\psi_{\nu},0}^{-2}n$, $b\lesssim (mn)^{1/2}$ and $m^{\alpha}\gg p\Psi_{2,\alpha} \Psi_{2,0}$.

Recall $G\sim N(0,\Xi)$ and $\tilde{S}_{n,y}^{(m)}\sim N(0,\tilde{\Xi})$. Let $D_{n2}=Cm^{-\alpha}\Phi_{\psi_{\nu},\alpha}(\log n)^{(1+2\nu)/2}$ for some sufficiently large constant $C>0$. Consider the event $\mathcal{E}_2:=\{|S_{n,x}-S_{n,x}^{(m)}|_{\infty}\leq D_{n2}\}$. Let $R_2=p^{1/2}D_{n2}$. By triangle inequality and Lemma \ref{la:comparison_convex}, we have
\begin{align*}
\zeta_{n}&\leq \sup_{U\in\mathcal{U}}|\mathbb{P}(S_{n,x}\in U)-\mathbb{P}\{\tilde S_{n,y}^{(m)}\in U\}|+\sup_{U\in\mathcal{U}}|\mathbb{P}\{\tilde S_{n,y}^{(m)}\in U\}-\mathbb{P}(G\in U)|\\
&\lesssim\zeta_{n,2}^{(m)}+\mathbb{P}(\mathcal{E}_2^c)+\sup_{U\in\mathcal{U}}\mathbb{P}\{\tilde S_{n,y}^{(m)}\in U^{R_2}\backslash U\}+\sup_{U\in\mathcal{U}}\mathbb{P}\{\tilde S_{n,y}^{(m)}\in U\backslash U^{-R_2}\}+p^2|\tilde\Xi-\Xi|_{\infty}^{1/2}\,,
\end{align*}
where the last step is identical to \eqref{eq:substituteblock}. Notice that $\sup_{U\in\mathcal{U}}\mathbb{P}\{\tilde S_{n,y}^{(m)}\in U^{R_2}\backslash U\}+\sup_{U\in\mathcal{U}}\mathbb{P}\{\tilde S_{n,y}^{(m)}\in U\backslash U^{-R_2}\}\lesssim p^{3/4}D_{n2}\lesssim p^{3/4}m^{-\alpha}\Phi_{\psi_{\nu},\alpha}(\log n)^{(1+2\nu)/2}$. Due to $b\lesssim (mn)^{1/2}$, it holds that 
$$|\tilde\Xi-\Xi|_\infty\lesssim \Phi_{\psi_{\nu},0}^2mb^{-1}+m^{-\alpha} \Psi_{2,\alpha} \Psi_{2,0}\,.$$ It follows from Lemma \ref{lem:tail_prob_bound_m-dep_approx} that $
\mathbb{P}(\mathcal{E}_2^c)\lesssim 
n^{-1/2}b^{1/2}p^{17/4}\Phi_{\psi_\nu,0}^3$. Then 
\begin{align*}
\zeta_{n}\lesssim&\,\, n^{-1/2}b^{1/2}p^{17/4}\Phi_{\psi_\nu,0}^3+b^{-1/2}m^{1/2}\Phi_{\psi_\nu,0}\{p^{3/4}(\log n)^{(1+2\nu)/2}+p^2\}\\
&+p^2m^{-\alpha/2} \Psi_{2,\alpha}^{1/2} \Psi_{2,0}^{1/2}+p^{3/4}m^{-\alpha}\Phi_{\psi_{\nu},\alpha}(\log n)^{(1+2\nu)/2}
\end{align*}
provided that $p\Phi_{\psi_{\nu},0}^2m \ll b\ll p^{-1}\Phi_{\psi_{\nu},0}^{-2}n$, $b\lesssim (mn)^{1/2}$ and $m^\alpha\gg p\Psi_{2,\alpha}\Psi_{2,0}$. Selecting $b\asymp (nm)^{1/2}$, we have 
\begin{align*}
\zeta_{n}\lesssim&\,\, n^{-1/4}m^{1/4}p^{3/4}\Phi_{\psi_\nu,0}\{p^{7/2}\Phi_{\psi_\nu,0}^2+(\log n)^{(1+2\nu)/2}\}+p^2m^{-\alpha/2} \Psi_{2,\alpha}^{1/2} \Psi_{2,0}^{1/2}\\
&+p^{3/4}m^{-\alpha}\Phi_{\psi_{\nu},\alpha}(\log n)^{(1+2\nu)/2}
\end{align*}
provided that $(p\Psi_{2,\alpha} \Psi_{2,0})^{1/\alpha}\ll m\ll p^{-2}\Phi_{\psi_\nu,0}^{-4}n$. Recall $p^{17/2}\Phi_{\psi_{\nu},0}^6+p^4\Psi_{2,\alpha} \Psi_{2,0}=o\{n^{\alpha/(1+2\alpha)}\}$. Selecting $m\asymp n^{1/(1+2\alpha)}$, then 
\begin{align*}\zeta_{n}\lesssim&\,\, n^{-\alpha/(2+4\alpha)}p^{3/4}\Phi_{\psi_\nu,0}\{p^{7/2}\Phi_{\psi_\nu,0}^2+(\log n)^{(1+2\nu)/2}\}\\
&+n^{-\alpha/(1+2\alpha)}p^{3/4}\Phi_{\psi_{\nu},\alpha}(\log n)^{(1+2\nu)/2}+n^{-\alpha/(2+4\alpha)}p^2\Psi_{2,\alpha}^{1/2} \Psi_{2,0}^{1/2}\,.
\end{align*}
We complete the proof of part {\rm(i)} of Proposition {\rm\ref{pn:BE_bound_physical}}. $\hfill\Box$

\subsection{Proof of Part {\rm(ii)} of Proposition {\rm\ref{pn:BE_bound_physical}}} Without loss of generality, we assume $p^4\Psi_{2,\alpha} \Psi_{2,0}=o\{n^{\alpha/(4+2\alpha)}\}$, since otherwise the assertions hold trivially. In this proof, we still follow the notation used in Section \ref{sec:pfphysical}.
Recall $\Xi=\cov(S_{n,x})$ with $S_{n,x}=n^{-1/2}\sum_{t=1}^{n}X_{t}$ and $\Xi^{(m)}=\cov\{S_{n,x}^{(m)}\}$ with $S_{n,x}^{(m)}=n^{-1/2}\sum_{t=1}^{n}X_{t}^{(m)}$. As shown in Section \ref{sec:pfphysical}, $|\Xi^{(m)}-\Xi|_{\infty}\lesssim m^{-\alpha}\Psi_{2,\alpha} \Psi_{2,0}$. For any $v\in \mathbb{S}^{p-1}$, we have
$$
|v^{\T}\{\Xi^{(m)}-\Xi\}v|\leq|\Xi^{(m)}-\Xi|_{\infty}|v|_{1}^2\lesssim pm^{-\alpha}\Psi_{2,\alpha} \Psi_{2,0}\,,$$
which implies $\|\Xi^{(m)}-\Xi\|_{2}\lesssim pm^{-\alpha}\Psi_{2,\alpha} \Psi_{2,0}$. Thus, by Weyl's inequality,
$
|\lambda_{\min}\{\Xi^{(m)}\}-\lambda_{\min}(\Xi)|\leq\|\Xi^{(m)}-\Xi\|_{2}\lesssim pm^{-\alpha}\Psi_{2,\alpha} \Psi_{2,0}$.
With choosing $m\gg p^{1/\alpha}\Psi_{2,\alpha}^{1/\alpha} \Psi_{2,0}^{1/\alpha}$, we have $\lambda_{\min}\{\Xi^{(m)}\}=\lambda_{\min}(\Xi)+o(1)$. By Condition \ref{as:eigenvalue1}, we know $\lambda_{\min}\{\Xi^{(m)}\}$ is uniformly bounded away from zero. Notice that $\{X_{t}^{(m)}\}_{t=1}^{n}$ is an $m$-dependent sequence with mean zero and $\max_{t\in [n],j\in [p]}\|X_{t,j}^{(m)}\|_{\psi_{\gamma_{1}}}\leq B_{n}$.
Applying Proposition \ref{pn:BE_bound_graph} with $D_n=2m+1$ and $D_n^*=4m+1$, we have
$$
\zeta_{n}^{(m)}:=\sup_{U\in\mathcal{U}}|\mathbb{P}\{S_{n,x}^{(m)}\in U\}-\mathbb{P}\{S_{n,y}^{(m)}\in U\}|\lesssim n^{-1/4}B_{n}^{3/2}mp^{13/4}\,,$$
where $S_{n,y}^{(m)}=n^{-1/2}\sum_{t=1}^{n}Y_{t}^{(m)}$ and $\{Y_{t}^{(m)}\}$ is a sequence of centered Gaussian random vectors such that $\Cov\{Y_{t}^{(m)}, Y_{s}^{(m)}\}  = \Cov\{X_{t}^{(m)}, X_{s}^{(m)}\}$ for all $t,s\in[n]$. Recall $G\sim N(0,\Xi)$. By triangle inequality and Lemma \ref{la:comparison_convex}, it holds that
\begin{align*}
\zeta_{n}&\leq \sup_{U\in\mathcal{U}}|\mathbb{P}(S_{n,x}\in U)-\mathbb{P}\{S_{n,y}^{(m)}\in U\}|+\sup_{U\in\mathcal{U}}|\mathbb{P}\{S_{n,y}^{(m)}\in U\}-\mathbb{P}(G\in U)|\\
&\lesssim\sup_{U\in\mathcal{U}}|\mathbb{P}(S_{n,x}\in U)-\mathbb{P}\{S_{n,y}^{(m)}\in U\}|+\frac{p^{2}\Psi_{2,\alpha}^{1/2} \Psi_{2,0}^{1/2}}{m^{\alpha/2}}\,.
\end{align*}
Let $
\tilde D_{n}=Cm^{-\alpha}\Phi_{\psi_{\nu},\alpha}(\log n)^{(1+2\nu)/2}$ for some sufficiently large constant $C>0$. Consider the event $\mathcal{E}:=\{|S_{n,x}-S_{n,x}^{(m)}|_{\infty}\leq \tilde{D}_{n}\}$. Let $R=p^{1/2}\tilde{D}_{n}$. Identical  to \eqref{eq:substituteblock}, for any $U\in\mathcal{U}$, it holds that
\begin{align*}
|\mathbb{P}(S_{n,x}\in U)-\mathbb{P}\{S_{n,y}^{(m)}\in U\}|\leq&\,\, \zeta_{n}^{(m)}+\mathbb{P}(\mathcal{E}^c)+\mathbb{P}\{S_{n,y}^{(m)}\in U^{R}\backslash U\}+\mathbb{P}\{S_{n,y}^{(m)}\in U\backslash U^{-R}\}\,.
\end{align*}
By the same arguments in Section \ref{sec:f.5} for deriving the bounds of $\mathbb{P}\{G^{(1)}\in U^{R}\backslash U\}$ and $\mathbb{P}\{G^{(1)}\in U\backslash U^{-R}\}$, it holds that 
\begin{align*}
&\sup_{U\in\mathcal{U}}\mathbb{P}\{S_{n,y}^{(m)}\in U^{R}\backslash U\}+\sup_{U\in\mathcal{U}}\mathbb{P}\{S_{n,y}^{(m)}\in U\backslash U^{-R}\}\\
&~~~~~\lesssim p^{3/4}\tilde{D}_{n}\lesssim p^{3/4}m^{-\alpha}\Phi_{\psi_{\nu},\alpha}(\log n)^{(1+2\nu)/2}\,.
\end{align*}
By Lemma \ref{lem:tail_prob_bound_m-dep_approx}, $
\mathbb{P}(\mathcal{E}^c)\lesssim 
B_n^{3/2}mp^{13/4}n^{-1/4}$. Hence, we have
$$
\zeta_{n}\lesssim B_{n}^{3/2}mp^{13/4}n^{-1/4}+p^{3/4}m^{-\alpha}\Phi_{\psi_{\nu},\alpha}(\log n)^{(1+2\nu)/2}+p^{2}\Psi_{2,\alpha}^{1/2} \Psi_{2,0}^{1/2}m^{-\alpha/2}$$ provided that $m\gg p^{1/\alpha}\Psi_{2,\alpha}^{1/\alpha} \Psi_{2,0}^{1/\alpha}$. Recall $p^4\Psi_{2,\alpha} \Psi_{2,0}=o\{n^{\alpha/(4+2\alpha)}\}$. With selecting $m\asymp n^{1/(4+2\alpha)}$, it holds that 
\begin{align*}
\zeta_{n}\lesssim n^{-\alpha/(8+4\alpha)}(B_{n}^{3/2}p^{13/4}+p^{2}\Psi_{2,\alpha}^{1/2} \Psi_{2,0}^{1/2})+n^{-\alpha/(4+2\alpha)}p^{3/4}\Phi_{\psi_{\nu},\alpha}(\log n)^{(1+2\nu)/2}\,.
\end{align*}
We complete the proof of part {\rm (ii)}  of Proposition {\rm\ref{pn:BE_bound_physical}}. $\hfill\Box$

\section{Proofs of Theorem \ref{tm:pb_rectangle}}\label{pf:pd_rectangle}

Part (i) of Theorem \ref{tm:pb_rectangle} can be derived immediately from Lemma \ref{lem:gaussian_comparison}. We now consider the proof of Part (ii). For any $A\in\mathcal{A}\subset\mathcal{A}^{\rm si}(a,d)$, there exsits a $K$-generated $A^{K}$ such that $A^{K}\subset A\subset A^{K,\epsilon}$ with $K\leq (pn)^{d}$ and $\epsilon=a/n$. Define
$
\bar{\rho}^{(K)}=|\mathbb{P}(\hat{G}\in A^{K}\,|\,\mathcal{X}_{n})-\mathbb{P}(G\in A^{K})|$ and $\bar{\rho}^{(K),\epsilon}=|\mathbb{P}(\hat{G}\in A^{K,\epsilon}\,|\,\mathcal{X}_{n})-\mathbb{P}(G\in A^{K,\epsilon})|$. Following the same arguments for deriving \eqref{eq:simple_upper_bound}, it holds that
$
|\mathbb{P}(\hat{G}\in A\,|\,\mathcal{X}_{n})-\mathbb{P}(G\in A)|\lesssim \epsilon(\log K)^{1/2}+\bar{\rho}^{(K)}+\bar{\rho}^{(K),\epsilon}$.
Together with Part (i) of Theorem \ref{tm:pb_rectangle}, repeating the arguments used in Section \ref{sec:pf_simple_convex} for bounding $\rho^{(K)}$ and $\rho^{(K),\epsilon}$, we have
$
\bar{\rho}^{(K)}+\bar{\rho}^{(K),\epsilon}\lesssim \Delta_{n}^{1/3}(\mathcal{A})(\log K)^{2/3}$. Due to $K\leq (pn)^{d}$, $\epsilon=a/n$ and $p\geq n^{\kappa}$, it holds that
$$
|\mathbb{P}(\hat{G}\in A\,|\,\mathcal{X}_{n})-\mathbb{P}(G\in A)|\lesssim an^{-1}(d\log p)^{1/2}+\Delta_{n}^{1/3}(\mathcal{A})(d\log p)^{2/3}$$ 
for any $A\in \mathcal{A}$. We complete the proof of Part (ii).

In the sequel, we consider the proof of Part (iii) which is similar to that of Theorem \ref{pn:GA_s-sparsely} stated in Section \ref{sec:pf_s-sparsely}. We only point out the differences here. Without loss of generality, we assume $\Delta_{n,r}\leq \min\{c(2s^{2})^{-1},1\}$ for some sufficiently small constant $c>0$. Otherwise, the assertions hold trivially. Let $B^{c}=\{w=(w_1,\ldots,w_p)^\T\in\mathbb{R}^{p}:\max_{j\in[p]}|w_{j}|> pn^{5/2}\}$.  Recall $\hat{G}=(\hat{G}_{1},\ldots,\hat{G}_{p})^{\T}$ and $G=({G}_{1},\ldots,{G}_{p})^{\T}$. As shown in Section \ref{sec:pf_s-sparsely}, $\mathbb{P}(G\in B^{c})\lesssim B_{n}n^{-1}$. Due to $\hat{G}\,|\,\mathcal{X}_n\sim N(0,\hat{\Xi}_{n})$, we have
\begin{align*}
\mathbb{E}\bigg(\max_{j\in [p]}|\hat{G}_{j}|\,\bigg|\,\mathcal{X}_n\bigg)\lesssim&\,\,(\log p)^{1/2}\max_{j\in [p]}\{\mathbb{E}(\hat{G}_{j}^{2}\,|\,\mathcal{X}_{n})\}^{1/2}\\
\lesssim&\,\, (\log p)^{1/2}\bigg[\max_{j\in [p]}\{\mathbb{E}(G_{j}^{2})\}^{1/2}+\Delta_{n,r}^{1/2}\bigg]\\
\lesssim&\,\, B_{n}(n\log p)^{1/2}+(\Delta_{n,r}\log p)^{1/2}\,,
\end{align*}
where the second step is due to $\max_{j\in [p]}|\mathbb{E}(\hat{G}_{j}^{2}\,|\,\mathcal{X}_{n})-\mathbb{E}(G_{j}^{2})|\leq \Delta_{n,r}$ and the third step is due to $\mathbb{E}(G_{j}^2)\lesssim B_{n}^2 n$. Then
\begin{align*}
\mathbb{P}(\hat{G}\in B^{c}\,|\,\mathcal{X}_{n})=&\,\mathbb{P}\bigg(\max_{j\in [p]}|\hat{G}_{j}|> pn^{5/2}\,\bigg|\,\mathcal{X}_{n}\bigg)\\
\leq&\, p^{-1}n^{-5/2}\mathbb{E}\bigg(\max_{j\in [p]}|\hat{G}_{j}|\,\bigg|\,\mathcal{X}_{n}\bigg)\lesssim B_{n}n^{-2}+n^{-5/2}\Delta_{n,r}^{1/2}\,.
\end{align*}
Hence, 
$$
|\mathbb{P}(\hat{G}\in A\cap B^{c}\,|\,\mathcal{X}_n)-\mathbb{P}(G\in A\cap B^{c})|
\lesssim B_{n}n^{-1}+n^{-5/2}\Delta_{n,r}^{1/2}\,.$$ 
Let $A^{*}=A\cap B=\cap_{q=1}^{K_{*}}A_{q}^{*}$. We next bound $|\mathbb{P}(\hat{G}\in A^{*}\,|\,\mathcal{X}_n)-\mathbb{P}(G\in A^{*})|$ in three cases.

{\bf Case 1.} ({\it $A^{*}$ contains a ball with radius $\epsilon=n^{-1}$}) As shown in Section \ref{sec:pf_s-sparsely}, $A^{*}$ must be a simple convex set with $a=1$ and $d=Cs^2$ for some positive constant $C$, and the corresponding $K$-generated $A^{K}$ satisfies $\max_{v\in\mathcal{V}(A^{K})}|v|_{0}\leq s$. Hence, $\max_{v_{1},v_{2}\in\mathcal{V}(A^{K})}|v_{1}(\hat{\Xi}_{n}-\Xi)v_{2}|\leq s^{2}\Delta_{n,r}$. By Part (ii), we have
$|\mathbb{P}(\hat{G}\in A^{*}\,|\,\mathcal{X}_n)-\mathbb{P}(G\in A^{*})| \lesssim sn^{-1}(\log p)^{1/2}+s^{2}\Delta_{n,r}^{1/3}(\log p)^{2/3}$.

{\bf Case 2.} ({\it $A^{*}=\cap_{q=1}^{K_{*}}A_{q}^*$ doesn't contain a ball with radius $\epsilon=n^{-1}$, and there exists one $A_{q}^{*}$ containing no ball with radius $\epsilon$}) As shown in Section \ref{sec:pf_s-sparsely}, $\mathbb{P}(G\in A^{*})\lesssim s^{1/4}n^{-1}$. Due to $|v^{\T}(\hat{\Xi}_{n}-\Xi)v|\leq s^{2}\Delta_{n,r}$ for any $v\in\mathbb{S}^{p-1}$ with $|v|_{0}\leq s$, we have
$v^{\T}\hat{\Xi}_{n}v\geq v^{\T}\Xi v-s^{2}\Delta_{n,r}\geq C$. Analogously, we also have  $\mathbb{P}(\hat{G}\in A^{*}\,|\,\mathcal{X}_{n})\lesssim s^{1/4}n^{-1}$. It then holds that
$|\mathbb{P}(\hat{G}\in A^{*}\,|\,\mathcal{X}_{n})-\mathbb{P}(G\in A^{*})| \lesssim s^{1/4}n^{-1}$.

{\bf Case 3.} ({\it each $A_{q}^{*}$ contains a ball with radius $\epsilon=n^{-1}$, but they are not all the same}) As shown in Section \ref{sec:pf_s-sparsely}, $
\mathbb{P}(G\in A^{0})\lesssim  sn^{-1}(\log p)^{1/2}$.
Notice that $\var(v^{\T}\hat{G}\,|\,\mathcal{X}_n)=v^{\T}\hat{\Xi}_{n}v\geq C$ holds for all $v\in\mathbb{S}^{p-1}$ with $|v|_{0}\leq s$. Analogously,
we also have $\mathbb{P}(\hat{G}\in A^{0}\,|\,\mathcal{X}_{n})\lesssim sn^{-1}(\log p)^{1/2}$. Since $A^{*}\subset A^{0}$, then $|\mathbb{P}(\hat{G}\in A^{*}\,|\,\mathcal{X}_{n})-\mathbb{P}(G\in A^{*})| \lesssim s n^{-1}(\log p)^{1/2}$. $\hfill\Box$

\section{Proof of Theorem {\rm\ref{tm:cov_compare}}}\label{sec:cov_compare}

To prove Theorem \ref{tm:cov_compare}, we need the following lemma.

\begin{lemma}
	\label{lem:tail_prob_bound_m-dep_approx_2}
	Let  $\{X_{t}\}$ be a sequence of $p$-dimentional centered random vectors generated from the model {\rm(\ref{eqn:weakly_dependent_ts_fdp})} such that $\Phi_{\psi_{\nu},0} < \infty$ for some $\nu \in (0, \infty)$. Then there exists a universal constant $C > 0$ depending only on $\nu$ such that
	$
	\max_{j\in[p]}\mathbb{P}(|n^{-1/2}\sum_{t=1}^{n}X_{t,j}| > u) \leq C \exp\{-(4e)^{-1}(1+2\nu)(u\Phi_{\psi_{\nu,0}}^{-1})^{2/(1+2\nu)}\}$
	for any $u>0$.
\end{lemma}

The proof of Lemma \ref{lem:tail_prob_bound_m-dep_approx_2} essentially follows from the arguments in proving Lemma C.4 of \cite{ZhangWu_2017} with the necessary modification using the uniform functional dependence measure to non-stationarity of the sequence $\{X_{t}\}_{t=1}^n$. Details are omitted. 

Define
$\Xi^* = \sum_{j=-n+1}^{n-1} \mathcal K(j/b_n)H_j$, where $H_{j}=n^{-1}\sum_{t=j+1}^{n}\mathbb{E}(X_{t}X_{t-j}^{\T})$ if $j\geq 0$ and $H_{j}=n^{-1}\sum_{t=-j+1}^{n}\mathbb{E}(X_{t+j}X_{t}^{\T})$ otherwise.
By triangle inequality, we have
\begin{align}\label{eq:cov_compare}
\Delta_{n,r}=|\hat{\Xi}_{n}-\Xi|_{\infty}&\leq |\hat{\Xi}_{n}-\Xi^{*}|_{\infty}+|\Xi^{*}-\Xi|_{\infty}\,.
\end{align}
To bound $\Delta_{n,r}$, it suffices to bound $|\hat{\Xi}_{n}-\Xi^{*}|_{\infty}$ and $|\Xi^{*}-\Xi|_{\infty}$, respectively.

\underline{\it Proofs of Parts {\rm(i)} and {\rm(ii)}.} As suggested in \cite{ChangHuLiuTang_2021}, we can formulate the $m$-dependent sequence and the $\alpha$-mixing sequence with $\alpha$-mixing coefficients $\{\alpha_n(k)\}_{k\geq1}$ in a unified framework. More specifically, let $\{X_t\}_{t=1}^n$ be an $\alpha$-mixing sequence of real-valued and centered random variables with $\alpha$-mixing coefficients $\{\tilde{\alpha}_n(k)\}_{k\geq1}$ satisfying $
\tilde{\alpha}_n(k)\leq a_1\exp(-a_2{L}_{{n}}^{-r_2}k^{r_2})$
for any integer $k\geq 1$, where ${L}_{{n}}>0$ may diverge with ${n}$. If $\{X_t\}_{t=1}^n$ satisfies Condition \ref{as:alpha-mixing}, we can set $L_n=1$ and $r_2=\gamma_2$. If $\{X_t\}_{t=1}^n$ is an $m$-dependent sequence, we can set $L_n=m$ and $r_2=\infty$. We assume $L_n(\log p)^{(2r_2+1)/r_2}=o(n)$. In the sequel, we will use this general setting of $\{X_t\}_{t=1}^n$ to show Parts (i) and (ii).

We first consider $|\Xi^{*}-\Xi|_{\infty}$. Write $\Xi=(\sigma_{\ell_1,\ell_2})_{p\times p}$ and $ H_{j}=(H_{j,\ell_{1},\ell_{2}})_{p\times p}$. For any $\ell_{1},\ell_{2}\in [p]$, we have
$
\sigma_{\ell_{1},\ell_{2}}=H_{0,\ell_1,\ell_2} +  \sum_{j=1}^{n-1} H_{-j,\ell_{1},\ell_{2}} +  \sum_{j=1}^{n-1} H_{j,\ell_{1},\ell_{2}}$. 
By triangle inequality and Davydov's inequality, $|H_{j,\ell_{1},\ell_{2}}| \leqslant n^{-1}\sum_{t=j+1}^{ n}|\mathbb{E}(X_{t,\ell_1}X_{t-j,\ell_2})| \lesssim  B_{n}^{2}\exp(-CL_{n}^{-r_2}j^{r_2})$ for any $j\geqslant 1$.
This bound also holds for $|H_{-j,\ell_1,\ell_2}|$ with $j\geqslant 1$. Write $\Xi^{*}=(\sigma_{\ell_{1},\ell_{2}}^*)_{p\times p}$, where $\sigma_{\ell_{1},\ell_{2}}^*= \sum_{j=-n+1}^{n-1} \mathcal K(j/b_n) H_{j,\ell_{1},\ell_{2}}$. By triangle inequality and Condition \ref{as:kernel}, 
\begin{align}\label{eq:covbd1}
|\Xi^{*}-\Xi|_{\infty}
\leqslant&\, \max_{\ell_1,\ell_2\in[p]}\sum_{j=1}^{n-1}\bigg|\mathcal K\bigg(\frac{j}{b_n}\bigg)-1\bigg|\big(|H_{j,\ell_{1},\ell_{2}}|
+|H_{-j,\ell_{1},\ell_{2}}|\big)
\lesssim n^{-\rho}B_{n}^{2}L_{n}^2 \,.
\end{align}

Next, we consider the convergence rate of $|\hat{\Xi}_{n}-\Xi^{*}|_{\infty}$. Recall $\hat{\Xi}_{n}=\sum_{j=-n+1}^{n-1}\mathcal{K}(j/b_{n})\hat{{H}}_{j}$ with $\hat{H}_{j}=n^{-1}\sum_{t=j+1}^{n}(X_{t}-\bar{X})(X_{t-j}-\bar{X})^{\T}$ if $j\geq 0$ and $\hat{{H}}_{j}=n^{-1}\sum_{t=-j+1}^{n}(X_{t+j}-\bar{X})(X_{t}-\bar{X})^{\T}$ otherwise. By triangle inequality, 
\begin{align}
\bigg|\sum_{j=0}^{n-1}\mathcal K\bigg(\frac{j}{b_n}\bigg)
(\hat{{H}}_j-{H}_j)\bigg|_\infty\leq&\,I_1+I_2+I_3+I_4\,,\label{eq:lrcovexp}
\end{align}
where 
\begin{align*}
I_{1}=&\,\,\bigg|\sum_{j=0}^{n-1}\mathcal K\bigg(\frac{j}{b_n}\bigg)\bigg[\frac{1}{n}\sum_{t=j+1}^{n}\{X_tX_{t-j}^\T - \mathbb{E}(X_tX_{t-j}^\T)\}\bigg]\bigg|_\infty\,,\\
I_2=&\,\,\bigg|\sum_{j=0}^{ n-1}\mathcal K\bigg(\frac{j}{b_n}\bigg)
	\bigg(\frac{1}{n}\sum_{t=j+1}^{n}X_t\bigg)\bar{X}^{\T}\bigg|_\infty\,,\\
	I_3=&\,\,\bigg|\sum_{j=0}^{ n-1}\mathcal K\bigg(\frac{j}{b_n}\bigg)
	\bar{X}\bigg(\frac{1}{n}\sum_{t=j+1}^{n}X_{t-j}\bigg)^{\T}\bigg|_\infty\,,~I_4=\bigg|\sum_{j=0}^{n-1}\frac{n-j}{n}\mathcal K\bigg(\frac{j}{b_n}\bigg)\bar{X}\bar{X}^{\T}\bigg|_\infty\,.
	\end{align*}
We will specify the convergence rates of $I_1, I_2, I_3$ and $I_4$, respectively. Given $M>0$ such that $b_n\ll M\ll n$, for any $u>0$, we have
\begin{align}\label{eq:I1}
\mathbb{P}(I_{1}>u)\leq&\,\,\mathbb{P}\bigg(\bigg|\sum_{j=0}^{M}\mathcal K\bigg(\frac{j}{b_n}\bigg) \bigg[\frac{1}{n}\sum_{t=j+1}^{n}\{X_tX_{t-j}^\T - \mathbb{E}(X_tX_{t-j}^\T)\}\bigg] \bigg|_\infty>\frac{u}{2}\bigg)\notag\\
&+\mathbb{P}\bigg(\bigg|\sum_{j=M+1}^{n-1}\mathcal K\bigg(\frac{j}{b_n}\bigg) \bigg[\frac{1}{n}\sum_{t=j+1}^{n}\{X_tX_{t-j}^\T - \mathbb{E}(X_tX_{t-j}^\T)\}\bigg] \bigg|_\infty>\frac{u}{2}\bigg)\\
=&:\,T_1+T_2\,.\notag
\end{align}
Write $Z_{t,\ell_{1},\ell_{2}}^{(j)}=X_{t,\ell_{1}}X_{t-j,\ell_{2}}-\mathbb{E}(X_{t,\ell_{1}}X_{t-j,\ell_{2}})$. For given $(j, \ell_{1},\ell_{2})$, $\{Z_{t,\ell_{1},\ell_{2}}^{(j)}\}$ is an $\alpha$-mixing sequence with $\alpha$-coefficients $\check{\alpha}_{n,\ell_1,\ell_2}^{(j)}(k)\leq 
a_1\exp(-a_{2}L_{n}^{-r_2}|k-j|_+^{r_2})$. By Condition \ref{as:tail} and Lemma 2 of \cite{CTW_2013}, $\max_{j,t,\ell_{1},\ell_{2}}	\mathbb{P}\{|Z_{t,\ell_{1},\ell_{2}}^{(j)}|>u\}\lesssim\exp(-Cu^{\gamma_{1}/2}B_{n}^{-\gamma_{1}})$. By Lemma \ref{la:large_deviation} with $(\tilde{B}_{\tilde{n}},\tilde{L}_{\tilde{n}},\tilde{j}_{\tilde{n}},r_1)=(B_n^2,L_n,j,\gamma_1/2)$, \begin{align*}
    \mathbb{P}\bigg\{\bigg|\frac{1}{n}\sum_{t=j+1}^{n}Z_{t,\ell_{1},\ell_{2}}^{(j)}\bigg|\geq x\bigg\}\lesssim&\,\,\exp\{-C(M+L_{n})^{-1}B_n^{-4}nx^2\}\\
    &+\exp(-C(M+L_{n})^{-\gamma_*}n^{\gamma_*}x^{\gamma_*}B_n^{-2\gamma_*})
    \end{align*}
    for any $j\in\{0\}\cup [M]$ and $\ell_{1},\ell_2\in [p]$, where $\gamma_{*}=(2+|2\gamma_1^{-1}-1|_{+}+r_2^{-1})^{-1}$. Since $\sum_{j=0}^{M}|\mathcal{K}(j/b_{n})|\lesssim b_{n}$ and $b_n\asymp n^{\rho}$, we have
    \begin{align*}
T_1\leq&\,\,\sum_{j=0}^M\sum_{\ell_1,\ell_2=1}^p\mathbb{P}\bigg\{\bigg| \frac{1}{n}\sum_{t=j+1}^{n}Z_{t,\ell_{1},\ell_{2}}^{(j)} \bigg|>Cub_{n}^{-1}\bigg\}\\
\lesssim&\,\, Mp^2\bigg[\exp\bigg\{-\frac{Cn^{1-2\rho}u^2}{B_n^{4}(M+L_{n})}\bigg\}+\exp\bigg\{-\frac{Cn^{\gamma_*(1-\rho)}u^{\gamma_*}}{(M+L_{n})^{\gamma_*}B_n^{2\gamma_*}}\bigg\}\bigg]\,.
\end{align*}
Since $\sum_{j=M+1}^{n-1}|\mathcal{K}(j/b_{n})|\lesssim n^{\vartheta\rho}M^{1-\vartheta}$, we have 
\begin{align*}
T_2
\leq&\,\,\sum_{j=M+1}^{n-1}\sum_{\ell_{1},\ell_{2}=1}^{p}\sum_{t=j+1}^{n}\mathbb{P}\{|Z_{t,\ell_{1},\ell_{2}}^{(j)}|>Cn^{-\vartheta\rho}M^{\vartheta-1}u\}\\
\lesssim&\,\, n^2p^2\exp\{-Cn^{-\vartheta\rho\gamma_{1}/2}B_{n}^{-\gamma_{1}}M^{(\vartheta-1)\gamma_{1}/2}u^{\gamma_{1}/2}\}\,.
\end{align*}
Note that $p\geq n^{\kappa}$ for some $\kappa>0$. By \eqref{eq:I1}, it holds that
\begin{align*}
I_{1}=&\,\,O_{\rm p}\{n^{-(1-2\rho)/2}B_{n}^2(M+L_{n})^{1/2}(\log p)^{1/2}\}+O_{\rm p}\{n^{\rho-1}B_{n}^2(M+L_{n})(\log p)^{1/\gamma_{*}}\}\\
&+O_{\rm p}\{M^{1-\vartheta}B_{n}^2n^{\vartheta\rho}(\log p)^{2/\gamma_{1}}\}\,.
\end{align*}
Applying  Lemma \ref{la:large_deviation} to the sequence $\{X_{t,\ell}\}$ with $(\tilde{B}_{\tilde{n}},\tilde{L}_{\tilde{n}},\tilde{j}_{\tilde{n}},r_1)=(B_n,L_n,0,\gamma_1)$, we have  $\max_{j,\ell}\mathbb{P}(
|n^{-1}\sum_{t=j+1}^{n}X_{t,\ell} |\geq u)\lesssim \exp\{-C(L_{n}^{-1}B_n^{-1}nu)^{\gamma_{**}}\}+\exp(-CL_{n}^{-1}nu^2B_n^{-2})$ with $\gamma_{**}=r_2/(2r_2+1)$. Then 
\begin{align*}
&\max_{j}\mathbb{P}\bigg\{
\bigg|\bigg(\frac{1}{n}\sum_{t=j+1}^{n}X_{t}\bigg)\bar{X}^{\T}\bigg|_\infty>u\bigg\}\lesssim p^2\max_{j,\ell}\mathbb{P}\bigg(
\bigg|\frac{1}{n}\sum_{t=j+1}^{n}X_{t,\ell}\bigg|>\sqrt{u}\bigg)\\
&~~~~\lesssim p^2\exp\{-C(L_{n}^{-1}B_{n}^{-1}n\sqrt{u})^{\gamma_{**}}\}+p^2\exp(-CL_{n}^{-1}B_n^{-2}nu)\,.
\end{align*}
Thus, 
\begin{align*}
\mathbb{P}(I_{2}> u)\leq&\,\,\sum_{j=0}^{n-1}\mathbb{P}\bigg\{
\bigg|\bigg(\frac{1}{n}\sum_{t=j+1}^{n}X_t\bigg)\bar{X}^{\T}\bigg|_\infty>Cub_{n}^{-1}\bigg\}\\
\lesssim&\,\, np^{2}\exp(-CL_{n}^{-1}B_n^{-2}n^{1-\rho}u)+np^{2}\exp\{-C(L_n^{-1}B_n^{-1}n^{(2-\rho)/2}u^{1/2})^{\gamma_{**}}\}\,,
\end{align*}
which implies that 
\begin{align*}
I_{2}=&\,\,O_{\rm p}\{n^{\rho-2}B_n^2L_{n}^2(\log p)^{2/\gamma_{**}}\}+O_{\rm p}(n^{\rho-1}B_{n}^2L_{n}\log p)\\
=&\,\,o_{\rm p}\{n^{(2\rho-1)/2}B_{n}^2L_{n}^{1/2}(\log p)^{1/2}\}+o_{\rm p}\{n^{\rho-1}B_{n}^2L_{n}(\log p)^{1/\gamma_{*}}\}\,,
\end{align*}
where the second step is due to $L_n(\log p)^{(2r_2+1)/r_2}=o(n)$. Similarly, 
$$
I_{3}=O_{\rm p}\{n^{\rho-2}B_n^2L_{n}^2(\log p)^{2/\gamma_{*}}\}+O_{\rm p}(n^{\rho-1}B_{n}^2L_{n}\log p)=I_4\,.$$ 
Hence, 
\begin{align*}
\bigg|\sum_{j=0}^{n-1}\mathcal K\bigg(\frac{j}{b_n}\bigg)
(\hat{H}_j-H_j)\bigg|_\infty=&\,\,O_{\rm p}\bigg\{\frac{B_{n}^2(M+L_{n})(\log p)^{1/\gamma_{*}}}{n^{1-\rho}}\bigg\}+O_{\rm p}\bigg\{\frac{B_{n}^2n^{\vartheta\rho}(\log p)^{2/\gamma_{1}}}{M^{\vartheta-1}}\bigg\}\\
&+O_{\rm p}\bigg\{\frac{B_{n}^2(M+L_{n})^{1/2}(\log p)^{1/2}}{n^{(1-2\rho)/2}}\bigg\}\,.
\end{align*}
Such convergence rate also holds for $|\sum_{j=-n+1}^{-1}\mathcal K(j/b_n)
(\hat{H}_j-H_j)|_\infty$. Together with \eqref{eq:covbd1}, by \eqref{eq:cov_compare}, it holds that 
\begin{align*}
\Delta_{n,r}=&\,\,O(n^{-\rho}B_{n}^{2}L_{n}^2)+O_{\rm p}\{n^{(2\rho-1)/2}B_{n}^2(M+L_{n})^{1/2}(\log p)^{1/2}\}\\
&+O_{\rm p}\{n^{\rho-1}B_{n}^2(M+L_{n})(\log p)^{1/\gamma_{*}}\}+O_{\rm p}\{M^{1-\vartheta}B_{n}^2n^{\vartheta\rho}(\log p)^{2/\gamma_{1}}\}\,.
\end{align*}
Now, we simplify the convergence rate of $\Delta_{n,r}$ under two cases by selecting optimal $M$. 

When $\{X_t\}_{t=1}^n$ satisfies Condition \ref{as:alpha-mixing}, we can set $L_n=1$ and $r_2=\gamma_2$. Correspondingly, we have
$
\Delta_{n,r}=O(n^{-\rho}B_{n}^{2})+O_{\rm p}\{M^{1-\vartheta}B_{n}^2n^{\vartheta\rho}(\log p)^{2/\gamma_{1}}\}+O_{\rm p}\{n^{(2\rho-1)/2}B_{n}^2M^{1/2}(\log p)^{1/2}\}+O_{\rm p}\{n^{\rho-1}B_{n}^2M(\log p)^{(2\gamma_2+1)/\gamma_2}\}$.
Select $M\asymp n^{(1-2\rho+2\vartheta\rho)/(2\vartheta-1)}(\log p)^{(4-\gamma_1)/\{\gamma_{1}(2\vartheta-1)\}} $. 
Due to $\rho<1<\vartheta$, we have $\rho<(1-2\rho+2\vartheta\rho)/(2\vartheta-1)<1$ which implies there exists a constant $\tilde{c}>0$ depending only on $(\gamma_{1},\rho,\vartheta)$ such that $b_n\ll M\ll n$ holds for such selected $M$ when $\log p=o(n^{\tilde{c}})$. Furthermore, if we restrict $0<\rho<(\vartheta-1)/(3\vartheta-2)$, we have $2\rho+\vartheta-3\vartheta\rho-1>0$ and $3\rho+2\vartheta-4\vartheta\rho-2>0$. Therefore, there exist two constants $c_1>0$ depending only on $(\rho,\vartheta)$ and $c_2>0$ depending only on $(\gamma_1,\gamma_2,\vartheta)$ such that $\Delta_{n,r}=O_{\rm p}\{B_n^2n^{-c_1}(\log p)^{c_2}\}+O(B_n^2n^{-\rho})$. We complete the proof of Part (i).

When $\{X_t\}_{t=1}^n$ is $m$-dependent, we can set $L_n=m$ and $r_2=\infty$. In this case, to make $\Delta_{n,r}=o_{\rm p}(1)$, we need to require $B_n^2m^2=o(n^{\rho})$. Due to $M\gg b_n\asymp n^\rho$, we have
$
\Delta_{n,r}=O(n^{-\rho}B_{n}^{2}m^2)+O_{\rm p}\{n^{(2\rho-1)/2}B_{n}^2M^{1/2}(\log p)^{1/2}\}+O_{\rm p}\{n^{\rho-1}B_{n}^2M(\log p)^{2}\}+O_{\rm p}\{M^{1-\vartheta}B_{n}^2n^{\vartheta\rho}(\log p)^{2/\gamma_{1}}\}$. 
Select $M\asymp n^{(1-2\rho+2\vartheta\rho)/(2\vartheta-1)}(\log p)^{(4-\gamma_{1})/\{\gamma_{1}(2\vartheta-1)\}}$. 
Following the same arguments for $\alpha$-mixing sequence $\{X_t\}$, we know there exist two constants $c_1>0$ depending only on $(\rho,\vartheta)$ and $c_2>0$ depending only on $(\gamma_1,\vartheta)$ such that $\Delta_{n,r}=O_{\rm p}\{B_n^2n^{-c_1}(\log p)^{c_2}\}+O(B_n^2m^2n^{-\rho})$. We complete the proof of Part (ii).

\underline{{\it Proof of Part {\rm(iii)}}.} Write $\calF_{t} = \sigma(\varepsilon_{t},\varepsilon_{t-1},\dots)$ and let $\calP_{t}(\cdot) = \mathbb{E}(\cdot\,|\,\calF_{t}) - \mathbb{E}(\cdot\,|\,\calF_{t-1})$ be a projection operator. As shown in Section \ref{se:pflemmoment_bound_cov_mat_compare} of the supplementary material, $\|\calP_{t-l}(X_{t,\ell})\|_{q} \leq \theta_{l,q,\ell}$ for any $\ell\in[p]$. Since $X_{t} = \sum_{l=0}^\infty \calP_{t-l}(X_{t})$, by triangle inequality and Cauchy-Schwarz inequality, for any $j\geq1$ and $t\geq j+1$, 
\begin{align*}
|\mathbb{E}(X_{t,\ell_{1}}X_{t-j,\ell_{2}})|=&\,\,\bigg|\sum_{l=0}^{\infty}\mathbb{E}\{\calP_{t-l}(X_{t,\ell_1})\calP_{t-l}(X_{t-j,\ell_2})\}\bigg|\leq\sum_{l=0}^{\infty}|\mathbb{E}\{\calP_{t-l}(X_{t,\ell_1})\calP_{t-l}(X_{t-j,\ell_2})\}|\\
\leq&\,\,\sum_{l=0}^{\infty}\|\calP_{t-l}(X_{t,\ell_1})\|_2\|\calP_{t-l}(X_{t-j,\ell_2})\|_2\leq\sum_{l=0}^{\infty}\theta_{l,2,\ell_1}\theta_{l-j,2,\ell_2}\,,
\end{align*}
which implies $|H_{j,\ell_{1},\ell_{2}}| \leqslant n^{-1}\sum_{t=j+1}^{ n}|\mathbb{E}(X_{t,\ell_1}X_{t-j,\ell_2})| \leq  \sum_{l=0}^{\infty}\theta_{l,2,\ell_1}\theta_{l-j,2,\ell_2}$ 
for any $j\geqslant 1$.
Analogously, $|H_{-j,\ell_1,\ell_2}|\leq  \sum_{l=0}^{\infty}\theta_{l,2,\ell_2}\theta_{l-j,2,\ell_1}$ for any $j\geqslant 1$. Based on the first inequality given in \eqref{eq:covbd1}, we have
$$
|\Xi^{*}-\Xi|_{\infty}\lesssim b_n^{-1}\max_{\ell_1,\ell_2\in[p]}\sum_{j=1}^{n-1}j\sum_{l=0}^{\infty}\theta_{l,2,\ell_1}\theta_{l-j,2,\ell_2}\,.$$
Due to $\theta_{j,2,\ell}=0$ for any $j<0$, 
$
\sum_{j=1}^{n-1}\sum_{l=0}^{\infty}j\theta_{l,2,\ell_1}\theta_{l-j,2,\ell_2}=\sum_{j=1}^{n-1}\sum_{l=0}^{\infty}j\theta_{l+j,2,\ell_1}\theta_{l,2,\ell_2}=\sum_{l=0}^{\infty}(\theta_{l,2,\ell_2}\sum_{j=1}^{n-1}\sum_{k=j}^{n-1} \theta_{l+k,2,\ell_1})\leq\Theta_{0,2,\ell_2}\sum_{j=1}^{n-1} \Theta_{j,2,\ell_1}$.
Furthermore, let $\varpi_n=(\log n)I(\alpha=1)+n^{1-\alpha}I(\alpha\neq 1)$. Since $\Theta_{0,2,\ell_2}\leq \Psi_{2,0}$ and $\Theta_{j,2,\ell_1}\leq \Psi_{2,\alpha}j^{-\alpha}$, we have $$\max_{\ell_{1},\ell_{2}\in [p]}\sum_{j=1}^{n-1}\sum_{l=0}^{\infty}j\theta_{l,2,\ell_1}\theta_{l-j,2,\ell_2}\lesssim \Psi_{2,0}\Psi_{2,\alpha}\varpi_{n}\,.$$ Notice that $b_n\asymp n^{\rho}$. Then
$|\Xi^{*}-\Xi|_{\infty}\lesssim n^{-\rho}\Psi_{2,0}\Psi_{2,\alpha}\varpi_{n}$.

Now we specify the convergence rates of $I_1, I_2, I_3$ and $I_4$ in \eqref{eq:lrcovexp}.
Recall $Z_{t,\ell_{1},\ell_{2}}^{(j)}=X_{t,\ell_{1}}X_{t-j,\ell_{2}}-\mathbb{E}(X_{t,\ell_{1}}X_{t-j,\ell_{2}})$. We define a $p^2$-dimensional vector $Z_{t}^{(j)}$ with components being $Z_{t,\ell_{1},\ell_{2}}^{(j)}$'s. Such defined $Z_{t}^{(j)}$ satisfies model \eqref{eqn:weakly_dependent_ts_fdp} with some jointly measurable function $f_t^{(j)}(\cdot)$. Let $Z_{t,\{m\}}'^{(j)}$ be the coupled version of $Z_{t}^{(j)}$ at the time lag $m$ with $\varepsilon_{t-m}$ replacing by 
$\varepsilon_{t-m}'$. By H\"older's inequality and Condition \ref{as:tail}, 
\begin{align*}
&\|Z_{t,\ell_1,\ell_2}^{(j)}-Z_{t,\ell_1,\ell_2,\{m\}}'^{(j)}\|_q=\|X_{t,\ell_{1}}X_{t-j,\ell_{2}}-X_{t,\ell_{1},\{m\}}X_{t-j,\ell_{2},\{m-j\}}\|_q\\
&~~~\leq\|X_{t-j,\ell_{2}}\|_{q}\|X_{t,\ell_{1}}-X_{t,\ell_{1},\{m\}}\|_q+\|X_{t,\ell_{1},\{m\}}\|_q\|X_{t-j,\ell_{2}}-X_{t-j,\ell_{2},\{m-j\}}\|_q\\
&~~~\lesssim B_n(\theta_{m,q,\ell_1}+\theta_{m-j,q,\ell_2})
\end{align*}
for any $m\geq j$.  
We define $\Phi_{\psi_{\nu},\alpha}^{(j)}$ in the same manner as $\Phi_{\psi_{\nu},\alpha}$ in \eqref{eq:aggnorm} by replacing $\{X_{t}\}$ with $\{Z_{t}^{(j)}\}$. Then $\Phi_{\psi_{\nu},\alpha}^{(j)}\lesssim B_nj^{\alpha}\Phi_{\psi_{\nu},\alpha}$. By Lemma \ref{lem:tail_prob_bound_m-dep_approx_2}, for $M>0$ such that $b_n\ll M\ll n$, 
$\max_{j,\ell_1,\ell_2}\mathbb{P}\{|(n-j)^{-1/2}\sum_{t=j+1}^{n}Z_{t,\ell_1,\ell_2}^{(j)}| > u\} \lesssim \exp\{-C(uM^{-\alpha}B_n^{-1}\Phi_{\psi_{\nu,0}}^{-1})^{2/(1+2\nu)}\}$. 
Since $\sum_{j=0}^{M}|\mathcal{K}(j/b_{n})|\lesssim b_{n}\asymp n^{\rho}$, $T_1$ specified in \eqref{eq:I1} in this setting satisfies 
\begin{align*}
T_1
\leq np^2\max_{j,\ell_1,\ell_2}\mathbb{P}\bigg\{\bigg|\frac{1}{n}\sum_{t=j+1}^{n}Z_{t,\ell_{1},\ell_{2}}^{(j)}\bigg|>Cub_{n}^{-1}\bigg\}\leq np^2\exp\bigg[-\bigg\{\frac{Cn^{(1-2\rho)/2}u}{M^{\alpha}B_n\Phi_{\psi_{\nu,0}}}\bigg\}^{2/(1+2\nu)}\bigg]\,.
\end{align*}
The upper bound for $T_2$ specified in \eqref{eq:I1} remains the same in this setting. Thus, by \eqref{eq:I1}, 
$
I_{1}=O_{\rm p}\{n^{(2\rho-1)/2}M^{\alpha}B_n\Phi_{\psi_{\nu,0}}(\log p)^{(1+2\nu)/2}\}+O_{\rm p}\{M^{1-\vartheta}B_{n}^2n^{\vartheta\rho}(\log p)^{2/\gamma_{1}}\}$. By Lemma \ref{lem:tail_prob_bound_m-dep_approx_2}, 
\begin{align*}
\max_{j,\ell}\mathbb{P}\bigg(\bigg|\frac{1}{\sqrt{n-j}}\sum_{t=j+1}^{n}X_{t,\ell}\bigg| > u\bigg) \lesssim \exp\{-C(u\Phi_{\psi_{\nu,0}}^{-1})^{2/(1+2\nu)}\}\,,
\end{align*}
which implies 
$
\max_{j}\mathbb{P}\{
|(n^{-1}\sum_{t=j+1}^{n}X_{t})\bar{X}^{\T}|_\infty>u\}\lesssim p^2\max_{j,\ell}\mathbb{P}(
|n^{-1}\sum_{t=j+1}^{n}X_{t,\ell}|>\sqrt{u})
\lesssim p^2\exp\{-C(\Phi_{\psi_{\nu,0}}^{-2}nu)^{1/(1+2\nu)}\}$.
Due to $\sum_{j=0}^{n-1}|\mathcal{K}(j/b_{n})|\lesssim b_{n}\asymp n^\rho$, we have
\begin{align*}
\mathbb{P}(I_{2}> u)\leq\sum_{j=0}^{n-1}\mathbb{P}\bigg\{
\bigg|\bigg(\frac{1}{n}\sum_{t=j+1}^{n}X_t\bigg)\bar{X}^{\T}\bigg|_\infty>Cub_{n}^{-1}\bigg\}\lesssim np^2\exp\{-C(\Phi_{\psi_{\nu,0}}^{-2}n^{1-\rho}u)^{1/(1+2\nu)}\}\,,
\end{align*}
which implies 
$
I_{2}=O_{\rm p}\{\Phi_{\psi_{\nu,0}}^2n^{\rho-1}(\log p)^{1+2\nu}\}$. Similarly, 
$
I_{3}=O_{\rm p}\{\Phi_{\psi_{\nu,0}}^2 n^{\rho-1}(\log p)^{1+2\nu}\}=I_4$. By \eqref{eq:lrcovexp}, 
\begin{align*}
\bigg|\sum_{j=0}^{n-1}\mathcal K\bigg(\frac{j}{b_n}\bigg)
(\hat{H}_j-H_j)\bigg|_\infty=&\,\,O_{\rm p}\bigg\{\frac{B_{n}^2n^{\vartheta\rho}(\log p)^{2/\gamma_{1}}}{M^{\vartheta-1}}\bigg\}+O_{\rm p}\bigg\{\frac{\Phi_{\psi_{\nu,0}}^2(\log p)^{1+2\nu}}{n^{1-\rho}}\bigg\}\\
&+O_{\rm p}\bigg\{\frac{M^{\alpha}B_n\Phi_{\psi_{\nu,0}}(\log p)^{(1+2\nu)/2}}{n^{(1-2\rho)/2}}\bigg\}\,.
\end{align*}
Same convergence rate also holds for $|\sum_{j=-n+1}^{-1}\mathcal K(j/b_n)
(\hat{H}_j-H_j)|_\infty$. By \eqref{eq:cov_compare}, 
\begin{align*}
\Delta_{n,r}=&\,\,O_{\rm p}\bigg\{\frac{B_{n}^2n^{\vartheta\rho}(\log p)^{2/\gamma_{1}}}{M^{\vartheta-1}}\bigg\}+O_{\rm p}\bigg\{\frac{\Phi_{\psi_{\nu,0}}^2(\log p)^{1+2\nu}}{n^{1-\rho}}\bigg\}\\
&+O_{\rm p}\bigg\{\frac{M^{\alpha}B_n\Phi_{\psi_{\nu,0}}(\log p)^{(1+2\nu)/2}}{n^{(1-2\rho)/2}}\bigg\}+O\bigg(\frac{\Psi_{2,0}\Psi_{2,\alpha}\varpi_{n}}{n^{\rho}}\bigg)\,.
\end{align*}
Select $M\asymp \{n^{(1-2\rho+2\vartheta\rho)/2}(\log p)^{(4-\gamma_{1}-2\gamma_{1}\nu)/(2\gamma_{1})}\}^{1/(\alpha+\vartheta-1)}$. 
Due to $\rho<\min\{(2\alpha+2\vartheta-3)/(2\vartheta-2),1/(2\alpha)\}$, we have $\rho<(1-2\rho+2\vartheta\rho)/\{2(\alpha+\vartheta-1)\}<1$ which implies there exists a constant $\tilde{c}>0$ depending only on $(\rho,\vartheta,\alpha)$ such that $b_n\ll M\ll n$ holds for such selected $M$ when $\log p=o(n^{\tilde{c}})$. 
Therefore, there exist two constants $c_1>0$ depending only on $(\rho,\vartheta,\alpha)$ and $c_2>0$ depending only on $(\alpha,\vartheta,\gamma_{1},\nu)$ such that $\Delta_{n,r}=O_{\rm p}\{(B_n\Phi_{\psi_{\nu,0}}+B_n^2+\Phi_{\psi_{\nu,0}}^2)n^{-c_1}(\log p)^{c_2}\}+O(n^{-\rho}\Psi_{2,0}\Psi_{2,\alpha}\varpi_{n})$. We complete the proof of Part (iii). $\hfill\Box$






\begin{thebibliography}{}


\bibitem[Andrews(1991)]{Andrews_1991}
Andrews, D. W. K. (1991). Heteroskedasticity and autocorrelation consistent covariance matrix
estimation. {\sl Econometrica} {\bf59} 817--858.


\bibitem[Baldi and Rinott(1989)]{BaldiRinott_1989}
Baldi, P. and Rinott, Y. (1989). On normal approximations of distributions in terms of dependency graphs. {\sl Ann. Probab.} {\bf 17} 1646--1650.


\bibitem[Bentkus(2003)]{Bentkus_2003}
Bentkus, V. (2003). On the dependence of the Berry-Esseen bound on dimension. {\sl J. Statist. Plann. Inference} {\bf113} 385--402.


\bibitem[Berkes, Liu and Wu(2014)]{BerkesLiuWu_2014}
Berkes, I., Liu, W. and Wu, W. B. (2014). Koml\'os-Major-Tusn\'ady approximation under dependence.
{\sl Ann. Probab.} {\bf42} 794--817.





\bibitem[Bhattacharya and Holmes(2010)]{BhattacharyaHolmes_2010}
Bhattacharya, R. N. and Holmes, S. (2010). An exposition of G\"otze’s estimation of the rate of convergence in the multivariate central limit theorem. {\sl arXiv:1003.4254}.

\bibitem[Bhattacharya and Rao(2010)]{BhattacharyaRao_2010}
Bhattacharya, R. N. and Rao, R. R. (2010). {\it Normal Approximation and Asymptotic Expansions}. Society for Industrial and Applied Mathematics.

\bibitem[Boussama, Fuchs and Stelzer(2011)]{BFS11}
Boussama, F., Fuchs, F. and Stelzer, R. (2011). Stationary and geometric ergodicity of BEKK multivariate GARCH models. {\sl Stochastic Process. Appl.} {\bf121} 2331--2360.

\bibitem[Bradley(2005)]{Bradley_2005}
Bradley, R. C.  (2005). Basic properties of strong mixing conditions. A survey and some open questions. {\sl Probab. Surv.} {\bf 2} 107--144.

\bibitem[Bradley(2007)]{Bradley_2007}
Bradley, R. C.  (2007). {\it Introduction to Strong Mixing Conditions}, Vols. 1, 2 and 3. Kendrick Press, Heber City (Utah).


\bibitem[Cai, Liu and Xia(2014)]{CaiLiuXia_2014}
Cai, T. T., Liu, W. and Xia, Y. (2014). Two-sample test of high dimensional means under dependence. {\sl J. R. Stat. Soc. Ser. B. Stat. Methodol.} {\bf76} 349--372.

\bibitem[Chang et al.(2023)]{ChangHuLiuTang_2021}
Chang, J., Hu, Q., Liu, C. and Tang, C. Y. (2023). Optimal covariance matrix estimation for high-dimensional noise in high-frequency data. {\sl J. Econometrics}, in press.

\bibitem[Chang, Jiang and Shao(2023)]{ChangJiangShao_2021}
Chang, J., Jiang, Q. and Shao, X. (2023). Testing the martingale difference hypothesis in high dimension. {\sl J. Econometrics}, in press.

\bibitem[Chang et al.(2018)]{ChangQiuYaoZou_2018}
Chang, J., Qiu, Y., Yao, Q. and Zou, T. (2018). Confidence regions for entries of a large precision matrix. {\sl J. Econometrics} {\bf206} 57--82.


\bibitem[Chang, Tang and Wu(2013)]{CTW_2013}
Chang, J., Tang, C. Y. and Wu, Y. (2013). Marginal empirical likelihood and sure independence
feature screening. {\sl Ann. Statist.} {\bf41} 2123--2148.


\bibitem[Chang, Yao and Zhou(2017)]{ChangYaoZhou_2017}
Chang, J., Yao, Q. and Zhou, W. (2017). Testing for high-dimensional white noise using maximum cross-correlations.
{\sl Biometrika} {\bf104} 111--127.


\bibitem[Chang et al.(2017a)]{ChangZhengZhouZhou_2017}
Chang, J., Zheng, C., Zhou, W. X. and Zhou, W. (2017). Simulation-based hypothesis testing of high dimensional means under covariance heterogeneity. {\sl Biometrics} {\bf73} 1300--1310.


\bibitem[Chang et al.(2017b)]{ChangZhouZhouWang_2017}
Chang, J., Zhou, W., Zhou, W. X. and Wang, L. (2017). Comparing large covariance matrices under weak conditions on the dependence structure and its application to gene clustering. {\sl Biometrics} {\bf73} 31--41.


\bibitem[Chen and Qin(2010)]{ChenQin_2010}
Chen, S. X. and Qin, Y. L. (2010). A two sample test for high dimensional data with applications to gene-set testing. {\sl Ann. Statist.} {\bf 38} 808--835.


\bibitem[Chen(2018)]{Chen_2018}
Chen, X. (2018). Gaussian and bootstrap approximations for high-dimensional U-statistics and their applications.
{\sl Ann. Statist.} {\bf46} 642--678.


\bibitem[Chen and Kato(2019)]{ChenKato_2019}
Chen, X. and Kato, K. (2019). Randomized incomplete $U$-statistics in high dimensions. {\sl Ann. Statist.} {\bf 47} 3127--3156.


\bibitem[Chen and Kato(2020)]{ChenKato_2020}
Chen, X. and Kato, K. (2020). Jackknife multiplier bootstrap: finite sample approximations to the $U$-process supremum with applications. {\sl Probab. Theory Related Fields} {\bf 176} 1097--1163.

\bibitem[Chen and Shao(2004)]{ChenShao_2004}
Chen, L. H. Y. and Shao, Q.-M. (2004). Normal approximation under local dependence. {\sl Ann. Probab.} {\bf32} 1985--2028.


\bibitem[Chernozhukov, Chetverikov and Kato(2013)]{CCK_2013}
Chernozhukov, V., Chetverikov, D. and Kato, K. (2013). Gaussian approximations and multiplier bootstrap for maxima of sums of high-dimensional random vectors. {\sl Ann. Statist.} {\bf41} 2786--2819.


\bibitem[Chernozhukov, Chetverikov and Kato(2015)]{CCK_2015}
Chernozhukov, V., Chetverikov, D. and Kato, K. (2015). Comparison and anti-concentration bounds for maxima of Gaussian random vectors. {\sl Probab. Theory Related Fields} {\bf162} 47--70.



\bibitem[Chernozhukov, Chetverikov and Kato(2017)]{CCK_2017}
Chernozhukov, V., Chetverikov, D. and Kato, K. (2017). Central limit theorems and bootstrap in high dimensions. {\sl Ann. Probab.} {\bf45} 2309--2352.


\bibitem[Chernozhukov, Chetverikov and Kato(2019)]{CCK_2014}
Chernozhukov, V., Chetverikov, D. and Kato, K. (2019). Inference on causal and structural parameters using many moment inequalities. {\sl Rev. Econ. Stud.} {\bf 86} 1867--1900.


\bibitem[Chernozhukov et al.(2022)]{CCKK_2019}
Chernozhukov, V., Chetverikov, D., Kato, K. and Koike, Y. (2022). Improved central limit theorems and bootstrap approximation in high dimensions. {\sl Ann. Statist.} {\bf50} 2562--2586.


\bibitem[Chernozhukov, Chetverikov and Koike(2023)]{CCK_2020}
Chernozhukov, V., Chetverikov, D. and Koike, Y. (2023). Nearly optimal central limit theorem and bootstrap approximations in high dimensions. {\sl Ann. Appl. Probab.}, in press.


\bibitem[Das and Lahiri(2021)]{DasLahiri_2020}
Das, D. and Lahiri, S. (2021). Central limit theorem in high dimensions : The optimal bound on dimension growth rate. {\sl Trans. Amer. Math. Soc.} {\bf374} 6991--7009.


\bibitem[Davydov(1968)]{Davydov_1968}
Davydov, Y. A. (1968). Convergence of distributions generated by stationary stochastic
processes. {\sl Theory Probab. Appl.} {\bf13} 691--696.



\bibitem[Den Haan and Levin(1997)]{DenHanLevin_1997}
Den Haan, W. J. and Levin, A. (1997). A practitioner's guide to robust covariance matrix estimation. {\sl Handbook of Statist.} {\bf15} 299--342.


\bibitem[Deng(2020)]{Deng_2020}
Deng, H. (2020). Slightly conservative bootstrap for maxima of sums. {\sl arXiv:2007.15877}.


\bibitem[Deng and Zhang(2020)]{DengZhang_2020}
Deng, H. and Zhang, C. H. (2020). Beyond Gaussian approximation: Bootstrap for maxima of sums of independent random vectors. {\sl Ann. Statist.} {\bf48} 3643--3671.


\bibitem[Doukhan, Massart and Rio(1994)]{DoukhanMasssartRio_1994}
Doukhan, P., Massart, P. and Rio, E. (1994). The functional central limit theorem for strongly mixing processes. {\sl Ann. Inst. Henri Poincar'e Probab. Stat.} {\bf30} 63--82.

\bibitem[Fan and Yao(2003)]{FanYao_2003}
Fan, J. and Yao, Q. (2003). {\it Nonlinear Time Series: Nonparametric and Parametric Methods}. Springer, New York.


\bibitem[Fang and Koike(2021)]{FangKoike_2020}
Fang, X. and Koike, Y. (2021). High-dimensional central limit theorems by stein's method. {\sl Ann. Appl. Probab.} {\bf31} 1660--1686.

\bibitem[Hafner and Preminger(2009)]{HP09}
Hafner, C. M. and Preminger, A. (2009). On asymptotic theory for multivariate GARCH models.  {\sl J. Multivariate Anal.} {\bf100} 2044--2054.


\bibitem[Hoeffding and Robbins(1948)]{HoeffdingRobbins_1948}
Hoeffding, W. and Robbins, H. (1948). The central limit theorem for dependent random variables. {\sl Duke Math. J.} {\bf15} 773--780.

\bibitem[H\"{o}rmann(2009)]{Hormann_2009}
H\"{o}rmann, S. (2009). Berry-Esseen bounds for econometric time series. {\sl ALEA Lat. Am. J. Probab. Math. Stat.} {\bf6} 377--397.


\bibitem[Jirak(2016)]{Jirak_2016}
Jirak, M. (2016). Berry-Esseen theorems under weak dependence. {\sl Ann. Probab.} {\bf44} 2024--2063.

\bibitem[Kiefer, Vogelsang and Bunzel(2000)]{Kieferetaj_2000}
Kiefer, N. M., Vogelsang, T. J. and Bunzel, H. (2000). Simple robust testing of regression hypothesis. {\sl Econometrica} {\bf68} 695--714.





\bibitem[Koike(2021)]{Koike_2019b}
Koike, Y. (2021). Notes on the dimension dependence in high-dimensional central limit theorems for hyperrectangles. {\sl Jpn. J. Stat. Data Sci.} {\bf4} 257--297.


\bibitem[Koike(2023)]{Koike_2019c}
Koike, Y. (2023). High-dimensional central limit theorems for homogeneous sums. {\sl J. Theoret. Probab.}, in press.


\bibitem[Kuchibhotla and Rinaldo(2020)]{KuchibhotlaRinaldo_2020}
Kuchibhotla, A. K. and Rinaldo, A. (2020). High-dimensional CLT for sums of non-degenerate random vectors: $n^{-1/2}$-rate. {\sl  arXiv:2009.13673}.


\bibitem[Lahiri(2003)]{Lahiri_2003}
Lahiri, S. N. (2003). {\it Resampling Methods for Dependent Data}. Springer, New York.


\bibitem[Lopes(2022)]{Lopes_2020}
Lopes, M. E. (2022). Central limit theorem and bootstrap approximation in high dimensions: Near $1/\sqrt{n}$ rates via implicit smoothing. {\sl Ann. Statist.} {\bf 50} 2492--2513.


\bibitem[Lopes, Lin and M\"uller(2020)]{LopesLinMueller_2020}
Lopes, M. E., Lin, Z. and M\"uller, H. G. (2020). Bootstrapping max statistics in high dimensions: Near-parametric rates under weak variance decay and application to functional and multinomial data. {\sl Ann. Statist.} {\bf48} 1214--1229.













\bibitem[Petrov(1995)]{Petrov_1995}
Petrov, V. V. (1995). {\it Limit Theorems of Probability Theory: Sequences of Independent Random Variables}. Clarendon Press, Oxford.

\bibitem[Rai\v{c}(2019)]{Raic_2019}
Rai\v{c}, M. (2019). A multivariate Berry-Esseen theorem with explicit constants. {\sl Bernoulli} {\bf 25} 2824--2853.



\bibitem[Rio(2009)]{Rio_2009}
Rio, E. (2009). Moment inequalities for sums of dependent random variables under projective conditions. {\sl J. Theoret. Probab.} {\bf22} 146--163.


\bibitem[Rio(2017)]{Rio_2013}
Rio, E. (2017). {\it Asymptotic Theory of Weakly Dependent Random Processes}. Springer, Berlin.



\bibitem[Rosenblatt(1956)]{Rosenblatt_1956}
Rosenblatt, M. (1956). A central limit theorem and a strong mixing condition. {\sl Proc. Natl. Acad. Sci. USA} {\bf42} 43--47.


\bibitem[Saulis and Statulevi\v{c}ius(1991)]{SS_1991}
Saulis, L. and Statulevi\v{c}ius, V. A. (1991). {\it Limit Theorems for Large Deviations}. Springer, Netherlands.



\bibitem[Song, Chen and Kato(2019)]{SongChenKato_2019}
Song, Y., Chen, X. and Kato, K. (2019). Approximating high-dimensional infinite-order $U$-statistics: Statistical and computational guarantees. {\sl Electron. J. Stat.} {\bf 13} 4794--4848.

\bibitem[Sunklodas(1984)]{Sunklodas_1984}
Sunklodas, J. (1984). Rate of convergence in the central limit theorem for random variables with strong mixing. {\sl Lith. Math. J.} {\bf24} 182--190.


\bibitem[van der Vaart(1998)]{vanderVaart_1998}
van der Vaart, A. W. (1998). {\it Asymptotic Statistics.} Cambridge University Press.




\bibitem[Wainwright(2019)]{wainwright2019high}
Wainwright, M. J. (2019). \emph{High-dimensional Statistics: A Non-asymptotic Viewpoint}. Cambridge University Press.


\bibitem[Wong, Li and Tewari(2020)]{Wong2020}
Wong, K. C., Li, Z. and Tewari, A. (2020). Lasso guarantees for $\beta$-mixing heavy-tailed time series.  {\sl Ann. Statist.} {\bf48} 1124--1142.

\bibitem[Wu(2005)]{Wu_2005}
Wu, W. B. (2005). Nonlinear system theory: Another look at dependence. {\sl Proc. Natl. Acad. Sci. USA} {\bf102} 14150--14154.


\bibitem[Wu(2007)]{Wu_2007}
Wu, W. B. (2007). Strong invariance principles for dependent random variables. {\sl Ann. Probab.} {\bf35} 2294--2320.


\bibitem[Wu and Shao(2004)]{WuShao_2004}
Wu, W. B. and Shao, X. (2004). Limit theorems for iterated random functions. {\sl J. Appl. Probab.} {\bf41} 425--436.


\bibitem[Wu and Wu(2016)]{WuWu_2016}
Wu, W. B. and Wu, Y. N. (2016). {Performance bounds for parameter estimates of high-dimensional linear models with correlated errors}. {\sl Electron. J. Stat.} {\bf10} 352--379.




\bibitem[Yu and Chen(2021)]{YuChen_2020}
Yu, M. and Chen. X. (2021). {Finite sample change point inference and identification for high-dimensional mean vectors}. {\sl J. R. Stat. Soc. Ser. B. Stat. Methodol.} {\bf83} 247--270.


\bibitem[Yu and Chen(2022)]{YuChen_2019}
Yu, M. and Chen. X. (2022). {A robust bootstrap change point test for high-dimensional location parameter}. {\sl Electron. J. Stat.} {\bf16} 1096--1152.


\bibitem[Zhang and Wu(2017)]{ZhangWu_2017}
Zhang, D. and Wu, W. B. (2017). Gaussian approximation for high dimensional time series. {\sl Ann. Statist.} {\bf45} 1895--1919.

\bibitem[Zhang(2015)]{Zhang_2015}
Zhang, X. (2015). Testing high dimensional mean under sparsity. {\sl arXiv:1509.08444}.


\bibitem[Zhang and Cheng(2018)]{ZhangCheng_2018}
Zhang, X. and Cheng, G. (2018). Gaussian approximation for high dimensional vector under physical dependence. {\sl Bernoulli} {\bf 24} 2640--2675.















\end{thebibliography}
\end{document}